\documentclass[a4,10pt]{article}
\usepackage{amsmath,amscd,amssymb}

\newcommand{\nbigd}{\mathcal{D}}
\newcommand{\nbige}{\mathcal{E}}
\newcommand{\nbigf}{\mathcal{F}}
\newcommand{\nbigg}{\mathcal{G}}
\newcommand{\nbigh}{\mathcal{H}}
\newcommand{\nbigi}{\mathcal{I}}

\newcommand{\nbigk}{\mathcal{K}}
\newcommand{\nbigl}{\mathcal{L}}
\newcommand{\nbigm}{\mathcal{M}}
\newcommand{\nbign}{\mathcal{N}}
\newcommand{\nbigo}{\mathcal{O}}
\newcommand{\nbigp}{\mathcal{P}}
\newcommand{\nbigq}{\mathcal{Q}}
\newcommand{\nbigr}{\mathcal{R}}
\newcommand{\nbigs}{\mathcal{S}}
\newcommand{\nbigt}{\mathcal{T}}
\newcommand{\nbigu}{\mathcal{U}}
\newcommand{\nbigv}{\mathcal{V}}

\newcommand{\nbigx}{\mathcal{X}}
\newcommand{\nbigy}{\mathcal{Y}}
\newcommand{\nbigz}{\mathcal{Z}}

\newcommand{\proj}{\mathbb{P}}
\newcommand{\seisuu}{\mathbb{Z}}

\newcommand{\cnum}{{\boldsymbol C}}
\newcommand{\real}{{\boldsymbol R}}

\newcommand{\Tate}{\mathbb{T}}

\newcommand{\DD}{\mathbb{D}}
\newcommand{\EE}{\mathbb{E}}

\newcommand{\gbigv}{\mathfrak V}

\newcommand{\gbigz}{\mathfrak Z}

\newcommand{\gminia}{\mathfrak a}
\newcommand{\gminib}{\mathfrak b}
\newcommand{\gminic}{\mathfrak c}

\newcommand{\gminie}{\mathfrak e}

\newcommand{\gminih}{\mathfrak h}

\newcommand{\gminik}{\mathfrak k}

\newcommand{\gminip}{\mathfrak p}

\newcommand{\vecv}{{\boldsymbol v}}
\newcommand{\vecu}{{\boldsymbol u}}
\newcommand{\vecw}{{\boldsymbol w}}

\newcommand{\veczero}{{\boldsymbol 0}}
\newcommand{\vecalpha}{{\boldsymbol \alpha}}
\newcommand{\veca}{{\boldsymbol a}}
\newcommand{\vecb}{{\boldsymbol b}}

\newcommand{\vecdelta}{{\boldsymbol \delta}}

\newcommand{\vecc}{{\boldsymbol c}}

\newcommand{\vech}{{\boldsymbol h}}
\newcommand{\veck}{{\boldsymbol k}}
\newcommand{\vecm}{{\boldsymbol m}}

\newcommand{\vecN}{{\boldsymbol N}}

\newcommand{\vecf}{{\boldsymbol f}}

\newcommand{\vecE}{{\boldsymbol E}}

\newcommand{\vecn}{{\boldsymbol n}}

\newcommand{\vecz}{{\boldsymbol z}}

\newcommand{\lrarr}{\longrightarrow}

\newcommand{\pf}{{\bf Proof}\hspace{.1in}}
\newcommand{\qed}{\mbox{\rule{1.2mm}{3mm}}}

\def\Hom{\mathop{\rm Hom}\nolimits}

\def\End{\mathop{\rm End}\nolimits}

\def\Image{\mathop{\rm Im}\nolimits}

\def\Re{\mathop{\rm Re}\nolimits}

\def\Gr{\mathop{\rm Gr}\nolimits}

\def\rank{\mathop{\rm rank}\nolimits}

\def\Ker{\mathop{\rm Ker}\nolimits}

\def\Gr{\mathop{\rm Gr}\nolimits}
\def\Sym{\mathop{\rm Sym}\nolimits}

\def\Res{\mathop{\rm Res}\nolimits}

\def\ord{\mathop{\rm ord}\nolimits}

\def\tr{\mathop{\rm tr}\nolimits}

\def\can{\mathop{\rm can}\nolimits}

\def\id{\mathop{\rm id}\nolimits}

\def\Irr{\mathop{\rm Irr}\nolimits}

\newcommand{\del}{\partial}
\newcommand{\delbar}{\overline{\del}}

\newcommand{\nbar}{\underline{n}}

\newcommand{\jbar}{\underline{j}}
\newcommand{\mbar}{\underline{m}}

\newcommand{\ibar}{\underline{i}}

\newcommand{\pbar}{\underline{p}}
\newcommand{\itibar}{\underline{1}}
\newcommand{\nibar}{\underline{2}}

\newcommand{\sankaku}{\triangle}

\newcommand{\harmonicbundle}{(E,\delbar_E,\theta,h)}

\newcommand{\barz}{\overline{z}}
\newcommand{\zbar}{\barz}

\newcommand{\baralpha}{\overline{\alpha}}
\newcommand{\alphabar}{\baralpha}
\newcommand{\barlambda}{\overline{\lambda}}
\newcommand{\lambdabar}{\barlambda}

\newcommand{\etabar}{\overline{\eta}}

\newcommand{\Hbar}{\overline{H}}

\newcommand{\Poin}{{\bf p}}
\newcommand{\poin}{\Poin}

\newcommand{\prolongg}[2]{{}_{#1}{#2}}

\newcommand{\DDlambda}{\DD^{\lambda}}

\newcommand{\nbigelambda}{\nbige^{\lambda}}

\newcommand{\nbigxlambda}{\nbigx^{\lambda}}

\newcommand{\KMS}{{\mathcal{KMS}}}

\newcommand{\Par}{{\mathcal Par}}
\newcommand{\Sp}{{\mathcal Sp}}

\newcommand{\kmsmap}{\gminik}
\newcommand{\paramap}{\gminip}
\newcommand{\eigenmap}{\gminie}

\newcommand{\lefttop}[1]{{}^{#1}}

\def\irr{\mathop{\rm irr}\nolimits}

\newcommand{\Fzero}{F^{(\lambda_0)}}

\newcommand{\nbigqzero}{\nbigq^{(\lambda_0)}}

\newcommand{\lamda}{\lambda}

\newcommand{\openclosed}[2]{]#1,#2]}

\newcommand{\Ehat}{\widehat{E}}

\newcommand{\vecvhat}{\widehat{\vecv}}

\newcommand{\Gtilde}{\widetilde{G}}

\newcommand{\Etilde}{\widetilde{E}}
\newcommand{\vecEtilde}{\widetilde{\vecE}}
\newcommand{\thetatilde}{\widetilde{\theta}}

\newcommand{\nablahat}{\widehat{\nabla}}
\newcommand{\Vtilde}{\widetilde{V}}

\newcommand{\vecvtilde}{\widetilde{\vecv}}

\newcommand{\vecwbar}{\overline{\vecw}}

\newcommand{\Phitilde}{\widetilde{\Phi}}
\newcommand{\vtilde}{\widetilde{v}}

\newcommand{\Ftilde}{\widetilde{F}}

\newcommand{\gminiabar}{\overline{\gminia}}

\newcommand{\stilde}{\widetilde{s}}

\newcommand{\nbigltilde}{\widetilde{\nbigl}}
\newcommand{\nbigdhat}{\widehat{\nbigd}}
\newcommand{\DDtilde}{\widetilde{\DD}}
\newcommand{\nbigxhat}{\widehat{\nbigx}}

\newcommand{\Dhat}{\widehat{D}}

\newcommand{\vecwhat}{\widehat{\vecw}}

\newcommand{\Sbar}{\overline{S}}
\newcommand{\nbigehat}{\widehat{\nbige}}

\newcommand{\nbigvtilde}{\widetilde{\nbigv}}

\newcommand{\ellsitabar}{\underline{\ell}}

\newcommand{\Dtilde}{\widetilde{D}}
\newcommand{\Xtilde}{\widetilde{X}}

\newcommand{\Zhat}{\widehat{Z}}
\newcommand{\What}{\widehat{W}}
\newcommand{\DDhat}{\widehat{\DD}}

\def\ord{\mathop{\rm ord}\nolimits}

\def\full{\mathop{\rm full}\nolimits}

\def\Gal{\mathop{\rm Gal}\nolimits}

\def\exchange{\mathop{\rm exchange}\nolimits}
\def\HS{\mathop{\rm HS}\nolimits}
\def\TNIL{\mathop{\rm TNIL}\nolimits}
\def\Glue{\mathop{\rm Glue}\nolimits}

\newcommand{\Wtilde}{\widetilde{W}}

\newcommand{\nbigutilde}{\widetilde{\nbigu}}
\newcommand{\nbigstilde}{\widetilde{\nbigs}}

\newcommand{\gtilde}{\widetilde{g}}

\newcommand{\Stilde}{\widetilde{S}}

\newcommand{\nbigftilde}{\widetilde{\nbigf}}

\newcommand{\nbigpzero}{\nbigp^{(\lambda_0)}}

\newcommand{\wtilde}{\widetilde{w}}

\newcommand{\nbigxtilde}{\widetilde{\nbigx}}

\newcommand{\nbigdtilde}{\widetilde{\nbigd}}

\newcommand{\Ttilde}{\widetilde{T}}

\newcommand{\nbigsbar}{\overline{\nbigs}}

\newcommand{\nbigitilde}{\widetilde{\nbigi}}

\newcommand{\Multisector}{\nbigm\nbigs}

\newcommand{\Irrbar}{\overline{\Irr}}

\newcommand{\nbigktilde}{\widetilde{\nbigk}}

\newcommand{\nbiggzero}{\nbigg^{(\lambda_0)}}

\newcommand{\TTtilde}{{T}\widetilde{T}}
\newcommand{\Hhat}{\widehat{H}}
\newcommand{\Htilde}{\widetilde{H}}
\newcommand{\Omegatilde}{\widetilde{\Omega}}

\newcommand{\nbigxzero}{\nbigx^{(\lambda_0)}}
\newcommand{\nbigdzero}{\nbigd^{(\lambda_0)}}
\newcommand{\nutilde}{\widetilde{\nu}}
\newcommand{\vecmtilde}{\widetilde{\vecm}}
\newcommand{\gminiatilde}{\widetilde{\gminia}}

\newcommand{\cnumtilde}{\widetilde{\cnum}}

\newcommand{\HSbar}{\overline{\HS}}
\newcommand{\phitilde}{\widetilde{\phi}}
\newcommand{\nbightilde}{\widetilde{\nbigh}}

\newcommand{\Vbar}{\overline{V}}
\newcommand{\Wbar}{\overline{W}}
\newcommand{\Nbar}{\overline{N}}
\newcommand{\DDbar}{\overline{\DD}}
\newcommand{\Ebar}{\overline{E}}

\newcommand{\ahat}{\widehat{a}}
\newcommand{\bhat}{\widehat{b}}
\newtheorem{thm}{Theorem}[section]
\newtheorem{cor}[thm]{Corollary}

\newtheorem{rem}[thm]{Remark}
\newtheorem{lem}[thm]{Lemma}
\newtheorem{prop}[thm]{Proposition}
\newtheorem{df}[thm]{Definition}

\newtheorem{notation}[thm]{Notation}

\begin{document}

\title{Asymptotic behaviour of
 variation of pure polarized TERP structures}
\author{Takuro Mochizuki}
\date{}
\maketitle

\begin{abstract}
The purpose of this paper is twofold.
One is to give a survey of our study
on the reductions of harmonic bundles,
and the other is to explain a simple application
in the study of TERP structure.
In particular, we investigate
the asymptotic behaviour of
the ``new supersymmetric index''
for variation of pure polarized TERP structures.

\vspace{.1in}

\end{abstract}

\section{Introduction}

In our previous papers
\cite{mochi},
\cite{mochi2} and \cite{mochi7},
we studied asymptotic behaviour of
tame and wild harmonic bundles.
Briefly,
one of the main results
is the following sequence of 
reductions of harmonic bundles:
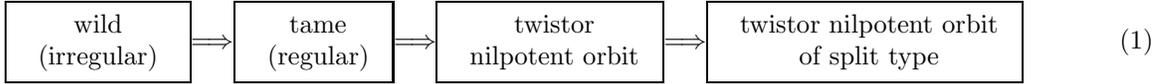
\begin{equation}
 \label{eq;08.7.27.1}
 \fbox{
 \begin{tabular}{c}
 wild \\
(irregular)
 \end{tabular}
 }
{\Longrightarrow}
 \fbox{
 \begin{tabular}{c}
 tame\\
(regular)
 \end{tabular}}
{\Longrightarrow}
  \fbox{
 \begin{tabular}{c}
 twistor \\
 nilpotent orbit
 \end{tabular}}
{\Longrightarrow}
  \fbox{
 \begin{tabular}{c}
 twistor 
 nilpotent orbit\\
 of split type
 \end{tabular}}
\end{equation}
A reduced object is simpler than the original one,
but it still gives a good approximation
of the original one.
And, a twistor nilpotent orbit of split type
comes from a variation of polarized pure Hodge structures,
whose asymptotic behaviour was deeply studied by
E. Cattani, A. Kaplan, M. Kashiwara, 
T. Kawai and W. Schmid.
Thus, 
we can say that 
the asymptotic behaviour of wild
harmonic bundles is understood pretty well.

The main purpose of this paper is twofold.
One is to give a survey of these reductions,
and the other is to explain a simple application
in the study of TERP structure.

C. Hertling \cite{Hertling} initiated the study of
TERP structures inspired by mathematical physics
and singularity theory.
The study was further developed by
Hertling and C. Sevenheck.
For example, they investigated
``nilpotent orbit'' \cite{Hertling-Sevenheck},
asymptotic behaviour of 
tame variation of TERP structures
and classifying spaces \cite{Hertling-Sevenheck3}.
We refer to the above papers and 
a survey \cite{Hertling-Sevenheck2}
for more details and precise.

\begin{rem}
Their ``nilpotent orbit'' is called
``HS-orbit'' 
(Hertling-Sevenheck orbit) in this paper.
We can consider several kinds of
generalization of ``nilpotent orbit''
in the theory of TERP structures
and twistor structures.
HS-orbit is the one.
Another one is twistor nilpotent orbit
studied in {\rm \cite{mochi2}},
which we will mainly use in this paper.
\hfill\qed
\end{rem}

\begin{rem}
We prefer to regard
TERP structure as
integrable twistor structure
with a real structure and a pairing
studied by C. Sabbah.
It is called twistor-TERP structure
in this paper.
\hfill\qed
\end{rem}

We will give an enrichment
of the sequence (\ref{eq;08.7.27.1})
with TERP structures
or integrable twistor structures.
As an application,
we will study the behaviour of 
``new supersymmetric index''
of variation of pure polarized TERP structures.
Let $\nabla$ be a meromorphic connection
of $V=\nbigo_{\proj^1}^{\oplus\,r}$
admitting a pole at $\{0,\infty\}$
of at most order two.
Let $d$ be the natural connection of $V$.
Then, we have the expression
$ \nabla=d+\bigl(
 \lambda^{-1}\cdot\nbigu_1
-\nbigq-\lambda\cdot\nbigu_2
 \bigr)\cdot d\lambda/\lambda$,
where
$\nbigu_i,\nbigq\in \End(V)$.
If $(V,\nabla)$ is equipped with
a real structure and a polarization
(see Subsection \ref{subsubsection;08.8.8.1}),
there is some more restriction on them.
Anyway, $\nbigq$ is called 
the supersymmetric index of
$(V,\nabla)$.
We set $X:=\bigl\{(z_1,\ldots,z_n)\,\big|\,
 |z_i|<1 \bigr\}$
and $D:=\bigcup_{i=1}^n\{z_i=0\}$.
Let $(\nbigv,\DDtilde^{\sankaku},\nbigs,\kappa)$
be a variation of pure polarized 
twistor-TERP structures of weight $0$
on $\proj^1\times(X-D)$.
(See Subsection \ref{subsection;08.8.20.10}.)
It is called unramifiedly good wild (resp. tame),
if the underlying harmonic bundle
$(E,\delbar_E,\theta,h)$ is so.
(See Subsection \ref{subsection;08.10.27.1}.)
For each point $P\in X-D$,
we have the new supersymmetric index
$\nbigq_P
 \in \End(\nbigv^{\sankaku}_{|\proj^1\times P})
\simeq\End(E_{|P})$ of 
$(\nbigv^{\sankaku},
 \DDtilde^{\sankaku})_{|\proj^1\times P}$,
and thus we obtain a $C^{\infty}$-section
$\nbigq$ of $\End(E)$.
We are interested in the behaviour of
$\nbigq$ around $(0,\ldots,0)$.
The result is the following:
\begin{itemize}
\item
In the case of twistor-TERP nilpotent orbit of split type,
the new supersymmetric index 
can be easily computed from the data of 
the corresponding polarized mixed twistor-TERP structure.
In particular, their eigenvalues are constant.
(See Section \ref{section;08.8.9.1}.)
\item
From a twistor-TERP nilpotent orbit
$(\nbigv,\DDtilde^{\sankaku},\nbigs,\kappa)$,
we obtain a twistor-TERP nilpotent orbit 
of split type
$(\nbigv_0,\DDtilde^{\sankaku}_0,\nbigs_0,\kappa_0)$, 
by taking Gr with respect to the weight filtration.
(Precisely, Gr is taken for
 the corresponding polarized mixed twistor-TERP structure.)
The new supersymmetric index
$\nbigq$ of $(\nbigv,\DDtilde^{\sankaku})$
can be approximated
by the new supersymmetric index
$\nbigq_0$ of $(\nbigv_0,\DDtilde^{\sankaku}_0)$
up to $O\Bigl(\sum(-\log |z_i|)^{-1/2}\Bigr)$.
In particular, the eigenvalues of $\nbigq$
are constant up to 
$O\Bigl(\sum(-\log |z_i|)^{-\delta}\Bigr)$
for some $\delta>0$.
(See Section \ref{section;08.8.9.2}.)
\item
From a tame variation of 
polarized pure twistor-TERP-structures
$(\nbigv,\DDtilde^{\sankaku},\nbigs,\kappa)$,
we obtain a twistor-TERP nilpotent orbit 
$(\nbigv_0,\DDtilde^{\sankaku}_0,\nbigs_0,\kappa_0)$
associated to the limit mixed twistor-TERP structure
which was essentially considered 
in \cite{Hertling-Sevenheck3}
as an enrichment of the limit mixed twistor structure
in \cite{mochi2}.
We can approximate the new supersymmetric index
$\nbigq$ of $(\nbigv,\DDtilde^{\sankaku})$
by the new supersymmetric index
$\nbigq_0$ of $(\nbigv_0,\DDtilde^{\sankaku}_0)$
up to $O\Bigl(\sum |z_i|^{\epsilon}\Bigr)$
for some $\epsilon>0$.
In particular, 
the eigenvalues of $\nbigq_0$
approximate those of $\nbigq$
 up to $O\Bigl(\sum |z_i|^{\epsilon'}\Bigr)$
for some $\epsilon'>0$.
(See Subsection \ref{subsection;08.8.9.3}
 for more precise statements.)
\item
From a wild variation of 
polarized pure twistor-TERP structures
$(\nbigv,\DDtilde^{\sankaku},\nbigs,\kappa)$,
we obtain a tame variation of 
polarized pure twistor-TERP structures
$(\nbigv_0,\DDtilde^{\sankaku}_0,\nbigs_0,\kappa_0)$,
by taking Gr with respect to Stokes filtrations.
We can approximate the new supersymmetric index
$\nbigq$ of $(\nbigv,\DDtilde^{\sankaku})$
by the new supersymmetric index
$\nbigq_0$ of $(\nbigv_0,\DDtilde^{\sankaku}_0)$
up to a term with exponential decay.
In particular, the eigenvalues of $\nbigq_0$
approximate those of $\nbigq$
up to exponential decay.
(See Subsection \ref{subsection;08.8.9.4}
for more precise statements.)
\end{itemize}

In each case,
we will construct a $C^{\infty}$-map
$\nbigv_0\lrarr \nbigv$,
which does not preserve 
but approximate the additional structures.
(More precisely,
 $\nbigv_0$ should be twisted.)
It would be interesting to clarify the precise relation
between these results and the celebrated 
nilpotent orbit theorem for Hodge structures 
due to W. Schmid \cite{sch}.
(See also \cite{Hertling-Sevenheck3}.)

As a corollary,
we obtain the convergence of 
the eigenvalues of new supersymmetric indices
of wild harmonic bundles on a punctured disc.
In his recent work 
(Section {\rm 3} of {\rm\cite{sabbah10}}),
Sabbah studied
the eigenvalues of new supersymmetric indices
for polarized wild pure integrable twistor $D$-modules
on curves.
Since wild harmonic bundles
are prolonged to polarized wild pure twistor
$D$-modules {\rm\cite{mochi7}},
we can also deduce the above convergence
in the curve case from his results.

We also show that
if a TERP-structure induces an HS-orbit,
then it is a mixed-TERP structure
in the sense of \cite{Hertling-Sevenheck}
by using the reduction from wild to tame,
which was conjectured by Hertling and Sevenheck.

\paragraph{Outline of this paper}

In Subsection \ref{subsection;08.8.20.10}
we recall integrable pure twistor structure
and TERP structure
and their variations in our convenient way,
which were originally studied by 
Hertling, Sabbah and Sevenheck.
We look at some basic examples 
in Subsection \ref{subsection;08.9.10.20}.
In particular,
we introduce the notions of 
integrable twistor nilpotent orbit
and twistor-TERP nilpotent orbit.
In Subsection \ref{subsection;08.8.20.11},
we argue a convergence of 
integrable pure twistor structures
and new supersymmetric indices.
The result will be used in many times.
In Subsection \ref{subsection;08.8.20.13},
we consider a variation of
polarized mixed twistor structures.
In Subsection \ref{subsubsection;08.8.11.2},
we explain the reduction
from polarized mixed twistor structure
to polarized mixed twistor structure 
{\em of split type}.
In Subsection \ref{subsubsection;08.8.20.14},
we give a $C^{\infty}$-splitting
of weight filtrations
compatible with nilpotent maps,
which is a preparation for 
Section \ref{section;08.8.9.2}.

In Section \ref{section;08.8.9.1},
we study polarized mixed twistor structure
{\em of split type}
with some additional structures.
It is quite easy to handle.
In Section \ref{section;08.8.9.2},
we show the correspondence 
between twistor nilpotent orbits
and polarized mixed twistor structures.
We have already established
the way from twistor nilpotent orbits
to polarized mixed twistor structures
in \cite{mochi2}.
The converse was also established
in the curve case.
The higher dimensional case is new.
The correspondence is easily enriched
with integrability and real structures.
We also show that a twistor nilpotent orbit
is approximated with
a twistor nilpotent orbit {\em of split type}.

In Section \ref{section;08.8.20.40},
we give a review on 
Stokes structure and reductions
for a family of meromorphic $\lambda$-flat bundles,
studied in Sections 7 and 8 in \cite{mochi7}.
We give some minor complementary results
on connections along the $\lambda$-direction
and pseudo-good lattices.

In Section \ref{section;08.8.20.41},
we explain the reduction
from unramifiedly good wild harmonic bundles
to polarized mixed twistor structures,
studied in \cite{mochi2} and \cite{mochi7}.
We give a review on the prolongation of
harmonic bundles in
Subsection \ref{subsection;08.8.3.6}.
Then, in Subsection \ref{subsection;08.8.3.15},
we review the reduction
from unramifiedly good wild bundles
to tame harmonic bundles as the Gr
with respect to Stokes filtrations,
which is one of the main results in
\cite{mochi7},
and in Subsection \ref{subsection;08.8.4.10},
we review the reduction
from tame harmonic bundles
to polarized mixed twistor structure
as the Gr with respect to KMS-structure,
which is one of the main results in
\cite{mochi2}.
Together with the result
in Section \ref{section;08.8.9.2},
we can regard it as the reduction
to nilpotent orbits.

In Section \ref{section;08.8.20.42},
we argue an enrichment of the reductions
with integrability and real structure.
One of the main issues is to obtain
a meromorphic extension of
the connection along the $\lambda$-direction.
For that purpose, we prepare some estimate
in Subsection \ref{subsection;08.8.20.44}.
Then, it is easy to obtain
the meromorphic prolongment of
variations of integrable twistor structures
and the enrichment of
the sequence of reductions
as in (\ref{eq;08.7.27.1}).
We also show that 
the reduced one gives a good approximation
of the original one.
In particular, we obtain the results
on approximation of
the new supersymmetric indices
of wild or tame variation of
integrable twistor structures.

In Section \ref{section;08.8.20.45},
we study the reduction of HS-orbit.

\paragraph{Acknowledgement}

This paper grew out of my effort
to understand the work
due to Claus Hertling, Claude Sabbah 
and Christian Sevenheck
on TERP structure and integrable twistor structure.
I am grateful to them who attracted
my attention to this subject.
I also thank their comments
on the earlier versions of this paper.
In particular,
Hertling kindly sent a surprisingly detailed 
and careful report,
which was quite helpful for improving this paper
and correcting some errors in earlier versions.

I wish to express my thanks
to Yoshifumi Tsuchimoto and Akira Ishii
for their constant encouragement.

I gave talks on the sequence (\ref{eq;08.7.27.1})
at the conferences 
``From tQFT to tt$^{\ast}$ and integrability'' in Augsburg
and 
``New developments in Algebraic Geometry, 
 Integrable Systems and Mirror symmetry'' in Kyoto.
This paper is an enhancement of the talks.
I would like to express my gratitude
to the organizers of the conferences
on this occasion.

I am grateful to the partial financial support
by Ministry of Education, Culture, 
Sports, Science and Technology.

\section{Preliminary}

\subsection{Integrable twistor structure}

\label{subsection;08.8.20.10}

We recall the notion of integrable twistor structures
and TERP structures in our convenient way
just for our understanding.
See 
\cite{Hertling},
\cite{Hertling-Sevenheck} and \cite{sabbah2}
for the original definitions 
and for more details.
We also recall twistor structures
introduced in \cite{s3}.
See also \cite{mochi} and \cite{mochi2}.

\subsubsection{Some sheaves and differential operators
 on $\proj^1\times X$}

Let $\proj^1$ denote a one dimensional
complex projective space.
We regard it as the gluing of 
two complex lines
$\cnum_{\lambda}$ and $\cnum_{\mu}$
by $\lambda=\mu^{-1}$.
We set $\cnum_{\lambda}^{\ast}:=
 \cnum_{\lambda}-\{0\}$.

Let $X$ be a complex manifold.
We set $\nbigx:=\cnum_{\lambda}\times X$
and $\nbigx^0:=\{0\}\times X$.
Let $\Omegatilde^{1,0}_{\nbigx}$ be 
the $C^{\infty}$-bundle
associated to
$\Omega^{1,0}_{\nbigx}(\log \nbigx^0)
\otimes \nbigo_{\nbigx}(\nbigx^0)$.
We put $\Omegatilde^{0,1}_{\nbigx}:=
\Omega^{0,1}_{\nbigx}$,
and we define
\[
 \Omegatilde^1_{\nbigx}:=
 \Omegatilde^{1,0}_{\nbigx}
\oplus
 \Omegatilde^{0,1}_{\nbigx},
\quad
 \Omegatilde^{\cdot}_{\nbigx}:=
 \bigwedge^{\cdot}
 \Omegatilde^{1}_{\nbigx}
\]
The associated sheaves of $C^{\infty}$-sections
are denoted by the same symbols.
Let $\DDtilde_X^f: 
 \Omegatilde^{\cdot}_{\nbigx}
\lrarr
 \Omegatilde^{\cdot+1}_{\nbigx}$
denote the differential operator induced by
the exterior differential $d$.

Let $X^{\dagger}$ denote the conjugate of $X$.
We set $\nbigx^{\dagger}:=
 \cnum_{\mu}\times X^{\dagger}$.
By the same procedure,
we obtain the $C^{\infty}$-bundles
$\Omegatilde^{\cdot}_{\nbigx^{\dagger}}$
with the differential operator
$\DDtilde^{\dagger\,f}_X$.

Their restrictions to 
$\cnum_{\lambda}^{\ast}\times X
=\cnum_{\mu}^{\ast}\times X^{\dagger}$
are naturally isomorphic:
\[
 \bigl(\Omegatilde^{\cdot}_{\nbigx},
 \DDtilde_X^f\bigr)
 _{|\cnum_{\lambda}^{\ast}\times X}
=
 \bigl(\Omega^{\cdot}
 _{\cnum_{\lambda}^{\ast}\times X},d
 \bigr)
=\bigl(\Omegatilde^{\cdot}_{\nbigx^{\dagger}},
 \DDtilde_X^{\dagger\,f}\bigr)
 _{|\cnum_{\mu}^{\ast}\times X^{\dagger}}
\]
By gluing them,
we obtain the $C^{\infty}$-bundles
$\Omegatilde^{\cdot}_{\proj^1\times X}$
with a differential operator
$\DDtilde^{\sankaku}_X$.

\begin{rem}
$\DDtilde^f_X$ and 
$\DDtilde^{\dagger\,f}_X$
are denoted also by $d$,
if there is no risk of confusion.
\hfill\qed
\end{rem}

We have the decomposition
$\Omegatilde^1_{\proj^1\times X}
=\xi\Omega^1_{X}
\oplus
 \Omegatilde^1_{\proj^1}$
into the $X$-direction and the $\proj^1$-direction.
The restriction of $\DDtilde^{\sankaku}_X$
to the $X$-direction is denoted by
$\DD^{\sankaku}_X$.
The restriction to the $\proj^1$-direction
is denoted by $d_{\proj^1}$.
We have the decomposition
\[
 \Omegatilde^{1}_{\proj^1}
=\pi^{\ast}
 \Omega^{1,0}_{\proj^1}(2\cdot\{0,\infty\})
\oplus 
 \pi^{\ast}\Omega^{0,1}_{\proj^1},
\]
into the $(1,0)$-part and the $(0,1)$-part,
where $\pi$ denotes the projection
$\proj^1\times X\lrarr \proj^1$.
We have the corresponding decomposition
$d_{\proj^1}=\del_{\proj^1}+\delbar_{\proj^1}$.

\vspace{.1in}

Let $\nu:\proj^1\lrarr\proj^1$
be a diffeomorphism.
Assume $\nu$ satisfies one of the following:
\begin{description}
\item[(A1)]
$\nu$ is holomorphic with
$\nu(0)=0$ and $\nu(\infty)=\infty$.
\item[(A2)]
$\nu$ is anti-holomorphic with
$\nu(0)=\infty$ and $\nu(\infty)=0$.
\end{description}
In particular,
we will often use the maps
$\sigma$, $\gamma$ and $j$:
\[
 \sigma([z_0:z_1])
=[-\zbar_1:\zbar_0], 
\quad
 \gamma([z_0:z_1])
=[\zbar_1:\zbar_0],
\quad
 j([z_0:z_1])=[-z_0:z_1]
\]
The induced diffeomorphism
$\proj^1\times X\lrarr\proj^1\times X$
is also denoted by $\nu$.
In the case (A1),
we have the natural isomorphism
$\Phi_{\nu}:
 \nu^{\ast}\Omegatilde^{\cdot}_{\proj^1\times X}
\simeq
 \Omegatilde^{\cdot}_{\proj^1\times X}$
of $C^{\infty}$-vector bundles
given by the ordinary pull back.
In the case (A2),
the multiplication of $C^{\infty}$-functions on 
$\nu^{\ast}\Omegatilde^{\cdot}_{\proj^1\times X}$
is twisted as 
$g\cdot \nu^{\ast}(\omega)
=\nu^{\ast}\bigl(
 \overline{\nu^{\ast}(g)}\cdot\omega
 \bigr)$
for a function $g$ and 
a section $\omega$ of 
$\Omegatilde^{\cdot}_{\proj^1\times X}$.
Then, we have the $C^{\infty}$-isomorphism
$\Phi_{\nu}:
 \nu^{\ast}\Omegatilde^{\cdot}_{\proj^1\times X}
 \simeq
 \Omegatilde^{\cdot}_{\proj^1\times X}$
given by the complex conjugate
and the ordinary pull back
\[
 \Phi_{\nu}(\nu^{\ast}\omega)
=\overline{\nu^{\ast}(\omega)}.
\]
It is easy to check that 
$\Phi_{\nu}\circ
 \nu^{\ast}(\DDtilde_X^{\sankaku})
=\DDtilde_X^{\sankaku}\circ\Phi_{\nu}$.
Similar relations hold for
$\DD^{\sankaku}_X$ and $d_{\proj^1}$.
If we are given an additional bundle $\nbigf$,
the induced isomorphism
$\nbigf\otimes \nu^{\ast}\bigl(
 \Omegatilde^{\cdot}_{\proj^1\times X}
 \bigr)
\simeq
 \nbigf\otimes\Omegatilde^{\cdot}_{\proj^1\times X}$
is also denoted by $\Phi_{\nu}$.

\subsubsection{Definitions and some remarks}

\paragraph{Variation of twistor structures}

Let $V$ be a $C^{\infty}$-vector bundle
on $\proj^1\times X$.
We use the same symbol
to denote the associated sheaf of
$C^{\infty}$-sections.
A $\proj^1$-holomorphic structure of
$V$ is defined to be a differential operator
\[
 d''_{\proj^1,V}:
 V\lrarr V\otimes
 \pi^{\ast}\Omega^{0,1}_{\proj^1}
\]
satisfying
(i) $d''_{\proj^1,V}(f\cdot s)
=f\cdot d''_{\proj^1,V}(s)
+\delbar_{\proj^1}(f)\cdot s$ 
for a $C^{\infty}$-function $f$
and a section $s$ of $V$,
(ii) $d''_{\proj^1,V}\circ d''_{\proj^1,V}=0$.
Such a tuple $(V,d''_{\proj^1,V})$
is called a $\proj^1$-holomorphic vector bundle.

A $\TTtilde$-structure of
$(V,d''_{\proj^1,V})$ is a differential operator
\[
 \DD^{\sankaku}_V:
 V\lrarr V\otimes\xi\Omega^1_X
\]
such that
(i) $\DD^{\sankaku}_{V}(f\cdot s)
=f\cdot \DD^{\sankaku}_{V}(s)
+\DD^{\sankaku}_X(f)\cdot s$ 
for a $C^{\infty}$-function $f$
and a section $s$ of $V$,
(ii) $(d''_{\proj^1,V}+\DD^{\sankaku}_{V})^2=0$.
Such a tuple 
$(V,d''_{\proj^1,V},\DD^{\sankaku}_V)$
is called {\em a $\TTtilde$-structure} in \cite{Hertling},
or {\em a variation of 
$\proj^1$-holomorphic vector bundles}
in \cite{mochi2}.
In this section, we prefer to call it
{\em variation of twistor structures.}

If $X$ is a point,
it is just a holomorphic vector bundle
on $\proj^1$.

\begin{rem}
We will often omit to specify $d''_{\proj^1,V}$
when we consider 
$\proj^1$-holomorphic bundles
or variations of twistor structures
(variations of $\proj^1$-holomorphic bundles).
\hfill\qed
\end{rem}

\paragraph{Variation of integrable twistor
 structures}

A $\TTtilde E$-structure of $V$ is 
a differential operator
\[
 \DDtilde^{\sankaku}_V:
 V\lrarr 
 V\otimes
 \Omegatilde^1_{\proj^1\times X}
\]
satisfying
(i) $ \DDtilde^{\sankaku}_V(f\cdot s)
=\DDtilde^{\sankaku}_X(f)\cdot s
+f\cdot\DDtilde^{\sankaku}_V(s)$
for a $C^{\infty}$-function $f$
and a section $s$ of $V$,
(ii) $ \DDtilde^{\sankaku}_V\circ
 \DDtilde^{\sankaku}_V=0$.
Such a tuple
$(V,\DDtilde^{\sankaku}_V)$
is called a variation of integrable twistor structures.

If $X$ is a point, it is equivalent to
a holomorphic vector bundle $V$ on $\proj^1$
with a meromorphic connection $\nabla$
which admits a pole at $\{0,\infty\}$
with at most order $2$,
 i.e.,
\[
 \nabla(V)\subset 
 V\otimes\Omega^1\bigl(2\cdot \{0,\infty\}\bigr).
\]
In this case, it is simply called
an integrable twistor structure.

\paragraph{Morphisms}

A morphism of variation of twistor structures
$F:(V_1,d''_{\proj^1,V_1}, \DD^{\sankaku}_{V_1})
\lrarr (V_2,d''_{\proj^1,V_2},
 \DD^{\sankaku}_{V_2})$
is defined to be 
a morphism of the associated sheaves
of $C^{\infty}$-sections,
compatible with the differential operators.
If $X$ is a point,
it is equivalent to an $\nbigo_{\proj^1}$-morphism.

A morphism of variation of
integrable twistor structures
$F:(V_1,\DDtilde^{\sankaku}_{V_1})
\lrarr (V_2,\DDtilde^{\sankaku}_{V_2})$
is defined to be 
a morphism of the associated sheaves
of $C^{\infty}$-sections,
compatible with the differential operators.
If $X$ is a point,
it is equivalent to an $\nbigo_{\proj^1}$-morphism
compatible with the meromorphic connections.

\paragraph{Some functoriality}

Let $(V,\DDtilde^{\sankaku}_V)$
be a variation of integrable twistor structures.
Let $f:Y\lrarr X$ be a holomorphic map
of complex manifolds.
Then, we have the naturally induced
variation of integrable twistor structures
$f^{\ast}(V,\DDtilde^{\sankaku}_V)$
as in the case of ordinary connections.

Let $\nu:\proj^1\lrarr\proj^1$ be 
a diffeomorphism
satisfying one of (A1) or (A2) above.
Then, $\nu^{\ast}V$
is naturally equipped with
a $\TTtilde E$-structure
$\DDtilde^{\sankaku}_{\nu^{\ast}V}$
given as follows:
\[
 \DDtilde^{\sankaku}_{\nu^{\ast}V}
\bigl(
 \Phi_{\nu}(\nu^{\ast}s)
\bigr)
=\Phi_{\nu}\Bigl(
 \nu^{\ast}\bigl(
 \DDtilde^{\sankaku}_V(s)\bigr)
 \Bigr)
\]
Here, $s$ denotes a section of
$V\otimes\Omegatilde_X^{\cdot}$.

We also have the pull back of
variation of twistor structures
via $f$ and $\nu$ as above.

\paragraph{Pure and mixed}

Let $(V,d''_{\proj^1,V})$ be 
a $\proj^1$-holomorphic vector bundle
on $\proj^1\times X$.
It is called pure of weight $w$
if the restrictions 
$V_P:=(V,d''_{\proj^1,V})_{|\proj^1\times\{P\}}$
are pure twistor structures of weight $w$
for any $P\in X$,
i.e., $V_P$ are isomorphic to
direct sums of $\nbigo_{\proj^1}(w)$.
A variation of (integrable) twistor structures
is called pure of weight $w$,
if the underlying $\proj^1$-holomorphic
vector bundle is pure of weight $w$.

\vspace{.1in}
Let $W$ be a filtration of $V$ by vector subbundles
indexed by integers.
We say that $W$ is $\proj^1$-holomorphic,
if each $W_n$ are preserved by $d''_{\proj^1,V}$.
We have induced $\proj^1$-holomorphic 
vector bundles $\Gr^W_n(V,d''_{\proj^1,V})$.
Then, $(V,d''_{\proj^1},W)$ is called mixed,
if each $\Gr^W_n(V,d''_{\proj^1,V})$
is pure of weight $n$.
When $(V,d''_{\proj^1})$ is equipped with
$\TTtilde$-structure $\DD^{\sankaku}_V$
(resp. $\TTtilde E$-structure $\DDtilde^{\sankaku}_V$),
we say that $W$ is
$\DD^{\sankaku}_V$-flat
(resp. $\DDtilde^{\sankaku}_V$-flat)
or more simply flat,
if each $W_n$ is preserved by the operator.
In that case,
$(V,d''_{\proj^1,V},\DD^{\sankaku}_V,W)$
(resp. $(V,\DDtilde^{\sankaku}_V,W)$)
is called mixed,
if $(V,d''_{\proj^1},W)$ is mixed.

\paragraph{New supersymmetric index}

Let $(V,\nabla)$ be 
a pure integrable twistor structure
of weight $0$.
We have a global trivialization
$V\simeq \nbigo_{\proj^1}^{\oplus r}$,
which is uniquely determined up to obvious ambiguity.
Let $d$ denote the natural connection
of $\nbigo_{\proj^1}^{\oplus r}$.
Then, we have the decomposition
\begin{equation}
 \label{eq;08.7.29.15}
 \nabla
=d+\Bigl(
\lambda^{-1}\nbigu_1
-\nbigq-\lambda\cdot\nbigu_2
\Bigr)\frac{d\lambda}{\lambda},
\end{equation}
where 
$\nbigu_1,\nbigu_2,\nbigq\in
 H^0\bigl(\proj^1,\End(V)\bigr)$,
The operator $\nbigq$ is called 
the new supersymmetric index.
If $(V,\nabla)$ is equipped with a polarization
(Subsection \ref{subsubsection;08.9.12.1}),
$\nbigu_2$ and $\nbigu_1$ are adjoint
with respect to the induced hermitian metric,
as observed by Hertling and Sabbah.

If we are given a variation of
polarized pure integrable twistor structures,
we obtain such operators in family.

\subsubsection{Simple examples}
\label{subsubsection;08.9.11.1}

We recall some simplest examples
of integrable pure twistor structures.

\paragraph{Example (Tate object)}
Let $\Tate(w)$ be a Tate object 
in the theory of twistor structures.
(See \cite{s3} and 
 Subsection 3.3.1 of \cite{mochi2}.)
It is isomorphic to $\nbigo_{\proj^1}(-2w)$,
and it is equipped with the distinguished frames
\[
 \Tate(w)_{|\cnum_{\lambda}}
=\nbigo_{\cnum_{\lambda}}\cdot
 t^{(w)}_0,
\quad
  \Tate(w)_{|\cnum_{\mu}}
=\nbigo_{\cnum_{\mu}}\cdot
 t^{(w)}_{\infty},
\quad
   \Tate(w)_{|\cnum_{\lambda}^{\ast}}
=\nbigo_{\cnum_{\lambda}^{\ast}}\cdot
 t^{(w)}_{1}.
\]
The transformation is given by
\[
  t_0^{(w)}=(\sqrt{-1}\lambda)^{w}\cdot t_1^{(w)},
\quad
 t_{\infty}^{(w)}=(-\sqrt{-1}\mu)^w\cdot t_1^{(w)}.
\]
In particular,
$(\sqrt{-1}\lambda)^{-2w}t_0^{(w)}=
 t_{\infty}^{(w)}$.
We have the meromorphic connection 
$\nabla_{\Tate(w)}$ on $\Tate(w)$ determined by
\[
 \nabla_{\Tate(w)}t_1^{(w)}=0,
\quad
 \nabla_{\Tate(w)}t_0^{(w)}= 
 t_0^{(w)}\cdot \left(
 w\cdot \frac{d\lambda}{\lambda}
 \right),
\quad
 \nabla_{\Tate(w)}t_{\infty}^{(w)}= 
 t_{\infty}^{(w)}\cdot \left(
 w\cdot \frac{d\mu}{\mu}
 \right).
\]
In the following,
the connection of $\Tate(w)$
is always given as above,
and hence we often omit to specify it explicitly.

We may identify $\Tate(w)$ with
$\nbigo_{\proj^1}\bigl(-w\cdot 0-w\cdot\infty\bigr)$
by the correspondence $t_1^{(w)}\longleftrightarrow 1$,
up to constant multiplication.
In particular,
we implicitly use the identification
of $\Tate(0)$ with $\nbigo_{\proj^1}$
by $t_1^{(0)}\longleftrightarrow 1$.
We will also implicitly use 
the identification
$\Tate(m)\otimes\Tate(n)
\simeq \Tate(m+n)$
given by 
$t_a^{(m)}\otimes
 t_a^{(n)}
\longleftrightarrow
 t_a^{(m+n)}$.

\paragraph{Example}

In Subsection 3.3.2 of \cite{mochi2},
we considered 
a line bundle $\nbigo(p,q)$
on $\proj^1$,
which is isomorphic to $\nbigo_{\proj^1}(p+q)$
and equipped with the distinguished frames:
\[
 \nbigo(p,q)_{|\cnum_{\lambda}}
=\nbigo_{\cnum_{\lambda}}
 \cdot f_0^{(p,q)},
\quad
 \nbigo(p,q)_{|\cnum_{\mu}}
=\nbigo_{\cnum_{\mu}}
 \cdot f_{\infty}^{(p,q)},
\quad
 \nbigo(p,q)_{|\cnum_{\lambda}^{\ast}}
=\nbigo_{\cnum_{\lambda}^{\ast}}
 \cdot f_{1}^{(p,q)}.
\]
The transformation is given by 
\[
 f_0^{(p,q)}
=(\sqrt{-1}\lambda)^{-p}\cdot
 f_1^{(p,q)},
 \quad
  f_{\infty}^{(p,q)}
=(-\sqrt{-1}\mu)^{-q}\cdot
 f_1^{(p,q)}.
\]
In particular, 
$(\sqrt{-1}\lambda)^{p+q}f_0^{(p,q)}
=f_{\infty}^{(p,q)}$.
We have the meromorphic connection
$\nabla_{\nbigo(p,q)}$ on $\nbigo(p,q)$
determined by
\[
 \nabla_{\nbigo(p,q)} f_1^{(p,q)}=0,
\quad
 \nabla_{\nbigo(p,q)}f_0^{(p,q)}
=f_0^{(p,q)}\cdot
 \left(
 -p\frac{d\lambda}{\lambda}
 \right),
\quad
  \nabla_{\nbigo(p,q)}f_{\infty}^{(p,q)}
=f_{\infty}^{(p,q)}\cdot
 \left(
 -q\frac{d\mu}{\mu}
 \right)
\]
In the following,
the connection of $\nbigo(p,q)$
is always given as above,
and hence we will often omit to  specify it
explicitly.

We may naturally identify
$\nbigo(p,q)$ with
$\nbigo_{\proj^1}(p\cdot 0+q\cdot\infty)$
by the correspondence
$f_1^{(p,q)}\longleftrightarrow 1$,
up to constant multiplication.
We will implicitly use the identification
$\nbigo(p,q)\otimes\nbigo(p',q')
\simeq
 \nbigo(p+p',q+q')$
given by
$f_a^{(p,q)}\otimes
 f_a^{(p',q')}
\longleftrightarrow
 f_a^{(p+p',q+q')}$.
We will also implicitly identify 
$\Tate(w)$ with $\nbigo(-w,-w)$
by $t_a^{(w)}=f_a^{(-w,-w)}$
for $a=0,1,\infty$.

\vspace{.1in}

Let $X$ be a complex manifold.
We have the pull back of
$\Tate(w)$ and $\nbigo(p,q)$
via the map from $X$ to a point.
They are denoted by
$\Tate(w)_X$ and $\nbigo(p,q)_X$,
respectively.
We will often omit to denote $X$,
if there is no risk of confusion.

\subsubsection{Polarization}
\label{subsubsection;08.9.12.1}

Recall that we have the isomorphism (\cite{mochi2})
\[
 \iota_{\Tate(w)}:
 \sigma^{\ast}\Tate(w)\simeq
 \Tate(w),
\]
given by the natural identification
$\sigma^{\ast}\nbigo\bigl(-w\cdot 0-w\cdot\infty\bigr)
\simeq
 \nbigo\bigl(-w\cdot 0-w\cdot\infty\bigr)$
via $\sigma^{\ast}(1)\longleftrightarrow 1$,
or equivalently,
\[
 \sigma^{\ast}t_1^{(w)}\longleftrightarrow
 t_1^{(w)},
\quad
 \sigma^{\ast}t_0^{(w)}\longleftrightarrow
 (-1)^w\cdot t_0^{(w)},
\quad
  \sigma^{\ast}t_{\infty}^{(w)}\longleftrightarrow
 (-1)^w\cdot t_{\infty}^{(w)}.
\]
It preserves the flat connections,
i.e.,
$\iota_{\Tate(w)}:
 \sigma^{\ast}(\Tate(w),\nabla_{\Tate(w)})
\simeq
 (\Tate(w),\nabla_{\Tate(w)})$.

For a variation of integrable twistor 
structures $(V,\DDtilde^{\sankaku}_V)$
on $\proj^1\times X$,
a morphism
\[
 \nbigs:(V,\DDtilde^{\sankaku}_V)
\otimes
 \sigma^{\ast}(V,\DDtilde^{\sankaku}_V)
\lrarr \Tate(-w)_X
\]
is called a pairing of weight $w$,
if it is $(-1)^w$-symmetric
in the following sense:
\[
  \iota_{\Tate(-w)}\circ
 \sigma^{\ast}\nbigs
=(-1)^{w}
 \nbigs\circ
 \exchange:
 \sigma^{\ast}V\otimes V\lrarr \Tate(-w)_X
\]
Here, $\exchange$ denotes the natural morphism
$\sigma^{\ast}V\otimes V\lrarr
 V\otimes\sigma^{\ast}V$
induced by the exchange of the components.
Similarly,
we have the notion of pairing
for variations of twistor structures.

\begin{df}
Let $(V,\DDtilde^{\sankaku}_V)$ be 
a variation of integrable pure twistor structure
of weight $w$ on $\proj^1\times X$.
Let $\nbigs:(V,\DDtilde^{\sankaku}_V)\otimes
 \sigma^{\ast}(V,\DDtilde^{\sankaku}_V)\lrarr\Tate(-w)_X$
be a pairing of weight $w$.
We say that $\nbigs$ is a polarization of 
$(V,\DDtilde^{\sankaku}_V)$,
if $\nbigs_{P}:=\nbigs_{|\proj^1\times\{P\}}$
is a polarizations of
$V_{P}:=
(V,d''_{\proj^1})_{|\proj^1\times \{P\}}$
for each $P\in X$.
Namely, the following holds:
\begin{itemize}
\item
 If $w=0$,
 the induced Hermitian pairing
 $H^0(\nbigs_P)$ of $H^0(\proj^1,V_{P})$
 is positive definite.
\item
 In the general case,
 the induced pairing
 $\nbigs_P\otimes\nbigs_{0,-w}$
 of $V_P\otimes\nbigo(0,-w)$
 is a polarization of
 the pure twistor structure. 
 (See Example $2$ below for
 $\nbigs_{0,-w}$.)
\end{itemize}
The notion of polarization
for variation of pure twistor structures
is defined in a similar way.
\hfill\qed
\end{df}

\paragraph{Example 1}

The identification $\iota_{\Tate(w)}$
induces the flat morphism
$\nbigs_{\Tate(w)}:
 \Tate(w)\otimes \sigma^{\ast}\Tate(w)
 \lrarr  \Tate(2w)$,
which is a polarization of $\Tate(w)$ of weight $-2w$.

\paragraph{Example 2}
The flat isomorphism
$\iota_{(p,q)}:
 \sigma^{\ast}\nbigo(p,q)
\simeq
 \nbigo(q,p)$ in \cite{mochi2} is given by
\[
 \sigma^{\ast}f_0^{(p,q)}
 \longmapsto
 (\sqrt{-1})^{p+q}f_{\infty}^{(q,p)},
\quad
 \sigma^{\ast}f_{\infty}^{(p,q)}
\longmapsto
 (-\sqrt{-1})^{p+q}
 f_0^{(q,p)},
\quad
 \sigma^{\ast}f_1^{(p,q)}
 \longmapsto
 (\sqrt{-1})^{q-p}f_1^{(q,p)}.
\]
Hence, we obtain the morphism
$\nbigs_{p,q}:
 \nbigo(p,q)\otimes\sigma^{\ast}\nbigo(p,q)
\lrarr \Tate(-p-q)$,
which is a polarization of weight $p+q$.

\subsubsection{Real structure and twistor-TERP structure}
\label{subsubsection;08.8.8.1}

\begin{df}
A real structure of a variation of
 integrable twistor structure 
$(V,\DDtilde^{\sankaku}_V)$
is defined to be an isomorphism
\[
 \kappa:\gamma^{\ast}
 (V,\DDtilde^{\sankaku}_V)
\simeq
 (V,\DDtilde^{\sankaku}_V)
\]
such that
$\gamma^{\ast}(\kappa)\circ\kappa=\id$.
\hfill\qed
\end{df}

We fix the real structure
$\kappa_{\Tate(w)}$ of $\Tate(w)$
given by the correspondence
\[
 \gamma^{\ast}t_1^{(w)}
\longleftrightarrow
 t_1^{(w)},
\quad
 \gamma^{\ast}t_0^{(w)}
\longleftrightarrow
 t_{\infty}^{(w)},
\quad
  \gamma^{\ast}t_{\infty}^{(w)}
\longleftrightarrow
 t_{0}^{(w)}.
\]

\begin{df}
Let $(V,\DDtilde^{\sankaku}_V)$ 
be a variation of integrable twistor structures
equipped with a pairing $\nbigs$ of weight $w$
and a real structure $\kappa$.
We say that $\kappa$ and $\nbigs$
are compatible,
if the following diagram is commutative:
\[
 \begin{CD}
 \gamma^{\ast}V\otimes
 \gamma^{\ast}\sigma^{\ast}V 
 @>{\gamma^{\ast}\nbigs}>> 
 \gamma^{\ast}\Tate(-w)\\
 @V{\kappa\otimes\sigma^{\ast}\kappa}VV 
 @V{\kappa_{\Tate(-w)}}VV\\
 V\otimes\sigma^{\ast}V
 @>{\nbigs}>> \Tate(-w)
 \end{CD}
\]
Namely, 
 $\kappa_{\Tate(-w)}\circ\gamma^{\ast}\nbigs
=\nbigs\circ(\kappa\otimes\sigma^{\ast}\kappa)$
holds.
In that case,
we also say
that $\kappa$ is a real structure of
$(V,\DDtilde^{\sankaku}_V,\nbigs)$,
or that $\nbigs$ is a pairing of
$(V,\DDtilde^{\sankaku}_V,\kappa)$ 
with weight $w$.
\hfill\qed
\end{df}

\begin{df}
Let $(V,\DDtilde^{\sankaku}_V)$ be 
a variation of integrable twistor structure
equipped with a pairing $\nbigs$ of weight $w$
and a real structure $\kappa$.
The tuple
$(V,\DDtilde^{\sankaku}_V,
 \nbigs,\kappa,-w)$
is called a variation of twistor-TERP structure,
if (i) $\nbigs$ is perfect,
(ii) $\nbigs$ and $\kappa$ are compatible.

If $X$ is a point,
it is called a twistor-TERP structure.
\hfill\qed
\end{df}

It is easy to observe that
twistor-TERP structure
is just an expression of
TERP structure \cite{Hertling}
in terms of twistor structures,
which we will explain later.

\begin{df}
A variation of twistor-TERP structures
$(V,\DDtilde^{\sankaku}_V,
 \nbigs,\kappa,-w)$ is called pure,
if $(V,\DDtilde^{\sankaku}_V)$ is pure 
with weight $w$.
It is called polarized,
if $(V,\DDtilde^{\sankaku}_V,\nbigs)$ 
is polarized.
\hfill\qed
\end{df}

\begin{rem}
If a variation of twistor-TERP structure 
$(V,\DDtilde^{\sankaku}_V,
 \nbigs,\kappa,-w)$ is pure,
we also say that
``$(V,\DDtilde^{\sankaku}_V,
 \nbigs,\kappa)$ is a variation of
pure twistor-TERP structure of weight $w$.''
\hfill\qed
\end{rem}

\paragraph{Example}

A Tate object 
$\bigl(\Tate(w),\nabla_{\Tate(w)},\nbigs_{\Tate(w)},
 \kappa_{\Tate(w)},2w\bigr)$
is a pure polarized twistor-TERP structure.

\subsubsection{Gluing construction}
\label{subsubsection;08.7.29.10}

\paragraph{Variation of
 integrable twistor structures}

We can describe a variation of integrable
twistor structures as gluing.
We set $\nbigx:=\cnum_{\lambda}\times X$,
$\nbigx^{0}:=\{0\}\times X$,
$\nbigx^{\dagger}:=
 \cnum_{\mu}\times X^{\dagger}$
and 
$\nbigx^{\dagger 0}:=
 \{0\}\times X^{\dagger}$.

Let $V_0$ be a holomorphic vector bundle on $\nbigx$
with a meromorphic flat connection
($TE$-structure \cite{Hertling})
\[
 \nabla_{V_0}:V_0\lrarr 
 V_0\otimes
 \Omega^{1,0}_{\nbigx}(\log\nbigx^0)
 \otimes\nbigo_{\nbigx}(\nbigx^0).
\]
We use the same symbol to denote
the associated differential operator
$V_0\lrarr V_0\otimes \Omegatilde^1_{\nbigx}$
in the $C^{\infty}$-category.
(The holomorphic structure 
 $d''_{V_0}$ is also included.)
Let $V_{\infty}$ be a holomorphic vector bundle
on $\nbigx^{\dagger}$
with a meromorphic flat connection
($\Ttilde E$-structure \cite{Hertling})
\[
 \nabla_{V_{\infty}}:
 V_{\infty} \lrarr 
 V_{\infty}\otimes
 \Omega^{1,0}_{\nbigx^{\dagger}}
 (\log\nbigx^{\dagger\,0})
 \otimes\nbigo_{\nbigx^{\dagger}}
 (\nbigx^{\dagger\,0}).
\]
We use the same symbol to denote
the associated differential operator
$V_{\infty}\lrarr 
 V_{\infty}\otimes
 \Omegatilde^{1}_{\nbigx^{\dagger}}$
in the $C^{\infty}$-category.
Assume that we are given 
an isomorphism $\Phi$ of
$C^{\infty}$-flat bundles:
\[
 \Phi:
 (V_0,\nabla_{V_0})
 _{|\cnum_{\lambda}^{\ast}\times X}
\simeq
 (V_{\infty},\nabla_{V_{\infty}})_{
 |\cnum_{\mu}^{\ast}\times X^{\dagger}}
\]
We obtain the $C^{\infty}$-vector bundle $V$
on $\proj^1\times X$
by gluing $V_0$ and $V_{\infty}$ via $\Phi$.
Since $\Phi$ is flat,
$\nabla_{V_0}$ and 
$\nabla_{V_{\infty}}$
induce the $\TTtilde E$-structure
$ \DDtilde^{\sankaku}_V:
 V\lrarr 
 V\otimes\Omegatilde^{1}_{\proj^1\times X}$.
Thus, we obtain a variation of
integrable twistor structures
$(V,\DDtilde^{\sankaku}_V)$.

Conversely,
we naturally obtain a tuple of
$(V_0,\nabla_{V_0})$,
$(V_{\infty},\nabla_{V_{\infty}})$
and $\Phi$ as above
from a variation of integrable twistor structures
$(V,\DDtilde^{\sankaku}_V)$
as the restriction to
$\nbigx$ and $\nbigx^{\dagger}$,
respectively.
In this situation,
we set 
\[
 \Glue\bigl(
 (V_0,\nabla_{V_0}),
 (V_{\infty},\nabla_{V_{\infty}}),
 \Phi
 \bigr):=(V,\DDtilde_V^{\sankaku}).
\]

\paragraph{Pairing and real structure}

Note that
we have the natural isomorphisms
$\nu^{\ast}
 \Omegatilde^1_{\nbigx^{\dagger}}
\simeq \Omegatilde^1_{\nbigx}$
and 
$\nu^{\ast}\Omegatilde^1_{\nbigx}
\simeq \Omegatilde^1_{\nbigx^{\dagger}}$
for anti-holomorphic diffeomorphism
$\nu:\cnum_{\lambda}\lrarr \cnum_{\mu}$
or $\cnum_{\mu}\lrarr\cnum_{\lambda}$,
as in the case of
$\Omegatilde^1_{\proj^1\times X}$.
Let $V_0$ be a holomorphic vector bundle on $\nbigx$
with a $TE$-structure $\nabla_{V_0}$.
By the above isomorphisms,
$\gamma^{\ast}V_0$
and $\sigma^{\ast}V_0$
are naturally equipped with
$\Ttilde E$-structure
$\nabla_{\gamma^{\ast}V_0}$
and $\nabla_{\sigma^{\ast}V_0}$.
Similarly,
if we are given a holomorphic vector bundle
$V_{\infty}$ on $\nbigx^{\dagger}$
with $\Ttilde E$-structure,
$\sigma^{\ast}V_{\infty}$
and $\gamma^{\ast}V_{\infty}$
are naturally equipped with $TE$-structures.
We remark that there exist
the natural isomorphisms:
\[
 \Glue\bigl(
 \gamma^{\ast}(V_{\infty},
 \nabla_{V_{\infty}}),
 \gamma^{\ast}(V_0,\nabla_{V_0}),
 \gamma^{\ast}\Phi^{-1}
 \bigr)
\simeq
 \gamma^{\ast}
 \Glue\bigl(
  (V_0,\nabla_{V_0}),
  (V_{\infty},\nabla_{V_{\infty}}),
  \Phi
 \bigr)
\]
\[
  \Glue\bigl(
 \sigma^{\ast}(V_{\infty},
 \nabla_{V_{\infty}}),
 \sigma^{\ast}(V_0,
 \nabla_{V_0}),
 \sigma^{\ast}\Phi^{-1}
 \bigr)
\simeq
 \sigma^{\ast}
 \Glue\bigl(
  (V_0,\nabla_{V_0}),
  (V_{\infty},\nabla_{V_{\infty}}),
  \Phi
 \bigr)
\]

A real structure of variation of integrable twistor structure 
corresponds to a pair of isomorphisms
\[
 \kappa_0:
 \gamma^{\ast}(V_{\infty},
 \nabla_{V_{\infty}})
\simeq
 (V_0,\nabla_{V_0}),
\quad
 \kappa_{\infty}:
  \gamma^{\ast}
 (V_0,\nabla_{V_0})
\simeq
(V_{\infty},\nabla_{V_{\infty}})
\]
such that 
(i) $\gamma^{\ast}\kappa_{0}=\kappa_{\infty}^{-1}$,
(ii) the following commutativity holds
on $\cnum_{\lambda}^{\ast}\times X$:
\[
 \begin{CD}
 \gamma^{\ast}V_{\infty}
 @>{\kappa_0}>> V_0 \\
 @V{\gamma^{\ast}\Phi^{-1}}VV @V{\Phi}VV \\
 \gamma^{\ast}V_0
 @>{\kappa_{\infty}}>>
 V_{\infty}
 \end{CD}
\]
A pairing of weight $w$ corresponds to 
\[
\nbigs_0:
 (V_0,\nabla_{V_0})
\otimes
 \sigma^{\ast}
 (V_{\infty},\nabla_{V_{\infty}})
\lrarr
 \Tate(-w)_{|\nbigx},
\quad
 \nbigs_{\infty}:
 (V_{\infty},\nabla_{V_{\infty}})
\otimes
 \sigma^{\ast}
 (V_{0},\nabla_{V_{0}})
\lrarr
 \Tate(-w)_{|\nbigx^{\dagger}}
\]
such that 
(i) $\iota_{\Tate(-w)}\circ\sigma^{\ast}\nbigs_{\infty}
=(-1)^w\nbigs_0\circ\exchange$,
(ii) it is compatible with the gluing.
Compatibility of $\nbigs$ and $\kappa$
is $\kappa_{\Tate(-w)}\circ
 \gamma^{\ast}\nbigs_{\infty}
=\nbigs_0\circ\bigl(
 \kappa_0\otimes\sigma^{\ast}\kappa_{\infty}\bigr)$.

\paragraph{Variation of twistor structures}

The above gluing description is 
essentially the same as that for a variation of 
twistor structures in \cite{s3},
which we recall in the following.
See also \cite{mochi2}.
We have the decomposition
$ \Omegatilde^1_{\nbigx}
=\xi\Omegatilde^1_{X|\nbigx}
\oplus
 \Omegatilde_{\cnum_{\lambda}}$
into the $X$-direction
and the $\cnum_{\lambda}$-direction.
Let $d_X$ denote the restriction of 
the exterior differential to the $X$-direction.
Similarly,
we have the decomposition
$ \Omegatilde^1_{\nbigx^{\dagger}}
=\xi\Omegatilde^1_{X|\nbigx^{\dagger}}
\oplus
 \Omegatilde_{\cnum_{\mu}}$,
and the restriction of $\DDtilde^{\dagger\,f}_X$
to the $X$-direction is denoted by $d_{X^{\dagger}}$.
The notions of 
$\cnum_{\lambda}$-holomorphic bundles
or $\cnum_{\mu}$-holomorphic bundles
are defined
as in the case of $\proj^1$-holomorphic bundles.

Let $(V_0,d''_{\cnum_{\lambda},V_0})$ 
be a $\cnum_{\lambda}$-holomorphic 
bundle on $\nbigx$.
A $T$-structure \cite{Hertling} of $V_0$
is a differential operator
\[
 \DD^f_{V_0}:V_0\lrarr 
 V_0\otimes
 \xi\Omega^1_{X|\nbigx}
\]
satisfying 
(i) $\DD^f_{V_0}(f\cdot s)
=d_Xf\cdot s+f\cdot \DD^f_{V_0}(s)$
for a function $f$ and a section $s$ of $V$,
(ii) 
 $\bigl(
 d''_{\cnum_{\lambda},V_0}+\DD^f_{V_0}
\bigr)^2=0$.

Let $(V_{\infty},d''_{\cnum_{\mu},V_{\infty}})$
be a $\cnum_{\mu}$-holomorphic vector bundle
on $\nbigx^{\dagger}$.
A $\Ttilde$-structure \cite{Hertling}
is defined to be a differential operator
\[
 \DD^{\dagger\,f}_{V_{\infty}}:
 V_{\infty} \lrarr 
 V_{\infty}\otimes
 \xi\Omega^1_{X|\nbigx^{\dagger}}
\]
satisfying conditions
similar to (i) and (ii) above.

Assume that we are given 
an isomorphism $\Phi$:
\begin{equation}
\label{eq;08.9.10.30}
 \Phi:
 (V_0,d''_{\cnum_{\lambda},V_0},
 \DD^f_{V_0})
 _{|\cnum_{\lambda}^{\ast}\times X}
\simeq
 (V_{\infty},d''_{\cnum_{\mu},V_{\infty}},
 \DD^{\dagger\,f}_{V_{\infty}})_{
 |\cnum_{\mu}^{\ast}\times X^{\dagger}}
\end{equation}
We obtain the $C^{\infty}$-vector bundle $V$
on $\proj^1\times X$
by gluing $V_0$ and $V_{\infty}$ via $\Phi$.
By the condition (\ref{eq;08.9.10.30}),
$d''_{\cnum_{\lambda},V_0}$
and $d''_{\cnum_{\mu},V_{\infty}}$
give $\proj^1$-holomorphic structure 
$d''_{\proj^1,V}$,
and 
$\DD^f_{V_0}$
and $\DD^{\dagger\,f}_{V_{\infty}}$
induce the $\TTtilde$-structure
$\DD^{\sankaku}_V$.
Thus, we obtain a variation of
twistor structures
$(V,d''_{\proj^1,V},\DD^{\sankaku}_V)$.

Conversely,
we naturally obtain such a tuple of
$ (V_0,d''_{\cnum_{\lambda},V_0},
 \DD^f_{V_0})$,
$(V_{\infty},d''_{\cnum_{\mu},V_{\infty}},
 \DD^{\dagger\,f}_{V_{\infty}})$
and 
$\Phi$
from a variation of twistor structures
$(V,d''_{\proj^1,V},\DD^{\sankaku}_V)$
as the restriction to
$\nbigx$ and $\nbigx^{\dagger}$,
respectively.
In this situation,
we set 
\[
 \Glue\bigl(
 (V_0,d''_{\cnum_{\lambda,,V_0},}
 \DD^f_{V_0}),
 (V_{\infty},d''_{\cnum_{\mu},V_{\infty}},
 \DD^{\dagger\,f}_{V_{\infty}}),
 \Phi
 \bigr):=(V,\DDtilde_V^{\sankaku})
\]

\begin{rem}
Let $p_{\lambda}$ be the projection
$\nbigx\lrarr X$.
Under the natural isomorphism
\[
 \xi\Omega^1_{X|\nbigx}
=\lambda^{-1}\cdot
 p_{\lambda}^{-1}\Omega^{1,0}_X
 \oplus
 p_{\lambda}^{-1}\Omega^{0,1}_X
\simeq
 p_{\lambda}^{-1}\Omega^{1,0}_X
 \oplus
 p_{\lambda}^{-1}\Omega^{0,1}_X
=
 p_{\lambda}^{-1}\Omega^1_X,
\]
a $T$-structure $\DD^f_{V_0}$ induces
a holomorphic family of flat $\lambda$-connections
$\DD_{V_0}$.
Similarly,
a $\Ttilde$-structure of 
$\DD^{\dagger\,f}_{V_{\infty}}$
naturally induces 
a holomorphic family of flat $\mu$-connections
$\DD^{\dagger}_{V_{\infty}}$.
Hence, 
a variation of twistor structure
is regarded as the gluing of
families of $\lambda$-flat bundles
and $\mu$-flat bundles.
\hfill\qed
\end{rem}

\subsubsection{Relation with harmonic bundles}
\label{subsubsection;08.7.29.30}

We recall a fundamental equivalence 
due to Hertling and Sabbah.
Let $X$ be a complex manifold.
Let $(\nbige^{\sankaku},
 \DDtilde^{\sankaku},\nbigs)$
be a variation of pure polarized integrable twistor structures
of weight $0$ on $\proj^1\times X$.
By the equivalence between
harmonic bundles and 
variations of pure polarized twistor structures
due to Simpson,
we have the underlying harmonic bundle
$(E,\delbar_E,\theta,h)$ on $X$.
Moreover, it is equipped with 
$C^{\infty}$-sections
$\nbigu$ and $\nbigq$ of $\End(E)$
satisfying the following equations:
\begin{equation}
 \label{eq;08.9.12.2}
 \delbar_E\nbigu=0,\quad
 [\nbigu,\theta]=0,\quad
 \nbigq=\nbigq^{\dagger}
\end{equation}
\begin{equation}
 \label{eq;08.7.29.16}
  \del_E\nbigu-[\theta,\nbigq]+\theta=0,
\quad
 \del_E\nbigq+[\theta,\nbigu^{\dagger}]=0
\end{equation}
Here, $\nbigu_{|Q}$ and $\nbigq_{|Q}$ $(Q\in X)$
are obtained as in (\ref{eq;08.7.29.15}),
and $\nbigu^{\dagger}$ and $\nbigq^{\dagger}$
denote the adjoint of $\nbigu$ and $\nbigq$
with respect to $h$, respectively.
Conversely,
we obtain a variation of 
polarized pure integrable twistor structures
$(\nbige^{\sankaku},\DDtilde^{\sankaku},
 \nbigs)$
from a harmonic bundle
$(E,\delbar_E,\theta,h)$ with
$\nbigu$ and $\nbigq$
satisfying (\ref{eq;08.9.12.2})
and (\ref{eq;08.7.29.16}).
Let $p:\proj^1\times X\lrarr X$
be the projection.
We set
$\nbige^{\sankaku}:=p^{-1}E$
on which we have the natural connection 
$d_{\proj^1}$ along the $\proj^1$-direction.
We set
\[
 \nabla_{\lambda}:=d_{\proj^1}
+\bigl(
 \lambda^{-1}\cdot\nbigu
-\nbigq-\lambda\cdot\nbigu^{\dagger}
 \bigr)\frac{d\lambda}{\lambda}
\]
It gives a flat connection of $\nbige^{\sankaku}$
along the $\proj^1$-direction.
Then, we obtain a $\TTtilde E$-structure
\[
 \DDtilde^{\sankaku}:=
 \bigl( \delbar_E+\lambda\theta^{\dagger}\bigr)
+\bigl(
\del_E+\lambda^{-1}\theta
\bigr)
+\nabla_{\lambda}:
 \nbige^{\sankaku}\lrarr 
 \nbige^{\sankaku}\otimes
 \Omegatilde^1_{\proj^1\times X}.
\]
The pairing $\nbigs$ is induced by
$\nbigs(u\otimes\sigma^{\ast}v)
=h(u,\sigma^{\ast}v)$.

\vspace{.1in}

Let us also see the gluing construction
of the above
$(\nbige^{\sankaku},\DDtilde^{\sankaku},\nbigs)$.
Let $(E,\delbar_E,\theta,h,\nbigu,\nbigq)$ be as above.
Let $p_{\lambda}$ be the projection
$\nbigx\lrarr X$.
Let $\nbige$ be the holomorphic vector bundle
$\bigl(p_{\lambda}^{-1}E,\,\delbar_E
+\lambda\theta^{\dagger}+\delbar_{\lambda}\bigr)$,
where $\delbar_{\lambda}$
denotes the natural $\lambda$-holomorphic structure
of $\nbige$.
We have the family of flat $\lambda$-connections
$ \DD=\delbar_E+\lambda\theta^{\dagger}
+\lambda\del_E+\theta$
of $\nbige$.
The associated family of flat connections
is given by 
$\DD^f=\delbar_E+\lambda\theta^{\dagger}
+\del_E+\lambda^{-1}\theta$.
Then,
$\DDtilde^f:=
 \DD^f+\nabla_{\lambda}$
gives a meromorphic flat connection of $\nbige$.

Let $p_{\mu}$ be the projection
$\nbigx^{\dagger}\lrarr X^{\dagger}$.
Let $\nbige^{\dagger}$ be the holomorphic vector bundle
$\bigl(p_{\mu}^{-1}E,\,\del_E
+\mu\theta+\delbar_{\mu}\bigr)$,
where $\delbar_{\mu}$
denotes the natural $\mu$-holomorphic structure
of $\nbige^{\dagger}$.
We have the family of flat $\mu$-connections
$ \DD^{\dagger}=
\del_E+\mu\theta
+\mu\delbar_E+\theta^{\dagger}$
of $\nbige^{\dagger}$.
The associated family of flat connections
is given by 
$\DD^{\dagger\,f}=
 \del_E+\mu\theta
+\delbar_E+\mu^{-1}\theta^{\dagger}$.
Then,
$\DDtilde^{\dagger\,f}:=
 \DD^{\dagger\,f}+\nabla_{\lambda}$
gives a meromorphic flat connection 
of $\nbige^{\dagger}$.

We have the induced pairings
$ \nbigs_0:
 \nbige\otimes\sigma^{\ast}\nbige^{\dagger}
\lrarr \nbigo_{\nbigx}$
and 
$\nbigs_{\infty}:
 \nbige^{\dagger}\otimes\sigma^{\ast}\nbige
\lrarr \nbigo_{\nbigx^{\dagger}}$
induced by $h$.
Then,
$(\nbige^{\sankaku},
 \DDtilde^{\sankaku},\nbigs)$
is obtained as the gluing of
$(\nbige,\DDtilde^f)$,
$(\nbige^{\dagger},\DDtilde^{\dagger\,f})$
and $(\nbigs_0,\nbigs_{\infty})$
by the procedure in 
Subsection \ref{subsubsection;08.7.29.10}.

\subsubsection{TERP and twistor-TERP}
\label{subsubsection;08.9.11.30}

Let us observe that
the notions of TERP-structure and
twistor-TERP structure
are equivalent.
First, let us introduce a pairing $P$
induced by $\kappa$ and $\nbigs$.
Then, we argue the equivalence
in the case that $X$ is a point,
for simplicity.
We give a remark for the family case
in the end.

\paragraph{The induced pairing $P$}

We set $j:=\gamma\circ\sigma
=\sigma\circ\gamma$,
which is a holomorphic involution of $\proj^1$.
We have the induced isomorphisms
\[
 \sigma^{\ast}\kappa:
 j^{\ast}\Tate(w)\simeq
 \sigma^{\ast}\Tate(w),
\quad
 j^{\ast}\kappa:
 \sigma^{\ast}\Tate(w)\simeq
 j^{\ast}\Tate(w).
\]
We have the following equality:
\[
 \sigma^{\ast}\kappa
\circ
 j^{\ast}\kappa
=j^{\ast}\bigl(
 \gamma^{\ast}\kappa\circ\kappa
 \bigr)
=j^{\ast}(\id)=\id
\]
We will use similar relations implicitly.
We also remark the commutativity
of the following diagram,
which can be checked by a direct calculation:
\[
 \begin{CD}
 j^{\ast}\Tate(w) @>{\gamma^{\ast}\iota_{\Tate(w)}}>>
 \gamma^{\ast}\Tate(w) \\
 @V{\sigma^{\ast}\kappa_{\Tate(w)}}VV 
 @V{\kappa_{\Tate(w)}}VV\\
 \sigma^{\ast}\Tate(w)
 @>{\iota_{\Tate(w)}}>>  \Tate(w)
 \end{CD}
\]
The composite
$j^{\ast}\Tate(w)\lrarr \Tate(w)$
is denoted by $\rho_{\Tate(w)}$.

Let $(V,\DDtilde^{\sankaku}_V,
 \nbigs,\kappa,-w)$ be a variation of
twistor-TERP structure.
We define a pairing
$P:V\otimes j^{\ast}V\lrarr \Tate(-w)$
by 
\begin{equation}
 \label{eq;08.7.21.12}
 P:=(\sqrt{-1})^w\cdot
 \nbigs\circ(1\otimes\sigma^{\ast}\kappa).
\end{equation}

\begin{lem}
\label{lem;08.7.21.10}
$P$ is $(-1)^w$-symmetric
in the sense that the following diagram
is commutative:
\[
\begin{CD}
 j^{\ast}V\otimes V
 @>{j^{\ast}P}>>
 j^{\ast}\Tate(-w) \\
 @V{\exchange}VV @V{\rho_{\Tate(-w)}}VV \\
 V\otimes j^{\ast}V
 @>{(-1)^wP}>>
 \Tate(-w)
\end{CD}
\]
Namely,
$\rho_{\Tate(-w)}\circ j^{\ast}P
=(-1)^w\cdot P\circ\exchange$.
Here,
$\exchange$ denotes the natural morphism
exchanging the components.
\end{lem}
\pf
We have the following equality:
\begin{multline}
\label{eq;08.7.28.2}
\rho_{\Tate(-w)}\circ
 j^{\ast}P
=(\sqrt{-1})^w
 \kappa_{\Tate(-w)}\circ
 \gamma^{\ast}\iota_{\Tate(-w)}\circ
 j^{\ast}\nbigs\circ
 \bigl(1\otimes j^{\ast}\sigma^{\ast}\kappa\bigr)
 \\
=(\sqrt{-1})^w
 \kappa_{\Tate(-w)}\circ
 \gamma^{\ast}\iota_{\Tate(-w)}\circ
 \bigl(\gamma^{\ast}\sigma^{\ast}\nbigs\bigr)
 \circ
 \bigl(1\otimes\gamma^{\ast}\kappa\bigr)
=(\sqrt{-1})^w\kappa_{\Tate(-w)}
 \circ
 \gamma^{\ast}\Bigl(
 \iota_{\Tate(-w)}\circ\sigma^{\ast}\nbigs
 \Bigr) 
 \circ \bigl(1\otimes \gamma^{\ast}\kappa\bigr)
\end{multline}
By using the compatibility
of $\nbigs$ and $\kappa$,
we obtain
\begin{multline}
\label{eq;08.7.28.3}
(-1)^w
 P\circ\exchange
=(\sqrt{-1})^w
 (-1)^w\nbigs\circ(1\otimes\sigma^{\ast}\kappa)
 \circ\exchange
=(\sqrt{-1})^w
 (-1)^w\nbigs\circ(\kappa\otimes\sigma^{\ast}\kappa) 
 \circ(\gamma^{\ast}\kappa\otimes 1)
 \circ\exchange \\
=(\sqrt{-1})^w
 \kappa_{\Tate(-w)}\circ
 \gamma^{\ast}\Bigl(
 (-1)^w\nbigs\circ\exchange
\Bigr)
 \circ \bigl(1\otimes\gamma^{\ast}\kappa\bigr)
\end{multline}
Thus, we are done.
\hfill\qed

\begin{lem}
\label{lem;08.7.21.11}
The following diagram is commutative:
\[
 \begin{CD}
 \gamma^{\ast}V\otimes\sigma^{\ast}V
 @>{\gamma^{\ast}P}>>
 \gamma^{\ast}\Tate(-w)\\
 @V{\kappa\otimes j^{\ast}\kappa}VV 
 @V{\kappa_{\Tate(-w)}}VV \\
 V\otimes j^{\ast}V
 @>{(-1)^wP}>>
 \Tate(-w)
 \end{CD}
\]
Namely,
$(-1)^wP\circ\bigl(\kappa\otimes j^{\ast}\kappa\bigr)
=\kappa_{\Tate(-w)}\circ\gamma^{\ast}P$.
\end{lem}
\pf
We have the following equalities:
\begin{equation}
 \label{eq;08.7.28.4}
(\sqrt{-1})^{-w}P\circ(\kappa\otimes j^{\ast}\kappa)
=\nbigs\circ
 (1\otimes\sigma^{\ast}\kappa)
\circ
 (\kappa\otimes j^{\ast}\kappa)
=\nbigs\circ(\kappa\otimes\sigma^{\ast}\kappa)
 \circ (1\otimes j^{\ast}\kappa)
\end{equation}
\begin{equation}
 \label{eq;08.7.28.5}
\kappa_{\Tate(-w)}\circ\gamma^{\ast}\bigl(
 (\sqrt{-1})^{-w}P\bigr)
=\kappa_{\Tate(-w)}\circ
 \gamma^{\ast}\bigl(\nbigs\circ
 (1\otimes\sigma^{\ast}\kappa)\bigr)
=\kappa_{\Tate(-w)}\circ
 (\gamma^{\ast}\nbigs)\circ
 (1\otimes j^{\ast}\kappa)
\end{equation}
Then, the claim of the lemma follows from
the compatibility of $\nbigs$ and $\kappa$.
\hfill\qed

\paragraph{From twistor-TERP to TERP}

Let $(V,\nabla,\nbigs,\kappa,-w)$
be a twistor-TERP structure.
Let us explain how to associate
a TERP structure
$(H,H'_{\real},\nabla,P',-w)$
in the sense of Hertling.
We set $H:=V_{|\cnum_{\lambda}}$
and $H':=V_{|\cnum_{\lambda}^{\ast}}$.
In general,
for a $\cnum$-vector bundle $U$,
let $\overline{U}$ denote the conjugate 
of $U$,
i.e., $\overline{U}=U$ as an $\real$-vector bundle,
and the multiplication of $\sqrt{-1}$
on $\overline{U}$ is given by
the multiplication of $-\sqrt{-1}$ on $U$.
Note that
$\gamma^{\ast}(H)_{|\lambda}$
for $\lambda\neq 0$
is naturally identified with
$\Hbar_{|\lambdabar^{-1}}$.

The following diagram for $\lambda\neq 0$
is commutative by the flatness of $\kappa$:
\begin{equation}
 \label{eq;08.9.10.1}
 \begin{CD}
 \Hbar_{|\lambda}
@>{\Pi_{\lambda}}>>
 \Hbar_{|\lambdabar^{-1}}\\
@V{\kappa_{|\lambdabar^{-1}}}VV
@VV{\kappa_{|\lambda}}V\\
 H_{|\lambdabar^{-1}}
@>{\Pi_{\lambda}^{-1}}>>
 H_{|\lambda}  
 \end{CD}
\end{equation}
Here, $\Pi_{\lambda}$ denotes
the parallel transform
along the segment connecting
$\lambda$ and $\lambdabar^{-1}$,
as often used in \cite{Hertling}.
A flat isomorphism
$\kappa':
 \overline{H}_{|\cnum^{\ast}_{\lambda}}\simeq
 H_{|\cnum^{\ast}_{\lambda}}$ is given by
the composite of the morphisms,
i.e.,
$\kappa'_{|\lambda}
:=\kappa_{|\lambda}\circ \Pi_{\lambda}$.
Because $\gamma^{\ast}\kappa\circ\kappa=\id$,
the composite
\[
\begin{CD}
 \Hbar_{|\lambda}
 @>{\kappa_{|\lambdabar^{-1}}}>>
 H_{|\lambdabar^{-1}}
 @>{\kappa_{\lambda}}>>
 \Hbar_{|\lambda}
\end{CD}
\]
is the identity.
Let us check 
$\kappa'\circ\kappa'=\id$
by using the commutativity 
(\ref{eq;08.9.10.1}):
\[
 \kappa'_{\lambda}\circ
 \kappa'_{\lambda}
=\bigl(
 \kappa_{\lambda}\circ\Pi_{\lambda}
 \bigr)
\circ
 \bigl(
 \Pi_{\lambda}^{-1}\circ
 \kappa_{|\lambdabar^{-1}}
 \bigr)
=\kappa_{\lambda}\circ
 \kappa_{|\lambdabar^{-1}}=\id
\]
Hence, 
$\kappa'$ gives a flat real structure of $H'$.
Thus, we obtain
a real flat subbundle 
$H'_{\real}$ of $H_{|\cnum^{\ast}_{\lambda}}$.
By restricting $P$,
we obtain a pairing:
\[
P_{|\cnum_{\lambda}}:
 H\otimes j^{\ast}H\lrarr
 \Tate(-w)_{|\cnum_{\lambda}}
=\nbigo_{\cnum_{\lambda}}\cdot
 (\sqrt{-1}\lambda)^{-w}
 t_1^{(-w)}
\]
By taking the coefficients of $t_1^{(-w)}$,
we obtain a flat morphism
\[
P':
 H'\otimes j^{\ast}H'
 \lrarr
 \nbigo_{\cnum_{\lambda}^{\ast}}
\]
such that
$\lambda^w\cdot P'$
induces a perfect pairing
$H\otimes j^{\ast}H\lrarr \nbigo_{\cnum_{\lambda}}$.
By Lemma \ref{lem;08.7.21.10},
$P'$ is $(-1)^w$-symmetric.

\begin{lem}
$P'\bigl(
 H'_{\real}\otimes_{\real}
 j^{\ast}H'_{\real}
 \bigr)
\subset 
 (\sqrt{-1})^w\real$.
\end{lem}
\pf
Note that $\kappa$ gives real structures
$\kappa_{|a}:
 \overline{H_{|a}}\simeq H_{|a}$ for $a=1,-1$.
By Lemma \ref{lem;08.7.21.11},
we have
\begin{equation}
 \label{eq;08.7.28.1}
(\sqrt{-1})^w\cdot
 P_{|1}\circ
 \bigl(\kappa_{|1}\otimes \kappa_{|-1}\bigr)
=
 (\kappa_{\Tate(-w)})_{|1}
\circ
 \bigl((\sqrt{-1})^wP_{|1}\bigr).
\end{equation}
We obtain 
$ P'_{|1}\bigl(
 H_{|1}\otimes
 H_{|-1}
 \bigr)
\subset
 (\sqrt{-1})^w\real$.
Then, the claim of the lemma follows from
the flatness of $P'$.
\hfill\qed

\vspace{.1in}
Thus, we obtain a TERP-structure
$(H,H'_{\real},\nabla,P',-w)$.

\paragraph{From TERP to twistor-TERP}

Conversely,
we obtain a twistor-TERP structure
$(V,\nabla,\kappa,\nbigs,-w)$ 
from a TERP structure
$(H,H'_{\real},\nabla,P',-w)$.
We set $V_0:=H$
and $V_{\infty}:=\gamma^{\ast}H$.
We have the flat isomorphism
\[
 \tau_{\rm real}:H_{|\cnum^{\ast}_{\lambda}}
\simeq
 \gamma^{\ast}(H_{|\cnum^{\ast}_{\lambda}}),
\]
obtained as the composite of the conjugate
with respect to the real structure
and the parallel transform
along the segment connecting
$\lambda$ and $\lambdabar^{-1}$.
By gluing $(H,\nabla)$ and 
$\gamma^{\ast}(H,\nabla)$ via 
$\tau_{\rm real}$,
we obtain an integrable twistor structure
$(V,\nabla)$.

By construction,
we have
$\gamma^{\ast}(\tau_{\rm real})
=\tau_{\rm real}^{-1}$,
and the following diagram is commutative:
\[
 \begin{CD}
 H_{|\cnum_{\lambda}^{\ast}}
 @>{\tau_{\rm real}}>>
 \gamma^{\ast}(H_{|\cnum_{\lambda}^{\ast}})
 \\ 
 @V{=}VV @V{=}VV \\
 \gamma^{\ast}
\bigl(\gamma^{\ast}
 H_{|\cnum_{\lambda}^{\ast}}\bigr)
 @>{\gamma^{\ast}\tau_{\rm real}^{-1}}>>
 \gamma^{\ast}H_{|\cnum_{\lambda}^{\ast}}
 \end{CD}
\]
Hence, 
a morphism
$\kappa:\gamma^{\ast}(V,\nabla)
 \simeq (V,\nabla)$
is given by the gluing of
$\gamma^{\ast}V_{\infty}\simeq V_0$
and $\gamma^{\ast}V_0\simeq V_{\infty}$
induced by the identity.
Clearly it satisfies
$\gamma^{\ast}\kappa\circ\kappa=\id$.
The restriction
$\kappa_{|\cnum_{\lambda}^{\ast}}:
 \gamma^{\ast}(V)_{|\cnum_{\lambda}^{\ast}}
\lrarr
 V_{|\cnum_{\lambda}^{\ast}}$
is identified with
$\tau_{\rm real}^{-1}:
 \gamma^{\ast}H_{|\cnum_{\lambda}^{\ast}}
\simeq
 H_{|\cnum_{\lambda}^{\ast}}$.

Let $P_0:V_0\otimes j^{\ast}V_0
\lrarr 
 \nbigo_{\cnum_{\lambda}}\cdot t_0^{(-w)}$
be given by
\[
 P_0=P'\cdot t_1^{(-w)}
=P'\cdot \bigl(\sqrt{-1}\lambda\bigr)^w
 \cdot t_0^{(-w)}.
\]
We have the induced morphism
\[
\kappa_{\Tate(-w)}\circ \gamma^{\ast}P_0:
 V_{\infty}\otimes j^{\ast}V_{\infty}
\lrarr 
 \nbigo_{\cnum_{\mu}}\cdot t_{\infty}^{(-w)}.
\]
We obtain the following equalities for linear maps
$\overline{H_{|1}}
 \otimes
 \overline{H_{|-1}}
\lrarr \Tate(-w)_{|1}$
from
$P'\bigl(H'_{\real}\otimes_{\real}
 j^{\ast}H'_{\real}\bigr)
\subset
 (\sqrt{-1})^w\real$:
\[
(\sqrt{-1})^w\cdot
 P_{0|1}\circ
 \bigl(\kappa_{|1}
 \otimes \kappa_{|-1}
 \bigr)
= (\kappa_{\Tate(-w)})_{|1}
\circ
 \bigl((\sqrt{-1})^wP_{0|1}\bigr)
=(-\sqrt{-1})^w
(\kappa_{\Tate(-w)})_{|1}
\circ
 \bigl(\gamma^{\ast}P_{0|1}\bigr)
\]
Here, we have used the natural identification
$P_{0|1}=(\gamma^{\ast}P_0)_{|1}$.
The first and third terms are 
obtained as the restrictions of morphisms
$(V_{\infty} \otimes  j^{\ast} V_{\infty})
 _{|\cnum_{\lambda}^{\ast}}
\lrarr
 \nbigo_{\cnum_{\lambda}^{\ast}}\cdot t_{1}^{(-w)}$
to the fiber over $1$.
By flatness,
we obtain the following equality 
on $\cnum_{\lambda}^{\ast}$:
\begin{equation}
\label{eq;08.9.9.100}
(-1)^w\cdot
 P_0\circ\bigl(
 \kappa \otimes  j^{\ast}\kappa
 \bigr)
=
 \kappa_{\Tate(-w)}\circ\gamma^{\ast}P_0
\end{equation}
Hence, 
the pairings $P_0$ and 
$(-1)^w\kappa_{\Tate(-w)}
 \circ\gamma^{\ast}P_0$ 
induce 
$P:V\otimes j^{\ast}V\lrarr \Tate(-w)$.
Since $P'$ is $(-1)^w$-symmetric,
$P$ is also $(-1)^w$-symmetric
in the sense of Lemma \ref{lem;08.7.21.10}.
From (\ref{eq;08.9.9.100}),
we obtain 
\begin{equation}
  \label{eq;08.7.28.6}
 (-1)^w\cdot P\circ
 \bigl(\kappa\otimes j^{\ast}\kappa\bigr)
=\kappa_{\Tate(-w)}\circ \gamma^{\ast}P.
\end{equation}

The pairing $\nbigs$
is constructed from $P$ and $\kappa$
by the relation (\ref{eq;08.7.21.12}).
The compatibility of $\kappa$ and $\nbigs$
follows from (\ref{eq;08.7.28.4}),
(\ref{eq;08.7.28.5})
and (\ref{eq;08.7.28.6}).
The pairing $\nbigs$ is $(-1)^w$-symmetric,
which follows from 
(\ref{eq;08.7.28.2}), (\ref{eq;08.7.28.3})
and the compatibility with $\kappa$.
Thus, we obtain 
a twistor-TERP structure
$(V,\nabla,\nbigs,\kappa,-w)$.

\paragraph{Hertling's vector bundle}

Let $(H,H_{\real}',\nabla,P,-w)$
be a TERP-structure
corresponding to 
a twistor-TERP structure
$(V,\nabla,\nbigs,\kappa,-w)$.
Recall that 
Hertling constructed 
an integrable twistor structure
$(\Hhat,\nabla)$
from a TERP-structure
$(H,H_{\real}',\nabla,P,-w)$
by gluing $H$ and $\gamma^{\ast}H$
via a map $\tau$.
(See \cite{Hertling}.)
We do not recall $\tau$
and his construction here,
but $\Hhat$ is naturally isomorphic to
$V\otimes \nbigo(0,-w)$
by the following correspondence:
\[
 H=V_0\longleftrightarrow
 V_0\otimes\nbigo(0,-w)_0,
\quad
 a\longleftrightarrow
 a\otimes f_0^{(0,-w)}
\]
\[
 \gamma^{\ast}H
 \longleftrightarrow
 \gamma^{\ast}V_0\otimes
 \nbigo(0,-w)_{\infty},
\quad
 \gamma^{\ast}b
 \longleftrightarrow
 \gamma^{\ast}b\otimes
 (\sqrt{-1})^w f_{\infty}^{(0,-w)}
\]
According to \cite{Hertling} and
\cite{Hertling-Sevenheck},
$(H,H'_{\real},\nabla,P,-w)$
is defined to be pure 
if $(\Hhat,\nabla)$ is pure of weight $0$.
They consider
the hermitian pairing $h$ of 
$H^0(\proj^1,\Hhat)$
given by 
 $\lambda^w\cdot P'\circ
 (1\otimes\tau)$,
and $(H,H'_{\real},\nabla,P,-w)$
is defined to be polarized
if $h$ is positive definite.

\begin{lem}
$(H,H'_{\real},\nabla,P,-w)$ is pure (polarized),
if and only if 
$(V,\nabla,\nbigs,\kappa,-w)$ is pure (polarized).
\end{lem}
\pf
The claim for purity is obvious.
Let us consider polarizability.
We have only to show that
$h$ is the hermitian pairing induced by
$\nbigstilde:=\nbigs\otimes \nbigs_{0,-w}$,
under the identification of
$\Hhat$ and $V\otimes \nbigo(0,-w)$.

Let $\ahat,\bhat\in H^0(\proj^1,\Hhat)$.
Under the identification
$\Hhat_{|\cnum_{\lambda}}=H$
and $\Hhat_{|\cnum_{\mu}}=\gamma^{\ast}H$,
the sections $a$ and $b$ of $H$ are determined by
$ a:=\ahat_{|\cnum_{\lambda}}$
and
$\gamma^{\ast}b:=
 \bhat_{|\cnum_{\mu}}$.
By definition, we have 
\[
 h(\ahat,\bhat)
=\lambda^wP'(a,j^{\ast}b)
\]
Let us look at
$\nbigstilde_{|\cnum_{\lambda}}$.
Under the above identification,
the pairing of $\ahat$ and $\bhat$
is given by
\[
 \nbigstilde\Bigl(a\otimes f_0^{(0,-w)},
 \sigma^{\ast}\bigl(
 \gamma^{\ast}b\otimes
 (\sqrt{-1})^wf_{\infty}^{(0,-w)}
 \bigr)
\Bigr)
=\nbigs\bigl(a,\sigma^{\ast}(\gamma^{\ast}b)\bigr)
 \cdot t_0^{(w)}
=\nbigs\bigl(a,j^{\ast}b\bigr)
 \cdot t_0^{(w)}
=:\nbigs_0\bigl(a,j^{\ast}b\bigr)
\]

Let us compare $\lambda^wP'(a,j^{\ast}b)$
and $\nbigs_0\bigl(a,j^{\ast}b\bigr)$.
Since $\kappa_{|\cnum_{\mu}}$
is the same as the identity
$ V_{\infty}=
 \gamma^{\ast}H
\lrarr 
\gamma^{\ast}V_0
=\gamma^{\ast}H$,
we have
\[
 P_{|\cnum_{\lambda}}
 =(\sqrt{-1})^w\nbigs
 \circ(1\otimes\sigma^{\ast}\kappa)_{|\cnum_{\lambda}}
=(\sqrt{-1})^w\nbigs_{|\cnum_{\lambda}}
\]
Hence, we have the following equality:
\[
  P'(a,j^{\ast}b)\cdot t_1^{(-w)}
=P\bigl(a,j^{\ast}b\bigr)
=\nbigs\bigl(
 a,j^{\ast}b
 \bigr)\cdot (\sqrt{-1})^w
=(\sqrt{-1})^w
\cdot
 \nbigs_0\bigl(a,j^{\ast}b\bigr)
 \cdot t_0^{(-w)}
=\lambda^{-w}
 \nbigs_0\bigl(a,j^{\ast}b\bigr)
 \cdot t_1^{(-w)}
\]
Thus, we obtain
$ \lambda^w\cdot  P'(a,j^{\ast}b)
=\nbigs_0(a,j^{\ast}b)$.
Therefore, $\nbigstilde$ induces $h$.
\hfill\qed

\paragraph{Family version}

The correspondence
is generalized in the family case.
We set $H:=V_{|\cnum_{\lambda}\times X}$.
It is equipped with $TE$-structure $\nabla$
obtained as the restriction of
$\DDtilde^{\sankaku}_V$.
As in the previous case,
we obtain a flat $\cnum$-anti-linear isomorphism
$\kappa':
 H_{|\cnum_{\lambda}^{\ast}\times X}\simeq 
 H_{|\cnum_{\lambda}^{\ast}\times X}$
and a flat pairing 
$P:H'\otimes j^{\ast}H'\lrarr
 \nbigo_{\cnum^{\ast}\times X}$.
It is easy to check that
$(H,H'_{\real},\nabla,P,-w)$
is a variation of TERP structures.
The converse can be constructed 
similarly.
The correspondence preserves
``pure''  and ``polarized'',
for which we have only to check
the case in which $X$  is a point.

\subsection{Basic examples}
\label{subsection;08.9.10.20}

\subsubsection{Example associated to
 a holomorphic function}
\label{subsubsection;08.8.3.20}

Let $\gminia$ be a holomorphic function on 
a complex manifold $X$.
We set 
\[
 V_0:=\nbigo_{\cnum_{\lambda}\times X}\cdot e,
\quad
 \nabla_{V_0}(e)=e\cdot
 d\bigl(\lambda^{-1}\cdot \gminia\bigr),
\]
\[
 V_{\infty}:=
 \nbigo_{\cnum_{\mu}\times X^{\dagger}}
 \cdot e^{\dagger},
\quad
 \nabla_{V_{\infty}}(e^{\dagger})=
 e^{\dagger}\cdot
d\bigl(\mu^{-1}\cdot \gminiabar\bigr).
\]
We put
$s:=\exp\bigl(-\lambda^{-1}\gminia\bigr)\cdot e$
and 
$s^{\dagger}:=
 \exp\bigl(-\mu^{-1}\gminiabar\bigr)\cdot e^{\dagger}$,
which are flat sections of
$V_{0|\cnum_{\lambda}^{\ast}\times X}$
and $V_{\infty|\cnum_{\mu}^{\ast}\times X^{\dagger}}$,
respectively.
A gluing 
$\Phi:V_{0|\cnum_{\lambda}^{\ast}\times X}
\simeq V_{\infty|\cnum_{\mu}^{\ast}\times X^{\dagger}}$
is given by $\Phi(s)=s^{\dagger}$,
in other words,
\[
 \Phi(e)
=\exp\bigl(
 \lambda^{-1}\gminia
-\mu^{-1}\gminiabar
 \bigr)\cdot e^{\dagger}.
\]
Let $V$ be the $C^{\infty}$-bundle
obtained as the gluing of $V_0$ and $V_{\infty}$
via $\Phi$,
which is equipped with $\TTtilde E$-structure.
For each point $P\in X$,
the restriction $V_{|\proj^1\times\{P\}}$
is isomorphic to $\nbigo_{\proj^1}$,
and hence $(V,\DDtilde_V^{\sankaku})$
is pure of weight $0$.
A real structure $\kappa$
is given by
$\kappa(\gamma^{\ast}e^{\dagger})=e$
and $\kappa(\gamma^{\ast}e)=e^{\dagger}$.
We can check that $\kappa$ actually gives
a flat isomorphism $\gamma^{\ast}V\simeq V$.
A pairing $\nbigs$ of $V$
with weight $0$ is given by
$e\otimes\sigma^{\ast}e^{\dagger}
\longmapsto t_0^{(0)}$
and 
$e^{\dagger}\otimes\sigma^{\ast}e
\longmapsto t_{\infty}^{(0)}$.
It is easy to check that
$\nbigs$ actually gives a symmetric flat pairing
$V\otimes\sigma^{\ast}V\lrarr\Tate(0)_X$.
The compatibility of 
$\nbigs$ and $\kappa$ can be checked
by a direct calculation:
\[
 \kappa_{\Tate(0)}\circ
 \gamma^{\ast}
 \nbigs(\gamma^{\ast}e^{\dagger}
 \otimes \gamma^{\ast}\sigma^{\ast}e)
=
\kappa_{\Tate(0)}\bigl(
 \gamma^{\ast}\bigl(\nbigs(e^{\dagger}
 \otimes\sigma^{\ast}e)
 \bigr) \bigr)
=\kappa_{\Tate(0)}\gamma^{\ast}t_{\infty}^{(0)}
=t_{0}^{(0)}
\]
\[
  \nbigs\circ(\kappa\otimes\sigma^{\ast}\kappa)
 (\gamma^{\ast}e^{\dagger}
 \otimes\sigma^{\ast}\gamma^{\ast}e)
=\nbigs\bigl(
 e\otimes \sigma^{\ast}e^{\dagger}
 \bigr)
=t_0^{(0)}
\]
Hence, we obtain a variation of
twistor-TERP structure
denoted by $L(\gminia)$.
It is polarized.
The underlying harmonic bundle is
given by
the line bundle $\nbigo_X\cdot v$
with the Higgs field
$\theta\cdot v=v\cdot d\gminia$
and the hermitian metric
$h(v,v)=1$,
where $v:=e_{|\{0\}\times X}$.
The operators $\nbigu$ and $\nbigq$
are 
$\nbigu=-\gminia$
and $\nbigq=0$.

\subsubsection{Example associated to
 unitary flat bundles of rank one}
\label{subsubsection;08.8.10.30}

In general,
a variation of pure polarized Hodge structures
provides us with an example of
variation of pure polarized 
integrable twistor structures.
Any unitary flat bundle naturally gives
a variation of pure polarized Hodge structures,
and hence an integrable variation of pure polarized 
integrable twistor structure.

In particular, we will use the following example.
Let $X:=\cnum^n$ and
$D:=\bigcup_{i=1}^{\ell}\bigl\{z_i=0\bigr\}$.
For any $\veca\in\real^{\ell}$,
we have the unitary flat bundle
\[
 \nbigo_{X-D}\cdot e,
\quad
 \nabla e=e\cdot
 \Bigl(
 -\sum_{i=1}^{\ell}
 a_i\cdot \frac{dz_i}{z_i}
 \Bigr)
\]
The associated variation of integrable
polarized pure integrable structures
is denoted by $L(\veca)$.

More specifically,
it is obtained as the gluing of the following 
meromorphic flat bundles:
\[
 V_0=\nbigo_{\cnum_{\lambda}\times (X-D)}\cdot e,
\quad
 \nabla_{V_0}e
=e\cdot \Bigl(
 -\sum_{i=1}^{\ell}
 a_i\cdot \frac{dz_i}{z_i}
 \Bigr)
\]
\[
  V_{\infty}= 
\nbigo_{\cnum_{\mu}\times (X^{\dagger}-D^{\dagger})}
 \cdot e^{\dagger},
\quad
 \nabla_{V_{\infty}}e^{\dagger}
=e^{\dagger}\cdot \Bigl(
 \sum_{i=1}^{\ell}
 a_i\cdot \frac{d\zbar_i}{\zbar_i}
 \Bigr)
\]
The gluing is given by
$\Phi(e)=
 \prod_{i=1}^{\ell}|z_i|^{-2a_i}\cdot e^{\dagger}$.
The pairing is given by
$\nbigs(e,\sigma^{\ast}e^{\dagger})=1$.
The underlying harmonic bundle is 
the line bundle
$\nbigo_{X-D}\cdot v$
with $\theta\cdot v=0$
and $h(v,v)=\prod_{i=1}^{\ell}|z_i|^{-2a_i}$,
where $v=e_{|\{0\}\times (X-D)}$.
The operators $\nbigu$ and $\nbigq$ are $0$.

\subsubsection{Example induced by nilpotent maps}
\label{subsubsection;08.7.26.1}

Let $Y$ be a complex manifold.
We set $X:=\cnum^{\ell}\times Y$,
$D=\bigcup_{i=1}^{\ell}\{z_i=0\}\times Y$.
We put $\nbigx:=\cnum_{\lambda}\times X$
and $\nbigx^{\dagger}:=
 \cnum_{\mu}\times X^{\dagger}$.
We use the symbols
$\nbigd$, $\nbigy$,
$\nbigd^{\dagger}$ and $\nbigy^{\dagger}$
in similar meanings.
Let $q_0:\nbigx\lrarr \nbigy$
and 
$q_{\infty}:
 \nbigx^{\dagger}\lrarr \nbigy^{\dagger}$
denote the naturally defined projections.

Let $(V,\DD^{\sankaku})$ be a variation of
$\proj^1$-holomorphic vector bundles
on $\proj^1\times Y$
with a tuple $\vecf$ of nilpotent morphisms
\[
 f_i:V\lrarr V\otimes \Tate(-1),
 \quad
 i=1,\ldots,\ell
\]
such that
(i) $[f_i,f_j]=0$,
(ii) they are $\proj^1$-holomorphic
and $\DD^{\sankaku}$-flat.
We recall a construction of
the variation of $\proj^1$-holomorphic vector bundles
on $\proj^1\times (X-D)$
associated to $(V,\vecf)$
given in Subsection 3.5.3 of \cite{mochi2}
with a minor generalization.
(We considered the case in which
$Y$ is a point in \cite{mochi2}.)

We regard 
$(V,\DD^{\sankaku}_{V})$
as the gluing of
a family of $\lambda$-flat bundles
$(V_0,\DD_{V_0})$ on $\nbigy$,
and a family of $\mu$-flat bundles
$(V_{\infty},\DD^{\dagger}_{V_{\infty}})$
on $\nbigy^{\dagger}$.
We obtain a holomorphic vector bundle
$\nbigv_0:=q_0^{\ast}V_0$ on $\nbigx-\nbigd$
with a family of flat $\lambda$-connections
$q_0^{\ast}\DD_{V_0}$.
We naturally identify
$\Tate(0)_{|\nbigx-\nbigd}
\simeq\nbigo_{\nbigx-\nbigd}$
by the trivialization $t_0^{(0)}$.
We also use the natural identification
$\Tate(-1)\otimes\Tate(1)\simeq\Tate(0)$.
We have the $q_0^{\ast}\DD_{V_0}$-flat
endomorphisms
$q_0^{\ast}f_i\otimes t_0^{(1)}
 \in \End(\nbigv_0)$.
We obtain the family of flat $\lambda$-connections
on $\nbigv_0$ given as follows:
\[
 \DD_{\nbigv_0}:=
 q_0^{\ast}\DD_{V_0}
+\sum_{i=1}^{\ell}
 q_0^{\ast}
 f_i\otimes t_0^{(1)}
 \frac{dz_i}{z_i}
\]
Similarly, we obtain a holomorphic vector bundle
$\nbigv_{\infty}:=
 q_{\infty}^{\ast}V_{\infty}$
on $\nbigx^{\dagger}-\nbigd^{\dagger}$
with a family of flat $\mu$-connections
$q_{\infty}^{\ast}\DD^{\dagger}_{V_{\infty}}$.
We have the 
$q_{\infty}^{\ast}
 \DD^{\dagger}_{V_{\infty}}$-flat
endomorphisms
$q_{\infty}^{\ast}
 f_i\otimes t_{\infty}^{(1)}
\in \End(\nbigv_{\infty})$.
Hence, we obtain the following 
family of flat $\mu$-connections:
\[
 \DD^{\dagger}_{\nbigv_{\infty}}:=
 q_{\infty}^{\ast}\DD^{\dagger}_{V_{\infty}}
+\sum_{i=1}^{\ell}
 q_{\infty}^{\ast}
 f_i\otimes t_{\infty}^{(1)}
 \frac{d\zbar_i}{\zbar_i}
\]

Let $\Psi_{V}:
 V_{0|\cnum_{\lambda}^{\ast}\times Y}
\simeq
 V_{\infty|\cnum_{\mu}^{\ast}\times Y}$
denote the gluing.
An isomorphism
$\Psi:\nbigv_{0|\cnum_{\lambda}^{\ast}\times(X-D)}
 \lrarr
 \nbigv_{\infty|
 \cnum_{\mu}^{\ast}\times (X^{\dagger}-D^{\dagger})}$
is given as follows:
\begin{equation}
 \label{eq;08.10.1.1}
\Psi:=
 \Psi_V\circ
 \exp\Bigl(
 \sum_{i=1}^{\ell}
 \log|z_i|^2\cdot 
 q_0^{\ast}
 f_i\otimes \sqrt{-1}t_1^{(1)}
 \Bigr)
\end{equation}
By construction,
$\Psi$ is holomorphic with respect to $\lambda$.
\begin{lem}
\label{lem;08.8.11.1}
$\Psi\circ\DD^{f}_{\nbigv_0}
=\DD^{\dagger\,f}_{\nbigv_{\infty}}\circ\Psi$.
\end{lem}
\pf
We have the following expressions:
\begin{equation}
 \label{eq;08.8.18.15}
  \DD_{\nbigv_0}^f=
 q_0^{\ast}\DD_{V_0}^f
+\sum_{i=1}^{\ell}
 q_0^{\ast}
 f_i\otimes (\sqrt{-1}t_1^{(1)})
 \frac{dz_i}{z_i},
\quad
 \DD^{\dagger\,f}_{\nbigv_{\infty}}=
 q_{\infty}^{\ast}\DD^{\dagger\,f}_{V_{\infty}}
+\sum_{i=1}^{\ell}
 q_{\infty}^{\ast}
 f_i\otimes (-\sqrt{-1}t_{1}^{(1)})
 \frac{d\zbar_i}{\zbar_i}
\end{equation}
Because
$\Psi_V\circ\DD^f_{V_0}
=\DD^{\dagger\,f}_{V_{\infty}}
 \circ \Psi_V$,
we have the following:
\begin{multline*}
q_{\infty}^{\ast}\DD^{\dagger\,f}_{V_{\infty}}
 \circ\Psi
-\Psi\circ q_0^{\ast}\DD^f_{V_0}
=\Psi_V\circ
 q_0^{\ast}\DD^f_{V_0}
\Bigl(
\exp\Bigl(
\sum_{i=1}^{\ell}
 \log|z_i|^2\cdot
 q_0^{\ast}f_i\otimes\sqrt{-1}t_1^{(1)}
 \Bigr)
\Bigr)\\
=\Psi\circ
 \Bigl(
 \sum_{i=1}^{\ell}
 \Bigl(\frac{dz_i}{z_i}
 +\frac{d\zbar_i}{\zbar_i}\Bigr)
 \cdot
 q_0^{\ast}f_i\otimes \sqrt{-1}t_1^{(1)}
 \Bigr)
\end{multline*}
Then, the claim of the lemma follows.
\hfill\qed

\vspace{.1in}

Let $\TNIL(V,\DD^{\sankaku}_V,\vecf)$
denote 
the variation of $\proj^1$-holomorphic bundles
on $(X-D)\times \proj^1$
obtained as the gluing of
$(\nbigv_0,\DD_{\nbigv_0})$ and
$(\nbigv_{\infty},
 \DD^{\dagger}_{\nbigv_{\infty}})$
via $\Psi$.

Assume that
 $(V,\DD^{\sankaku}_V)$ is equipped with
a $(-1)^w$-symmetric pairing
$\nbigs:(V,\DD^{\sankaku}_V)
\otimes
 \sigma^{\ast}(V,\DD^{\sankaku}_V)
\lrarr \Tate(-w)$
such that 
$\nbigs(f_i\otimes \id)
+\nbigs(\id\otimes \sigma^{\ast}f_i)=0$
for any $i$.
Then, we have the induced $(-1)^w$-symmetric
pairing
\[
 \TNIL(\nbigs):
 \TNIL(V,\DD^{\sankaku}_V,\vecf)
\otimes
 \sigma^{\ast}
 \TNIL(V,\DD^{\sankaku}_V,\vecf)
\lrarr \Tate(-w).
\]
It is obtained as the gluing of the pairings
\[
 \nbigs_0:
 \nbigv_0\otimes\sigma^{\ast}\nbigv_{\infty}
\lrarr \Tate(-w)_{|\nbigx-\nbigd},
\quad
  \nbigs_{\infty}:
 \nbigv_{\infty}\otimes\sigma^{\ast}\nbigv_{0}
\lrarr 
\Tate(-w)_{|\nbigx^{\dagger}-\nbigd^{\dagger}},
\]
which are the pull backs of
$V_0\otimes\sigma^{\ast}V_{\infty}
\lrarr \Tate(-w)_{|\cnum_{\lambda}}$
and $V_{\infty}\otimes
 \sigma^{\ast}V_0\lrarr\Tate(-w)_{|\cnum_{\mu}}$.
(See Subsection 3.6.1 of \cite{mochi2}.)

\paragraph{Enrichment}

Assume that 
$(V,\DD^{\sankaku}_V)$ is enriched
to a variation of integrable twistor structures
$(V,\DDtilde^{\sankaku}_V)$
such that $f_j$ are $\DDtilde^{\sankaku}_V$-flat,
which is obtained as the gluing of
$(V_0,\nabla_{V_0})$
and $(V_{\infty},\nabla_{V_{\infty}})$
via $\Psi_V$.
Then $\TNIL(V,\DD^{\sankaku}_V,\vecf)$
is also enriched to the variation of
integrable twistor structures
$\TNIL(V,\DDtilde^{\sankaku}_V,\vecf)$,
which can be checked by an obvious enrichment
of the argument in the proof of 
Lemma \ref{lem;08.8.11.1}.
The $TE$-structure $\nabla_{\nbigv_0}$
and the $\Ttilde E$-structure
$\nabla_{\nbigv_{\infty}}$
are given by essentially the same formula
as (\ref{eq;08.8.18.15}):
\[
  \nabla_{\nbigv_0}=
 q_0^{\ast}\nabla_{V_0}
+\sum_{i=1}^{\ell}
 q_0^{\ast}
 f_i\otimes (\sqrt{-1}t_1^{(1)})
 \frac{dz_i}{z_i},
\quad
 \nabla_{\nbigv_{\infty}}=
 q_{\infty}^{\ast}\nabla_{V_{\infty}}
+\sum_{i=1}^{\ell}
 q_{\infty}^{\ast}
 f_i\otimes (-\sqrt{-1}t_{1}^{(1)})
 \frac{d\zbar_i}{\zbar_i}
\]
If we are given a pairing $\nbigs$
of $(V,\DDtilde_V^{\sankaku})$ with weight $w$
such that 
$\nbigs\circ(f_j\otimes\id)
+\nbigs\circ(\id\otimes \sigma^{\ast}f_j)=0$,
we have a naturally induced pairing
$\TNIL(\nbigs)$ of 
$\TNIL(V,\DDtilde_V^{\sankaku},\vecf)$ with weight $w$.
Assume that we are given a real structure $\kappa$ 
of $(V,\DDtilde^{\sankaku}_V,\nbigs)$
such that 
$\kappa\circ \gamma^{\ast}f_i=f_i\circ\kappa$.
Because
$\kappa_0\circ\gamma^{\ast}(f_i\otimes t_1^{(1)})
=(f_i\otimes t_1^{(1)})\circ \kappa_0$,
we obtain isomorphisms:
\[
 \kappa_0:
 \gamma^{\ast}
 (\nbigv_{\infty},\nabla_{\nbigv_{\infty}})
\simeq 
 (\nbigv_0,\nabla_{\nbigv_0}),
\quad
 \kappa_{\infty}:
  \gamma^{\ast}(\nbigv_{0},\nabla_{\nbigv_0})
\simeq 
 (\nbigv_{\infty},\nabla_{\nbigv_0})
\]
The following diagram
on $(X-D)\times\cnum_{\lambda}^{\ast}$
is commutative:
\[
 \begin{CD}
 \gamma^{\ast}
 \nbigv_{\infty}
 @>{\kappa_0}>>
 \nbigv_0\\
 @V{\gamma^{\ast}\Psi^{-1}}VV
 @V{\Psi}VV \\
 \gamma^{\ast}\nbigv_0
 @>{\kappa_{\infty}}>> \nbigv_{\infty}
 \end{CD}
\]
To see it, we have only to remark
\begin{multline}
 \Psi\circ\kappa
=\Psi_V\circ\exp\Bigl(
 \sum_{i=1}^n\log|z_i(P)|^2\cdot
 f_i\otimes \sqrt{-1}t_1^{(1)}
 \Bigr)\circ\kappa \\
=\kappa\circ\gamma^{\ast}\Psi^{-1}_{V}\circ
 \exp\Bigl(
 -\sum_{i=1}^n\log|z_i(P)|^2
 \cdot \gamma^{\ast}\bigl(
 f_i\otimes\sqrt{-1}t_1^{(1)}
 \bigr)
 \Bigr)
=\kappa\circ\gamma^{\ast}\Psi^{-1}
\end{multline}
Hence, we obtain the isomorphism
$\TNIL(\kappa):
 \gamma^{\ast}\TNIL(V,\DDtilde^{\sankaku}_V,\vecf)
\simeq
 \TNIL(V,\DDtilde^{\sankaku}_V,\vecf)$.
By construction,
it is easy to check
\[
\gamma^{\ast}\TNIL(\kappa)\circ
 \TNIL(\kappa)=\id.
\]
It is also easy to check the compatibility condition,
if the original $\nbigs$ and $\kappa$ are compatible.
Therefore,
we obtain a variation of twistor-TERP structures 
$\TNIL(V,\DDtilde^{\sankaku}_V,\vecf,\nbigs,\kappa,-w)$
on $X-D$
from a variation of twistor-TERP structures
$(V,\DDtilde^{\sankaku}_V,\nbigs,\kappa,-w)$
with $\vecf=(f_i)$ as above.

\begin{df}
\mbox{{}}
Let $(V,\DD^{\sankaku}_V,\vecf,\nbigs)$
be as above.
We set $X^{\ast}(R):=
Y\times \bigl\{
 (z_1,\ldots,z_n)\,\big|\,
 0<|z_i|<R\bigr\}$ for $R>0$.
\begin{itemize}
\item
If there exists $R>0$ such that
$\TNIL(V,\DD^{\sankaku}_V,\vecf,\nbigs)_{|
 \proj^1\times X^{\ast}(R)}$ 
is pure and polarized,
it is called a twistor nilpotent orbit of weight $w$.
\item
If moreover $(V,\DD_V^{\sankaku})$ is enriched to
a variation of integrable twistor structures
$(V,\DDtilde_V^{\sankaku})$
such that $f_j$ and $\nbigs$ are 
$\DDtilde_V^{\sankaku}$-flat,
$\TNIL(V,\DDtilde^{\sankaku}_V,\vecf,\nbigs)_{|
 \proj^1\times X^{\ast}(R)}$ 
is called an integrable twistor nilpotent orbit of weight $w$.
(We often omit to distinguish ``integrable''
if there is no risk of confusion.)
\item
If moreover $(V,\DDtilde^{\sankaku}_V,\nbigs)$
is equipped with a real structure $\kappa$
such that
 $\kappa\circ\gamma^{\ast}f_i=f_i\circ\kappa$,
the variation
$\TNIL(V,\DDtilde^{\sankaku}_V,
 \nbigs,\kappa,-w)_{|
 \proj^1\times X^{\ast}(R)}$ 
is called a twistor-TERP nilpotent orbit.
\hfill\qed
\end{itemize}
\end{df}

\begin{rem}
The notion of a twistor-TERP nilpotent orbit
is different from ``nilpotent orbit'' defined by
Hertling and Sevenheck.
Their ``nilpotent orbit'' is called
HS-orbit in this paper.
\hfill\qed
\end{rem}

\subsection{Convergence}
\label{subsection;08.8.20.11}

\subsubsection{Complement on
 convergence of pure polarized twistor structures}
\label{subsection;08.7.25.3}

Let $(V^{(i)},\nbigs^{(i)})$ $(i=0,1)$
be polarized pure twistor structures with weight $0$
of rank $r$.
Let $h^{(i)}$ be the hermitian metrics of $V^{(i)}$
corresponding to $\nbigs^{(i)}$,
and let $d^{(i)}$ denote 
the associated flat unitary connections
of $V^{(i)}$,
which are equal to the natural
connection given by
holomorphic trivializations
$V^{(i)}\simeq\nbigo_{\proj^1}^{\oplus r}$.
Let $\delbar^{(i)}$ denote the $(0,1)$-part
of $d^{(i)}$,
which is the same as 
the holomorphic structures of $V^{(i)}$.
We fix a hermitian metric $g$ of
$\Omega^{1,0}_{\proj^1}
\oplus
 \Omega^{0,1}_{\proj^1}$.
Let $\Phi:V^{(0)}\lrarr V^{(1)}$ be
a $C^{\infty}$-isomorphism such that
the following holds for some $\epsilon>0$:
\begin{description}
\item[(A1)]
 $\bigl|
 \Phi^{\ast}\delbar^{(1)}
-\delbar^{(0)}\bigr|_{h^{(0)},g}
\leq \epsilon$
as a $C^{\infty}$-section of
$\End(V^{(0)})\otimes \Omega^{0,1}$.
\item[(A2)]
$\bigl|\Phi^{\ast}\nbigs^{(1)}
 -\nbigs^{(0)}\bigr|_{h^{(0)}}\leq \epsilon$
as a $C^{\infty}$-section of
$\Hom\bigl(V^{(0)}\otimes 
 \sigma^{\ast}V^{(0)},\Tate(0)\bigr)$.
\item[(A3)]
$
\bigl|
\delbar^{(0)}\bigl(
 \Phi^{\ast}\nbigs^{(1)}
-\nbigs^{(0)}
\bigr)
 \bigr|_{h^{(0)},g}
=
\bigl|
\delbar^{(0)}
 \Phi^{\ast}\nbigs^{(1)}
 \bigr|_{h^{(0)},g}
 \leq \epsilon$ as a $C^{\infty}$-section 
 of $Hom\bigl(V^{(0)}\otimes 
 \sigma^{\ast}V^{(0)},\Tate(0)\bigr)
 \otimes\Omega^{0,1}_{\proj^1}$,
 where
 $\delbar^{(0)}$ denotes
 the induced holomorphic structure on  
 $Hom\bigl(V^{(0)}\otimes 
 \sigma^{\ast}V^{(0)},\Tate(0)\bigr)$.
\end{description}

\begin{lem}
\label{lem;08.10.27.11}
There exists a constant $C_0>0$,
which is independent of $\epsilon$,
with the following property:
\begin{itemize}
\item
If 
$B^{-1}\cdot\Phi^{\ast}h^{(1)}
\leq
 h^{(0)}
\leq
 B\cdot \Phi^{\ast}h^{(1)}$
for some $B>1$,
the following holds:
\begin{equation}
 \label{eq;08.10.27.10}
 \bigl|
 \Phi^{\ast}d^{(1)}-d^{(0)}
 \bigr|_{h^{(0)},g}
\leq
 C_0\cdot B^2\cdot \epsilon
\end{equation}
\end{itemize}
\end{lem}
\pf
In the following argument,
$C_i$ denote positive constants
independent of $\epsilon$.
Let $\del^{(i)}_{h^{(i)}}$ denote the $(1,0)$-part
of $d^{(i)}$, which are determined by
$h^{(i)}$ and $\delbar^{(i)}$.
To show (\ref{eq;08.10.27.10}),
we have only to estimate
$\bigl|\del^{(0)}_{h^{(0)}}
-\Phi^{\ast}\del^{(1)}_{h^{(1)}}\bigr|_{h^{(0)}}$.

Let $e_1,\ldots,e_r$ be an orthogonal frame 
of $V^{(0)}$ with respect to $h^{(0)}$.
Because
$\Phi^{\ast}h^{(1)}(e_i,e_j)
=\Phi^{\ast}\nbigs^{(1)}\bigl(
 e_i\otimes\sigma^{\ast}e_j \bigr)$,
we have the following estimate
for any $i,j$:
\[
 \Bigl|
 \delbar \bigl(
 \Phi^{\ast}h^{(1)}(e_i,e_j)\bigr)
\Bigr|_{g}
=\Bigl|
 \delbar^{(0)}(\Phi^{\ast}\nbigs^{(1)})
 (e_i\otimes\sigma^{\ast}e_j)
\Bigr|_g
\leq C_1\cdot\epsilon
\]
Hence, we obtain
\begin{equation}
 \label{eq;08.7.24.2}
 \Bigl|
 \del \bigl(\Phi^{\ast}h^{(1)}(e_i,e_j)\bigr)
\Bigr|_{g}\leq C_1\cdot\epsilon
\end{equation}
Let $\del^{(0)}_{h^{(1)}}$ denote 
the $(1,0)$-operator determined by
$\Phi^{\ast}h^{(1)}$ and $\delbar^{(0)}$.
From $B^{-1}\cdot \Phi^{\ast}h^{(1)}
\leq h^{(0)}
\leq B\cdot \Phi^{\ast}h^{(1)}$
and (\ref{eq;08.7.24.2}),
we obtain 
\[
\bigl|
 \del^{(0)}_{h^{(1)}}-\del^{(0)}_{h^{(0)}}
\bigr|_{h^{(0)},g}
\leq C_2\cdot B\cdot\epsilon
\]
By
$\bigl|
\Phi^{\ast}\delbar^{(1)}-\delbar^{(0)}
\bigr|_{h^{(0)}}\leq C_3\cdot\epsilon$
and 
$B^{-1}\cdot \Phi^{\ast}h^{(1)}
\leq h^{(0)}
\leq B\cdot \Phi^{\ast}h^{(1)}$,
we have
\[
 \bigl|
\Phi^{\ast}\delbar^{(1)}-\delbar^{(0)}
\bigr|_{\Phi^{\ast}h^{(1)}}
 \leq C_4\cdot B\cdot\epsilon
\]
Hence, we obtain
$\bigl|\Phi^{\ast}\del^{(1)}_{h^{(1)}}
-\del^{(0)}_{h^{(1)}}\bigr|_{\Phi^{\ast}h^{(1)}}
\leq
 C_4\cdot B\cdot\epsilon$,
which implies
\[
 \bigl|\Phi^{\ast}\del^{(1)}_{h^{(1)}}
-\del^{(0)}_{h^{(1)}}\bigr|_{h^{(0)}}
\leq
 C_5\cdot B^2\cdot\epsilon
\]
Thus, we obtain (\ref{eq;08.10.27.10}).
\hfill\qed

\begin{lem}
\label{lem;08.7.25.5}
There exist constants $\epsilon_0>0$,
$C_{10}>0$ and $C_{11}>0$ 
such that the following holds
if $\epsilon\leq \epsilon_0$:
\begin{equation}
 \label{eq;08.8.9.10}
\bigl|\Phi^{\ast}h^{(1)}-h^{(0)}
\bigr|_{h^{(0)}}
\leq C_{10}\cdot\epsilon
\end{equation}
\begin{equation}
\label{eq;08.8.9.11}
\bigl|
 \Phi^{\ast}d^{(1)}-d^{(0)}
\bigr|_{h^{(0)}}
\leq C_{11}\cdot \epsilon
\end{equation}
\end{lem}
\pf
According to the result in Subsection 2.8
of \cite{mochi7},
if $\epsilon_0$ is sufficiently small,
(\ref{eq;08.8.9.10}) holds for some $C_{10}$.
Then, we obtain (\ref{eq;08.8.9.11})
from Lemma \ref{lem;08.10.27.11}.
\hfill\qed

\subsubsection{Approximation of 
 pure polarized integrable twistor structures}

Let $(V^{(i)},\nabla^{(i)},\nbigs^{(i)})$ $(i=1,2)$
be integrable polarized pure twistor structures.
Let $h^{(i)}$ be the hermitian metrics of
$V^{(i)}$
corresponding to $\nbigs^{(i)}$.
We fix a hermitian metric $\gtilde$ of
$\Omega^{1,0}_{\proj^1}
 \bigl(2\cdot 0+2\cdot\infty\bigr)
\oplus
 \Omega^{0,1}_{\proj^1}$.
Let $\Phi:V^{(0)}\lrarr V^{(1)}$ be
a $C^{\infty}$-isomorphism such that
the following holds for some $\epsilon>0$:
\begin{description}
\item[(B1)]
 $\bigl|\Phi^{\ast}\nabla^{(1)}
 -\nabla^{(0)}\bigr|_{h^{(0)},\gtilde}\leq \epsilon$
as a $C^{\infty}$-section of
$ \End(V^{(0)})\otimes\Bigl(
 \Omega^{1,0}_{\proj^1}
 \bigl(2\cdot 0+2\cdot\infty\bigr)
\oplus
 \Omega^{0,1}_{\proj^1}\Bigr)$.
Note that it implies 
(A1) in Subsection \ref{subsection;08.7.25.3}.
\item[(B2)]
Conditions (A2) and (A3) are satisfied.
\end{description}

\begin{lem}
\label{lem;08.10.27.15}
There exists a constant $C_{20}>0$,
which is independent of $\epsilon$,
with the following property:
\begin{itemize}
\item
 If 
$B^{-1}\cdot\Phi^{\ast}h^{(1)}
\leq
 h^{(0)}
\leq
 B\cdot \Phi^{\ast}h^{(1)}$
for some $B>1$,
the following holds:
\[
\bigl|
 \Phi^{\ast}\nbigu^{(1)}-\nbigu^{(0)}
\bigr|_{h^{(0)}}
\leq C_{20}\cdot B^2\cdot\epsilon,
\quad\quad
\bigl|
 \Phi^{\ast}\nbigq^{(1)}-\nbigq^{(0)}
\bigr|_{h^{(0)}}
\leq C_{20}\cdot B^2\cdot\epsilon.
\]
\end{itemize}
\end{lem}
\pf
In the following argument,
$C_i$ denote positive constants
independent of $\epsilon$.
By Lemma \ref{lem;08.7.25.5},
we have
$\bigl|
 \Phi^{\ast}d^{(1)}-d^{(0)}
 \bigr|_{h^{(0)},\gtilde}
\leq C_{21}\cdot B^2\cdot \epsilon$.
We obtain the following estimate:
\[
 \left|
\Bigl(
 \lambda^{-1}
 \cdot(\Phi^{\ast}\nbigu^{(1)}-\nbigu^{(0)})
-(\Phi^{\ast}\nbigq^{(1)}-\nbigq^{(0)})
-\lambda
 \cdot(\Phi^{\ast}\nbigu^{(1)\dagger}
-\nbigu^{(0)\dagger})
\Bigr)\cdot d\lambda/\lambda
\right|_{h^{(0)},\gtilde}
\leq C_{22}\cdot B^2\cdot\epsilon
\]
Then, the claim of the lemma follows.
\hfill\qed

\begin{lem}
\label{lem;08.7.26.11}
There exist  constants 
$\epsilon_0>0$ and $C_{30}$,
such that the following holds
for any $0<\epsilon\leq \epsilon_0$:
\[
\bigl|\Phi^{\ast}h^{(1)}-h^{(0)}
\bigr|_{h^{(0)}}
\leq C_{30}\cdot\epsilon,
\quad
\bigl|
 \Phi^{\ast}\nbigu^{(1)}-\nbigu^{(0)}
\bigr|_{h^{(0)}}
\leq C_{30}\cdot\epsilon,
\quad\quad
\bigl|
 \Phi^{\ast}\nbigq^{(1)}-\nbigq^{(0)}
\bigr|_{h^{(0)}}
\leq C_{30}\cdot\epsilon.
\]
\end{lem}
\pf
It can be shown by the argument
in the proof of Lemma \ref{lem;08.10.27.15}.
\hfill\qed

\subsection{Variation of 
polarized mixed twistor structures
and its enrichment}
\label{subsection;08.8.20.13}

\subsubsection{Definitions}

\label{subsubsection;08.8.10.40}

\paragraph{Variation of
polarized mixed twistor structures}

Let $X$ be a complex manifold.
Let $(V,\DD^{\sankaku})$ be 
a variation of $\proj^1$-holomorphic 
vector bundles on $\proj^1\times X$
equipped with an increasing filtration $W$
indexed by $\seisuu$
in the category of vector bundles,
which is $\proj^1$-holomorphic
and $\DD^{\sankaku}$-flat.
If each $\Gr^W_n(V)$ is 
a variation of pure twistor structure
of weight $n$, 
$(V,W,\DD^{\sankaku})$ is called
a variation of mixed twistor structures.
Assume we are given
the following data on $(V,W,\DD^{\sankaku})$,
which are $\proj^1$-holomorphic
and $\DD^{\sankaku}$-flat:
\begin{itemize}
\item
 A tuple $\vecf$ of 
 nilpotent morphisms
 $f_j:V\lrarr V\otimes\Tate(-1)$
 $(j=1,\ldots,n)$,
 which are mutually commutative.
\item
 A $(-1)^w$-symmetric pairing 
 $\nbigs:V\otimes\sigma^{\ast}V\lrarr \Tate(-w)$.
\item
 For each $P\in X$,
 the restriction
 $(V,W,\vecf,\nbigs)_{|\proj^1\times \{P\}}$
 is a polarized mixed twistor structure
 of weight $w$ in $n$-variables.
 (See Subsection 3.48 of \cite{mochi2}.)
\end{itemize}
Then, such a tuple
$(V,\DD^{\sankaku},
 W,\vecf,\nbigs)$ is called 
a variation of polarized mixed twistor structures.
Since $W$ is determined by $\vecf$
as the weight filtration of 
$f(\nbar):=\sum_{j=1}^{n} f_j$ up to shift by $w$,
we sometimes omit to denote $W$.

\paragraph{Enrichment}

If $\DD^{\sankaku}$ and the $\proj^1$-holomorphic
structure are extended to
$\TTtilde E$-structure $\DDtilde^{\sankaku}$
for which $\vecf$ and $\nbigs$ are flat,
$(V,\DDtilde^{\sankaku},
 W,\vecf,\nbigs)$ is called 
a variation of polarized mixed 
integrable twistor structures
of weight $w$ in $n$-variables.
Note that $W$ is automatically
$\DDtilde^{\sankaku}$-flat.

If moreover $(V,\DDtilde^{\sankaku},\nbigs)$
is equipped with real structure $\kappa$
such that
$\kappa\circ\gamma^{\ast}f_j=f_j\circ\kappa$,
then such a tuple
$(V,\DDtilde^{\sankaku},
 W,\vecf,\nbigs,\kappa,-w)$ is called 
a variation of polarized mixed 
twistor-TERP structures in $n$-variables.

\begin{rem}
The notion of polarized mixed twister-TERP structure
 is different from
``mixed TERP structure''
defined by 
Hertling and Sevenheck
(Section {\rm 9} of 
{\rm \cite{Hertling-Sevenheck}}).
\hfill\qed
\end{rem}

\paragraph{Split type}

Let $(V,W,\DD^{\sankaku})$
be a variation of mixed twistor structures.
It is called of split type,
if it is equipped with a grading
$V=\bigoplus V_m$ such that
(i) it is $\proj^1$-holomorphic and $\DD^{\sankaku}$-flat,
(ii) $W_m=\bigoplus_{p\leq m}V_p$.
Each $(V_m,\DD^{\sankaku})$ 
is a variation of pure twistor structures
of weight $m$.

A variation of polarized mixed twistor structures
of weight $w$ in $n$-variables
$(V,W,\DD^{\sankaku},\vecf,\nbigs)$
is called of split type,
if the underlying variation of mixed twistor structures
$(V,W,\DD^{\sankaku})$ is of split type
with a grading $V=\bigoplus V_m$.
By using
$H^0\bigl(\proj^1,\nbigo_{\proj^1}(m)\bigr)=0$
for any $m<0$,
we can show that the following:
\begin{itemize}
\item
$f_j(V_p)\subset V_{p-2}\otimes\Tate(-1)$.
\item
The restriction of $\nbigs$ to
$V_p\otimes \sigma^{\ast}V_q$ is $0$ 
unless $p+q=2w$.
\end{itemize}
Similarly,
a variation of
polarized mixed integrable twistor structures
of weight $w$ in $n$-variables
$(V,W,\DDtilde^{\sankaku},\vecf,\nbigs)$
is called of split type,
if the underlying variation of
polarized mixed twistor structure
is of split type with a $\DDtilde^{\sankaku}$-flat
grading.

A polarized mixed twistor-TERP structures
$(V,W,\nabla,\vecf,\nbigs,\kappa,-w)$
in $n$-variables is called of split type,
if the underlying variation of mixed integrable
twistor structures is of split type
with a grading $V=\bigoplus V_m$
such that $\kappa(\gamma^{\ast}V_m)=V_m$.

\subsubsection{Reduction}
\label{subsubsection;08.8.11.2}

Let $(V,W,\DD^{\sankaku},\vecf,\nbigs)$
be a variation of polarized mixed twistor
structures of weight $w$ in $n$-variables.
We obtain a variation of $\proj^1$-holomorphic
vector bundles
$(V^{(0)},\DD^{(0)\sankaku}):=
 \Gr^W(V,\DD^{\sankaku})$.
It is naturally equipped with a grading
$V^{(0)}=\bigoplus \Gr^W_m(V)$
and a filtration
$W_m^{(0)}=\bigoplus_{p\leq m}\Gr^W_p(V)$.
We have induced morphisms
$f_j^{(0)}:\Gr^W_m(V)\lrarr 
 \Gr^W_{m-2}(V)\otimes\Tate(-1)$,
and hence
$ f_j^{(0)}:V^{(0)}\lrarr V^{(0)}\otimes\Tate(-1)$.
We also obtain induced morphisms
$\nbigs^{(0)}:
 \Gr^W_{w-m}(V)\otimes
 \sigma^{\ast}\Gr^W_{w+m}(V)
\lrarr \Tate(-w)$,
and hence
$ \nbigs^{(0)}:
 V^{(0)}\otimes\sigma^{\ast}(V^{(0)})
\lrarr
 \Tate(-w)$.
It is known that
$(V^{(0)},W^{(0)},\vecf^{(0)},\nbigs^{(0)})
_{|\proj^1\times\{P\}}$
are polarized mixed twistor structures
of split type with weight $w$
in $n$-variables. 
(See \cite{mochi2}.
It can be shown directly and easily.)
Hence, 
$(V^{(0)},W^{(0)},\DD^{(0)\sankaku},
 \vecf^{(0)},\nbigs^{(0)})$
is a variation of polarized mixed twistor structures
of split type with weight $w$ in $n$-variables.
It is denoted by
$\Gr^W(V,W,\DD^{\sankaku},\vecf,\nbigs)$.

If $(V,W,\DD^{\sankaku},\vecf,\nbigs)$
is enriched to 
a variation of polarized mixed integrable
twistor structures
with weight $w$ in $n$-variables,
$\Gr^W(V,W,\DD^{\sankaku},\vecf,\nbigs)$
is also integrable.
If moreover the variation of polarized mixed
integrable twistor structures is enriched to
a variation of polarized mixed twistor-TERP structures,
$\Gr^W$ is also enriched to 
a variation of polarized mixed twistor-TERP
structures of split type.

\subsubsection{Splittings}
\label{subsubsection;08.8.20.14}

\paragraph{Preliminary}

Let $(V_i,W,\DD_i^{\sankaku})$ $(i=1,2)$ 
be variations of mixed twistor structures
on $\proj^1\times X$
with a morphism 
$F:(V_1,W,\DD_1^{\sankaku})
 \lrarr
 (V_2,W,\DD_2^{\sankaku})$.
We set 
$(V_i^{(0)},\DD^{(0)\sankaku}_i)
 :=\Gr^W(V_i,\DD^{\sankaku}_i)$
on which we have the naturally induced filtrations
$W^{(0)}$.
We also obtain induced morphism
$F^{(0)}:(V_1^{(0)},W^{(0)},\DD^{(0)\sankaku}_1)
\lrarr (V_2^{(0)},W^{(0)},\DD^{(0)\sankaku}_2)$.
The following lemma is standard.
\begin{lem}
The rank of $F_{|(\lambda,P)}$ is independent of
$(\lambda,P)\in\proj^1\times X$.
The morphism
$F$ is strict with respect to
the weight filtration.
Hence,
$\Ker F$ with the induced filtration $W(\Ker F)$
is a mixed twistor structure,
and we have the isomorphism
$ \Ker F^{(0)}\simeq
 \Gr^W\bigl(\Ker F\bigr)$.
\end{lem}
\pf
If $X$ is a point,
the claims are well known
and easy to show.
Namely, it is shown in Lemma 2.20 of \cite{mochi}
that 
(i) $\Ker(F)$ is a subbundle of $V_1$,
(ii) $F$ is strict with respect to the weight filtrations,
 i.e., $F(W_l(V_1))=F(V_1)\cap W_l(V_2)$,
(iii) $\Ker(F)$ with the induced weight filtration
 is a mixed twistor structure.
We obtain the isomorphism
$ \Ker F^{(0)}\simeq
 \Gr^W\bigl(\Ker F\bigr)$ 
from the strictness.

Let us consider the general case.
By using the flatness,
it is easy to show that
$\rank F_{|(1,P)}$ and
$\rank F^{(0)}_{|(1,P)}$
are independent of the choice of
a point $P\in X$.
Then, the claim of the lemma follows.
\hfill\qed

\begin{cor}
\label{cor;08.10.27.20}
Let $(V_i,W,\DD^{\sankaku}_i)$ $(i=0,1,\ldots,m)$
be variations of mixed twistor structures
with morphisms 
$F_i:(V_0,W,\DD^{\sankaku}_0)
 \lrarr 
 (V_i,W,\DD^{\sankaku}_i)$
$(i=1,\ldots,m)$.
Then,
we have the following natural isomorphism
of variations of mixed twistor structures:
\[
 \Gr^W\Bigl(
 \bigcap_{i=1}^m\Ker F_i
 \Bigr)
\simeq
 \bigcap_{i=1}^m\Ker F_i^{(0)}
\]
Here, $F_i^{(0)}$ denote induced morphisms
$(V_{0}^{(0)},W^{(0)},\DD_0^{(0)\sankaku})
 \lrarr (V_i^{(0)},W^{(0)},\DD_i^{(0)\sankaku})$.
\hfill\qed
\end{cor}

\paragraph{Local splitting}

Let $(V,W,\DD^{\sankaku})$ 
be a variation of mixed twistor structures.
Let $\vecN=(N_j\,|\,j=1,\ldots,\ell)$ 
be a tuple of morphisms
$N_j:(V,W,\DD^{\sankaku})
 \lrarr 
 (V,W,\DD^{\sankaku})\otimes \Tate(-1)$
which are mutually commutative.
Let $(V^{(0)},W^{(0)},\DD^{(0)\sankaku})$ 
be as above.
Let $\vecN^{(0)}=(N_j^{(0)}\,|\,j=1,\ldots,\ell)$
be the induced commuting tuple of morphisms
$N_j^{(0)}:(V^{(0)},W^{(0)},\DD^{(0)\sankaku})
\lrarr (V^{(0)},W^{(0)},\DD^{(0)\sankaku})
 \otimes\Tate(-1)$.

We set
$\Vbar:=\Hom(V^{(0)},V)$,
which is naturally equipped with
the operator $\DDbar^{\sankaku}$
and an induced filtration $\Wbar$.
Let $\Nbar_j:(\Vbar,\Wbar,\DDbar)
\lrarr (\Vbar,\Wbar,\DDbar)\otimes\Tate(-1)$
be  the morphisms of
mixed twistor structures
given by
$\Nbar_j(f)=N_j\circ f-f\circ N_j^{(0)}$.
Similarly, we set 
$\Vbar^{(0)}:=\Hom(V^{(0)},V^{(0)})$
on which we have the naturally induced
operator $\DDbar^{(0)\sankaku}$,
filtration $\Wbar^{(0)}$
and morphisms of
mixed twistor structures
$\Nbar^{(0)}_j:
 (\Vbar^{(0)},\Wbar^{(0)},\DDbar^{(0)})
\lrarr 
 (\Vbar^{(0)},\Wbar^{(0)},\DDbar^{(0)})
 \otimes\Tate(-1)$.

We have the natural isomorphism
$\Gr^{\Wbar}(\Vbar)
\simeq \Vbar^{(0)}$.
The induced filtrations and the morphisms
coincide.
According to Corollary \ref{cor;08.10.27.20},
we have the following isomorphism
of variations of mixed twistor structures:
\[
 \Gr^{\Wbar}\Bigl(
 \bigcap\Ker\Nbar_j
 \Bigr)
\simeq
 \bigcap
 \Ker\Nbar^{(0)}_j.
\]
Then, we obtain the following corollary.
\begin{cor}
\label{cor;08.7.28.15}
Let $(\lambda,P)$ be any point of 
$\cnum_{\lambda}\times X$,
and let $U$ be a small neighbourhood of 
$(\lambda,P)$.
There exists a $C^{\infty}$-morphism
$F:V^{(0)}_{|U} \lrarr V_{|U}$
with the following property:
\begin{itemize}
\item
 It preserves the weight filtration,
 and the induced morphism
 on $\Gr^{W}(V^{(0)}_{|U})
\lrarr
 \Gr^W(V_{|U})$ is the identity.
\item
$F\circ N_j^{(0)}=N_j\circ F$
for $j=1,\ldots,\ell$.
\hfill\qed
\end{itemize}
\end{cor}

\paragraph{$C^{\infty}$-splitting}

Let $(V,W,\DD^{\sankaku},\vecN)$
and $(V^{(0)},W^{(0)},\DD^{(0)\sankaku},
 \vecN^{(0)})$
be as above.

\begin{lem}
\label{lem;08.7.29.1}
There exists a $C^{\infty}$-isomorphism
$\Phi:V^{(0)}\lrarr V$ with the following property:
\begin{itemize}
\item
 $\Phi$ preserves the weight filtration $W$,
 and
 $\Gr^W\Phi$ is the identity 
 $\Gr^W(V^{(0)})=\Gr^W(V)$.
\item
 $\Phi\circ N_j^{(0)}=N_j\circ\Phi$
 for $j=1,\ldots,\ell$.
\end{itemize}
\end{lem}
\pf
Let $U\subset\cnum_{\lambda}$ be 
a compact region with
$U\cup\sigma(U)=\proj^1$.
We take a locally finite open covering
$U\times X\subset 
 \bigcup_{p\in I} \nbigu_p$
such that 
we have $C^{\infty}$-isomorphisms
$\Phi_{\nbigu_p}:
 V^{(0)}_{|\nbigu_p}
\lrarr
 V_{|\nbigu_p}$
as in  Corollary \ref{cor;08.7.28.15},
i.e.,
$\Phi_{\nbigu_p}\circ 
 N^{(0)}_j=N_j\circ\Phi_{\nbigu_p}$
for any $j$.
Similarly,
we take a locally finite open covering
$\sigma(U)\times X^{\dagger}
\subset
 \bigcup_{q\in J} \nbigu_q^{\dagger}$
such that 
we have $C^{\infty}$-isomorphisms
$\Phi_{\nbigu_q^{\dagger}}:
 V^{(0)}_{|\nbigu_q^{\dagger}}
 \simeq
 V_{|\nbigu_q^{\dagger}}$
as in Corollary \ref{cor;08.7.28.15}.
We take a partition of unity
$\bigl\{
 \chi_{\nbigu_p},\chi_{\nbigu_q^{\dagger}}\,
 \big|\,p\in I,\,\,q\in J \bigr\}$
subordinated to the covering
$\bigl\{\nbigu_p,\nbigu_q^{\dagger}\,\big|\,
 p\in I,\,\,q\in J \bigr\}$ of 
$\proj^1\times X$.
We obtain the $C^{\infty}$-isomorphism
\[
 \Phi:=
 \sum_{p\in I}\chi_{\nbigu_p}\cdot \Phi_{\nbigu_p}
+\sum_{q\in J}\chi_{\nbigu_q^{\dagger}}
 \cdot \Phi_{\nbigu_q^{\dagger}}
 :V^{(0)}\lrarr V.
\]
By construction,
it has the desired property.
\hfill\qed

\section{Polarized mixed integrable twistor structure
of split type}
\label{section;08.8.9.1}

\subsection{Basic examples in one variable}
\label{subsection;08.8.9.22}

\subsubsection{Rank two}

Let us recall a basic example 
studied in Subsection 3.7.2 of \cite{mochi2}
with a minor enrichment.
We set $V^{[2]}:=
 \nbigo(0,-1)\oplus
 \nbigo(1,0)$.
(See Subsection \ref{subsubsection;08.9.11.1}
 for $\nbigo(p,q)$.)
It is naturally equipped with 
a meromorphic connection $\nabla^{[2]}$,
and $(V^{[2]},\nabla^{[2]})$ is
an integrable twistor structure.
We put
\[
 W_{-2}(V^{[2]}):=0,
\quad
 W_{-1}(V^{[2]})=W_{0}(V^{[2]}):=\nbigo(0,-1),
\quad
W_{1}(V^{[2]}):=V^{[2]}.
\]
Let $F^{[2]}:V^{[2]}\lrarr V^{[2]}\otimes\Tate(-1)$
be given by
\[
 f_a^{(1,0)}\longmapsto
 f_a^{(0,-1)}\otimes t^{(-1)}_a,
\quad (a=0,1,\infty),
\quad
 f_a^{(0,-1)}\longmapsto 0.
\]
A flat morphism
$S^{[2]}:V^{[2]}\otimes\sigma^{\ast}V^{[2]}
\lrarr \Tate(0)$ is given by the following correspondence:
\[
 S^{[2]}\bigl(
 f_1^{(1,0)}\otimes\sigma^{\ast}f_1^{(0,-1)}
 \bigr)
=\sqrt{-1}t_1^{(0)},
\quad
  S^{[2]}\bigl(
 f_1^{(0,-1)}\otimes\sigma^{\ast}f_1^{(1,0)}
 \bigr)
=-\sqrt{-1}t_1^{(0)},
\]
\[
  S^{[2]}\bigl(
 f_1^{(1,0)}\otimes\sigma^{\ast}f_1^{(1,0)}
 \bigr)
=0,
\quad
  S^{[2]}\bigl(
 f_1^{(0,-1)}\otimes\sigma^{\ast}f_1^{(0,-1)}
 \bigr)
=0.
\]
Recall that $(V^{[2]},W,F^{[2]},S^{[2]})$
is a polarized mixed twistor structure of split type
in one variable with weight $0$
(Lemma 3.90 of \cite{mochi2}).
Hence, $(V^{[2]},W,\nabla^{[2]},
F^{[2]},S^{[2]})$ is a polarized mixed integrable
twistor structure of split type.

\subsubsection{Twist}
\label{subsubsection;08.7.26.2}

The bundle $V^{[2]}$ is obtained as
the gluing of 
$V^{[2]}_0:=V^{[2]}_{|\cnum_{\lambda}}$
and $V^{[2]}_{\infty}:=V^{[2]}_{|\cnum_{\mu}}$.
We would like to explain a twist of the gluing
given in Subsection 3.7.2 of \cite{mochi2},
related with the construction in 
Subsection \ref{subsubsection;08.7.26.1}.
Let $N:=F^{[2]}\otimes t_1^{(1)}$.
Let $v\in V^{[2]}_{|\lambda}$
for $\lambda\neq 0,\infty$.
The induced elements of $V^{[2]}_{0|\lambda}$ 
and $V^{[2]}_{\infty|\mu}$ are denoted by $v$
and $v^{\dagger}$,
respectively.
The gluing for $V^{[2]}$ is given by 
$v=v^{\dagger}$.
For $y\in\cnum$,
a vector bundle $\Vtilde^{[2]}_y$
is given by the following  twisted gluing:
\[
 \exp\bigl(\sqrt{-1}y\cdot N\bigr)\cdot v
=v^{\dagger}
\]
Since $N$ is flat,
we have the naturally induced flat connection
$\nabla^{[2]}_y$ of $\Vtilde^{[2]}_y$.
We also have the induced pairing
$\Stilde^{[2]}_y$ of $(\Vtilde^{[2]}_y,\nabla^{[2]}_y)$
of weight $0$.

For $y\neq 0$,
we have a frame of $\Vtilde^{[2]}_y$ given as follows:
\[
 \stilde_1:=
 \sqrt{-1}\lambda\cdot f^{(1,0)}_0
+\sqrt{-1}y\cdot f_0^{(0,-1)}
=f_{\infty}^{(1,0)}
\]
\[
 \stilde_2:=f_0^{(1,0)}
=-\sqrt{-1}\mu\cdot f_{\infty}^{(1,0)}
-\sqrt{-1}y\cdot f_{\infty}^{(0,-1)}
\]
In particular,
$(\Vtilde^{[2]}_y,\nabla^{[2]}_y)$ is 
a pure integrable twistor structure
of weight $0$ for any $y\neq 0$.
If $y$ is a positive real number,
$\Stilde^{[2]}_y$ gives a polarization
of $(\Vtilde^{[2]}_y,\nabla^{[2]}_y)$
(Lemma 3.91 of \cite{mochi2}).
Actually, 
$\stilde_i$ $(i=1,2)$ give
an orthogonal frame:
\[
 \Stilde^{[2]}_y(\stilde_i,\sigma^{\ast}\stilde_i)=y\,\,\,\,(i=1,2),
\quad
 \Stilde^{[2]}_y(\stilde_1,\sigma^{\ast}\stilde_2)=0.
\]

Note that $\nabla^{[2]}_y$ is logarithmic
with respect to the lattice $\Vtilde_y^{[2]}$.
For any $y\neq 0$,
we have the decomposition
\[
 \nabla^{[2]}_y=d^{[2]}_y
-\nbigq^{[2]}_y
 \frac{d\lambda}{\lambda}
\]
Here, $d^{[2]}_y$ is a natural flat connection
of $V^{[2]}_y\simeq\nbigo_{\proj^1}(0)^{\oplus\,2}$.
Let us calculate $\nbigq^{[2]}$.
By easy calculations,
\[
 \nabla^{[2]}_y \stilde_1=0,
\quad
 \nabla^{[2]}_y\stilde_2
=\stilde_2\cdot \left(-\frac{d\lambda}{\lambda}\right).
\]
Hence, $\nbigq^{[2]}$ is expressed 
by the following matrix
with respect to the frame
$\stilde_1,\stilde_2$:
\[
 \left(
 \begin{array}{cc}
 0 & 0 \\
 0 & 1
 \end{array}
 \right)
\]
In particular,
the eigenvalues are independent of $y$.

\begin{rem}
For our application,
we essentially need only the case 
in which $y$ is a positive real number.
Recall that we have considered 
a twisted isomorphism {\rm (\ref{eq;08.10.1.1})}.
We will use the above consideration
by setting
$y=-\sum_{i=1}^{\ell}\log|z_i|^2$.
\hfill\qed
\end{rem}

\subsubsection{Rank $\ell$}
\label{subsubsection;08.7.26.5}

For any positive integer $\ell$,
we set 
$(V^{[\ell]},\nabla^{[\ell]}):=
\Sym^{\ell-1}(V^{[2]},\nabla^{[2]})$,
equipped with 
a morphism 
 $F^{[\ell]}:V^{[\ell]}\lrarr
 V^{[\ell]}\otimes \Tate(-1)$
and a pairing
$S^{[\ell]}:
 V^{[\ell]}\otimes\sigma^{\ast}V^{[\ell]}
 \lrarr \Tate(0)$.
For any $y\in\cnum$,
we obtain an integrable twistor structure
$(\Vtilde^{[\ell]}_y,\nabla^{[\ell]}_y,\Stilde^{[\ell]}_y)$
with a pairing of weight $0$,
by the procedure
in Subsection \ref{subsubsection;08.7.26.2}.
It is also obtained as the $(\ell-1)$-th symmetric product of
$(\Vtilde^{[2]}_y,\nabla^{[2]}_y,\Stilde^{[2]}_y)$.
Hence,
$(\Vtilde^{[\ell]}_y,\nabla^{[\ell]}_y)$ is pure
with weight $0$ for each $y\neq 0$,
and 
$\Stilde^{[\ell]}_y$ gives a polarization
for each $y>0$.
We have the decomposition
\[
 \nabla^{[\ell]}_y=
 d^{[\ell]}_y-\nbigq^{[\ell]}\frac{d\lambda}{\lambda}
\]
Let $y\neq 0$.
A frame of $\Vtilde^{[\ell]}_y$
is given by symmetric products
$\stilde^{[\ell]}_p:=
 \stilde^{\ell-1-p}_1\cdot
 \stilde^{p}_2$
$(p=0,1,\ldots,\ell-1)$,
for which $\nbigq^{[\ell]}$ is expressed by
the diagonal matrix 
whose $p$-th entry is $p$
$(p=0,1,\ldots,\ell-1)$.
In particular,
the eigenvalues are independent of $y$.

\subsection{Twistor nilpotent orbits of split type
 and their new supersymmetric indices}

\subsubsection{One variable case}
\label{subsubsection;08.7.26.6}
Let $Y$ be a complex manifold.
Let $(V,W,\DD^{\sankaku},N,S)$ be
a variation of polarized mixed twistor structures
of split type with weight $0$ in one variable
on $\proj^1\times Y$.
The following lemma is 
essentially the same as 
Corollary 3.97 of \cite{mochi2}.

\begin{prop}
\label{prop;08.8.9.21}
There exist variations of 
polarized pure twistor structures
$(U_{\ell},\DD^{\sankaku}_{\ell},S_{\ell})$
of weight $0$  on $\proj^1\times Y$
for $\ell\geq 1$,
such that 
(i) $(V,\DD^{\sankaku})
 \simeq\bigoplus_{\ell\geq 1} 
 (U_{\ell},\DD^{\sankaku}_{\ell})\otimes V^{[\ell]}$,
(ii)
 $N=\bigoplus \id_{U_{\ell}}\otimes F^{[\ell]}$
 and
 $S=\bigoplus S_{\ell}\otimes S^{[\ell]}$
  under the isomorphism.

If $(V,W,\DD^{\sankaku},N,S)$ is 
enriched to be integrable,
$(U_{\ell},\DD^{\sankaku},S_{\ell})$
are also integrable.
\end{prop}
\pf
We have the grading $V=\bigoplus_{j\in\seisuu} V_{j}$.
For each $j\geq 0$,
we set
$PV_{j}:=\Ker\Bigl(N^{j+1}:
 V_j\lrarr V_{-j-2}\otimes \Tate(j+1)\Bigr)$.
It is a variation of pure twistor structure of weight $j$,
and equipped with the induced polarization $S_j$.
For $\ell\geq 1$,
we set
\[
 U_{\ell}:=
 PV_{\ell-1}
\otimes
 \nbigo(0,-\ell+1).
\]
which are naturally 
variations of polarized pure twistor
structures.
Then, it is easy to observe that
$V$ has the desired decomposition.
The integrable case is also easy.
\hfill\qed

\vspace{.1in}

Let $q:Y\times \cnum^{\ast}\lrarr Y$
denote the projection.
We obtain the variations of polarized pure twistor
structures on 
$\proj^1\times (Y\times\cnum^{\ast})$
obtained as the pull back of
$(U_{\ell},\DD_{\ell}^{\sankaku},S_{\ell})$,
denoted by
$q^{\ast}
 (U_{\ell},\DD_{\ell}^{\sankaku},S_{\ell})$.
Recall the construction in
Subsection \ref{subsubsection;08.7.26.1}.
We obtain the following isomorphism
from Proposition \ref{prop;08.8.9.21}:
\begin{equation}
 \label{eq;08.7.26.3}
 \TNIL(V,\DD^{\sankaku},N,S)
\simeq
\bigoplus_{\ell}
 q^{\ast}(U_{\ell},\DD^{\sankaku}_{\ell},S_{\ell})
\otimes
 \TNIL(V^{[\ell]},\nabla^{[\ell]},F^{[\ell]},
 S^{[\ell]})
\end{equation}
By using the result in Subsection 
\ref{subsection;08.8.9.22},
we can conclude the following:
\begin{prop}
We set $X_+:=
Y\times\bigl\{z\in\cnum\,\big|\,0<|z|<1\bigr\}$
and $X_-:=
Y\times \bigl\{z\in\cnum\,\big|\,|z|>1
 \bigr\}$.
Then,
$\TNIL(V,\DD^{\sankaku},N)$ 
is a variation of pure integrable twistor structures
on $\proj^1\times (X_+\cup X_-)$,
and the restriction 
$\TNIL(V,\DD^{\sankaku},N,S)_{|\proj^1\times X_+}$
is a twistor nilpotent orbit.
\hfill\qed
\end{prop}

Assume that
$(V,\DD^{\sankaku})$ is enriched to
integrable $(V,\DDtilde^{\sankaku})$
such that $S$ and $N$ are 
$\DDtilde^{\sankaku}$-flat.
Let $\nbigq$ and $\nbigq^{[\ell]}$
be the new supersymmetric indices of
$\TNIL(V,\DDtilde^{\sankaku},N)$
and $\TNIL(V^{[\ell]},\nabla^{[\ell]},F^{[\ell]})$,
respectively.
We also have the new supersymmetric index
$\nbigq_{\ell}$ of $(U_{\ell},\DDtilde^{\sankaku}_{\ell})$.
By construction, we have the following equality,
under the isomorphism (\ref{eq;08.7.26.3}):
\[
 \nbigq=\bigoplus_{\ell\geq 1}
 \bigl(\nbigq_{\ell}\otimes\id+\id\otimes\nbigq^{[\ell]}\bigr)
\]
The eigenvalues of $\nbigq$ are easily calculable,
once we know those of $\nbigq_{\ell}$.
In particular, we obtain the following.
\begin{cor}
The eigenvalues of $\nbigq_{|q^{-1}(y)}$ are constant
for any $y\in Y$,
where $q:(X_+\cup X_-)\lrarr Y$ denotes the projection.
\hfill\qed
\end{cor}

\subsubsection{Several variable case}
\label{subsubsection;08.7.26.7}

Let $(V,W,\DD^{\sankaku}_V,\vecN,S)$ be
a variation of polarized mixed twistor structures
of split type with weight $0$ in $n$-variables
on $\proj^1\times Y$.
We have the associated variation of
twistor structures
$\TNIL(V,\DD^{\sankaku}_V,\vecN,S)$
with a pairing of weight $0$
on $(\cnum^{\ast})^n\times Y$.
We set $X^{\ast}=\bigl\{
 (z_1,\ldots,z_n)\in\cnum^n\,\big|\,
 0<|z_i|<1 \bigr\}\times Y$.

\begin{prop}
$\TNIL(V,\DD_V^{\sankaku},\vecN,S)_{|
 \proj^1\times X^{\ast}}$
is a twistor nilpotent orbit.
\end{prop}
\pf
For any $\veca\in\real^{n}_{>0}$,
we set $N(\veca):=\sum_{i=1}^n a_i\cdot N_i$.
We obtain 
a variation of mixed polarized twistor structures
$(V,W,\DD^{\sankaku},
 N(\veca),S)$ of split type
with weight $0$ in one variable
on $\proj^1\times Y$.
Applying the result in 
Subsection \ref{subsubsection;08.7.26.5}
to $(V,W,\DD^{\sankaku}_V,N(\veca),S)$,
we obtain the desired property of
$(V,W,\DD^{\sankaku}_V,\vecN,S)$.
\hfill\qed

\begin{df}
An (integrable) twistor nilpotent orbit 
is called of split type,
if it is associated to
(integrable) polarized mixed twistor structures
of split type.
\hfill\qed
\end{df}

If $(V,W,\DD^{\sankaku}_V,\vecN,S)$
is enriched to integrable
$(V,W,\DDtilde^{\sankaku}_V,\vecN,S)$,
the associated twistor nilpotent orbit
is also enriched to integrable
$\TNIL(V,\DDtilde^{\sankaku}_V,\vecN,S)$.
Let us consider its new supersymmetric index $\nbigq$.
For any $\veca\in\real^{n}_{>0}$,
we set $N(\veca):=\sum_{i=1}^n a_i\cdot N_i$.
According to Proposition \ref{prop;08.8.9.21},
there exist variations of polarized pure integrable
twistor structures
$(U_{\veca,\ell},\DDtilde^{\sankaku}_{\veca,\ell})$
for $\ell\geq 1$
such that 
\[
 \bigl(V,\DDtilde^{\sankaku}_V,N(\veca)\bigr)
\simeq \bigoplus_{\ell\geq 1}
 \bigl(U_{\veca,\ell},\DDtilde^{\sankaku}_{\veca,\ell}
\bigr)
\otimes
 \bigl(
 V^{[\ell]},\nabla^{[\ell]},
 F^{[\ell]}\bigr).
\]
\begin{lem}
\label{lem;08.8.9.30}
For any $\veca,\vecb\in\real_{>0}^n$,
we have an isomorphism
$\bigl(U_{\veca,\ell},
 \DDtilde^{\sankaku}_{\veca,\ell}\bigr)
\simeq
 \bigl(U_{\vecb,\ell},
 \DDtilde^{\sankaku}_{\vecb,\ell}\bigr)$.
\end{lem}
\pf
Let $V=\bigoplus V_j$ be the splitting.
For any $\veca\in\real^{n}_{>0}$
and $j\geq 0$,
we set
\[
 (PV_{j,\veca},\DDtilde^{\sankaku}):=\Ker\Bigl(
 N(\veca)^{j+1}:
 (V_j,\DDtilde^{\sankaku})\lrarr 
 (V_{-j-2},\DDtilde^{\sankaku})\otimes\Tate(-j-1)
 \Bigr)
\]
We have only to show that
$(PV_{j,\veca},\DDtilde^{\sankaku})$
and
 $(PV_{j,\vecb},\DDtilde^{\sankaku})$ 
are isomorphic,
if $\vecb$ is sufficiently close to $\veca$.

We  set
$ (Y_{j,\veca},\DDtilde^{\sankaku})
=\Image\Bigl(
 N(\veca):
 (V_{j+2},\DDtilde^{\sankaku})\otimes\Tate(1)
 \lrarr (V_j,\DDtilde^{\sankaku})
 \Bigr)$.
Then, we obtain the flat splittings
$(V_j,\DDtilde^{\sankaku})=
 (PV_{j,\veca},\DDtilde^{\sankaku})
\oplus
 (Y_{j,\veca},\DDtilde^{\sankaku})$.
If $\vecb$ is sufficiently close to $\veca$,
flat isomorphisms
$PV_{j,\veca}\lrarr PV_{j,\vecb}$ are
induced by inclusions and projections.
Thus, we are done.
\hfill\qed

\vspace{.1in}
By Lemma \ref{lem;08.8.9.30} and
the result in Subsection
\ref{subsubsection;08.7.26.6},
the eigenvalues of 
$\nbigq$ are easily calculable
once we know the new supersymmetric indices
of $(U_{\veca,\ell},\DDtilde^{\sankaku}_{\veca,\ell})$
for $\veca\in\real_{>0}^{\ell}$
and $\ell\geq 1$.
In particular, we obtain the following.
\begin{cor}
The eigenvalues of $\nbigq_{|q^{-1}(y)}$ 
are constant for any $y\in Y$,
where $q:X^{\ast}\lrarr Y$ denotes 
the natural projection.
\hfill\qed
\end{cor}

\section{Integrable twistor nilpotent orbit}
\label{section;08.8.9.2}

\subsection{Statements}

\subsubsection{Twistor nilpotent orbits
 and polarized mixed twistor structures}

Let $Y$ be a complex manifold.
Let $(V,\DD^{\sankaku}_V)$ be a 
variation of $\proj^1$-holomorphic vector bundles
on $\proj^1\times Y$
equipped with the following 
$\proj^1$-holomorphic 
$\DD_V^{\sankaku}$-flat data:
\begin{itemize}
\item
 A $(-1)^w$-symmetric pairing
 $S: V\otimes\sigma^{\ast}V
\lrarr \Tate(-w)$.
\item
 A tuple $\vecN$ of nilpotent morphisms
 $N_j:V\lrarr V\otimes\Tate(-1)$
 $(j=1,\ldots,n)$,
 which are mutually commutative.
\item
 $S(N_j\otimes\id)
+S(\id\otimes\sigma^{\ast}N_j)=0$
 for $j=1,\ldots,n$.
\end{itemize}
For simplicity of the statement,
we assume the following:
\begin{itemize}
\item $Y$ is contained in another complex manifold
 $Y'$ as a relatively compact subset,
 and $(V,\DD^{\sankaku}_V,S,\vecN)$
 is extended on $Y'$.
\end{itemize}
We set $X^{\ast}(R):=
 \bigl\{
 (z_1,\ldots,z_n)\,\big|\,
 |z_i|<R \bigr\}\times Y$.

\begin{thm}
\label{thm;08.7.26.8}
$(V,\DD^{\sankaku}_V,\vecN,S)$
is a variation of polarized mixed twistor structures
with weight $w$ in $n$-variables,
if and only if
$\TNIL(V,\DD^{\sankaku}_V,\vecN,S)_{|
 \proj^1\times X^{\ast}(R)}$
is a twistor nilpotent orbit with weight $w$
for some $R>0$.
\end{thm}

Note that the ``if'' part follows from 
Theorem 12.22 of \cite{mochi2}.
The ``only if'' part immediately follows from
Proposition \ref{prop;08.7.29.3} below
and a result in Subsection 2.8 of \cite{mochi7}.
(We apply Proposition \ref{prop;08.7.29.3}
 to each point of $Y'$.)
The one dimensional case was 
proved in Proposition 3.105 of \cite{mochi2}.
Such an equivalence for Hodge structure
was established by 
Cattani-Kaplan-Schmid and Kashiwara-Kawai.

\begin{cor}
\label{cor;08.7.26.9}
Let $(V,\DD^{\sankaku}_V,\vecN,S)$ be 
as above.
\begin{itemize}
\item
 Assume that $(V,\DD^{\sankaku}_V)$ 
 is enriched to integrable
 $(V,\DDtilde^{\sankaku}_V)$
 such that 
 $\vecN$ and $S$ are flat with respect to 
 $\DDtilde^{\sankaku}_V$.
  Then, 
 $(V,\DDtilde^{\sankaku}_V,\vecN,S)$ 
 is a variation of polarized mixed integrable twistor structures
 with weight $w$ in $n$ variables,
 if and only if 
 $\TNIL(V,\DDtilde^{\sankaku}_V,\vecN,S)_{|
 \proj^1\times X^{\ast}(R)}$
 is integrable twistor nilpotent orbit  for some $R>0$.
\item
Assume moreover that
$(V,\DDtilde^{\sankaku}_V,S)$ 
is equipped with a real structure $\kappa$
which is compatible with $\vecN$.
Then, 
$(V,\DDtilde^{\sankaku}_V,\vecN,S,\kappa,-w)$
is a variation of polarized mixed twistor-TERP structures
 if and only if 
 the associated 
$\TNIL(V,\DDtilde^{\sankaku}_V,\vecN,\nbigs,\kappa,-w)$
is a twistor-TERP nilpotent orbit 
on $X^{\ast}(R)$ for some $R>0$.
\hfill\qed
\end{itemize}
\end{cor}

\begin{rem}
As the one variable case of 
Corollary {\rm \ref{cor;08.7.26.9}},
we obtain the correspondence between
twistor-TERP nilpotent orbits and
polarized mixed twistor-TERP structures.
This is different from 
the correspondence
between mixed TERP structure and HS-orbit
in the regular singular case
established by Hertling and Sevenheck
{(\rm \cite{Hertling} and
 \cite{Hertling-Sevenheck})}.
\hfill\qed
\end{rem}

\subsubsection{Construction of 
 an approximating $C^{\infty}$-isomorphism}
\label{subsubsection;08.7.29.2}

Let $(V,W,\DD^{\sankaku},\vecN,S)$ 
be a variation of
polarized mixed twistor structures of weight $0$
in $n$-variables on $\proj^1\times Y$.
As explained in 
Subsection \ref{subsubsection;08.8.11.2},
we obtain a variation of polarized mixed twistor structure
of split type
$\bigl(
 V^{(0)},W^{(0)},\DD^{(0)\sankaku},
 \vecN^{(0)}, S^{(0)}
 \bigr)$
by taking Gr with respect to the weight filtration.
We obtain the families of $\proj^1$-holomorphic
vector bundles
$(\nbigv^{\sankaku},\DD^{\sankaku})
 :=\TNIL(V,\DD^{\sankaku},\vecN)$
and
$(\nbigv^{(0)\sankaku},\DD^{(0)\sankaku})
:=\TNIL(V^{(0)},\DD^{(0)\sankaku},\vecN^{(0)})$
on $(\cnum^{\ast})^n\times Y$.
They are equipped with the induced pairings
$\nbigs$ and $\nbigs^{(0)}$.
By the result in
Subsection \ref{subsubsection;08.7.26.7},
$(\nbigv^{(0)},\DD^{(0)\sankaku},
 \nbigs^{(0)})$ is 
a variation of polarized pure twistor structure
on $\proj^1\times X^{\ast}(1)$.
Let $h^{(0)}$ be the corresponding
pluri-harmonic metric.

We take a $C^{\infty}$-isomorphism
$\Phi:V^{(0)}\lrarr V$
as in Lemma \ref{lem;08.7.29.1},
i.e.,
it satisfies 
(i) $\Phi\circ N_i^{(0)}=N_i\circ\Phi$
 for $i=1,\ldots,n$,
(ii) $\Phi$ preserves the weight filtration $W$,
 and
 $\Gr^W\Phi$ is the identity of 
 $\Gr^W(V^{(0)})=\Gr^W(V)$.
By the property (i) for $\Phi$
and the construction of 
$\nbigv^{\sankaku}$ and 
$\nbigv^{(0)\sankaku}$,
we obtain a naturally induced
$C^{\infty}$-isomorphism
$\Phitilde:\nbigv^{(0)\sankaku}
\lrarr
 \nbigv^{\sankaku}$.

Let $\delbar_{\nbigv^{\sankaku},\proj^1}$ denote
the $\proj^1$-holomorphic structure of $\nbigv^{\sankaku}$.
We use the symbol
$\delbar_{\nbigv^{(0)\sankaku},\proj^1}$
in a similar meaning.
We obtain the following $C^{\infty}$-section
of $\End(\nbigv^{(0)\sankaku})\otimes 
\Omega_{\proj^1}^{0,1}$ 
on $\proj^1\times X^{\ast}(1)$:
\[
 F:=
 \delbar_{\nbigv^{(0)\sankaku},\proj^1}
-\Phitilde^{\ast}
 \bigl(\delbar_{\nbigv^{\sankaku},\proj^1}\bigr)
\]
We also obtain the following $C^{\infty}$-morphism:
\[
 G:=\nbigs^{(0)}-\Phitilde^{\ast}\nbigs:
 \nbigv^{(0)\sankaku}
\otimes
 \sigma^{\ast}\nbigv^{(0)\sankaku}
\lrarr
 \Tate(0)
\]
We fix a K\"ahler metric $g$ of $\proj^1$.
Although
the following proposition looks rather auxiliary,
it means that
$(\nbigv^{(0)\sankaku},
 \DD^{(0)\sankaku},\nbigs^{(0)})$
approximates
$(\nbigv^{\sankaku},
 \DD^{\sankaku},\nbigs)$
via $\Phitilde$
around $\proj^1\times\{0\}\times Y$.
We will prove it in Subsection 
\ref{subsubsection;08.7.29.4}.

\begin{prop}
\label{prop;08.7.29.3}
For any $P\in Y$,
there exist a positive constant $R_P>0$
and a neighbourhood $U_P$ of $P$ in $Y$
such that the following estimate holds
$\proj^1\times
 \bigl\{(z_1,\ldots,z_n)\,\big|\,
 0<|z_j|<R_P
 \bigr\}\times U_P$:
\[
 \bigl|F\bigr|_{h^{(0)},g}
=
 O\Bigl(
 \sum_{j=1}^n
\bigl(
 -\log|z_j|
\bigr)^{-1/2}
 \Bigr)
\]
\[
 \bigl|G\bigr|_{h^{(0)}}
=
 O\Bigl(
 \sum_{j=1}^n
 \bigl(-\log|z_j|\bigr)^{-1/2}
 \Bigr),
 \quad
 \bigl|\delbar_{\nbigv^{(0)\sankaku},\proj^1} 
 G\bigr|_{h^{(0)},g}
=O\Bigl(
 \sum_{j=1}^n
 \bigl(-\log|z_j|\bigr)^{-1/2}
 \Bigr)
\]
\end{prop}

\subsubsection{Estimate of 
the new supersymmetric index}

Assume that $(V,\DD^{\sankaku}_V,\vecN,S)$
is enriched to integrable 
$(V,\DDtilde_V^{\sankaku},\vecN,S)$.
By taking Gr with respect to the weight filtration,
we obtain a polarized mixed integrable twistor
structure of split type
$\bigl(V^{(0)},W^{(0)},
 \DD^{(0)\sankaku}_V,
 \vecN^{(0)},S^{(0)}\bigr)$.
Let $(\nbigv,\DDtilde^{\sankaku},\nbigs)
 =\TNIL(V,\DDtilde^{\sankaku}_V,
 \vecN,S)_{|\proj^1\times X^{\ast}(R)}$
and 
$(\nbigv^{(0)},\DDtilde^{(0)\sankaku},
 \nbigs^{(0)})
=\TNIL(V^{(0)},\DDtilde_{V}^{(0)},\vecN^{(0)},S^{(0)})
 _{|\proj^1\times X^{\ast}(R)}$
be the associated nilpotent orbits
(Corollary \ref{cor;08.7.26.9}).
Let $\nbigq$ and $h$
(resp. $h^{(0)}$ and $\nbigq^{(0)}$)
denote the new supersymmetric index
and the pluri-harmonic metric
of $(\nbigv,\DDtilde^{\sankaku},\nbigs)$
(resp. $(\nbigv^{(0)},\DDtilde^{(0)\sankaku},
 \nbigs^{(0)})$).
We will prove the following proposition
in Subsection \ref{subsubsection;08.8.18.20}.

\begin{prop}
\label{prop;08.7.29.5}
Let $\Phitilde:\nbigv^{(0)}\lrarr \nbigv$
be a $C^{\infty}$-isomorphism 
constructed in Subsection
{\rm\ref{subsubsection;08.7.29.2}}.
For any $P\in Y$,
there exist $R>0$
and a neighbourhood $U_P$ of $P$ in $Y$
such that the following estimate holds
with respect to $h^{(0)}$
on $\proj^1\times
 \bigl\{(z_1,\ldots,z_n)\,\big|\,
 0<|z_j|<R
 \bigr\}\times U_P$:
\[
 \Phitilde^{\ast}h-h^{(0)}
=O\Bigl(
 \sum_{i=1}^n \bigl(-\log|z_i|\bigr)^{-1/2}
 \Bigr),
\quad
 \Phitilde^{\ast}\nbigq-\nbigq^{(0)}
=O\Bigl(
 \sum_{i=1}^n \bigl(-\log|z_i|\bigr)^{-1/2}
 \Bigr)
\]
In particular,
the eigenvalues of
$\nbigq_{|q^{-1}(y)}$ are constant up to
$O\Bigl(
 \sum_{i=1}^n \bigl(-\log|z_i|\bigr)^{-\delta}
 \Bigr)$ for some $\delta>0$,
where $q:X^{\ast}(1)\lrarr Y$
denotes the natural projection.
\end{prop}

\subsection{Proof}

\subsubsection{Proof of Proposition 
 \ref{prop;08.7.29.3}}
\label{subsubsection;08.7.29.4}

Let $C>0$.
Fix $P\in Y$.
In the following,
we will shrink $Y$ 
instead of taking a neighbourhood $U_P$,
for simplicity of description.
We set
\[
 Z(C):=
 \Bigl\{
 (z_1,\ldots,z_n)\in (\cnum^{\ast})^n\,\Big|\,
 \bigl|z_{i}\bigr|^C
\leq
 \bigl|z_{i+1}\bigr|<1,\,\,\,
 i=1,\ldots,n-1
 \Bigr\}\times Y.
\]
It is easy to observe that
we have only to estimate $F$,
$G$ and $\delbar_{\nbigv^{(0)\sankaku},\proj^1}G$
on $\proj^1\times Z(C)$.
For $m=1,\ldots,n$,
we put $N^{(0)}(\mbar):=\sum_{i\leq m}N_i^{(0)}$.
Let $W(\mbar)$ denote the weight filtration
of $V^{(0)}$ induced by $N^{(0)}(\mbar)$.
Recall that the filtrations
$W(\itibar),
 W(\nibar),
 \ldots,
 W(\nbar)$
are compatible
(Lemma 3.116 of \cite{mochi2}).

We take a compact region
$\nbigu\subset\cnum_{\lambda}$
such that 
the union of the interior parts of
$\nbigu$ and $\sigma(\nbigu)$ cover $\proj^1$.
Let $\vecv=(v_i)$ be a frame of
$V^{(0)}_{|\nbigu\times Y}$ 
compatible with
$W(\itibar),
 W(\nibar),
 \ldots,
 W(\nbar)$.
For $m=1,\ldots,n$,
we set 
\[
 k_m(v_i):=\frac{1}{2}
 \deg^{W(\mbar)}(v_i).
\]
We formally put
$k_0(v_i)=0$.

\begin{lem}
\label{lem;08.7.28.120}
Let $A$ be determined by
$\bigl(
 \Phi^{-1}\circ\delbar\Phi
\bigr)
 \vecv=\vecv\cdot A$.
Then, $A_{i,j}=0$ unless
$k_m(v_i)\leq k_m(v_j)$ $(m=1,\ldots,n-1)$
and $k_n(v_i)<k_n(v_j)$.
\end{lem}
\pf
Because of our choice of $\Phi$,
it preserves the filtrations 
$W(\mbar)$ $(m=1,\ldots,n)$,
and $\Gr^{W(\nbar)}\Phi$ is holomorphic.
Then, the claim of Lemma 
\ref{lem;08.7.28.120} immediately follows.
\hfill\qed

\vspace{.1in}

Let $q_0:
 \cnum_{\lambda}\times X^{\ast}(1)
 \lrarr \cnum_{\lambda}\times Y$ be the projection.
Recall
$\nbigv^{(0)\sankaku}_{|
 \cnum_{\lambda}\times X^{\ast}(1)}
 = q_0^{\ast}V_0$.
Let $\vtilde_i$ be the section of
$\nbigv^{(0)\sankaku}_{|
 \nbigu\times X^{\ast}(1)}$
induced by $v_i$,
and we put
\[
  v_j':=\vtilde_j\cdot
 \prod_{m=1}^n
 \bigl(-\log|z_m|\bigr)^{-k_m(v_j)+k_{m-1}(v_j)}
=\vtilde_j\cdot 
 \prod_{m=1}^{n-1}
 \left(
 \frac{-\log|z_m|}{-\log|z_{m+1}|}
 \right)^{-k_m(v_j)}
\cdot
 \bigl(-\log|z_{n}|\bigr)^{-k_n(v_j)}.
\]
Due to the norm estimate
for tame harmonic bundles
(Theorem 13.25 of \cite{mochi2}),
the $C^{\infty}$-frame
$\vecv'=(v_j')$ is adapted to 
the metric $h^{(0)}$
on $Z(C)$,
i.e., the hermitian matrix-valued functions
$H=\bigl(h(v_i',v_j')\bigr)$ and $H^{-1}$
are bounded on $Z(C)$.
Let $A'$ be the matrix-valued function
determined by
$F\vecv'=\vecv'\cdot A'$.
Then, we have
\[
 A'_{i,j}=
 A_{i,j}\cdot 
 \prod_{m=1}^{n-1}
 \left(\frac{-\log|z_m|}{-\log|z_{m+1}|}
 \right)^{k_m(v_i)-k_m(v_j)}
 \cdot\bigl(-\log|z_n|
 \bigr)^{k_n(v_i)-k_n(v_j)}.
\]
Hence, we obtain
$A_{i,j}'=O\Bigl(
 \bigl(-\log|z_n|\bigr)^{-1/2}
 \Bigr)$.
It implies the desired estimate for $F$
on $\nbigu\times Z(C)$.
Similarly, we obtain the estimate
on $\sigma(\nbigu)\times Z(C)$,
and thus on $\proj^1\times Z(C)$.

\vspace{.1in}

Let $\vecw$ be a frame of
$V_{|\sigma(\nbigu)\times Y^{\dagger}}$
compatible with the filtrations
$W(\itibar),
W(\nibar),\ldots,
 W(\nbar)$.
For $m=1,\ldots,n$,
we set
\[
 k_m(w_i):=\frac{1}{2}
 \deg^{W(\mbar)}(w_i).
\]
We formally put $k_0(w_i)=0$.
We set 
$G_0:=S^{(0)}-\Phi^{\ast}S:
 V^{(0)}\otimes\sigma^{\ast}V^{(0)}
\lrarr \Tate(0)$.

\begin{lem}
\label{lem;08.7.28.31}
$G_0(v_i,\sigma^{\ast}w_j)=0$
unless the following holds:
\[
 k_m(v_i)+k_m(w_j)\geq 0
\quad
 (m=1,\ldots,n-1),
\quad
 k_n(v_i)+k_n(w_j)>0.
\]
\end{lem}
\pf
By the relation
$S(N_i\otimes \id)+
 S\bigl(\id\otimes\sigma^{\ast}(N_i)\bigr)=0$,
we have
$S\bigl(W_p(\mbar)\otimes 
 \sigma^{\ast}W_q(\mbar)\bigr)=0$
unless $p+q\geq 0$.
We have similar vanishings
for $S^{(0)}$.
Note that $\Phi$ preserves the filtrations
$W(\mbar)$ for $m=1,\ldots,n$,
and $\Gr^{W(\nbar)}\Phi$
is compatible with
$S$ and $S^{(0)}$.
Thus, we obtain the claim of
Lemma \ref{lem;08.7.28.31}.
\hfill\qed

\vspace{.1in}

Let $q_{\infty}:
 \cnum_{\mu}\times X^{\ast}(1)^{\dagger}
 \lrarr
 \cnum_{\mu}\times Y^{\dagger}$
 be the projection.
Recall 
$\nbigv^{(0)\sankaku}_{
 |\cnum_{\mu}\times X^{\ast}(1)^{\dagger}}
=q_{\infty}^{\ast}V_{\infty}$.
Let $\wtilde_j$ be the section of
$\nbigv^{(0)\sankaku}_{
 |\sigma(\nbigu)\times X^{\ast}(1)^{\dagger}}$
induced by $w_j$,
and we put
\[
  w_j':=\wtilde_j\cdot
 \prod_{m=1}^{n-1}
 \left(\frac{-\log|z_m|}{-\log|z_{m+1}|}\right)^{-k_m(w_j)}
\cdot 
  \bigl(-\log|z_{n}|\bigr)^{-k_n(w_j)}.
\]
Note the following:
\[
 G(v_i',\sigma^{\ast}w_j')=
 G_0(v_i,\sigma^{\ast}w_j)
\times
 \prod_{m=1}^{n-1}
 \left(\frac{-\log|z_m|}{-\log|z_{m+1}|}
 \right)^{-k_m(v_i)-k_m(w_j)}
 \cdot\bigl(-\log|z_n| \bigr)^{-k_n(v_i)-k_n(v_j)}
\]
Hence, we obtain 
$ \bigl| G \bigr|_{h^{(0)}}
=O\Bigl(
 \bigl(-\log|z_n|\bigr)^{-1/2}
 \Bigr)$.
Similarly, we obtain the estimate for
$\bigl|\delbar_{\nbigv^{(0)},\proj^1}G\bigr|$.
Thus, the proof of Proposition 
\ref{prop;08.7.29.3} is finished.
The proof of Theorem \ref{thm;08.7.26.8}
is also finished.

\subsubsection{Proof of 
Proposition \ref{prop;08.7.29.5}}
\label{subsubsection;08.8.18.20}

We have the decompositions
$\DDtilde^{\sankaku}
=\DD^{\sankaku}_{\nbigv_0}
+\nabla_{\lambda}$
and 
$\DDtilde^{(0)\sankaku}
=\DD^{(0)\sankaku}_{\nbigv_0}
+\nabla_{\lambda}^{(0)}$.
By an argument used in the proof of
Proposition \ref{prop;08.7.29.3},
we obtain the following estimate
with respect to $h^{(0)}$:
\[
 \Phitilde^{\ast}\nabla_{\lambda}
-\nabla_{\lambda}^{(0)}
=O\Bigl(
 \sum_{i=1}^n \bigl(-\log|z_i|\bigr)^{-1/2}
 \Bigr)
\]
Then, Proposition \ref{prop;08.7.29.5}
follows from Lemma \ref{lem;08.7.26.11}
with Proposition \ref{prop;08.7.29.3}.
\hfill\qed

\section{Family of 
meromorphic $\lambda$-flat bundles}
\label{section;08.8.20.40}

We will review some results on
family of meromorphic $\lambda$-flat bundles
mainly explained in Sections 7 and 8 of \cite{mochi7}.
See also \cite{majima} and \cite{sabbah4}
for the earlier works on asymptotic analysis
of meromorphic flat bundles.

\subsection{Good lattice in the level $\vecm$}

\subsubsection{Preliminary}

\paragraph{Good set of irregular values
 in the level $\vecm$}

Let $\Delta^{\ell}:=
\bigl\{(z_1,\ldots,z_{\ell})\,\big|\,
 |z_i|<1,\,i=1,\ldots,\ell \bigr\}$
denote the $\ell$-dimensional 
multi-disc.
Let $X:=\Delta^{\ell}\times Y$
for some complex manifold $Y$.
Let $D_i:=\{z_i=0\}$ and 
$D:=\bigcup_{i=1}^{\ell}D_i$
be hypersurfaces of $X$.
Let $M(X,D)$ (resp. $H(X)$)
denote the space of meromorphic (resp. holomorphic)
functions on $X$ whose poles are contained in $D$.
For $\vecm=(m_1,\ldots,m_{\ell})\in\seisuu^{\ell}$,
we put $\vecz^{\vecm}:=\prod_{i=1}^{\ell}z_i^{m_i}$.

Let $\vecm\in\seisuu_{\leq 0}^{\ell}-\{\veczero\}$.
A finite set $\nbigi$ of meromorphic functions
 $\bigl\{\gminia=
 \gminia_{\vecm}\cdot\vecz^{\vecm}\bigr\}
 \subset M(X,D)$
is called 
a good set of irregular values on $(X,D)$
in the level $\vecm$,
if the following holds:
\begin{itemize}
\item
$\gminia_{\vecm}$ are 
 holomorphic functions on $X$.
\item
 $\gminia_{\vecm}-\gminib_{\vecm}$
 are nowhere vanishing holomorphic functions on $X$
 for any two distinct $\gminia,\gminib\in \nbigi$.
\end{itemize}
Let $i(0)$ be the integer such that
$m_{i(0)}<0$.
If moreover the following condition holds,
$\nbigi$ is called a good set of irregular values
on $(X,D)$
in the level $(\vecm,i(0))$.
\begin{itemize}
\item
 $\gminia_{\vecm}$ 
 are independent of the variable $z_{i(0)}$
 for any $\gminia\in\nbigi$.
\end{itemize}

\begin{rem}
The first condition is not essential.
If we do not impose it,
the third condition should be replaced with
that $\gminia_{\vecm}-\gminib_{\vecm}$
are independent of $z_{i(0)}$
for any $\gminia,\gminib\in\nbigi$.
\hfill\qed
\end{rem}

\paragraph{Multi-sectors and
orders on good sets of irregular 
values in the level $\vecm$}

Let $X:=\Delta^{\ell}\times Y$
for some complex manifold $Y$.
Let $D_i:=\{z_i=0\}$ and 
$D:=\bigcup_{i=1}^{\ell}D_i$
be hypersurfaces of $X$.
Let $\nbigk$ be a region of $\cnum_{\lambda}$
or a point in $\cnum_{\lambda}^{\ast}$.
(For Definition \ref{df;08.11.6.10},
 we may admit $\nbigk=\{0\}$.
 Since we do not have to consider Stokes
 structure in this case,
 we exclude it in the following.)
The product $\nbigk\times X$
is expressed by $\nbigx$.
We use the symbols like
$\nbigy$ and $\nbigd$ 
in similar meanings.
We put $W:=\nbigd\cup \bigl(\{0\}\times X\bigr)$
in the case $0\in \nbigk$,
and $W:=\nbigd$ otherwise.
Let $\pi:\nbigxtilde(W)\lrarr \nbigx$ 
denote the real blow up
of $\nbigx$ along $W$.

In this paper,
a sector of a punctured disc $\Delta^{\ast}$
means a subset of the form
$ \bigl\{
 z\,\big|\,
 0<|z|<R,\,\,
 \theta_{0}\leq \arg(z)\leq \theta_1
 \bigr\}$
for some $\theta_0<\theta_1$.
It may be standard to admit
the case $|\theta_1-\theta_0|\geq 2\pi$,
but we do not care about it.

By a ``multi-sector of $\nbigx-W$'',
we mean a subset of the following form
\[
 U\times
 \prod_{i=1}^{\ell} S_i
\times V,
\quad\mbox{\rm or }\,\,
 S_{\lambda}
\times
\prod_{i=1}^{\ell}
 S_i
\times V.
\]
\begin{itemize}
\item
 $U$ denotes a 
 compact region in $\nbigk$.
 (If $\nbigk$ is a point, $U=\nbigk$.)
\item
 $S_{\lambda}$ denotes 
 a sector of $\nbigk-\{0\}$.
 (If $0\not\in\nbigk$, we do not consider 
 the subsets of the second type.)
\item
 $S_i$ denote sectors of $\Delta_{z_i}^{\ast}$.
\item
 $V$ denotes a compact region in $Y$.
\end{itemize}
For a multi-sector $S$,
let $\Sbar$ denote the closure of $S$
in $\nbigxtilde(W)$.

\begin{notation}
Let $\Multisector(\nbigx-W)$ denote 
the set of multi-sectors in $\nbigxtilde(W)$.
For any point $P\in\nbigxtilde(W)$,
let $\Multisector(P,\nbigx-W)$ denote the set of
multi-sectors $S$ such that
$P$ is contained in the interior part of $\Sbar$.
\hfill\qed
\end{notation}

Let $\nbigi$ be 
a good set of irregular values on $(\nbigx,\nbigd)$
in the level $\vecm$.
We put $F_{\gminia,\gminib}:=
-\Re\bigl(\lambda^{-1}\cdot
 (\gminia-\gminib)
 \bigr)\cdot|\lambda|\cdot |\vecz^{-\vecm}|$
for any distinct $\gminia,\gminib\in \nbigi$.
They determine the $C^{\infty}$-functions
on $\nbigxtilde(W)$.

\begin{notation}
Let $A$ be any subset of $\nbigxtilde(W)$.
We say 
$\gminia<_A\gminib$ for 
$(\gminia,\gminib)\in \nbigi^2$
if $F_{\gminia,\gminib}(Q)<0$ for any $Q\in A$.
We say
$\gminia\leq_A\gminib$ for 
$(\gminia,\gminib)\in \nbigi^2$
if either $\gminia<_A\gminib$ or $\gminia=\gminib$ holds.
The relation $\leq_A$ gives the partial order of $\nbigi$.

We use the symbol
$\leq_P$ in the case $A=\{P\}$.
For a multi-sector $S$,
we prefer the symbol $\leq_S$ to $\leq_{\Sbar}$.
We also use 
$\leq^{\lambda}_S$  and $\leq^{\lambda}_P$
when we emphasize the twist by $\lambda^{-1}$.
\hfill\qed
\end{notation}

For any point $P\in \pi^{-1}(W)$,
there exists 
$S_P\in \Multisector(P,\nbigx-W)$
such that 
the relations $\leq_P$ and $\leq_{S_P}$ coincide.
Let $\Multisector(P,\nbigx-W,\nbigi)$
denote the set of such $S_P$.
(The definitions of $\Multisector(P,\nbigx-W,\nbigi)$
 is slightly different from
 that in {\rm\cite{mochi7}}.)

\subsubsection{Good lattice in the level $\vecm$}
\label{subsubsection;08.8.3.2}

Let $Y$ be a complex manifold
with a simple normal crossing divisor $D_Y'$.
Let $X:=\Delta_z^k\times Y$,
$D_{z,i}:=\{z_i=0\}$
and $D_z:=\bigcup_{i=1}^{k}D_{z,i}$.
We also put $D_Y:=\Delta_z^k\times D_Y'$
and $D:=D_z\cup D_Y$.
Let $\nbigk$ be a point of 
$\cnum_{\lambda}^{\ast}$
or a compact region in $\cnum_{\lambda}$.
We put $\nbigx:=\nbigk\times X$.
We use the symbols $\nbigy$,
$\nbigd_z$, $\nbigd$ in similar meanings.
Let $p_{\lambda}$ denote the projection
forgetting the $\nbigk$-component.
The completion of $\nbigx$ along $\nbigd_z$
is denoted by $\nbigdhat_z$.
(See \cite{banica}, \cite{bingener}
and \cite{krasnov} for completion
of complex analytic spaces.)
We use the symbol $\nbigdhat$
in a similar meaning.
Let $d_X$ denote the restriction of
the exterior derivative to the $X$-direction.

Let $E$ be a locally free $\nbigo_{\nbigx}$-module
with a family of meromorphic flat $\lambda$-connections
$\DD:E\lrarr E\otimes 
 p_{\lambda}^{\ast}\Omega^{1}_{X}(\ast D)$.
Let $\vecm\in\seisuu_{<0}^k$
and $i(0)\in[1,k]=\{1,\ldots,k\}$.
We put $\vecm(1):=\vecm+\vecdelta_{i(0)}$.

\begin{df}
\label{df;08.11.6.10}
We say that $(E,\DD)$ is an unramifiedly good lattice
of a family of meromorphic $\lambda$-flat bundles
in the level $(\vecm,i(0))$,
if there exists a good set of irregular values
$\nbigi$ in the level $(\vecm,i(0))$
on $(\nbigx,\nbigd_z)$,
and a decomposition
\begin{equation}
 \label{eq;07.12.15.1}
 (E,\DD)_{|\nbigdhat_z}
=\bigoplus_{\gminia\in\nbigi}
 (\Ehat_{\gminia},\DDhat_{\gminia})
\end{equation}
with $\ord(\DDhat_{\gminia}-d_X\gminia)\geq\vecm(1)$
in the sense
$(\DDhat_{\gminia}-d_X\gminia)\Ehat_{\gminia}
\subset
 \vecz^{\vecm(1)}\cdot\Ehat_{\gminia}\otimes
 p_{\lambda}^{\ast}\Omega^1_{X}\bigl(\log D\bigr)$.

The decomposition {\rm(\ref{eq;07.12.15.1})}
is called the irregular decomposition
in the level $(\vecm,i(0))$,
(or simply $\vecm$).
We also often say that $(E,\DD)$ is a good lattice
in the level $(\vecm,i(0))$ for simplicity.
\hfill\qed
\end{df}

In the case $0\in \nbigk$,
we put $\nbigx^0:=\{0\}\times X$
and $\nbigd_z^0:=\{0\}\times D_z$.
By shrinking $X$,
we obtain the irregular decomposition
$(E,\DD)_{|\nbigx^0}
=\bigoplus_{\gminia\in\nbigi}
 (E_{\gminia,\nbigx^0},\DD^0_{\gminia})$
whose completion along
$\nbigd_z^0$ is equal to
the one induced by  (\ref{eq;07.12.15.1}).
It is uniquely extended to the $\DD$-flat decomposition
on the completion $\nbigxhat^0$ of $\nbigx$
along $\nbigx^0$:
\[
 (E,\DD)_{|\nbigxhat^0}
=\bigoplus_{\gminia\in\nbigi}
 (\Ehat_{\gminia,\nbigxhat^0},
 \DDhat_{\gminia})
\]
We put $W:=\nbigx^0\cup\nbigd_z$.
Let $\What$ denote the completion along $W$.
We obtain the decomposition:
\begin{equation}
\label{eq;07.12.15.2}
 (E,\DD)_{|\What}
=\bigoplus_{\gminia\in \nbigi}
 (\Ehat_{\gminia,\What},\DDhat_{\gminia})
\end{equation}
The decomposition (\ref{eq;07.12.15.2})
is also called the irregular decomposition
in the level $(\vecm,i(0))$
if $0\in\nbigk$.

In the following,
we formally set $W:=\nbigd_z$
if $0\not\in\nbigk$.
Let $\pi:\nbigxtilde(W)\lrarr\nbigx$ denote 
the real blow up of $\nbigx$ along $W$.
Let $O_z$ be the origin of $\Delta_z^k$,
and we put $\gbigz:=\pi^{-1}(O_z\times \nbigy)$.
We consider the case
that $Y=\Delta^{n}_{\zeta}$
and 
$D_{Y}':=\bigcup_{j=1}^{\ell}D_{\zeta,j}$,
where 
$D_{\zeta,j}:=\{\zeta_j=0\}$.
The restriction of $\DD$
to the $\Delta_z^{k}$-direction
is denoted by $\DD_z$.

\paragraph{Stokes structure in the level $\vecm$}

For any multi-sector $S$ in $\nbigx-W$,
let $\Sbar$ denote the closure of $S$
in $\nbigxtilde(W)$,
and let $Z$ denote $\Sbar\cap \pi^{-1}(W)$.
The irregular decomposition (\ref{eq;07.12.15.2})
on $\What$ induces the decomposition on $\Zhat$:
\begin{equation}
 \label{eq;07.12.15.50}
 (E,\DD)_{|\Zhat}=
 \bigoplus_{\gminia\in\nbigi}
 (\Ehat_{\gminia},\DDhat_{\gminia})_{|\Zhat}
\end{equation}
We put
$\nbigf^Z_{\gminia}:=
 \bigoplus_{\gminib\leq_S\gminia}
 \Ehat_{\gminib|\Zhat}$,
and then we obtain the filtration
$\nbigf^Z$ of $E_{|\Zhat}$ indexed by
$\bigl(\nbigi,\leq_S\bigr)$.
We can show the following proposition.
(See Subsections 7.2.1 and 8.1.1 of \cite{mochi7}.)
\begin{prop}
\label{prop;07.6.16.6}
\mbox{{}}
For any point $P\in \gbigz$,
there exists $S\in \Multisector(P,\nbigx-W,\nbigi)$
such that the following holds:
\begin{itemize}
\item
There exists the unique $\DD$-flat filtration
$\nbigf^S$ of $E_{|\Sbar}$
indexed by $\bigl(\nbigi,\leq_S\bigr)$
such that
$\nbigf^S_{|\Zhat}=\nbigf^Z$.
Moreover, 
if a $\DD_{z}$-flat filtration $\nbigf^{\prime\,S}$
of $E_{|\Sbar}$ indexed by $(\nbigi,\leq_S)$
satisfies $\nbigf^{\prime\,S}_{|\Zhat}=\nbigf^Z$,
then $\nbigf^{\prime\,S}=\nbigf^S$.
\item
There exists a $\DD_{z}$-flat splitting of
$\nbigf^{S}$ on $\Sbar$.
Note that if we take such a splitting,
the restriction to $\Zhat$ is 
the same as {\rm(\ref{eq;07.12.15.50})}.
\end{itemize}
We call $\nbigf^S$ the Stokes filtration of $(E,\DD)$
in the level $\vecm$.
\hfill\qed
\end{prop}

\begin{notation}
For any $P\in \gbigz$,
let $\Multisector^{\ast}(P,\nbigx-W,\nbigi)$
denote the set of
$S\in \Multisector(P,\nbigx-W,\nbigi)$
as in Proposition 
{\rm\ref{prop;07.6.16.6}}.
Let $\Multisector^{\ast}(\nbigx-W,\nbigi)$
denote the union of
$\Multisector^{\ast}(P,\nbigx-W,\nbigi)$
for $P\in \gbigz$.
\hfill\qed
\end{notation}

The following lemma is clear.
\begin{lem}
Let $S,S'\in \Multisector(P,\nbigx-W,\nbigi)$.
Assume 
(i) $S'\subset S$,
(ii) $S\in\Multisector^{\ast}(P,\nbigx-W,\nbigi)$.
Then, $S'\in \Multisector^{\ast}(P,\nbigx-W,\nbigi)$.
The filtration $\nbigf^{S'}$ is the restriction of
$\nbigf^S$.
\hfill\qed
\end{lem}

\paragraph{Compatibility of the Stokes filtrations}

Let $S,S'\in \Multisector^{\ast}(\nbigx-W,\nbigi)$
such that $S'\subset S$.
The natural map
$\bigl(\nbigi,\leq_S\bigr)\lrarr 
 \bigl(\nbigi,\leq_{S'}\bigr)$ is order-preserving.
We can show the following lemma easily
by using Proposition \ref{prop;07.6.16.6}.
(See Subsections 7.2.2 and 8.1.2 of \cite{mochi7}.)
\begin{lem}
\label{lem;07.12.16.9}
The filtrations
$\nbigf^{S}$ and $\nbigf^{S'}$
are compatible with respect to
$\bigl(\nbigi,\leq_S\bigr)\lrarr 
 \bigl(\nbigi,\leq_{S'}\bigr)$
in the following sense:
\begin{itemize}
\item
$\nbigf^{S'}_{\gminia}(E_{|\Sbar'})
=\nbigf^{S'}_{<\gminia}(E_{|\Sbar'})
+\nbigf^{S}_{\gminia}(E_{|\Sbar})_{|\Sbar'}$.
\item
The induced morphisms
$\Gr^{\nbigf^{S}}_{\gminia}(E_{|S})_{|S'}
\lrarr
 \Gr^{\nbigf^{S'}}_{\gminia}(E_{|S'})$
are isomorphisms.
\end{itemize}
In particular,
we have 
$\nbigf^{S}(E_{|\Sbar})_{|\Sbar'}
=\nbigf^{S'}(E_{|\Sbar'})$,
if $\bigl(\nbigi,\leq_S\bigr)\lrarr 
 \bigl(\nbigi,\leq_{S'}\bigr)$ is isomorphic.
\hfill\qed
\end{lem}

\paragraph{Splitting with nice property}

We have the induced morphisms
$\Res_j(\DD):E_{|\nbigd_{\zeta,j}}
\lrarr
 \vecz^{\vecm(1)}\cdot E_{|\nbigd_{\zeta,j}}$
for $j=1,\ldots,\ell$.
Since $\nbigf^S$ is $\DD$-flat,
$\Res_j(\DD)$ preserves 
$\nbigf^{S}_{|\nbigd_{\zeta,j}}$.
If we fix the coordinate,
we have the induced family of flat $\lambda$-connections
of $E_{|\nbigd_{\zeta,j}}$
which is denoted by $\lefttop{j}\DD$.
It also preserves the filtration $\nbigf^S_{|\nbigd_{\zeta,j}}$.
Let 
$\lefttop{j}F$ $(j=1,\ldots,\ell)$ be filtrations
of $E_{|\nbigd_{\zeta,j}}$,
which are preserved by
the endomorphism $\Res_j(\DD)$
and the flat connection $\lefttop{j}\DD$
of $E_{|\nbigd_{\zeta,j}}$.
We can show the following
(Subsections 7.2.3 and 8.1.3 of \cite{mochi7}).

\begin{prop}
\label{prop;07.9.30.22}
Let $P\in \gbigz$.
There exist 
$S\in\Multisector^{\ast}(P,\nbigx-W,\nbigi)$
and a $\DD_{z}$-flat splitting
of the filtration $\nbigf^{S}$,
whose restriction to $\Sbar\cap \nbigd_{\zeta,j}$
is compatible with $\Res_{j}(\DD)$ 
and the filtrations $\lefttop{j}F$
for $j=1,\ldots,\ell$.
\hfill\qed
\end{prop}

Under some more assumption,
we can take a $\DD$-flat splitting.
(See Subsection 7.2.3 of \cite{mochi7}.)
\begin{prop}
\label{prop;07.9.30.5}
Assume that $\nbigk$ is a point or 
a compact region in $\cnum_{\lambda}^{\ast}$.
Assume that the eigenvalues $\alpha,\beta$
of $\Res_{j}(\DD^f)_{|D_j\times\{\lambda\}}$
satisfy $\alpha-\beta\not\in (\seisuu-\{0\})$
for any $j=1,\ldots,\ell$ and for any $\lambda\in\nbigk$.
Then, we have a $\DD$-flat splitting
of $\nbigf^S$,
whose restriction to $\nbigd_{\zeta,j}$ is compatible
with $\lefttop{j}F$ for each $j=1,\ldots,\ell$.
\hfill\qed
\end{prop}

\paragraph{Some functoriality of Stokes filtrations}

We explain functoriality of Stokes filtrations.
See Subsections 7.2.4 and 8.1.4 of \cite{mochi7}
for more details.

In general,
when we are given vector spaces
$U\subset V$,
let $U^{\bot}$ denote the subspace of
the dual $V^{\lor}$
given by $U^{\bot}=\bigl\{
 f\in V^{\lor}\,\big|\,f(U)=0 \bigr\}$.
It is naturally generalized for vector bundles.
Let $(E,\DD,\nbigi)$ be an unramifiedly good lattice 
of a family of meromorphic $\lambda$-flat bundles
in the level $(\vecm,i(0))$ on $(\nbigx,\nbigd_z)$.
Let $S\in\Multisector^{\ast}(\nbigx-W,\nbigi)$.
We have the following 
for any $\gminia\in\nbigi^{\lor}:=\bigl\{
 -\gminib\,\big|\,\gminib\in\nbigi
 \bigr\}$:
\[
 \nbigf^S_{\gminia}(E^{\lor}_{|\Sbar})
=\left(
 \sum_{\substack{\gminic\in\nbigi\\
  \gminic\not\geq_{S}-\gminia}}
 \nbigf^S_{\gminic}(E_{|\Sbar})
 \right)^{\bot}
\]

Let $(E_p,\DD_p,\nbigi_p)$ $(p=1,2)$ be 
good lattices of families of 
meromorphic $\lambda$-flat bundles
in the level $(\vecm,i(0))$.
We assume that
$\nbigi_1\otimes\nbigi_2:=
 \bigl\{\gminia_1+\gminia_2\,\big|\,
 \gminia_p\in\nbigi_p \bigr\}$ 
is a good set of irregular values in the level $(\vecm,i(0))$.
We put $(\Etilde,\DDtilde):=
 (E_1,\DD_1)\otimes (E_2,\DD_2)$.
Let $S\in\bigcap_{p=1,2}\Multisector^{\ast}
\bigl(\nbigx-W,\nbigi_p\bigr)$.
We have the following for each
$\gminia\in\nbigi_1\otimes\nbigi_2$:
\[
\nbigf^S_{\gminia}(\Etilde_{|\Sbar})
=
\sum_{\gminia_1+\gminia_2\leq_S\gminia}
 \nbigf^S_{\gminia_1}(E_{1|\Sbar})
\otimes
 \nbigf^S_{\gminia_2}(E_{2|\Sbar}).
\]

Assume that
$\nbigi_1\oplus\nbigi_2
:=\nbigi_1\cup\nbigi_2$ is a good set of
irregular values in the level $(\vecm,i(0))$.
Let $S\in\bigcap_{p=1,2}\Multisector^{\ast}\bigl(
 \nbigx-W,\nbigi_p\bigr)$.
We have the following 
for each $\gminia\in\nbigi_1\oplus\nbigi_2$:
\[
\nbigf^S_{\gminia}\bigl((E_1\oplus E_2)_{|\Sbar}\bigr)
=
\nbigf^S_{\gminia}(E_{1|\Sbar})
\oplus
 \nbigf^S_{\gminia}(E_{2|\Sbar}).
\]

Let $F:(E_1,\DD_1)\lrarr(E_2,\DD_2)$
be a flat morphism.
For simplicity,
we assume that $\nbigi_1\cup\nbigi_2$
is a good set of irregular values
in the level $(\vecm,i(0))$.

\begin{lem}
\label{lem;07.12.7.6}
Let $S\in\bigcap_{p=1,2}
  \Multisector^{\ast}(\nbigx-W,\nbigi_p)$.
The restriction $F_{|\Sbar}$ preserves
the Stokes filtrations.
As a result, we obtain the following.
\begin{itemize}
\item
If the restriction of $F$ to $\nbigx-\nbigd$
is isomorphic,
we have $\nbigi_1=\nbigi_2$
and 
$\nbigf^S_{\gminia}(E_{1|S\setminus \nbigd})
=\nbigf^S_{\gminia}(E_{2|S\setminus \nbigd})$.
\item
In particular, the Stokes filtration $\nbigf^S$
depends only on 
the family of meromorphic $\lambda$-flat bundles
$\bigl(E(\ast \nbigd),\DD\bigr)$
in the sense that
it is independent of
the choice of an unramifiedly good lattice 
$E\subset E(\ast \nbigd)$ in the level $(\vecm,i(0))$.
\hfill\qed
\end{itemize}
\end{lem}

\paragraph{The associated graded bundle
 in the level $\vecm$}

For sectors $S\in\Multisector^{\ast}(\nbigx-W,\nbigi)$
and each $\gminia\in\nbigi$,
we obtain the bundle 
$\Gr_{\gminia}^{\vecm}(E_{|\Sbar})$ on $\Sbar$
associated to the Stokes filtration $\nbigf^S$
in the level $\vecm$.
By varying $S$ and gluing 
$\Gr^{\vecm}_{\gminia}(E_{|\Sbar})$,
we obtain the bundle
$\Gr^{\vecm}_{\gminia}(E_{|\nbigvtilde(W)})$
on $\nbigvtilde(W)$
with the induced family of flat $\lambda$-connections 
$\DD_{\gminia}$,
where $\nbigv$ denotes some neighbourhood
of $O_z\times \nbigy$,
and $\nbigvtilde(W)$
denotes the real blow up of $\nbigv$
along $W\cap\nbigv$.
It is shown that 
we have the descent of
$\Gr^{\vecm}_{\gminia}(E_{|\nbigvtilde(W)})$
to $\nbigv$,
i.e.,
there exists a locally free sheaf 
$\Gr^{\vecm}_{\gminia}(E)$
on $\nbigv$
with a family of meromorphic flat
$\lambda$-connections
$\DD_{\gminia}$,
such that 
\[ 
\pi^{-1}
(\Gr^{\vecm}_{\gminia}(E),
 \DD_{\gminia})
\simeq 
 (\Gr^{\vecm}_{\gminia}(E_{|\nbigvtilde(W)}),
 \DD_{\gminia}),
\quad\quad
(\Gr^{\vecm}_{\gminia}(E),
 \DD_{\gminia})_{|\What\cap\nbigv}
\simeq 
\bigl(\Ehat_{\gminia},
 \DDhat_{\gminia}\bigr)_{|\What\cap\nbigv} 
\]
(See Subsection 7.3 and Subsection 8.1.5
of \cite{mochi7}.)
If we set
$\DD'_{\gminia}:=\DD_{\gminia}-d_X\gminia$,
we have
\[
 \DD'_{\gminia}
 E_{\gminia}
\subset 
 \vecz^{\vecm(1)}
 \cdot E_{\gminia}
 \otimes
 p_{\lambda}^{\ast}\Omega_X^1(\log D).
\]

We give some statements for functoriality.
See Subsections 7.3.2 and 8.1.6 of \cite{mochi7}
for more details.

By taking Gr of the Stokes filtrations
of $(E^{\lor},\DD^{\lor},\nbigi^{\lor})$,
we obtain the associated graded bundle
$\Gr^{\vecm}(E^{\lor})
=\bigoplus_{\gminia\in\nbigi^{\lor}}
 \Gr^{\vecm}_{\gminia}(E^{\lor})$.
We have the natural flat isomorphism
\begin{equation}
 \label{eq;08.11.5.1}
\Gr^{\vecm}_{\gminia}(E^{\lor})
\simeq
 \Gr^{\vecm}_{-\gminia}(E)^{\lor}.
\end{equation}
Actually, by construction,
we have such an isomorphism on the real blow up,
which induces (\ref{eq;08.11.5.1}).

Let $(E_p,\nabla_p,\nbigi_p)$ $(p=1,2)$ be 
unramifiedly good lattices of 
families of meromorphic $\lambda$-flat bundles.
Assume $\nbigi_1\otimes\nbigi_2$ is a good set of
irregular values in the level $\vecm$.
Let $(\Etilde,\DDtilde):=
 (E_p,\DD_1)\otimes (E_p,\DD_2)$.
We have the following natural isomorphism
for each $\gminia\in\nbigi_1\otimes\nbigi_2$:
\begin{equation}
 \label{eq;07.12.7.4}
 \Gr_{\gminia}^{\vecm}(\Etilde)
\simeq
\bigoplus_{\substack{
 (\gminia_1,\gminia_2)\in
 \nbigi_1\times
 \nbigi_2\\
 \gminia_1+\gminia_2=\gminia}}
 \Gr^{\vecm}_{\gminia_1}(E_1)
\otimes
 \Gr^{\vecm}_{\gminia_2}(E_2)
\end{equation}
Assume $\nbigi_1\oplus\nbigi_2$ is a good set of
irregular values in the level $(\vecm,i(0))$.
For each $\gminia\in\nbigi_1\oplus\nbigi_2$,
we obviously have 
\[
 \Gr^{\vecm}_{\gminia}(E_1\oplus E_2)
\simeq
 \Gr^{\vecm}_{\gminia}(E_1)
\oplus
 \Gr^{\vecm}_{\gminia}(E_2).
\]

\begin{lem}
\label{lem;07.12.7.17}
Let $F:(E_1,\DD_1)\lrarr (E_2,\DD_2)$
be a flat morphism.
Assume $\nbigi_1\oplus\nbigi_2$ is a good set of
irregular values in the level $(\vecm,i(0))$.
We have the naturally induced morphism
$\Gr^{\vecm}_{\gminia}(F):
 \Gr^{\vecm}_{\gminia}(E_1)
\lrarr \Gr^{\vecm}_{\gminia}(E_2)$.
If the restriction $E_{1|\nbigx-\nbigd}\lrarr 
 E_{2|\nbigx-\nbigd}$ is an isomorphism,
the induced morphism 
\[
 \Gr^{\vecm}_{\gminia}(E_1)\otimes\nbigo(\ast \nbigd)
\lrarr
 \Gr^{\vecm}_{\gminia}(E_2)\otimes\nbigo(\ast \nbigd)
\]
is an isomorphism.

Hence, the associated meromorphic flat bundles
$\bigl(\Gr^{\vecm}_{\gminia}(E)\otimes
 \nbigo(\ast \nbigd),\DD_{\gminia}\bigr)$
are well defined for the meromorphic flat bundle
$\bigl(E(\ast \nbigd),\DD\bigr)$.
\hfill\qed
\end{lem}

\paragraph{A characterization of
sections of $E$}

Let $\vecw_{\gminia}$ be a frame of
$\Gr^{\vecm}_{\gminia}(E)$.
Let $S\in\Multisector^{\ast}(\nbigx-W,\nbigi)$,
and let $E_{|\Sbar}=\bigoplus E_{\gminia,S}$
be a $\DD_{z}$-flat splitting
of the Stokes filtration $\nbigf^S$.
By the natural isomorphism
$E_{\gminia,S}\simeq
 \Gr_{\gminia}^{\vecm}(E)_{|\Sbar}$,
we take a lift $\vecw_{\gminia,S}$ of $\vecw_{\gminia}$.
Thus, we obtain a frame 
$\vecw_S=\bigl(\vecw_{\gminia,S}\bigr)$
of $E_{|\Sbar}$.
The following proposition is clear,
which implies a
characterization of sections of $E$
by growth order with respect
to the frames $\vecw_S$ 
($S\in\Multisector^{\ast}(\nbigx-W,\nbigi)$).

\begin{prop}
\label{prop;08.8.2.2}
Let $\vecv$ be a frame of $E$,
and let $G_S$ be determined by
$\vecv=\vecw_S\cdot G_S$.
Then, $G_S$ and $G_S^{-1}$ are bounded on $S$.
\hfill\qed
\end{prop}

\paragraph{Complement
on the induced flat connection along 
 the $\lambda$-direction}

Assume that we are given 
a connection along the $\lambda$-direction
$ \nabla_{\lambda}:
E\lrarr E\otimes\Omega_{\nbigk}^1(\ast W)$
such that $\DD^f+\nabla_{\lambda}$
is a meromorphic flat connection of $E$.
\begin{lem}
\label{lem;08.8.3.3}
The Stokes filtrations are flat 
with respect to $\nabla_{\lambda}$,
and we have the induced meromorphic
flat connection $\nabla_{\lambda}$
along the $\lambda$-direction
on $\Gr^{\vecm}_{\gminia}(E)$.
\end{lem}
\pf
Take $N$ such that
$\lambda^N\nabla_{\lambda}(\del_{\lambda})
 E\subset 
 E\otimes\nbigo_{\nbigx}(\ast\nbigd)$.
Let $\vecw_{S}=(\vecw_{\gminia,S})$
be a frame of $E_{|\Sbar}$ as above.
Let $A=(A_{\gminia,\gminib})$
be the matrix-valued holomorphic function
on $S$ determined by
$\lambda^N\nabla(\del_{\lambda})
 \vecw_S=\vecw_S\cdot A$.
By using Proposition \ref{prop;08.8.2.2},
we can show that
$A_{\gminia,\gminib}$ are
of polynomial order.

Let $B_{\gminia}$ be the matrix-valued 
meromorphic one-forms determined by
$\DD_{\gminia,z}
\vecw_{\gminia}
=\vecw_{\gminia}\cdot
 \bigl(d_z\gminia+B_{\gminia}\bigr)$.
Note that
$\vecz^{-\vecm(1)}B_{\gminia}$
is logarithmic.
By the commutativity
$[\DD^f,\nabla_{\lambda}]=0$,
we obtain the following relation
for $\gminia\neq\gminib$:
\begin{equation}
 \label{eq;08.8.3.1}
 \lambda\cdot d_zA_{\gminia,\gminib}
+\bigl(d_z(\gminia-\gminib)\bigr)
 \cdot A_{\gminia,\gminib}
+\bigl(
 A_{\gminia,\gminib}B_{\gminib}
-B_{\gminia}A_{\gminia,\gminib}
 \bigr)=0
\end{equation}
By applying the results in Subsection 4.3
of \cite{mochi7}
to (\ref{eq;08.8.3.1}),
we obtain $A_{\gminia,\gminib}=0$
unless $\gminia\leq_S\gminib$,
which implies the first claim.
Since $A_{\gminia,\gminia}$
is of polynomial order,
the induced connection along the $\lambda$-direction
is meromorphic.
\hfill\qed

\paragraph{Prolongment of morphisms}

Let $(E_p,\DD_p,\nbigi_p)$ $(p=1,2)$ be 
good lattices in the level $(\vecm,i(0))$.
Assume that $\nbigi_1\cup\nbigi_2$
is a good set of irregular values in the level
$(\vecm,i(0))$.
Assume that we are given a flat morphism
$F:(E_1,\DD_1)_{|\nbigx-\nbigd_z}\lrarr 
 (E_2,\DD_2)_{|\nbigx-\nbigd_z}$
with the following property:
\begin{itemize}
\item
 For each small sector 
 $S\in \Multisector(\nbigx-\nbigd_z,\nbigi_1\cup\nbigi_2)$,
 the Stokes filtrations are preserved by $F_{|S}$.
\item
 The induced maps
 $\Gr^{\vecm}_{\gminia}(F):
 \Gr^{\vecm}_{\gminia}(E_1)_{|\nbigx-\nbigd_z}
\lrarr
 \Gr^{\vecm}_{\gminia}(E_2)_{|\nbigx-\nbigd_z}$
 are extended to
 $\Gr^{\vecm}_{\gminia}(E_1)
\lrarr
 \Gr^{\vecm}_{\gminia}(E_2)$ 
for any $\gminia\in\nbigi_1\cup\nbigi_2$.
\end{itemize}

\begin{lem}
\label{lem;08.11.6.15}
$F$ is extended to a morphism
$E_{1}\lrarr E_2$.
\end{lem}
\pf
Let $\vecw_{p,S}=(\vecw_{p,\gminia,S})$ be frames of
$E_{p|\Sbar}$ as above.
Let $A=(A_{\gminia,\gminib})$ be determined by
$F(\vecw_{1,S})
=\vecw_{2,S}\cdot A$.
By the assumption,
$A_{\gminia,\gminib}=0$
unless $\gminia\leq_S\gminib$,
and $A_{\gminia,\gminia}$ is bounded.
By applying an argument in the proof of 
Lemma \ref{lem;08.8.3.3}
to $A_{\gminia,\gminib}$ 
for $\gminia<_S\gminib$,
and by shrinking $X$,
we obtain
$A_{\gminia,\gminib}
=O\bigl(
 \exp\bigl(-\epsilon|\lambda^{-1}\cdot\vecz^{\vecm}|\bigr)
 \bigr)$
on $S\cap (\nbigx-\nbigd_z)$.
Then, the claim follows from
Proposition \ref{prop;08.8.2.2}.
\hfill\qed

\subsubsection{Pseudo-good lattice
 in the level $\vecm$}
\label{subsubsection;08.8.12.4}

Let $Y$ be a complex manifold.
Let $X:=\Delta_z^k\times Y$,
$D_{z,i}:=\{z_i=0\}$
and $D:=\bigcup_{i=1}^{k}D_{z,i}$.
Let $E$ be a locally free $\nbigo_{X}$-module.
For simplicity,
we consider a meromorphic flat connection
$\nabla:
 E\lrarr E\otimes \Omega^{1}_{X}(\ast D)$
instead of a family of meromorphic
flat $\lambda$-connections.
Let $\vecm\in\seisuu_{<0}^k$
and $i(0)\in[1,k]$.
We put $\vecm(1):=\vecm+\vecdelta_{i(0)}$.

\begin{df}
We say that $(E,\nabla)$ is 
an unramifiedly pseudo-good lattice
in the level $(\vecm,i(0))$,
if there exists an unramifiedly good lattice
$E'\supset E$ of $\bigl(E(\ast D),\nabla\bigr)$
with the irregular decomposition
$(E',\DD)_{|\Dhat}=
 \bigoplus_{\gminia\in\nbigi}
 (\Ehat'_{\gminia},\nablahat_{\gminia})$
in the level $(\vecm,i(0))$,
such that 
\begin{equation}
\label{eq;08.8.12.1}
 E_{|\Dhat}
=\bigoplus_{\gminia\in\nbigi} \bigl(
 \Ehat'_{\gminia}\cap E_{|\Dhat}
\bigr)
\end{equation}
The decomposition {\rm(\ref{eq;08.8.12.1})}
is called the irregular decomposition of $(E,\DD)$
in the level $(\vecm,i(0))$.
\hfill\qed
\end{df}

It is easy to observe that
$\Ehat_{\gminia}:=\Ehat'_{\gminia}\cap E_{|\nbigdhat_z}$
in (\ref{eq;08.8.12.1})
is independent of the choice of
a good lattice $E'\supset E$ in the level $\vecm$.
We have straightforward generalizations
of the results in Subsection
\ref{subsubsection;08.8.3.2}.
We naturally identify $X$ with 
$\{1\}\times X\subset
  \cnum_{\lambda}\times X$
when we consider the order
$\leq_S$ for multi-sectors $S\subset X-D$.

\paragraph{Construction of Gr}

We take an unramifiedly good lattice
$E'\supset E$ in the level $(\vecm,i(0))$.
By shrinking $X$ around $O_z\times Y$,
we have the vector bundle
$\Gr^{\vecm}_{\gminia}(E')$ on $X$
with a meromorphic flat connection
$\nabla_{\gminia}$
for each $\gminia\in\nbigi$.
Recall that we have the natural isomorphism
$\Gr^{\vecm}_{\gminia}(E')_{|\Dhat}
\simeq
 \Ehat'_{\gminia}$.
Hence, we have the sub-lattice
of $\Gr^{\vecm}_{\gminia}(E')$
corresponding to
$\Ehat_{\gminia}
 \subset \Ehat'_{\gminia}$,
which is denoted by
$\Gr^{\vecm}_{\gminia}(E)$.
It is equipped with a meromorphic flat 
connection $\nabla_{\gminia}$.
By construction,
we have the isomorphism
\begin{equation}
 \label{eq;08.10.1.4}
 (\Gr^{\vecm}_{\gminia}(E),\nabla_{\gminia})
 _{|\Dhat}
\simeq
 (\Ehat_{\gminia},\nablahat_{\gminia}).
\end{equation}

\begin{lem}
\label{lem;08.8.12.10}
Let $(E_i,\nabla_i)$ $(i=1,2)$ be
pseudo-good lattices in the level $(\vecm,i(0))$.
Let $F:(E_1,\nabla_1)\lrarr (E_2,\nabla_2)$
be a flat morphism.
Assume $\nbigi_1\oplus\nbigi_2$ is a good set of
irregular values in the level $(\vecm,i(0))$.
We have the naturally induced morphism
$\Gr^{\vecm}_{\gminia}(F):
 \Gr^{\vecm}_{\gminia}(E_1)
\lrarr \Gr^{\vecm}_{\gminia}(E_2)$.
\end{lem}
\pf
We can take good lattices $(E_i',\nabla_i)$
in the level $(\vecm,i(0))$
such that
$E_i\subset E_i'$ and $F(E_1')\subset E_2'$.
By Lemma \ref{lem;07.12.7.17},
we have the induced morphism
$\Gr^{\vecm}_{\gminia}(F):
 \Gr^{\vecm}_{\gminia}(E_1')\lrarr
 \Gr^{\vecm}_{\gminia}(E_2')$.
By considering the completion,
it is easy to observe that
$\Gr^{\vecm}_{\gminia}(E_1)\lrarr
 \Gr^{\vecm}_{\gminia}(E_2)$
is induced.
\hfill\qed

\paragraph{Flat splitting and Stokes filtration}

Let  $\pi:\Xtilde(D)\lrarr X$ be the real blow up.
Let $S\in\Multisector^{\ast}(X-D,\nbigi)$.
Let $\Sbar$ denote the closure of $S$
in $\Xtilde(D)$,
and let $Z$ denote $\Sbar\cap \pi^{-1}(D)$.
We have the Stokes filtration
$\nbigf^S$ of $E'_{|\Sbar}$,
and we can take a flat splitting
$E'_{|\Sbar}=\bigoplus E'_{\gminia,S}$
such that
$E'_{\gminia,S|\Zhat}
=\pi^{-1}(\Ehat'_{\gminia})$.
Because
$E_{|X-D}=E'_{|X-D}$,
it induces the flat decomposition of $E_{|S}$.

\begin{lem}
\label{lem;08.8.19.20}
It is extended to the decomposition
$E_{|\Sbar}=\bigoplus E_{\gminia,S}$
such that
$E_{\gminia,S|\Zhat}
=\pi^{-1}(\Ehat_{\gminia})$.
\end{lem}
\pf
Let $\vecw_{\gminia}$ and $\vecw'_{\gminia}$
be frames of $\Gr^{\vecm}_{\gminia}(E)$
and $\Gr^{\vecm}_{\gminia}(E')$.
Let $G_{\gminia}$ be determined by
$\vecw_{\gminia}=
 \vecw'_{\gminia}\cdot G_{\gminia}$.
They induce the frames
$\vecwhat_{\gminia}$
and $\vecwhat'_{\gminia}$
of $\Ehat_{\gminia}$ and $\Ehat'_{\gminia}$,
respectively.

By the isomorphism
$E'_{\gminia,S}\simeq
 \Gr^{\vecm}_{\gminia}(E')_{|\Sbar}$,
we obtain the frames
$\vecw'_{\gminia,S}$ of $E'_{\gminia,S}$.
Then,
$\vecw_{\gminia,S}:=
 \vecw'_{\gminia,S}\cdot G_{\gminia}$
gives a tuple of sections of $E'_{\gminia,S}$,
and we can observe that
$\vecw_{\gminia|\Zhat}
=\pi^{-1}(\vecwhat_{\gminia})$.
Let $E_{\gminia,S}$ be generated by
$\vecw_{\gminia,S}$,
and then we obtain the desired decomposition
$E=\bigoplus E_{\gminia,S}$.
\hfill\qed

\vspace{.1in}
Let $\vecw_S=(\vecw_{\gminia,S})$  be as above.
Let $\vecv$ be a frame of $E$ on $X$.
Let $G_S$ be determined by
$\vecv=\vecw_S\cdot G_S$.
Both $\vecv_{|\Zhat}$
and $\vecw_{S|\Zhat}$ give
the frame of $E_{|\Zhat}$,
we obtain the following.
\begin{prop}
\label{prop;08.8.12.11}
$G_S$ and $G_S^{-1}$
are bounded on $S$.
\hfill\qed
\end{prop}

\begin{prop}
\label{prop;08.8.12.5}

The flat subbundle
$\nbigf^S_{\gminia}(E_{|\Sbar}):=
 \bigoplus_{\gminib\leq_S\gminia}
 \Ebar_{\gminib,S}$ is independent
of the choice of a flat decomposition
$E_{|\Sbar}=\bigoplus_{\gminia\in\nbigi}
 \Ebar_{\gminia,S}$
such that $\Ebar_{\gminia,S|\Zhat}
=\pi^{-1}\Ehat_{\gminia}$.
\end{prop}
\pf
Let 
$E_{|\Sbar}=\bigoplus_{\gminia\in\nbigi}
 \Ebar_{\gminia,S}$
be another  flat decomposition
such that $\Ebar_{\gminia,S|\Zhat}
=\pi^{-1}\Ehat_{\gminia}$.
We take a frame
$\vecwbar_{\gminia,S}$ of
$\Ebar_{\gminia,S}$ such that
$\vecwbar_{\gminia,S|\Zhat}
=\vecwhat_{\gminia}$.
We set 
$\vecwbar'_{\gminia}:=
 \vecwbar_{\gminia}\cdot G_{\gminia}^{-1}$.
Then, $\vecwbar'_{\gminia|\Zhat}=
 \pi^{-1}\vecwhat_{\gminia}'$.
Let $\Ebar'_{\gminia}$ be generated by
$\vecwbar'_{\gminia}$.
Then, we obtain a flat decomposition
$E'_{|\Sbar}
=\bigoplus \Ebar_{\gminia}'$,
which has to be a splitting of the Stokes filtration 
$\nbigf^S(E'_{|\Sbar})$.
Because $\Ebar'_{\gminia|S}=\Ebar_{\gminia|S}$,
we obtain the well definedness 
of the filtration.
\hfill\qed

\vspace{.1in}

Thus, we obtain the filtration
$\nbigf^S$ of $E_{|\Sbar}$,
which is called the Stokes filtration.

\begin{lem}
\label{lem;08.8.20.1}
We have the natural isomorphism
$\Gr^{\nbigf^S}_{\gminia}(E_{|\Sbar})
\simeq
 \Gr_{\gminia}(E)_{|\Sbar}$.
\end{lem}
\pf
We use the notation in
the proof of Lemma \ref{lem;08.8.19.20}.
By the comparison
of $\vecw_{\gminia}$ and $\vecw_{\gminia,S}$,
we obtain 
$E_{\gminia,S}\simeq
 \Gr_{\gminia}(E)_{|\Sbar}$.
By the construction of the Stokes filtration,
we have the natural isomorphism
$\Gr^{\nbigf^S}_{\gminia}(E_{|\Sbar})
\simeq
 E_{\gminia,S}$.
Then, the claim of 
Lemma \ref{lem;08.8.20.1} is clear.
\hfill\qed

\subsubsection{A comparison}

Let $Y$ be a complex manifold.
Let $X:=\Delta_z^k\times Y$,
$D_{z,i}:=\{z_i=0\}$
and $D:=\bigcup_{i=1}^{k}D_{z,i}$.
Let $\nbigktilde$ be a compact region in $\cnum_{z_{k+1}}$.
We set $\nbigxtilde:=\nbigktilde\times X$.
We use the symbol
$\nbigdtilde$ in a similar meaning.
We set $\Wtilde:=
 \nbigdtilde\cup (\{0\}\times X)$.

Let $\nbigi\subset M(X,D)$ be 
a good set of irregular values in the level $(\vecm,i(0))$.
We set 
$\vecmtilde:=(\vecm,-1)\in\seisuu_{<0}^{k+1}$.
We put $\gminiatilde:=z_{k+1}^{-1}\gminia$
for $\gminia\in \nbigi$,
and we set
\[
 \nbigitilde:=\bigl\{
 \gminiatilde\,\big|\,
 \gminia\in \nbigi \bigr\}
 \subset 
 M(\nbigxtilde,\Wtilde)/H(\nbigxtilde).
\]
Then, it is a good set of irregular values
in the level $(\vecmtilde,i(0))$.

Let $\Etilde$ be a holomorphic vector bundle on
$\nbigxtilde$
with a meromorphic flat connection
$\nabla:\Etilde\lrarr 
 \Etilde\otimes\Omega^1_{\nbigxtilde}(\ast \Wtilde)$
such that $(\Etilde,\nabla)$ is 
an unramifiedly good lattice
in the level $(\vecmtilde,i(0))$ on 
$(\nbigxtilde,\Wtilde)$
with the irregular decomposition:
\begin{equation}
 \label{eq;08.8.12.3}
 (\Etilde,\nabla)_{|\widehat{\Wtilde}}
=\bigoplus_{\gminiatilde\in\nbigitilde}
  (\widehat{\Etilde}_{\gminiatilde},
 \widehat{\nabla}_{\gminiatilde})
\end{equation}
Applying a general theory in Subsection
\ref{subsubsection;08.8.12.4},
we obtain a holomorphic vector bundle
$\Gr^{\vecmtilde}_{\gminiatilde}(\Etilde)$
on $\nbigxtilde$
with the induced meromorphic flat connection
$\nabla_{\gminiatilde}$
for each $\gminiatilde\in \nbigitilde$.

By setting $\lambda=z_{k+1}$,
we obtain the isomorphism
$\cnum_{z_{k+1}}\simeq\cnum_{\lambda}$.
Let $\nbigk\subset\cnum_{\lambda}$
be the image of $\nbigktilde$.
We put $\nbigx:=\nbigk\times X$
and we use the symbol $\nbigd$
in a similar meaning.
We set $W:=\nbigd\cup\bigl(\{0\}\times X\bigr)$.
We have the natural isomorphism
$\iota:(\nbigx,\nbigd)\lrarr (\nbigxtilde,\nbigdtilde)$.
The pull back of $\Etilde$ is denoted by $E$.
Let $\DD^f$ denote the restriction of
$\iota^{\ast}\nabla$ to the $X$-direction.
We set $\DD:=\lambda\cdot \DD^f$.
Note the following:
\begin{itemize}
\item
$\DD(E)\subset
 E\otimes p_{\lambda}^{\ast}
 \Omega^{1}_X(\ast D)$,
i.e.,
$\DD$ gives a family of 
meromorphic $\lambda$-connections of $E$.
\item
$(E,\DD)$ is a good lattice in the level 
$(\vecm,i(0))$ on $(\nbigx,W)$,
and (\ref{eq;08.8.12.3})
naturally induces the irregular decomposition
of $(E,\DD)_{|\widehat{W}}$.
\end{itemize}
By applying a general theory explained
in Subsection \ref{subsubsection;08.8.3.2},
for each $\gminia\in \nbigi$,
we obtain 
$\Gr^{\vecm}_{\gminia}(E,\DD)$.

\vspace{.1in}

Let $\Stilde$ be a small sector in $\nbigxtilde-\Wtilde$.
We have the Stokes filtration 
$\nbigf^{\Stilde}$
of $\Etilde_{|\overline{\Stilde}}$ 
in the level $\vecmtilde$
indexed by
$(\nbigitilde,\leq_{\Stilde})$
(Proposition \ref{prop;08.8.12.5}).
For $S:=\iota^{-1}(\Stilde)$,
we have the Stokes filtration
$\nbigf^S$ of $E_{|\Sbar}$
in the level $\vecm$
indexed by $(\nbigi,\leq_{S})$.
We remark the following.
\begin{lem}
\label{lem;08.8.12.6}
Under the natural identification
$\nbigitilde=\nbigi$,
the orders $\leq_{\Stilde}$ and $\leq_S$
are the same.
Under the natural isomorphism
$E\simeq \iota^{\ast}\Etilde$,
the filtrations
$\nbigf^S$ and $\nbigf^{\Stilde}$
are the same.
\end{lem}
\pf
For the order $\leq_{\Stilde}$,
we use the identification
$\nbigxtilde
=\{1\}\times \nbigxtilde\subset
 \cnum_{\lambda}\times\nbigxtilde$.
Then, the first claim is clear.
Note that
both $\iota^{\ast}\nbigftilde^{\Stilde}$
and $\nbigf^S$ satisfy the condition
in Proposition \ref{prop;07.6.16.6}.
Hence, they are the same.
\hfill\qed

\begin{cor}
\label{cor;08.8.7.5}
We have the natural isomorphism
$\iota^{\ast}\Gr^{\vecmtilde}_{\gminiatilde}(\Etilde)
\simeq
 \Gr^{\vecm}_{\gminia}(E)$,
and $\DD_{\gminia}$ is induced 
by $\iota^{\ast}\nabla_{\gminiatilde}$
via the above procedure.
\end{cor}
\pf
By Lemma \ref{lem;08.8.12.6},
we obtain the isomorphism
$j:\iota^{\ast}
 \Gr^{\vecmtilde}_{\gminiatilde}
 (\Etilde)_{|\nbigx-W}
\simeq
 \Gr^{\vecm}_{\gminia}
 (E)_{|\nbigx-W}$,
on which $\DD_{\gminia}$
is induced by $\nabla_{\gminia}$
via the above procedure.
Since $j$ is extended on $\nbigxtilde(W)$,
it is extended on $\nbigx$.
\hfill\qed

\subsubsection{Stokes filtration
of the associated flat bundle on the real blow up}
\label{subsubsection;08.8.8.3}

We use the setting in Subsection
\ref{subsubsection;08.8.12.4}.
Let $\nbigi\subset M(X,D)$ be 
a good set of irregular values
in the level $(\vecm,i(0))$.
Let $E$ be a holomorphic vector bundle on $X$
with a meromorphic flat connection
$\nabla:E\lrarr E\otimes\Omega^{1}_X(\ast D)$
such that $(E,\nabla)$ is a pseudo-good lattice 
in the level $(\vecm,i(0))$.
(In other words,
we consider a family of meromorphic $\lambda$-flat
bundles on $\{1\}\times (X,D)$.)
Let $\pi:\Xtilde(D)\lrarr X$ be a real blow up of $X$
along $D$.
The flat bundle $E_{|X-D}$ is naturally
extended to the flat bundle $\gbigv$
on $\Xtilde(D)$.

We set $\gbigz:=\pi^{-1}(O_z\times Y)$.
For each $P\in \gbigz$,
we take a small sector 
$S\in\Multisector(P,X-D,\nbigi)$
on which we have the Stokes filtration
$\nbigf^S$ of $E_{|S}$.
The filtration is naturally extended to
the flat filtration of
$\gbigv_{|\Sbar}$.
By restricting it to the fiber
$\gbigv_{|P}$,
we obtain the filtration $\nbigf^P$
indexed by $\bigl(\nbigi,\leq_P\bigr)$.
It is easy to observe that
$\nbigf^P$ is well defined.

If $Q\in\pi^{-1}(\gbigz)$ is sufficiently close to $P$,
the map
$\bigl(\nbigi,\leq_P\bigr)
\lrarr \bigl(\nbigi,\leq_Q\bigr)$ preserves the orders,
and the filtrations
$\nbigf^P$ and $\nbigf^Q$ are compatible
under the identification
$\gbigv_{|P}\simeq\gbigv_{|Q}$
given by the parallel transport in $\Sbar_P$.
In particular,
we have $\nbigf^P=\nbigf^Q$
if $\leq_P=\leq_Q$.

We have the functoriality of 
the filtrations $\nbigf^P$
for dual, tensor product and direct sum
as in the case of $\nbigf^S$.

\begin{lem}
Let $F:(E_1,\nabla_1)\lrarr(E_2,\nabla_2)$
be a flat morphism.
For simplicity,
we assume that $\nbigi_1\cup\nbigi_2$
is a good set of irregular values
in the level $(\vecm,i(0))$.
The induced morphism
$F_{|P}:\gbigv_{1|P}\lrarr\gbigv_{2|P}$ 
preserves the Stokes filtrations $\nbigf^P$.
\hfill\qed
\end{lem}

\begin{rem}
We considered two vector bundles 
on $\Xtilde(D)$.
One is $\pi^{-1}(E)$ 
and the other is $\gbigv$.
We should emphasize that they are different
in general.
The bundle $\gbigv$ depends only on
the flat bundle $(E,\nabla)_{|X-D}$,
and $\pi^{-1}(E)$ depends on
the prolongment $(E,\nabla)$.

Let us see the simplest example
$E=\nbigo\cdot e$
and $\nabla(e)=e\cdot d(z^{-1})$.
A trivialization of $\pi^{-1}(E)$ 
is given by $\pi^{-1}(e)$.
A trivialization of $\gbigv$ is induced by
$\exp(-z^{-1})\cdot e$.
\hfill\qed
\end{rem}

\subsection{Unramifiedly good lattices of
 a family of meromorphic $\lambda$-flat bundles}
 \label{subsection;07.12.14.30}

\subsubsection{Preliminary}
\label{subsubsection;08.7.30.40}
\paragraph{Good set of irregular values}

We use the partial order $\leq_{\seisuu^{n}}$
of $\seisuu^{n}$ given by
$\veca\leq_{\seisuu^n}\vecb
\Longleftrightarrow
 a_i\leq b_i,\,(\forall i)$.
We say $\veca<_{\seisuu^n}\vecb$
in the case $a_i<b_i$ for any $i$,
and 
we say $\veca\lneq_{\seisuu^n}\vecb$
in the case $\veca\leq_{\seisuu^n}\vecb$
and $\veca\neq\vecb$.
Let $\vecdelta_j$ denote the element
$(\overbrace{0,\ldots,0}^{j-1},1,0,\ldots,0)$,
and let $\veczero$
denote the zero in $\seisuu^n$.
We also use $\veczero_n$
when we distinguish the dependence on $n$.

\vspace{.1in}

Let $Y$ be a complex manifold.
Let $X:=\Delta^{\ell}\times Y$.
Let $D_i:=\{z_i=0\}\times Y$ and
$D:=\bigcup_{i=1}^{\ell}D_i$
be the hypersurfaces of $X$.
We also put
$D_{\ellsitabar}=\bigcap_{i=1}^{\ell}D_i$,
which is naturally identified with $Y$.

For any $f\in M(X,D)$,
we have the Laurent expansion:
\[
 f=\sum_{\vecm\in\seisuu^{\ell}}
f_{\vecm}\bigl(y\bigr)
 \cdot \vecz^{\vecm}.
\]
Here $f_{\vecm}$ are holomorphic functions
on $D_{\ellsitabar}$.
We often use the following identification
implicitly:
\begin{equation}
 \label{eq;08.10.1.5}
 M(X,D)\big/\vecz^{\vecn}\cdot H(X)\simeq
 \Bigl\{f\in M(X,D)\,\Big|\,
 f_{\vecm}=0,\,\,\forall\vecm\geq\vecn\Bigr\}
\end{equation}

For any $f\in M(X,D)$,
let $\ord(f)$ denote 
the minimum of the set
$\bigl\{\vecm\in\seisuu^{\ell}
 \,\big|\,f_{\vecm}\neq 0\bigr\}
 \cup \{\veczero\}$
with respect to
$\leq_{\seisuu^{\ell}}$,
if it exists.
It is always contained in $\seisuu_{\leq 0}^{\ell}$,
if it exists.

For any $\gminia\in M(X,D)/H(X)$,
we take any lift $\gminiatilde$ to $M(X,D)$,
and we set $\ord(\gminia):=\ord(\gminiatilde)$,
if the right hand side exists.
If $\ord(\gminia)$ exists in $\seisuu^{\ell}-\{\veczero\}$,
$\gminiatilde_{\ord(\gminia)}$ is independent
of the choice of a lift $\gminiatilde$,
which is denoted by $\gminia_{\ord(\gminia)}$.

\begin{df}
\label{df;07.6.1.15}
A finite subset $\nbigi\subset M(X,D)\big/H(X)$ 
is called a good set of irregular values on $(X,D)$,
if the following conditions are satisfied:
\begin{itemize}
\item
$\ord(\gminia)$ exists
for each $\gminia\in \nbigi$,
and
$\gminia_{\ord(\gminia)}$
is nowhere vanishing on $D_{\ellsitabar}$
for $\gminia\neq 0$.
\item
For any two distinct $\gminia,\gminib\in \nbigi$,
$\ord(\gminia-\gminib)$ exists in 
 $\seisuu_{\leq 0}^{\ell}-\{\veczero\}$,
and 
$(\gminia-\gminib)_{\ord(\gminia-\gminib)}$
is nowhere vanishing on $D_{\ellsitabar}$.
\item
The set 
$\nbigt(\nbigi):=\bigl\{\ord(\gminia-\gminib)\,\big|\,
 \gminia,\gminib\in \nbigi
 \bigr\}$
is totally ordered
with respect to the partial order on $\seisuu^{\ell}$.
\hfill\qed
\end{itemize}
\end{df}
The condition in Definition {\rm\ref{df;07.6.1.15}}
does not depend on the choice of a holomorphic coordinate
such that $D=\bigcup_{i=1}^{\ell}\{z_i=0\}$.

We will use the following lemma implicitly.
\begin{lem}
The set 
 $\bigl\{
 \ord(\gminia)\,\big|\,
 \gminia\in\nbigi
 \bigr\}$ is totally ordered.
In particular,
the minimum 
 \[
  \vecm(0):=\min\bigl\{\ord(\gminia)\,\big|\,
 \gminia\in \nbigi\bigr\}
 \]
exists.
Moreover,
$\vecm(0)\leq_{\seisuu^{\ell}}\vecm$
for any $\vecm\in\nbigt(\nbigi)$.
\end{lem}
\pf
Let $\gminia,\gminib\in\nbigi$.
Assume $\ord(\gminia)\not\leq\ord(\gminib)$
and $\ord(\gminia)\not\geq\ord(\gminib)$.
Then, $\ord(\gminia-\gminib)$
does not exist,
which contradicts the second condition.
Hence, we obtain the first claim of the lemma.
For any $\vecm\in\nbigt(\nbigi)$,
there exists $\gminia\in\nbigi$
such that $\gminia_{\vecm}\neq 0$.
Hence, $\vecm(0)\leq_{\seisuu^{\ell}}\vecm$.
\hfill\qed

\begin{rem}
\label{rem;07.11.15.10}
It is often convenient to
use a coordinate such that
$\nbigt(\nbigi)
 \cup\{\vecm(0)\}
 \subset \coprod_{i=0}^{\ell}
 \seisuu_{<0}^i\times
 \veczero_{\ell-i}$.
\hfill\qed
\end{rem}

\paragraph{Auxiliary sequence}

Let $\nbigi$ be a good set of irregular values on $(X,D)$.
Since the set $\nbigt(\nbigi)$
is totally ordered with respect to
the partial order $\leq_{\seisuu^{\ell}}$,
we can take a sequence
\[
 \nbigm:=\bigl(
 \vecm(0),\vecm(1),\vecm(2),\ldots,
 \vecm(L),\vecm(L+1)
 \bigr)
\subset\seisuu_{\leq 0}^{\ell}
\]
with the following property:
\begin{itemize}
\item $\nbigt(\nbigi)\subset \nbigm$
and $\vecm(L+1)=\veczero_{\ell}$.
\item
 We have
$1\leq \gminih(i)\leq \ell$ such that
 $\vecm(i+1)=\vecm(i)+\vecdelta_{\gminih(i)}$
 for each $i\leq L$.
\end{itemize}
Such a sequence is called an auxiliary sequence for $\nbigi$.
It is not uniquely determined for $\nbigi$.
It is convenient for an inductive argument.

\paragraph{Truncation}

Let $\nbigi$ be a good set of irregular values.
We take an auxiliary sequence 
for $\nbigi$,
and let $\etabar_{\vecm(0)}:
 \nbigi\lrarr M(X,D)/H(X)$ be given as follows:
\[
 \etabar_{\vecm(0)}(\gminia):=
  \sum_{\vecn\not\geq \vecm(1)}
 \gminia_{\vecn}\cdot \vecz^{\vecn}
\]
Then, the image is 
a good set of irregular values 
in the level $(\vecm(0),i(0))$.
More generally,
$\etabar_{\vecm(j)}$ is defined as follows:
\[
 \etabar_{\vecm(j)}(\gminia):=
 \sum_{\vecn\not\geq\vecm(j)}
 \gminia_{\vecn}\cdot\vecz^{\vecn}
\]
We have $\etabar_{\vecm(L)}(\gminia)=\gminia$.
We set $\zeta_{\vecm(0)}(\gminia):=
 \etabar_{\vecm(0)}(\gminia)$
and $\zeta_{\vecm(j)}(\gminia):=
 \etabar_{\vecm(j)}(\gminia)-\etabar_{\vecm(j-1)}(\gminia)$
for $j=1,\ldots,L$.
Then, we have the decomposition
$\etabar_{\vecm(i)}(\gminia)=
 \sum_{j\leq i} \zeta_{\vecm(j)}(\gminia)$.

\vspace{.1in}
Let $\nbigi(\vecm(i))$ denote 
the image of 
$\etabar_{\vecm(i)}:\nbigi\lrarr
 M(X,D)/H(X)$.

\begin{lem}
If we shrink $X$ appropriately,
$\nbigi(\vecm(0))$ is a good set of 
irregular values in the level
$\bigl(\vecm(0),\gminih(0)\bigr)$.
\end{lem}
\pf
If $\etabar_{\vecm(0)}(\gminia-\gminib)\neq 0$
for $\gminia,\gminib\in\nbigi$,
we have 
$\ord(\gminia-\gminib)=\vecm(0)$
and 
$\bigl(
 \vecz^{-\vecm(0)}
 \etabar_{\vecm(0)}(\gminia-\gminib)
 \bigr)_{|D_{\ellsitabar}}$
is nowhere vanishing.
Hence,
$\bigl(
 \vecz^{-\vecm(0)}
 \etabar_{\vecm(0)}(\gminia-\gminib)\bigr)$
is nowhere vanishing on $X$
after $X$ is shrinked appropriately.
Similarly,
we may have 
$\bigl(
 \vecz^{-\vecm(0)}
 \etabar_{\vecm(0)}(\gminia)\bigr)$
is nowhere vanishing on $X$
after $X$ is shrinked appropriately.
\hfill\qed

\vspace{.1in}
We can use the following lemma
for inductive arguments.
\begin{lem}
For any $\gminib\in \nbigi(\vecm(0))$,
we fix any element
$\gminia^{(0)}\in\etabar_{\vecm(0)}^{-1}(\gminib)$.
Then, the set 
\[
 \bigl\{
 \gminia-\gminia^{(0)}\,\big|\,
 \etabar_{\vecm(0)}(\gminia)=\gminib
 \bigr\}
\]
is also a good set of irregular values.
\hfill\qed
\end{lem}

\paragraph{Example}
We give some examples.
\[
\gminia^{(1)}:=z_1^{-1}\cdot z_2^{-1},
\quad
\gminia^{(2)}:=z_1^{-1},
\quad
\gminia^{(3)}:=0.
\]
An auxiliary sequence is unique in this case,
and given as follows:
\begin{equation}
\label{eq;08.10.1.2}
 \vecm(0)=(-1,-1),\,\,
 \gminih(0)=2,\quad
 \vecm(1)=(-1,0),\,\,
 \gminih(1)=1,\quad
 \vecm(2)=(0,0)
\end{equation}
The truncations are given as follows:
\[
 \etabar_{\vecm(0)}(\gminia^{(1)})=\gminia^{(1)},
\quad
 \etabar_{\vecm(0)}(\gminia^{(2)})=0,
\quad
 \etabar_{\vecm(0)}(\gminia^{(3)})=0
\]
\[
 \etabar_{\vecm(1)}(\gminia^{(1)})
=\gminia^{(1)},
\quad
 \etabar_{\vecm(1)}(\gminia^{(2)})
=\gminia^{(2)},
\quad
 \etabar_{\vecm(1)}(\gminia^{(3)})
=\gminia^{(3)}
\]
The image of $\nbigi$ via
$\etabar_{\vecm(0)}$ is
$\bigl\{ \gminia^{(1)},0 \bigr\}$.

Let us consider the following set:
\[
 \gminib^{(1)}=z_1^{-1}\cdot z_2^{-1}
 +a\cdot z_2^{-1}+b\cdot z_1^{-1},
\quad
 \gminib^{(2)}=z_1^{-1}
\]
An auxiliary sequence is given by (\ref{eq;08.10.1.2}).
The truncation is given as follows:
\[
 \etabar_{\vecm(0)}(\gminib^{(1)})
=z_1^{-1}\cdot z_2^{-1}+a\cdot z_2^{-1},
\quad
 \etabar_{\vecm(0)}(\gminib^{(2)})=0
\]

\vspace{.2in}
We have the following picture in our mind
for truncation.

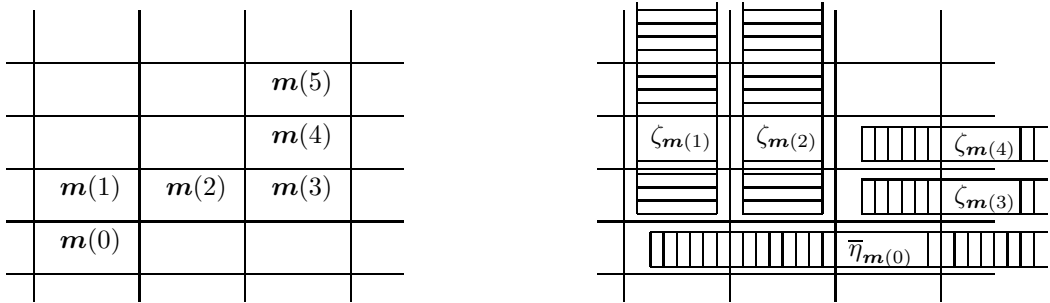
\begin{picture}(220,120)(-100,-80)
\put(0,0){$\vecm(5)$}
\put(0,-20){$\vecm(4)$}
\put(0,-40){$\vecm(3)$}
\put(-40,-40){$\vecm(2)$}
\put(-80,-40){$\vecm(1)$}
\put(-80,-60){$\vecm(0)$}

\put(-100,10){\line(1,0){150}}
\put(-100,-10){\line(1,0){150}}
\put(-100,-30){\line(1,0){150}}
\put(-100,-50){\line(1,0){150}}
\put(-100,-70){\line(1,0){150}}

\put(-90,-80){\line(0,1){110}}
\put(-50,-80){\line(0,1){110}}
\put(-10,-80){\line(0,1){110}}
\put(30,-80){\line(0,1){110}}

\end{picture}
\begin{picture}(150,120)(-100,-80)
\put(-100,10){\line(1,0){150}}
\put(-100,-10){\line(1,0){150}}
\put(-100,-30){\line(1,0){150}}
\put(-100,-50){\line(1,0){150}}
\put(-100,-70){\line(1,0){150}}

\put(-90,-80){\line(0,1){110}}
\put(-50,-80){\line(0,1){110}}
\put(-10,-80){\line(0,1){110}}
\put(30,-80){\line(0,1){110}}

\put(-80,-67){\line(1,0){150}}
\put(-80,-54){\line(1,0){150}}
\multiput(-80,-67)(5,0){15}{\line(0,1){13}}
\multiput(25,-67)(5,0){9}{\line(0,1){13}}
\put(-5,-63){$\etabar_{\vecm(0)}$}

\put(-85,-47){\line(0,1){80}}
\put(-55,-47){\line(0,1){80}}
\multiput(-85,-47)(0,5){5}{\line(1,0){30}}
\multiput(-85,-5)(0,5){8}{\line(1,0){30}}
\put(-80,-20){$\zeta_{\vecm(1)}$}

\put(-45,-47){\line(0,1){80}}
\put(-15,-47){\line(0,1){80}}
\multiput(-45,-47)(0,5){5}{\line(1,0){30}}
\multiput(-45,-5)(0,5){8}{\line(1,0){30}}
\put(-40,-20){$\zeta_{\vecm(2)}$}

\put(0,-47){\line(1,0){70}}
\put(0,-34){\line(1,0){70}}
\multiput(0,-47)(5,0){6}{\line(0,1){13}}
\multiput(60,-47)(5,0){2}{\line(0,1){13}}
\put(35,-43){$\zeta_{\vecm(3)}$}

\put(0,-27){\line(1,0){70}}
\put(0,-14){\line(1,0){70}}
\multiput(0,-27)(5,0){6}{\line(0,1){13}}
\multiput(60,-27)(5,0){2}{\line(0,1){13}}
\put(35,-23){$\zeta_{\vecm(4)}$}

\end{picture}

\begin{center}
$L=4$,
$\vecm(0)=(-2,-3)$,
$\vecm(1)=(-2,-2)$,
$\vecm(2)=(-1,-2)$, \\
$\vecm(3)=(0,-2)$,
$\vecm(4)=(0,-1)$,
$\vecm(5)=(0,0)$.
\end{center}

\subsubsection{Unramifiedly good lattices of
 a family of meromorphic $\lambda$-flat bundles}
\label{subsubsection;08.8.3.3}

Let $X$ be a complex manifold,
and let $D$ be a normal crossing divisor of $X$.
Let $\nbigk$ be a point or
a compact region in $\cnum_{\lambda}$.
Let $\nbigx$ and $\nbigd$ denote $\nbigk\times X$
and $\nbigk\times D$, respectively.
For $\lambda\in\nbigk$,
we set $\nbigx^{\lambda}:=\{\lambda\}\times X$
and $\nbigd^{\lambda}:=\{\lambda\}\times D$.
Let $(\nbige,\DD)$ be
a family of meromorphic $\lambda$-flat bundles
on $(\nbigx,\nbigd)$,
i.e.,
$\nbige$ is an $\nbigo_{\nbigx}(\ast \nbigd)$-coherent sheaf
with a holomorphic family of flat $\lambda$-connections
$\DD:\nbige\lrarr
 \nbige\otimes \Omega^1_{\nbigx/\nbigk}$.
The restriction to 
$(\nbigx^{\lambda},\nbigd^{\lambda})$
is denoted by
$(\nbigelambda,\DDlambda)$.

\begin{rem}
If $\nbigk$ is a point,
``family'' can be omitted.
\hfill\qed
\end{rem}

Let $E$ be an $\nbigo_{\nbigx}$-locally free lattice 
of $(\nbige,\DD)$.
Let $P$ be any point of $\nbigd$.
We can take a holomorphic coordinate
$(\nbigu,\lambda,z_1,\ldots,z_n)$ around $P$
such that
$\nbigd_{\nbigu}:=
 \nbigd\cap \nbigu
=\bigcup_{i=1}^{\ell}\nbigd_{\nbigu,i}$,
where $\nbigd_{\nbigu,i}:=\{z_i=0\}$.
We put 
 $\nbigd_{\nbigu,I}:=
 \bigcap_{i\in I}\nbigd_{\nbigu,i}$
and 
$\nbigd_{\nbigu}(I):=
 \bigcup_{i\in I}\nbigd_{\nbigu,i}$.
For any subset $I\subset\ellsitabar$,
we put $I^c:=\ellsitabar-I$.
The completion of $\nbigx$ along
$\nbigd_{\nbigu,I}$ (resp. $\nbigd_{\nbigu}(I)$)
is denoted by $\nbigdhat_{\nbigu,I}$
(resp. $\nbigdhat_{\nbigu}(I)$).

\begin{df}
\label{df;07.10.8.1}
We say that $E$ is unramifiedly good at $P$,
if the following holds:
\begin{itemize}
\item
We are given a good set of irregular values
$S\subset 
 M(\nbigu,\nbigd_{\nbigu})/H(\nbigu)$.
\item
For any $\emptyset\neq I\subset\ellsitabar$,
we have the decomposition:
\begin{equation}
 \label{eq;07.10.8.2}
 (E,\DD)_{|\nbigdhat_{\nbigu,I}}=
 \bigoplus_{\gminia\in S(I)}
 \bigl(\lefttop{I}\Ehat_{\gminia},
 \lefttop{I}\DDhat_{\gminia}\bigr)
\end{equation}
 Here $S(I)$ denotes the image of $S$
 via the map
 $M(\nbigu,\nbigd_{\nbigu})/H(\nbigu)
\lrarr M(\nbigu,\nbigd_{\nbigu})/
 M\bigl(\nbigu,\nbigd_{\nbigu}(I^c)\bigr)$.
\item
 $(\DD_{\gminia}-d\gminia)
 \bigl(\lefttop{I}\Ehat_{\gminia}\bigr)$
is contained in
 $\lefttop{I}\Ehat_{\gminia}
 \otimes 
 \Bigl(
 \Omega^1_{\nbigx/\nbigk}
 \bigl(\log \nbigd_{\nbigu}(I)\bigr)
+\Omega^1_{\nbigx/\nbigk}
 \bigl(\ast \nbigd_{\nbigu}(I^c)\bigr)
 \Bigr)$,
where $\gminia$ is lifted to
$M(\nbigu,\nbigd_{\nbigu})$.
This condition is independent of
the choice of a lift.
\end{itemize}
The property is independent of the choice of
the coordinate $(\nbigu,\lambda,z_1,\ldots,z_n)$.

We say that $(E,\DD)$ is unramifiedly good,
if it is unramifiedly good at any point.
\hfill\qed
\end{df}
See Subsection 5.7 of \cite{mochi7}
for another but equivalent formulation,
which seems easier to state.

The decomposition (\ref{eq;07.10.8.2})
is called the irregular decomposition
of $E_{|\nbigdhat_{\nbigu,I}}$.
The set $S$ is uniquely determined
if $\lefttop{\ellsitabar}E_{\gminia}\neq 0$
for each $\gminia\in S$.
So, it is denoted by $\Irr(\DD,P)$.
The restriction of $E$ to 
$\{\lambda\}\times X$
is denoted by $E^{\lambda}$.

If $E$ is an unramifiedly good lattice of 
$(\nbige,\DD)$,
we have the well defined endomorphism
$\Res_i(\DD)$ of $E_{|\nbigd_{i}}$ for 
each irreducible component $D_i$ of $D$.
It is called the residue of $\DD$ at $D_i$
with respect to the lattice $E$.
If $\nbigk\neq\{0\}$,
the eigenvalues of $\Res_{i}(\DD)$
are constant on $\nbigd_i^{\lambda}$
for each $\lambda\in\nbigk$.
(See Subsection 5.1.3 of \cite{mochi7},
for example.)

\begin{rem}
We have the notion of good lattice
which is locally a descent of
an unramifiedly good lattice.
See {\rm\cite{mochi7}}.
See also Definition
{\rm \ref{df;08.10.1.3}}
below.
\hfill\qed
\end{rem}

\paragraph{Irregular decompositions 
 in the level $\vecm(j)$}

In the following,
let $X:=\Delta^n$, $D_i:=\{z_i=0\}$
and $D:=\bigcup_{i=1}^{\ell}D_i$.
We set 
$\nbigd(\leq p):=\bigcup_{i\leq p} \nbigd_{i}$.
Let $(E,\DD)$ be an unramifiedly good lattice of
a family of meromorphic $\lambda$-flat bundles
on $(\nbigx,\nbigd)$
with the good set $\Irr(\DD)=\Irr(\DD,O)$.
We assume that the coordinate is 
as in Remark \ref{rem;07.11.15.10}
for $\Irr(\DD)$.
Let $\Irr(\DD,p)$ and $\Irr'(\DD,p)$
denote the image of $\Irr(\DD)$
by the natural maps
\[
\pi_p:M(\nbigx,\nbigd)/H(\nbigd)
 \lrarr M(\nbigx,\nbigd)/M(\nbigx,\nbigd(\leq p-1)),
\quad
\pi_p':M(\nbigx,\nbigd)/H(\nbigd)
 \lrarr M(\nbigx,\nbigd)/M(\nbigx,\nbigd(\neq p)).
\]
Note that the naturally induced map
$\Irr(\DD,p)\lrarr \Irr'(\DD,p)$
is bijective,
via which we identify them.

Take an auxiliary  sequence 
$\vecm(0),\ldots,\vecm(L)$
for the good set $\Irr(\DD)$.
Let $\Irrbar(\DD,\vecm(0))$ denote the image of
$\Irr(\DD)$ via $\etabar_{\vecm(0)}$.
Let $k(j)$ denote the number determined by
$\vecm(j)\in\seisuu_{<0}^{k(j)}
 \times\veczero_{\ell-k(j)}$.
For $p\leq k(j)$,
we have the map
$\Irr(\DD,p)\lrarr 
 M(\nbigx,\nbigd)/M(\nbigx,\nbigd(\leq p-1))$
induced by $\etabar_{\vecm(j)}$
which is denoted by
$\etabar_{\vecm(j),p}$.

By using a lemma in Subsection 5.1.2 of \cite{mochi7}
and the uniqueness of the decompositions,
we obtain the following decomposition 
on the completion $\nbigdhat(\leq k(j))$
along $\nbigd(\leq k(j))$:
\begin{equation}
\label{eq;07.12.17.3}
(E,\DD)_{|\nbigdhat(\leq k(j))}
=\bigoplus_{\gminib\in \Irrbar(\DD,\vecm(j))}
 \bigl(\Ehat^{\vecm(j)}_{\gminib},
 \DD_{\gminib}\bigr),
\quad
\mbox{\rm where }\,\,
 \Ehat^{\vecm(j)}_{\gminib|\nbigdhat_{p}}=
 \bigoplus_{\substack{
 \gminic\in \Irr(\DD,p)\\
 \etabar_{\vecm(j),p}(\gminic)=\pi_p(\gminib)
 }}
\lefttop{p}\Ehat_{\gminic},
\,\,\, \bigl(p\leq k(j)\bigr)
\end{equation}
The decomposition (\ref{eq;07.12.17.3})
is called the irregular decomposition
in the level $\vecm(j)$.

\begin{rem}
We do not have the irregular decomposition
in the level $\vecm(j)$ on $\Dhat$ in general,
which Sabbah remarked in {\rm\cite{sabbah4}}
for the surface case.
\hfill\qed
\end{rem}

\paragraph{The associated graded bundles
 with the family of meromorphic flat $\lambda$-connections}

Assume $\nbigk\neq\{0\}$.
We set $W:=\nbigx^0\cup\nbigd(\leq k(0))$.
It is easy to observe that
$(E,\DD)$ is an unramifiedly good lattice 
in the level $(\vecm(0),i(0))$
with the decomposition (\ref{eq;07.12.17.3})
for $j=0$.
The set of the irregular values in the level 
$\bigl(\vecm(0),i(0)\bigr)$
is $\Irrbar(\DD,\vecm(0))$.

As stated in Subsection \ref{subsubsection;08.8.3.2},
we obtain the holomorphic bundle
$\Gr^{\vecm(0)}_{\gminia}(E)$
with a family of meromorphic flat $\lambda$-connections
$\DD_{\gminia}^{\vecm(0)}$
on $(\nbigv,\nbigv\cap \nbigd)$ for each
$\gminia\in\Irrbar(\DD,\vecm(0))$,
where $\nbigv$ denotes a neighbourhood
of $\bigcap_{1\leq i\leq k(0)}\nbigd_i$.
Let 
$\Gr^{\vecm(0)}_{\gminia}(E,\DD):=
 \bigl(\Gr_{\gminia}^{\vecm(0)}(E),
 \DD^{\vecm(0)}_{\gminia}\bigr)$.
We obtain the following isomorphisms
for any $\gminia\in\Irrbar(\DD,\vecm(0))$
from (\ref{eq;08.10.1.4}):
\[
 \Gr^{\vecm(0)}_{\gminia}(E,\DD)_{|\What}
\simeq
  (\Ehat^{\vecm(0)}_{\gminia},\DD_{\gminia})
\]
In particular,
$\Gr^{\vecm(0)}_{\gminia}(E,\DD)$
are unramifiedly good lattices
whose set of irregular values is
$\Irr(\DD^{\vecm(0)}_{\gminia})=
 \etabar_{\vecm(0)}^{-1}(\gminia)$.

\vspace{.1in}

Let $\Irrbar(\DD,\vecm(j))$
denote the image of 
$\etabar_{\vecm(j)}:\Irr(\DD)
\lrarr M(X,D)/H(X)$ for any $j$.
Let us consider the case in which
$\Irrbar(\DD,\vecm(j-1))$
consists of a unique element.
We take any element 
$\gminia^{(1)}\in\Irr(\DD)$.
Let $\nbigl(\pm\gminia^{(1)})$
be a line bundle
$\nbigo_{\nbigx}\cdot e$
with a family of
meromorphic flat $\lambda$-connections
$\DD e=e\cdot (\pm d\gminia^{(1)})$.
Then,
$ (E',\DD'):=(E,\DD)\otimes\nbigl(-\gminia^{(1)})$
is an unramifiedly good lattice
with the good set
\[
 \Irr(\DD')=\bigl\{
 \gminia-\gminia^{(1)}\,\big|\,
 \gminia\in\Irr(\DD)
 \bigr\}.
\]
The sequence
$\vecm(j),\vecm(j+1),\ldots,\vecm(L)$
gives an auxiliary sequence
for $\Irr(\DD')$.
By applying the above procedure to $(E',\DD')$
and shrinking $X$,
we obtain
$\Gr^{\vecm(j)}_{\gminic}(E',\DD')$
for each $\gminic\in\Irrbar(\DD',\vecm(j))$.
For any $\gminib\in\Irrbar(\DD,\vecm(j))$,
we define
\[
 \Gr^{\vecm(j)}_{\gminib}(E,\DD):=
 \Gr^{\vecm(j)}_{
 \gminib-\etabar_{\vecm(j)}(\gminia^{(1)})}
 (E',\DD')
\otimes\nbigl(\gminia^{(1)})
\]
It is independent of the choice of $\gminia^{(1)}$
up to canonical isomorphisms.
(We may avoid tensor products.)
It is easy to observe that
$\Gr^{\vecm(j)}_{\gminib}(E,\DD)$
are also unramifiedly good lattices
with the good sets of irregular values
$\Irr(\DD^{\vecm(j)}_{\gminib})
=\etabar_{\vecm(j)}^{-1}(\gminib)$.
By construction,
$\Irrbar\bigl(\DD^{\vecm(j)}_{\gminib},\vecm(j)\bigr)$
consists of the unique element $\gminib$.

\vspace{.1in}

Let us consider the general case.
Let $\etabar_{\vecm(j-1),\vecm(j)}:
 \Irrbar(\DD,\vecm(j))\lrarr \Irrbar(\DD,\vecm(j-1))$
be the induced map.
For any $\gminia\in\Irrbar(\DD,\vecm(j))$,
we inductively define 
\[
 \Gr^{\vecm(j)}_{\gminia}(E,\DD):=
 \Gr^{\vecm(j)}_{\gminia}
 \Gr^{\vecm(j-1)}_{
 \etabar_{\vecm(j-1),\vecm(j)}(\gminia)}
 (E,\DD)
\]
For each $\gminia\in\Irr(\DD)$,
we set 
$\Gr^{\full}_{\gminia}(E,\DD):=
 \Gr^{\vecm(L)}_{\gminia}(E,\DD)$,
which is called the full reduction.
By construction,
$\Gr^{\full}_{\gminia}(E,\DD)
 \otimes\nbigl(-\gminia)$
is logarithmic.

We have the functoriality 
as in Subsection \ref{subsubsection;08.8.3.2}.

\paragraph{Deformation}

Assume $0\not\in\nbigk$.
We would like to regard $(E,\DD)$
as a prolongment of
$(E,\DD)_{|\nbigx-\nbigd(\leq k(0))}$.
For a given holomorphic function  $T(\lambda)$
with $\Re\bigl(T(\lambda)\bigr)>0$,
we have the other prolongment
$(E^{(T)},\DD^{(T)})$ of 
$(E,\DD)_{|\nbigx-\nbigd(\leq k(0))}$,
which is also an unramifiedly good lattice
with the set of irregular values
\[
 \Irr\bigl(E^{(T)},\DD^{(T)}\bigr):=
 \bigl\{
 T\cdot\gminia\,\big|\,
 \gminia\in\Irr(\DD)
 \bigr\}.
\]
We refer to Subsections 7.5, 7.8--7.9 of \cite{mochi7}
for the construction.
We mention some properties
(Subsection 7.8 of \cite{mochi7}):
\begin{description}
\item[(D1):]
$E^{(T_1\cdot T_2)}
\simeq
 \bigl(E^{(T_1)}\bigr)^{(T_2)}$,
if $\Re(T_i)>0$ and $\Re(T_1\cdot T_2)>0$.
\item[(D2):]
$(E^{(T)},\DD^{(T)})_{|\nbigdhat_{I_0}}
\simeq
 \bigoplus_{\gminia\in\nbigi}
 \bigl(
 \lefttop{I_0}
 \Ehat_{\gminia},
 \lefttop{I_0}
 \DDhat_{\gminia}+(T-1)d\gminia\bigr)$,
where $I_0:=\{1,\ldots,k(0)\}$.
In other brief words,
the deformation does not change the regular part.
\end{description}

We give some statements for functoriality.
See Subsection 7.8.1 of \cite{mochi7}
for more details.
Let $(E_p,\DD_p)$ $(p=1,2)$
be unramifiedly good.
We have the following natural isomorphisms:
\[
 (E_1\oplus E_2)^{(T)}
\simeq
 E_1^{(T)}\oplus E_2^{(T)},
\quad
  (E_1\otimes E_2)^{(T)}
\simeq
 E_1^{(T)}\otimes E_2^{(T)},
\quad
\bigl( E^{\lor}\bigr)^{(T)}
\simeq
 (E^{(T)})^{\lor}
\]
Here, we have assumed that
$(E_1,\DD_1)\oplus (E_2,\DD_2)$
and $(E_1,\DD_1)\otimes(E_2,\DD_2)$
are unramifiedly good.
Moreover,
let $F:(E_1,\DD_1)\lrarr (E_2,\DD_2)$
be a flat morphism.
Assume $\nbigi_1\cup\nbigi_2$ is a good set of
irregular values in the level $(\vecm,i(0))$.
Then, we have the naturally induced morphism
$(E_1^{(T)},\DD_1^{(T)})
\lrarr 
 (E_2^{(T)},\DD_2^{(T)})$.

\subsection{Smooth divisor case}
\label{subsubsection;08.8.10.1}

Let $X:=\Delta^n$
and $D:=\{z_1=0\}$.
Let $\nbigk\subset\cnum_{\lambda}$.
Let $(E,\DD)$ be an unramifiedly good lattice of
a family of meromorphic $\lambda$-bundles
on $(\nbigx,\nbigd)$ with a good set of irregular values
$\Irr(\DD)=\Irr(\DD,O)$.
We have the formal decomposition
$ (E,\DD)_{|\nbigdhat}
=\bigoplus_{\gminia\in \Irr(\DD)}
 (\Ehat_{\gminia},\DDhat_{\gminia})$,
where 
$\DDhat_{\gminia}-d\gminia\cdot\id_{\Ehat_{\gminia}}$
are logarithmic.
We set $W:=\nbigd\cup \nbigx^0$ 
in the case $0\in \nbigk$,
and $W:=\nbigd$ otherwise.
We obtain the decomposition on $\What$:
\begin{equation}
\label{eq;08.8.2.3}
 (E,\DD)_{|\What}
=\bigoplus_{\gminia\in \Irr(\DD)}
 (\Ehat_{\gminia},\DDhat_{\gminia})
\end{equation}

\paragraph{Full Stokes filtration}

In this case,
it is also easy and convenient 
to consider full Stokes filtration.
We explain it in the following.
Let $\pi:\nbigxtilde(W)\lrarr\nbigx$ denote 
the real blow up of $\nbigx$ along $W$.
We put $\gbigz:=\pi^{-1}(\nbigd)$.

For any multi-sector $S$ in $\nbigx-W$,
the order $\leq_S$ on $\Irr(\DD)$ is defined
as follows:
\begin{itemize}
\item
$\gminia\leq_S\gminib$
if and only if
$-\Re\bigl(\lambda^{-1}\gminia(\lambda,\vecz)\bigr)
 \leq_S
 -\Re\bigl(\lambda^{-1}\gminib(\lambda,\vecz)\bigr)$
for any $\vecz\in S$ such that $|z_1|$ is sufficiently small.
\end{itemize}
Let $\Sbar$ denote the closure of $S$
in $\nbigxtilde(W)$,
and let $Z$ denote $\Sbar\cap \pi^{-1}(W)$.
The irregular decomposition (\ref{eq;08.8.2.3})
on $\What$ induces the decomposition on $\Zhat$:
\begin{equation}
 \label{eq;08.8.2.4}
 (E,\DD)_{|\Zhat}=
 \bigoplus_{\gminia\in \Irr(\DD)}
 (\Ehat_{\gminia},\DDhat_{\gminia})_{|\Zhat}
\end{equation}
We put
$\nbigf^Z_{\gminia}:=
 \bigoplus_{\gminib\leq_S\gminia}
 \Ehat_{\gminib|\Zhat}$,
and then we obtain the filtration
$\nbigf^Z$ indexed by
$\bigl(\Irr(\DD),\leq_S\bigr)$.
By using Proposition \ref{prop;07.6.16.6}
and Lemma \ref{lem;07.12.16.9} successively
(or by using more classical results),
we obtain the following.

\begin{prop}
\mbox{{}}
For any point $P\in \gbigz$,
there exists $S\in \Multisector(P,\nbigx-W)$
such that the following holds:
\begin{itemize}
\item
There exists the unique $\DD$-flat filtration
$\nbigftilde^S$ of $E_{|\Sbar}$
on $\Sbar$
indexed by $\bigl(\Irr(\DD),\leq_S\bigr)$
such that
$\nbigftilde^S_{|\Zhat}=\nbigf^Z$.
\item
There exists a $\DD$-flat splitting of
$\nbigftilde^{S}$ on $\Sbar$.
\end{itemize}
We call $\nbigftilde^S$ the full Stokes filtration of $(E,\DD)$.

For $S'\subset S$,
the filtrations $\nbigftilde^{S'}$ and $\nbigftilde^S$
satisfy the compatibility condition
as in Lemma {\rm\ref{lem;07.12.16.9}}.
\hfill\qed
\end{prop}
The following lemma is clear from the definition 
of full Stokes filtrations.
\begin{lem}
Let $S,S'\in \Multisector(P,\nbigx-W)$.
Assume 
(i) $S'\subset S$,
(ii) $E_{|\Sbar}$ has 
 the full Stokes filtration $\nbigftilde^S$
 as above.
Then, the restriction of $\nbigftilde^S$
 to $\Sbar'$ is the full Stokes filtration 
 of $E_{|\Sbar'}$.
\hfill\qed
\end{lem}

We have functoriality of full Stokes filtrations
as in the case of Stokes filtrations
in the level $(\vecm,i(0))$.

\paragraph{The associated graded bundle}

For any sectors $S$ and each 
$\gminia\in\Irr(\DD)$,
we obtain the bundle 
$\Gr_{\gminia}^{\full}(E_{|S})$ on $S$
associated to the full Stokes filtration $\nbigftilde^S$.
By varying $S$ and gluing 
$\Gr^{\full}_{\gminia}(E_{|\Sbar})$,
we obtain the bundle
$\Gr^{\full}_{\gminia}(E_{|\nbigvtilde(W)})$
on $\nbigvtilde(W)$
with the induced family of flat $\lambda$-connections 
$\DD_{\gminia}$,
where $\nbigv$ denotes some neighbourhood of $\nbigd$,
and $\nbigvtilde(W)$ denote the real blow up
of $\nbigv$ along $W\cap\nbigv$.
As in Subsection \ref{subsubsection;08.8.3.2},
we can show that
$\Gr^{\full}_{\gminia}(E_{|\nbigvtilde(W)})$
has the descent to $\nbigv$,
i.e.,
there exists a locally free sheaf 
$\Gr^{\full}_{\gminia}(E)$ on $\nbigv$
with a family of meromorphic flat
$\lambda$-connections $\DD_{\gminia}$,
such that 
\[
\pi^{-1}\bigl(
 \Gr^{\full}_{\gminia}(E),\DD_{\gminia}
\bigr)
\simeq
\bigl(
\Gr^{\full}_{\gminia}(E_{|\nbigvtilde(W)}),
 \DD_{\gminia}
\bigr),
\quad\quad
(\Gr^{\full}_{\gminia}(E),
 \DD_{\gminia})_{|\What\cap\nbigv}
\simeq 
\bigl(\Ehat_{\gminia},
 \DD_{\gminia}\bigr)_{|\What\cap\nbigv}.
\]
By construction,
$\DD_{\gminia}-d\gminia$ is logarithmic
for each $\gminia\in\Irr(\DD)$.

As in the case of Gr with respect to Stokes filtrations
in the level $(\vecm,i(0))$,
we have the following isomorphisms:
\[
 \Gr^{\full}_{\gminia}(E^{\lor})
\simeq
 \Gr^{\full}_{-\gminia}(E)^{\lor},
\]
\[
   \Gr_{\gminia}^{\full}(E_1\otimes E_2)
\simeq
\bigoplus_{\substack{
 \gminia_i\in \Irr(\DD_i)
 \\
 \gminia_1+\gminia_2=\gminia}}
 \Gr^{\full}_{\gminia_1}(E_1)
\otimes
 \Gr^{\full}_{\gminia_2}(E_2),
\]
\[
  \Gr^{\full}_{\gminia}(E_1\oplus E_2)
\simeq
 \Gr^{\full}_{\gminia}(E_1)
\oplus
 \Gr^{\full}_{\gminia}(E_2).
\]
Here, we have assumed that
$(E_1,\DD_1)\otimes (E_2,\DD_2)$
and $(E_1,\DD_1)\oplus (E_2,\DD_2)$
are unramifiedly good lattices.

\begin{lem}
Let $(E_p,\DD_p)$ $(p=1,2)$ be unramifiedly
good lattices on $(\nbigx,\nbigd)$.
Assume $\nbigi_1\cup\nbigi_2$ is a good set of
irregular values.
Let $F:(E_1,\DD_1)\lrarr (E_2,\DD_2)$
be a flat morphism.
We have the naturally induced morphism
$\Gr^{\full}_{\gminia}(F):
 \Gr^{\full}_{\gminia}(E_1)
\lrarr \Gr^{\full}_{\gminia}(E_2)$.
\hfill\qed
\end{lem}

\paragraph{A characterization of sections of $E$}

Let $\vecw_{\gminia}$ be a frame of
$\Gr^{\full}_{\gminia}(E)$.
Let $S$ be a small multi-sector,
and let $E_{|\Sbar}=\bigoplus E_{\gminia,S}$
be a $\DD$-flat splitting
of the full Stokes filtration $\nbigftilde^S$.
By the natural isomorphism
$E_{\gminia,S}\simeq
 \Gr_{\gminia}^{\full}(E)_{|\Sbar}$,
we take a lift $\vecw_{\gminia,S}$ of $\vecw_{\gminia}$.
Thus, we obtain the frame 
$\vecw_S=\bigl(\vecw_{\gminia,S}\bigr)$
of $E_{|\Sbar}$.
The following proposition implies a
characterization of sections of $E$
by growth order with respect
to the frames $\vecw_S$
for small multi-sectors $S$.

\begin{prop}
\label{prop;08.9.11.3}
Let $\vecv$ be a frame of $E$,
and let $G_S$ be determined by
$\vecv=\vecw_S\cdot G_S$.
Then, $G_S$ and $G_S^{-1}$ are bounded on $S$.
\hfill\qed
\end{prop}

\paragraph{Deformation}

When $|\arg(T)|$ is sufficiently small,
we have a more direct local construction of the deformation
$(E,\DD)^{(T)}$.
We explain it in the smooth divisor case.

We take a covering 
$\nbigx-\nbigd=\bigcup_{i=1}^N S^{(i)}$
by sectors $S^{(i)}$
on which we have the full Stokes filtrations.
Assume that $\bigl|\arg(T)\bigr|$ is sufficiently small
such that the following holds:
\begin{itemize}
 \item
 $\gminia\leq_{S^{(i)}}\gminib$
$\Longleftrightarrow$
 $T\gminia\leq_{S^{(i)}}T\gminib$
 for any $\gminia,\gminib\in \Irr(\DD)$
 and for any $S^{(i)}$.
\end{itemize}
We take frames $\vecw_{\gminia}$
of $\Gr^{\full}_{\gminia}(E)$.
For each $S=S^{(i)}$,
we take a $\DD$-flat splitting
$E_{|S}=\bigoplus E_{\gminia,S}$
of the full Stokes filtration.
Let $\vecw_{S}=(\vecw_{\gminia,S})$
be as above.
We put 
$\vecw_{\gminia,S}^{(T)}:=
 \vecw_{\gminia,S}\cdot
 \exp\bigl((T-1)\cdot\lambda^{-1}\cdot\gminia\bigr)$
and
$\vecw_{S}^{(T)}:=
 \bigl(\vecw_{\gminia,S}^{(T)}\bigr)$.
Let $f$ be a holomorphic section of $E_{|\nbigx-\nbigd}$.
We have the corresponding decomposition
$f=\sum f_{\gminia,S}$ on each $S$.
We have the expression
$f_{\gminia,S}=
 \sum f^{(T)}_{\gminia,S,j}\cdot w^{(T)}_{\gminia,S,j}$.
We put $\vecf_{\gminia,S}:=\bigl(f^{(T)}_{\gminia,S,j}\bigr)$.
\begin{lem}
\label{lem;07.10.2.21}
$f$ gives a section of $E^{(T)}$
if and only if
$\vecf^{(T)}_{\gminia,S^{(i)}}$
is bounded for each $S^{(i)}$ and $\vecw_{S^{(i)}}$.
(See Subsection {\rm 7.9.1}
 of {\rm\cite{mochi7}}.)
\hfill\qed
\end{lem}

\paragraph{Prolongation of a flat morphism}

Let $(E_p,\DD_p)$ $(p=1,2)$ be unramifiedly
good lattices on $(\nbigx,\nbigd)$.
Assume $\Irr(\DD_1)\cup\Irr(\DD_2)$ is a good set of
irregular values.
Let $F:(E_1,\DD_1)_{|\nbigx-\nbigd}
 \lrarr (E_2,\DD_2)_{|\nbigx-\nbigd}$
be a flat morphism.

\begin{lem}
\label{lem;08.9.11.4}
If $F$ preserves the full Stokes filtrations
$\nbigftilde^S$ for each small sector $S$,
$F$ is extended to the meromorphic morphism
$F:E_1(\ast\nbigd)
 \lrarr E_2(\ast \nbigd)$.
\end{lem}
\pf
We have only to consider the case
$0\not\in\nbigk$
according to the Hartogs theorem.
Then, the claim follows from 
a result in Subsection 7.7.6 of \cite{mochi7}.
As another argument,
let $\vecw_S^{(i)}$ be frames of $E_{i|\Sbar}$
as in Proposition \ref{prop;08.9.11.3}.
We can directly show that
$F_{|S}$ is of polynomial order with respect to
the frames $\vecw_S^{(i)}$.
\hfill\qed

\paragraph{Complement on
a connection along the $\lambda$-direction}

let $X:=\Delta^n$, $D_i:=\{z_i=0\}$
and $D:=\bigcup_{i=1}^{\ell}D_i$.
Let $\nbigk\subset\cnum_{\lambda}^{\ast}$
be a compact region.
Let $(E,\DD)$ be an unramifiedly good lattice of
a family of meromorphic $\lambda$-flat bundles
on $(\nbigx,\nbigd)$
with a good set $\Irr(\DD)$.
Assume that $E$ is equipped with
a meromorphic connection along 
the $\lambda$-direction
$ \nabla_{\lambda}:
 E\lrarr E\otimes\Omega_{\nbigk}^1(\ast \nbigd)$,
such that
$\DD^f+\nabla_{\lambda}$ is flat.

\begin{lem}
\label{lem;08.8.2.5}
$\nabla_{\lambda}$ naturally 
induces a meromorphic
connection of $E^{(T)}$
along the $\lambda$-direction.
\end{lem}
\pf
It is easy to observe that
we have only to consider the case
in which $D$ is smooth
and $|\arg(T)|$ is sufficiently small.
Take $N$ such that
$\lambda^N\nabla_{\lambda}(\del_{\lambda})
 E\subset 
 E\otimes\nbigo_{\nbigx}(\ast\nbigd)$.
For $S=S^{(i)}$,
let $\vecw_{S}=(\vecw_{\gminia,S})$
be a frame of $E_{|\Sbar}$ as above. 
Let $A_S=(A_{S,\gminia,\gminia'})$
be the matrix-valued holomorphic function on $S$
determined by
$\lambda^N\nabla_{\lambda}(\del_{\lambda})
 \vecw_{S}
=\vecw_{S}\cdot A_S$.
Let $B_{\gminia}$ be the matrix-valued
holomorphic function on $\nbigx$ determined by
$\DD_{\gminia}(z_1\del_1)\vecw_{\gminia}
=\vecw_{\gminia}\cdot
 \bigl(z_1\del_1\gminia+B_{\gminia}\bigr)$.
Because $[\DD^f,\nabla_{\lambda}]=0$,
we have the following relation in the case
$\gminia\neq\gminib$:
\[
 \lambda\cdot z_1\del_1A_{S,\gminia,\gminib}
+\bigl(z_1\del_1(\gminia-\gminib)\bigr)
 \cdot A_{S,\gminia,\gminib}
+\bigl(
 A_{S,\gminia,\gminib}B_{\gminib}
-B_{\gminia}A_{S,\gminia,\gminib}
 \bigr)=0
\]
Hence, 
we have $A_{S,\gminia,\gminib}=0$
unless $\gminia\leq_S\gminib$,
and we obtain the estimate 
\[
 A_{S,\gminia,\gminib}\cdot
\exp\bigl(
 \lambda^{-1}(\gminia-\gminib)
 \bigr)
=O\Bigl(\exp\bigl(C|\lambda^{-1}|\cdot
 \log|z_1^{-1}|\bigr)
\Bigr)
\]
for some $C>0$
in the case $\gminia<_S\gminib$.

Let $A^{(T)}_S$ be the matrix-valued
holomorphic function on $S$
determined by 
$\lambda^N\nabla(\del_{\lambda})
 \vecw^{(T)}_S=
 \vecw^{(T)}_S\cdot
 A^{(T)}_S$.
We have $A^{(T)}_{S,\gminia,\gminib}=0$
unless $\gminia\leq_S\gminib$.
In the case $\gminia<_S\gminib$,
we have 
\[
 A^{(T)}_{S,\gminia,\gminib}
 \cdot
 \exp\bigl(\lambda^{-1}\cdot T\cdot(\gminia-\gminib)\bigr)
=A_{S,\gminia,\gminib}
\cdot
 \exp\bigl(\lambda^{-1}\cdot (\gminia-\gminib)\bigr)
=O\Bigl(\exp\bigl(C|\lambda^{-1}|\cdot
 \log|z_1^{-1}|\bigr)
 \Bigr).
\]
Therefore,
we obtain
$A_{S,\gminia,\gminib}^{(T)}
=O\bigl(
 \exp(-\epsilon|z_1^{-1}|)
 \bigr)$
for some $\epsilon>0$.
By a direct calculation,
we obtain
$ A^{(T)}_{S,\gminia,\gminia}
=A_{S,\gminia,\gminia}
+\lambda^N\cdot
 \del_{\lambda}\bigl(
\lambda^{-1}\cdot (1-T)\cdot\gminia
 \bigr)$,
which is of polynomial order.
Hence, the claim of the lemma follows from
Lemma \ref{lem;07.10.2.21}.
\hfill\qed

\subsection{Family of 
 good filtered $\lambda$-flat bundles}
\label{subsection;08.8.3.4}

\paragraph{Pull back of filtered bundle
 via a ramified covering}

The notion of filtered bundle
is introduced in \cite{s2} ($1$ dimension),
and  studied in \cite{mochi2} (arbitrary dimension).
Let $X$ be a complex manifold,
and let $D$ be a simple normal crossing hypersurface
with the irreducible decomposition
$D=\bigcup_{i\in I}D_i$.
A filtered bundle on $(X,D)$
is defined to be a sequence of locally free sheaves
$\vecE_{\ast}=\bigl(\prolongg{\veca}{E}
 \,\big|\,\veca\in\real^I \bigr)$
such that
(i) $\prolongg{\veca}{E}\subset
 \prolongg{\vecb}{E}$
 for $\veca\leq\vecb$ 
 and $\prolongg{\veca}{E}$
 is the intersection of  
 $\prolongg{\vecb}{E}$ for $\vecb>\veca$,
(ii) $\prolongg{\veca}{E}_{|X-D}
=\prolongg{\vecb}{E}_{|X-D}$,
(iii)
 $\prolongg{\veca}{E}\otimes
 \nbigo(\sum n_i\cdot D_i)
=\prolongg{\veca-\vecn}{E}$,
where $\vecn=(n_i)\in\seisuu^I$,
(iv) it satisfies some compatibility condition
at the intersection of the divisors.
The compatibility condition is
given in Definition 4.37 of \cite{mochi2}.
Although it is not difficult,
it is slightly complicated to state.
Later, Iyer and Simpson \cite{i-s} introduced
the notion of locally abelian condition,
which is equivalent to our compatibility condition.
Hertling and Sevenheck
(Chapter 4 of \cite{Hertling-Sevenheck3})
showed that it is equivalent to
another simple condition.
We refer to the above papers
for more details.

Let us recall the pull back of a filtered bundle
via a ramified covering.
See \cite{i-s} for more systematic treatment.
See also Subsection 2.9.1 of \cite{mochi7}.
Let $X:=\Delta^n_z$,
$D:=\bigcup_{i=1}^{\ell}\{z_i=0\}$,
$\Xtilde:=\Delta^n_{w}$ 
and $\Dtilde:=\bigcup_{j=1}^{\ell}\{w_j=0\}$.
Let $\varphi_e:\Xtilde\lrarr X$ be a ramified covering
$\varphi_e(w_1,\ldots,w_n)
=(w_1^{e},\ldots,w_{\ell}^{e},
 w_{\ell+1},\ldots,w_n)$.
For $\vecb\in\real^{\ell}$,
we put 
$ \nbigs(\vecb):=\bigl\{
 (\veca,\vecn)\in\real^{\ell}\times\seisuu_{\geq 0}^{\ell}
 \,\big|\,
 e\cdot\veca+\vecn\leq\vecb
 \bigr\}$.
For a given filtered bundle $\vecE_{\ast}$
on $(X,D)$,
we set
\[
 \prolongg{\vecb}{\Etilde}
=\sum_{(\veca,\vecn)\in\nbigs(\vecb)}
 \vecw^{-\vecn}\cdot
 \varphi_e^{\ast}\bigl(\prolongg{\veca}{E}\bigr).
\]
Then, it is easy to show that
$\vecEtilde_{\ast}$ is also a filtered bundle.
Let $\Gal(\Xtilde/X)$ denote the Galois group
of the ramified covering.
We can reconstruct $\vecE_{\ast}$
from $\vecEtilde_{\ast}$ with the natural
$\Gal(\Xtilde/X)$-action,
and hence $\vecE_{\ast}$ is called
the descent of $\vecEtilde_{\ast}$.
Since the construction is independent of 
the choice of coordinates,
it can be globalized.

\paragraph{Family of good filtered 
 $\lambda$-flat bundles}

We use the notation in 
Subsection \ref{subsection;07.12.14.30}.
A family of filtered $\lambda$-flat bundles
on $(\nbigx,\nbigd)$
is defined to be a filtered bundle
$\vecE_{\ast}$ on $(\nbigx,\nbigd)$
with a family of meromorphic flat $\lambda$-connections
$\DD$ of $\vecE=\bigcup \prolongg{\veca}{E}$.

\begin{df}
\mbox{{}}\label{df;08.10.1.3}
Let $(\vecE_{\ast},\DD)$ be a family of
filtered $\lambda$-flat bundles on $(\nbigx,\nbigd)$.
\begin{itemize}
\item
We say that $(\vecE_{\ast},\DD)$ is unramifiedly good,
if $\prolongg{\vecc}{E}$ are unramifiedly good lattices
for any $\vecc\in\real^{\ell}$.
\item 
Let $P\in \nbigd$.
We say that
$(\vecE_{\ast},\DD)$ is good at $P$,
if there exist a ramified covering
$\varphi_e:(\nbigutilde,\nbigdtilde_{\nbigu})
 \lrarr (\nbigu,\nbigd_{\nbigu})$
such that
$\bigl(\vecEtilde_{\ast},\,
 \varphi_e^{\ast}\DD\bigr)$
on $(\nbigutilde,\nbigdtilde_{\nbigu})$
is unramifiedly good.
Here, $\nbigu$ is a coordinate neighbourhood of $P$,
$\varphi_e$ is a ramified covering,
and $\vecEtilde_{\ast}$ is induced by
$\varphi$ and $\vecE_{\ast}$ as above.
\item
We say that $(\vecE_{\ast},\DD)$ is good,
if it is good at any point $P\in \nbigd$.
\hfill\qed
\end{itemize}
\end{df}

\paragraph{Induced filtrations}

Let $(\vecE_{\ast},\DD)$ be 
good family of filtered $\lambda$-flat bundles.
Let $\lefttop{i}F$ denote the induced filtration
of $\prolongg{\vecc}{E}_{|\nbigd_i}$.
We set 
$\lefttop{i}\Gr^F_a(\prolongg{\vecc}{E}):=
 \lefttop{i}F_a\big/\lefttop{i}F_{<a}$.
It can be shown that 
(i) we have the well defined residue endomorphism
$\Gr^F_a\Res_i(\DD)$
of $\lefttop{i}\Gr^F_a\bigl(\prolongg{\vecc}{E}\bigr)$
on $\nbigd_i$ for each $i\in \ellsitabar$,
(ii) it preserves the induced filtrations 
$\lefttop{j}F$ of 
$\lefttop{i}\Gr^F_a
 \bigl(\prolongg{\vecc}{E}\bigr)_{|\nbigd_i\cap \nbigd_j}$.
(See Subsection 6.1.3 of \cite{mochi7}.
 The residues are well defined 
 as endomorphisms of
 $\prolongg{\vecc}{E}_{|\nbigd_i}$
 in the non-ramified case,
 and as endomorphisms of
 $\lefttop{i}\Gr^F_a\bigl(\prolongg{\vecc}{E}\bigr)$
 even in the ramified case.)
In the following,
$\Gr^F_a\Res_{i}(\DD)$
are often denoted by
$\Res_i(\DD)$
for simplicity of the description.

Let $I$ be a subset of $\ellsitabar$.
We set $\nbigd_I:=\bigcap_{i\in I}\nbigd_i$.
For $\veca\in\real^I$,
we put
\[
 \lefttop{I}F_{\veca}\bigl(
 \prolongg{\vecc}{E}_{|\nbigd_I}
 \bigr):=
 \bigcap_{i\in I}
 \lefttop{i}F_{a_i}
 \bigl(\prolongg{\vecc}{E}_{|\nbigd_I}\bigr),
\quad\quad
 \lefttop{I}\Gr^F_{\veca}\bigl(
 \prolongg{\vecc}{E}
 \bigr):=
\frac{\lefttop{I}F_{\veca}\bigl(
 \prolongg{\vecc}{E}_{|\nbigd_I}\bigr)}
{\sum_{\vecb\lneq\veca}
 \lefttop{I}F_{\vecb}\bigl(
 \prolongg{\vecc}{E}_{|\nbigd_I}\bigr)}.
\]
We often consider the following sets:
\[
 \Par\bigl(\prolongg{\vecc}{E},I\bigr):=
 \bigl\{
 \veca\in\real^I\,\big|\,
\lefttop{I}\Gr^{F}_{\veca}(\prolongg{\vecc}{E})\neq 0
 \bigr\},
\quad
 \Par\bigl(\vecE_{\ast},I\bigr):=
 \bigcup_{\vecc\in\real^{\ell}}
 \Par\bigl(\prolongg{\vecc}{E},I\bigr)
\]
We have the induced endomorphisms
$\Res_i(\DD)$ ($i\in I$) of
$\lefttop{I}\Gr^F_{\veca}(\prolongg{\vecc}{E})$,
which are mutually commutative.

\paragraph{KMS structure for fixed $\lambda$}

Let us consider the case 
in which $\nbigk$ is a point $\{\lambda\}$.
In this case,
we prefer the symbol $\DDlambda$
to $\DD$.
If $\lambda\neq 0$,
the eigenvalues of $\Res_i(\DDlambda)$ are constant.
Hence, we have the generalized eigen decomposition
$\lefttop{I}\Gr^F_{\veca}\bigl(
 \prolongg{\vecc}{E}\bigr)
=\bigoplus_{\vecalpha}
\lefttop{I}\Gr^{F,\EE}_{(\veca,\vecalpha)}
 \bigl(\prolongg{\vecc}{E}\bigr)$,
where the eigenvalues of 
$\Gr^F\Res_i(\DDlambda)$ 
on $\lefttop{I}\Gr^{F,\EE}_{(\veca,\vecalpha)}
 \bigl(\prolongg{\vecc}{E}\bigr)$
are the $i$-th components of $\vecalpha$.
We put
\[
 \KMS(\prolongg{\vecc}{E},\DDlambda,I):=
 \bigl\{
 (\veca,\vecalpha)\,\big|\,
 \lefttop{I}\Gr^{F,\EE}_{(\veca,\vecalpha)}
 (\prolongg{\vecc}{E})\neq 0
 \bigr\},
\quad
 \KMS(\vecE_{\ast},\DDlambda,I):=
 \bigcup_{\vecc\in\real^S}
 \KMS(\prolongg{\vecc}{E},\DDlambda,I)
\]
\[
\Sp(\prolongg{\vecc}{E},\DDlambda,I):=
 \bigl\{
 \vecalpha\in\cnum^I\,\big|\,
 \exists\veca\in\real^I,\,\,
 (\veca,\vecalpha)\in 
 \KMS\bigl(\prolongg{\vecc}{E},\DDlambda,I\bigr)
 \bigr\},
\quad
 \Sp(\vecE_{\ast},\DDlambda,I):=
 \bigcup_{\vecc\in\real^S}
 \Sp(\prolongg{\vecc}{E},\DDlambda,I)
\]
Each element of $\KMS(\vecE_{\ast},\DDlambda,I)$
is called a KMS-spectrum of $(\vecE_{\ast},\DDlambda)$
at $D_I$.

Even in the case $\lambda=0$,
a similar definition makes sense
if the eigenvalues of $\Res_i(\DDlambda)$
are constant.
It is satisfied 
when we consider wild harmonic bundles.

\paragraph{KMS structure around $\lambda_0$}

Assume that
$\nbigk$ is a neighbourhood of $\lambda_0\in\cnum$,
and we regard that
$(\vecE_{\ast},\DD)$ is given
around $\{\lambda_0\}\times X$.
In this case, we prefer the symbols
$\lefttop{i}\Fzero$ to $\lefttop{i}F$.
Let $\paramap(\lambda):\real\times\cnum\lrarr\real$
and $\eigenmap(\lambda):\real\times\cnum\lrarr\cnum$
be given as follows:
\[
 \paramap\bigl(\lambda,(a,\alpha)\bigr)
=a+2\Re(\lambda\cdot\alphabar),
\quad
 \eigenmap\bigl(\lambda,(a,\alpha)\bigr)
=\alpha-a\cdot \lambda-\alphabar\cdot\lambda^2
\]
The induced map
$\real\times\cnum\lrarr\real\times\cnum$
is denoted by $\kmsmap(\lambda)$.

\begin{df}
\label{df;07.11.23.5}
We say that $(\vecE_{\ast},\DD)$ has the KMS-structure
at $\lambda_0$ indexed by
$T(i)\subset\real\times\cnum$ $(i\in S)$,
if the following holds:
\begin{itemize}
\item
 $\Par(\vecE_{\ast},i)$ is the image of
 $T(i)$ via the map $\paramap(\lambda_0)$.
\item 
 For each $a\in\Par(\vecE_{\ast},i)$,
 we put $\nbigk(a,i):=\bigl\{
 u\in T(i)\,\big|\,\paramap(\lambda_0,u)=a
 \bigr\}$.
 Then, the restrictions of
 $\Res_i(\DD)$
 to $\lefttop{i}\Gr^{\Fzero}_{a}
 \bigl(\prolongg{\vecc}{E}\bigr)
 _{|\nbigd^{\lambda}_i}$
 have the unique eigenvalue
 $\eigenmap(\lambda,u)$
 for any $u\in \nbigk(a,i)$.
\hfill\qed
\end{itemize}
\end{df}

Assume $(\vecE_{\ast},\DD)$
has the KMS-structure at $\lambda_0$.
We have the decomposition
\begin{equation}
 \label{eq;07.10.16.20}
 \lefttop{i}\Gr^{\Fzero}_{a}
 \bigl(\prolongg{\vecc}{E}\bigr)
=\bigoplus_{u\in\nbigk(a,i)}
 \lefttop{i}\nbigg^{(\lambda_0)}_{u}
\bigl(\prolongg{\vecc}{E}\bigr),
\end{equation}
such that
(i) it is preserved by $\Res_i(\DD)$,
(ii) the restriction of 
$\Res_i(\DD)-\eigenmap(\lambda,u)$ to
$\lefttop{i}\nbigg^{(\lambda_0)}_{u}
\bigl(\prolongg{\vecc}{E}\bigr)$ is nilpotent.
More generally, we have the decomposition
on $\nbigd_I$
\begin{equation}
 \label{eq;08.8.3.7}
  \lefttop{I}\Gr^{\Fzero}_{\veca}
 \bigl(\prolongg{\vecc}{E}\bigr)
=\bigoplus_{\vecu\in \prod\nbigk(a_i,i)}
 \lefttop{I}\nbigg^{(\lambda_0)}_{\vecu}
\bigl(\prolongg{\vecc}{E}\bigr),
\end{equation}
such that
(i) it is preserved by
$\Res_i(\DD)$ $(i\in I)$,
(ii) the restrictions of
 $\Res_i(\DD)-\eigenmap(\lambda,u_i)$ $(i\in I)$
 are nilpotent,
 where $u_i$ denotes the $i$-th component of $\vecu$.
Note $\lefttop{I}\nbigg_{\vecu}^{(\lambda_0)}
 \bigl(\prolongg{\vecc}{E}\bigr)$
can be $0$.

The following lemma is standard in our works.
(See Subsection 6.2.5 of \cite{mochi7}.)
\begin{lem}
\label{lem;07.11.23.5}
Let $(\vecE_{1\,\ast},\DD_1)$
and $(\vecE_{2\,\ast},\DD_2)$
be good filtered $\lambda$-flat bundles
on $(\nbigx,\nbigd)$
which have the $KMS$-structures at $\lambda_0$.
An isomorphism
$\varphi:(\vecE_{1},\DD_1)\simeq (\vecE_{2},\DD_2)$
of families of meromorphic $\lambda$-flat bundles
induces the isomorphism
$\varphi:(\vecE_{1\,\ast},\DD_1)\simeq 
 (\vecE_{2\,\ast},\DD_2)$
of families of filtered $\lambda$-flat bundles.
\hfill\qed
\end{lem}
We say that 
$(\vecE,\DD)$ has the KMS-structure at $\lambda_0$,
if there exists a good filtered $\lambda$-flat bundle
 $(\vecE_{\ast},\DD)$ which has the KMS-structure
at $\lambda_0$,
such that $\vecE=\bigcup \prolongg{\veca}{E}$.
It makes sense by the above lemma.

Pick $\vecc\in \real^S$
such that $c_i\not\in\Par\bigl(\vecE_{\ast},i\bigr)$
for each $i\in S$.
Assume that $\nbigk$ is a sufficiently small
neighbourhood of $\lambda_0$.
Take $\lambda_1\in \nbigk$,
and let $U(\lambda_1)\subset \nbigk$
be a neighbourhood of $\lambda_1$.
We set 
$\nbigx^{(\lambda_1)}:=U(\lambda_1)\times X$.
We use the symbols
$\nbigd^{(\lambda_1)}_i$ and
$\nbigd^{(\lambda_1)}$ in similar meanings.
Let $\pi_{i,a}$ denote the projection
$ \lefttop{i}\Fzero_a\bigl(
 \prolongg{\vecc}{E}_{|\nbigd_i}\bigr)
\lrarr
 \lefttop{i}\Gr^{\Fzero}_a
 \bigl(\prolongg{\vecc}{E}\bigr)$
for any $a\in\Par(\prolongg{\vecc}{E},i)$.
Let $b\in \openclosed{c_i-1}{c_i}$.
If $\paramap(\lambda_1,v)=b$ for some
$v\in \nbigk(a,i)$,
we put on $\nbigd_i^{(\lambda_1)}$
\[
 \lefttop{i}F^{(\lambda_1)}_{b}:=
 \bigoplus_{
 \substack{u\in\nbigk(a,i)\\
 \paramap(\lambda_1,u)\leq b}}
\pi_{i,a}^{-1}\bigl(
 \lefttop{i}\nbigg^{(\lambda_0)}_u
 (\prolongg{\vecc}{E})
\bigr).
\]
Otherwise,
let $b_0:=\max\bigl\{
 \paramap(\lambda_1,v)<b\,\big|\,
 v\in\nbigk(a,i)
 \bigr\}$,
and we set 
$\lefttop{i}F^{(\lambda_1)}_b:= 
 \lefttop{i}F^{(\lambda_1)}_{b_0}$.
Thus, we obtain the filtration 
$\lefttop{i}F^{(\lambda_1)}$ of
$\prolongg{\vecc}{E}_{|\nbigd^{(\lambda_1)}_i}$.
It induces the family of the filtered $\lambda$-flat bundles
$(\vecE^{(\lambda_1)}_{\ast},\DD)$
on $(\nbigx^{(\lambda_1)},\nbigd^{(\lambda_1)})$.
By construction,
$\Res_{i}(\DD)-\eigenmap(\lambda,u)$ 
are nilpotent on
$\lefttop{i}\Gr^{F^{(\lambda_1)}}_{
 \paramap(\lambda_1,u)}
 \bigl(\prolongg{\vecc}{E}\bigr)$.
Namely,
$(\vecE^{(\lambda_1)}_{\ast},\DD)$
has the KMS-structure at $\lambda_1$
indexed by $T(i)$.
Hence, if $(\vecE,\DD)$ has the KMS-structure
at $\lambda_0$,
it has the KMS-structure at any $\lambda$
sufficiently close to $\lambda_0$,
and the index set is independent of $\lambda$.
For each $\lambda\in \nbigk$,
we put $\vecE_{\ast}^{\lambda}:=
 (\vecE_{\ast}^{(\lambda)})_{|\nbigxlambda}$,
which is the good filtered $\lambda$-flat bundle.
The set
$\KMS(\vecE_{\ast}^{\lambda},i)$
is the image of $T(i)$ via the map $\kmsmap(\lambda)$.
Note $\KMS(\vecE_{\ast}^0,i)=T(i)$
if $0\in \nbigk$.
We often identify them.

\paragraph{Deformation}

Let $T(\lambda)$ be 
a holomorphic function
with $\Re\bigl(T(\lambda)\bigr)>0$.
We obtain the deformation
$(\vecE_{\ast}^{(T)},\DD)$.
If $(\vecE_{\ast},\DD)$ is unramified,
the set of irregular values is given by
\[
 \Irr\bigl(\DD,E^{(T)}\bigr):=
 \bigl\{
 T\cdot\gminia\,\big|\,
 \gminia\in\Irr(\DD)
 \bigr\}.
\]
Since the regular part of the completion
is unchanged,
the set of KMS-spectra is unchanged.

\section{Wild harmonic bundle}
\label{section;08.8.20.41}

\subsection{Definition of wild harmonic bundle}
\label{subsection;08.10.27.1}

\paragraph{Local condition for Higgs fields}

Let $(E,\delbar_E,\theta)$ be a Higgs bundle on $X-D$,
where $X$ is a complex manifold,
and $D$ is a normal crossing divisor of $X$.
We would like to explain some conditions
for the Higgs field $\theta$.
First, let us consider the case 
$X=\Delta^n=
 \bigl\{\vecz=(z_1,\ldots,z_n)\,\big|\,|z_i|<1\bigr\}$,
$D_i=\{z_i=0\}$
and $D=\bigcup_{i=1}^{\ell}D_i$.
We have the expression:
\[
  \theta=\sum_{j=1}^{\ell} F_j\cdot \frac{dz_j}{z_j}
+\sum_{j=\ell+1}^{n}G_j\cdot dz_j
\]
We have the characteristic polynomials
$\det\bigl(T-F_j(\vecz)\bigr)
=\sum A_{j,k}(\vecz)\cdot T^k$ and
$\det\bigl(T-G_j(\vecz)\bigr)=\sum B_{j,k}(\vecz)\cdot T^k$.
The coefficients $A_{j,k}$ and $B_{j,k}$
are holomorphic on $X-D$.

\begin{itemize}
\item
We say that $\theta$ is tame,
if the following conditions are satisfied:
\begin{description}
\item[(T1):]
$A_{j,k}$ and $B_{j,k}$ are holomorphic on $X$
for any $k$.
\item[(T2):]
The restriction of $A_{j,k}$ to $D_j$ are constant
for any $j=1,\ldots,\ell$
and any $k$.
In other words,
roots of
$\sum A_{j,k}(\vecz)\cdot T^k$
are independent of $\vecz\in D_j$.
\end{description}

\item
We say that $\theta$ is unramifiedly good,
if there exists a good set of irregular values
$\Irr(\theta)\subset M(X,D)\big/H(X)$
and a decomposition
$(E,\theta)=\bigoplus_{\gminia\in\Irr(\theta)}
 (E_{\gminia},\theta_{\gminia})$,
such that
$\theta_{\gminia}-d\gminia\cdot\pi_{\gminia}$
are tame,
where $\pi_{\gminia}$ denotes the projection
onto $E_{\gminia}$
with respect to the decomposition.
\item
We say that $\theta$ is good,
if $\varphi_e^{\ast}(\theta)$ is unramifiedly good
for some $e\in\seisuu_{>0}$,
where $\varphi_e$ is the covering
given by
$ \varphi_e(z_1,\ldots,z_n)=
 (z_1^{e},\ldots,z_{\ell}^{e},z_{\ell+1},\ldots,z_n)$.
\hfill\qed
\end{itemize}

\paragraph{Global condition for Higgs fields}

Let us consider the case in which
$X$ is a general complex manifold.
Let $D$ be a normal crossing hypersurface of $X$,
and let $(E,\theta)$ be a Higgs bundle on $X-D$.
\begin{itemize}
\item
We say that $\theta$ is 
(unramifiedly) good at $P\in D$,
if it is (unramifiedly) good 
on some holomorphic
coordinate neighbourhood of $P$.
\item
We say that $\theta$ is 
(unramifiedly) good,
if it is (unramifiedly) good 
at any point $P\in D$.
\hfill\qed
\end{itemize}

Let $Z$ be a closed analytic subset of $X$,
and let $(E,\theta)$ be a Higgs bundle 
on $X-Z$.
The Higgs field $\theta$ is called wild,
if there exists a regular birational map
$\varphi:X'\lrarr X$
such that (i) $\varphi^{-1}(D)$ is normal crossing,
(ii) $\varphi^{-1}\theta$ is good.

\begin{rem}
\label{rem;08.11.6.11}
Even if $Z$ is a normal crossing divisor,
wild $\theta$ is not necessarily good.
\hfill\qed
\end{rem}

\paragraph{Conditions for harmonic bundles}

Let $X$ be a complex manifold.
Let $D$ be a normal crossing hypersurface of $X$,
and let $\harmonicbundle$ be a harmonic bundle
on $X-D$.
\begin{itemize}
\item
 It is called tame,
 if $\theta$ is tame.
\item
 It is called (unramifiedly) good 
 wild harmonic bundle,
 if $\theta$ is 
 (unramifiedly) good.
\hfill\qed
\end{itemize}
Let $Z$ be a closed analytic subset of $X$.
A harmonic bundle $\harmonicbundle$ on $X-Z$
is called wild, if $\theta$ is wild.

\paragraph{Remark}
We give some remarks on the condition {\bf (T2)}
for tameness.
\begin{enumerate}
\item
If $\theta$ comes from a harmonic bundle
$(E,\delbar_E,\theta,h)$,
{\bf (T2)} is implied by {\bf (T1)}.
(See Lemma {\rm 8.2} of {\rm\cite{mochi2}}.)
\item
Let $(E,\delbar_E,\theta,h)$ be 
a harmonic bundle
with a good set of irregular values $\Irr(\theta)$
and a decomposition
$(E,\delbar_E,\theta)
=\bigoplus_{\gminia\in\Irr(\theta)}
 (E_{\gminia},\delbar_{E_{\gminia}},\theta_{\gminia})$
such that
$\thetatilde_{\gminia}:=\theta_{\gminia}
-d\gminia\cdot\pi_{\gminia}$ satisfy
the condition {\bf (T1)}.
The author does not know
whether {\bf (T2)} for $\thetatilde_{\gminia}$ is 
automatically satisfied or not.
But, if moreover $(E,\delbar_E,\theta,h)$ underlies
a variation of polarized pure integrable structures,
{\bf (T2)} is satisfied.
Actually, the roots of the polynomials are $0$.
(See Lemma  \ref{lem;08.7.26.110} below.)
\end{enumerate}

\subsection{Simpson's main estimate}
\label{subsubsection;08.7.30.3}

The first fundamental result 
is an estimate of Higgs field,
so called Simpson's main estimate.
For later use,
we recall it in the case that $D$ is smooth.
Let $X:=\Delta^n$ and $D:=\{z_1=0\}$.
Let $(E,\delbar_E,\theta,h)$
be an unramifiedly good wild harmonic bundle
on $X-D$.
(See Subsections 11.2 and 11.3 of \cite{mochi7}
for more details.)
We will be interested in the behaviour around $O$.
Hence, by shrinking $X$,
we assume that there exists 
 a holomorphic decomposition
$ (E,\theta)=\bigoplus_{
 (\gminia,\alpha)\in\Irr(\theta)\times\cnum}
 (E_{\gminia,\alpha},\theta_{\gminia,\alpha})$
satisfying the following conditions:
\begin{itemize}
\item
For each $(\gminia,\alpha)$,
let $\pi_{\gminia,\alpha}$ denote the projection
onto $E_{\gminia,\alpha}$
with respect to the decomposition.
We have the expression
\[
 \theta_{\gminia,\alpha}
 -\bigl(
 \alpha\cdot dz_1/z_1+d\gminia\bigr)
 \cdot \pi_{\gminia,\alpha}
=F_1\cdot \frac{dz_1}{z_1}
+\sum_{j=2}^nG_j\cdot dz_j.
\]
Then, 
the coefficients of 
$\det(T-F_1)$ and $\det(T-G_j)$
are holomorphic on $X$,
and $\det(T-F_1)_{|D}=T^{\rank E_{\gminia,\alpha}}$.
\end{itemize}
We also set 
$E_{\gminia}:=
 \bigoplus_{\alpha\in\cnum} E_{\gminia,\alpha}$,
and let $\pi_{\gminia}$ denote the projection
onto $E_{\gminia}$ with respect to the decomposition
$E=\bigoplus_{\gminia\in\Irr(\theta)} E_{\gminia}$.

\paragraph{Truncation}

For any $\gminia\in \Irr(\theta)$,
we have the expression
$\gminia=\sum_{j\leq -1} \gminia_j\cdot z_1^j$.
We put 
$\eta_{p}(\gminia):=
 \sum_{j\leq p}\gminia_j\cdot z_1^j$
and $\Irr(\theta,p):=\bigl\{
 \eta_p(\gminia)\,\big|\,
 \gminia\in\Irr(\theta) \bigr\}$.
For each $\gminib\in\Irr(\theta,p)$,
let $E_{\gminib}^{(p)}$ denote the direct sum of
$E_{\gminia}$ 
$(\gminia\in\Irr(\theta),\eta_p(\gminia)=\gminib)$,
and let $\pi^{(p)}_{\gminib}$
denote the projection onto
$E^{(p)}_{\gminib}$
with respect to the decomposition
$E=\bigoplus_{\gminib\in \Irr(\theta,p)} E^{(p)}_{\gminib}$.
We have 
$\Irr(\theta,-1)=\Irr(\theta)$ and
$E_{\gminia}=E^{(-1)}_{\gminia}$.
We have the induced maps
$\eta_{q,p}:\Irr(\theta,p)\lrarr\Irr(\theta,q)$
for $q\leq p$.

\paragraph{Asymptotic orthogonality}

We take total orders
$\leq'$ on $\Irr(\theta,p)$ $(p\leq -1)$
which are preserved by $\eta_{q,p}$.
For each $\gminib\in\Irr(\theta,p)$,
we set 
$F^{(p)}_{\gminib}(E):=
 \bigoplus_{\gminia\leq'\gminib}E^{(p)}_{\gminia}$.
Let $E^{(p)\prime}_{\gminib}$
be the orthogonal complement
of $F^{(p)}_{<\gminib}(E)$ 
in $F^{(p)}_{\gminib}(E)$.
We obtain an orthogonal decomposition
$E=\bigoplus_{\gminia\in\Irr(\theta,p)}
 E^{(p)\prime}_{\gminia}$.
Let $\pi^{(p)\prime}_{\gminia}$ 
denote the orthogonal projection
onto $E^{(p)\prime}_{\gminia}$.

We take a total order $\leq'$ on $\cnum$.
Then, we obtain the lexicographic order
on $\Irr(\theta)\times\cnum$.
We obtain the orthogonal decomposition
$E=\bigoplus E'_{\gminia,\alpha}$
by the procedure as above,
and let $\pi'_{\gminia,\alpha}$
denote the orthogonal projection onto 
$E'_{\gminia,\alpha}$.

\begin{prop}
\label{prop;08.7.30.2}
We have the following estimates
with respect to $h$.
\begin{itemize}
\item
$\pi^{(p)}_{\gminia}-\pi^{(p)\prime}_{\gminia}
=O\Bigl(\exp\bigl(-\epsilon|z_1^{p}|\bigr)\Bigr)$
for some $\epsilon>0$.
In particular, the decomposition
$E=\bigoplus E^{(p)}_{\gminib}$ is
$O\bigl(\exp(-\epsilon|z_1^p|)\bigr)$-asymptotically
orthogonal
in the sense that
there exists $A>0$ such that
\[
 \bigl|h(u,v)\bigr|
\leq
 A\cdot |u|_h\cdot |v|_h\cdot
 \exp\bigl(-\epsilon|z_1(Q)|^p\bigr)
\]
for any $Q\in X-D$,
$u\in E_{\gminia|Q}$ and
$v\in E_{\gminib|Q}$ $(\gminia\neq\gminib)$.
\item
$\pi_{\gminia,\alpha}-\pi'_{\gminia,\alpha}
=O\bigl(|z_1|^{\epsilon}\bigr)$
for some $\epsilon>0$.
In particular,
the decomposition
$E=\bigoplus E_{\gminia,\alpha}$ is
$O\bigl(|z_1|^{\epsilon}\bigr)$-asymptotically
orthogonal.
\hfill\qed
\end{itemize}
\end{prop}

\paragraph{Estimate of Higgs field}

We set 
$\thetatilde:=\theta
-\bigoplus_{\gminia,\alpha} 
 \bigl(d\gminia+\alpha\cdot dz_1/z_1\bigr)
 \pi_{\gminia,\alpha}$.
Let $g_{\poin}$ denote the Poincar\'e metric
of $X-D$.
The estimates in Subsection 11.2 of \cite{mochi7}
implies the following.
\begin{prop}
\label{prop;08.8.3.10}
$\thetatilde$ is bounded
with respect to $h$ and $g_{\poin}$.
\hfill\qed
\end{prop}

\paragraph{Estimate of curvatures}

As mentioned in Subsection
\ref{subsubsection;08.7.29.30},
we obtain a holomorphic vector bundle
$\nbigelambda=
 (E,\delbar_E+\lambda\theta^{\dagger})$
on $X-D$.
The curvature of the unitary connection
associated to $(\nbigelambda,h)$
equals to
$-(1+|\lambda|^2)\cdot
 \bigl[\theta,\theta^{\dagger}\bigr]$.

\begin{prop}
\label{prop;08.7.30.5}
$[\theta,\theta^{\dagger}]$ is 
bounded with respect to $h$ and $g_{\poin}$.
In particular,
$(\nbigelambda,h)$ is acceptable,
i.e.,
the curvature of
$(\nbigelambda,h)$ is bounded
with respect to $h$ and $g_{\poin}$.
\hfill\qed
\end{prop}

\subsection{Prolongation of unramifiedly good 
wild harmonic bundles}
\label{subsection;08.8.3.6}

\subsubsection{Prolongment $\nbigp\nbigelambda$}

Let $\harmonicbundle$ be 
a good wild harmonic bundle on $X-D$,
where $X$ is a complex manifold
and $D$ is a normal crossing divisor.
As mentioned in Subsection
\ref{subsubsection;08.7.29.30},
we obtain a holomorphic vector bundle
$\nbigelambda=
 (E,\delbar_E+\lambda\theta^{\dagger})$
on $X-D$ for each complex number $\lambda$.
It is important to prolong it
to a good filtered $\lambda$-flat bundle
on $(X,D)$.
For simplicity,
we explain it assuming the following.
(The general case can be easily
reduced to this case.)
\begin{itemize}
\item
$X=\Delta^n$ and
$D=\bigcup_{i=1}^{\ell}\{z_i=0\}$.
\item
$\harmonicbundle$ is unramifiedly good wild,
and the underlying Higgs bundle
has the following decomposition
\begin{equation}
 \label{eq;08.8.2.10}
 (E,\theta)
=\bigoplus_{\substack{\gminia\in\Irr(\theta)\\
 \vecalpha\in\cnum^{\ell}}}
 (E_{\gminia,\vecalpha},\theta_{\gminia,\vecalpha}),
\end{equation}
such that
(i) $\thetatilde_{\gminia}=\theta_{\gminia}
-\bigl(
d\gminia
+\sum_{j=1}^{\ell} \alpha_j\cdot dz_j/z_j
 \bigr)
 \cdot \pi_{\gminia,\vecalpha} $
are tame,
where $\pi_{\gminia,\vecalpha}$
denote the projections onto $E_{\gminia,\vecalpha}$,
(ii) 
$\det(T-F_j)_{|D_j}=T^{\rank E_{\gminia,\vecalpha}}$
for the expression
$\thetatilde_{\gminia}=
 \sum_{j=1}^{\ell} F_j\cdot dz_j/z_j
+\sum_{j=\ell+1}^nG_j\cdot dz_j$.
\end{itemize}

For any open subset $U\subset X$ and 
$\veca\in\real^{\ell}$,
we set 
\[
 \nbigp_{\veca}\nbigelambda(U):=
 \Bigl\{
 f\in \nbigelambda(U\setminus D)\,\big|\,
 |f|_h=O\Bigl(\prod_{i=1}^{\ell}|z_i|^{-a_i-\epsilon}\Bigr)\,\,
 \forall \epsilon>0
 \Bigr\}
\]
Thus, we obtain an increasing sequence of
$\nbigo_X$-modules
$\nbigp_{\ast}\nbigelambda:=
 \bigl(\nbigp_{\veca}\nbigelambda\,\big|\,
 \veca\in\real^{\ell}\bigr)$.
We obtain an $\nbigo_X(\ast D)$-module
$\nbigp\nbigelambda:=
 \bigcup_{\veca}\nbigp_{\veca}\nbigelambda$.
\begin{prop}
\mbox{{}}
\begin{itemize}
\item
(Subsection {\rm 11.4} of {\rm\cite{mochi7}})
$(\nbigp_{\ast}\nbigelambda,\DDlambda)$ is 
an unramifiedly good filtered
$\lambda$-flat bundle.
The set of irregular values
is given by
\[
 \Irr(\DDlambda,\nbigp\nbigelambda)
=\bigl\{
 (1+|\lambda|^2)\cdot\gminia\,\big|\,
 \gminia\in\Irr(\theta)
 \bigr\}.
\]
\item
(Subsection {\rm 12.2} of {\rm \cite{mochi7}})
$\kmsmap(\lambda)$ induces 
the bijection
$\KMS(\nbige^0,i)
\lrarr
 \KMS(\nbigelambda,i)$
for each $i$.
We also have
$ \dim\lefttop{i}\Gr^{F,\EE}_{a,\alpha}(\nbigp\nbige^0)
=\dim\lefttop{i}\Gr^{F,\EE}_{\kmsmap(\lambda,(a,\alpha))}
 \bigl(\nbigp\nbigelambda\bigr)$.
\hfill\qed
\end{itemize}
\end{prop}

Take an auxiliary sequence for $\Irr(\theta)$.
Let $\Irr(\theta,\vecm(0))$
denote the image of
$\Irr(\theta)$ via $\etabar_{\vecm(0)}$.
If $\lambda\neq 0$,
for each small sector $S$ in $\{\lambda\}\times(X-D)$,
we have the Stokes filtration
$\nbigf^S$ in the level $\vecm(0)$,
indexed by the ordered set
$\bigl\{
 (1+|\lambda|^2)\cdot\gminia\,\big|\,
 \gminia\in\Irr(\theta,\vecm(0))
 \bigr\}$ with $\leq_S$.
We have the following characterization
of the filtration by
the growth order of the norms of flat sections
with respect to $h$.
(See Subsection 11.4.1 of \cite{mochi7}
for more details.)

\begin{prop}
Assume $\lambda\neq 0$.
Let $f$ be a flat section of 
$\nbigelambda_{|S}$.
We have 
$f\in \nbigf^{S}_{(1+|\lambda|^2)\gminib}$
for $\gminib\in\Irr(\theta,\vecm(0))$,
if and only if 
\[
\Bigl|
 f\cdot\exp\bigl(
 (\lambda^{-1}+\lambdabar)
 \cdot\gminib
 \bigr)
\Bigr|_h
=O\Bigl(
 \exp\bigl(C\cdot|\vecz^{\vecm(1)}|
\bigr)
 \cdot \prod_{k(1)<j\leq \ell}
 |z_j|^{-N}
 \Bigr)
\]
holds for some $C>0$ and $N>0$,
where $k(1)$ is determined by
$\vecm(1)\in\seisuu_{<0}^{k(1)}\times
 \veczero_{\ell-k(1)}$.
\hfill\qed
\end{prop}

\subsubsection{Prolongment
 $\nbigpzero_{\ast}\nbige$}

It is important to consider families for $\lambda$.
In the tame case,
the family $\bigcup_{\lambda}\nbigp\nbigelambda$
gives a regular family of 
meromorphic $\lambda$-flat bundles.
More precisely,
if we consider the sheaf of
holomorphic sections of $\nbige$
of polynomial growth,
then (i) it is a locally free
$\nbigo_{\nbigx}(\ast\nbigd)$-module,
(ii) the specialization at each $\{\lambda\}\times X$
is naturally isomorphic to
$\nbigp\nbigelambda$.
(We need some more consideration
to take nice lattices.)

However, it does not in the non-tame case,
as suggested by the fact that the sets
\[
 \Irr(\nbigp\nbigelambda,\DDlambda)
=\bigl\{
 (1+|\lambda|^2)\gminia\,\big|\,
 \gminia\in\Irr(\theta)
 \bigr\}
\]
depend on $\lambda$
in a non-holomorphic way.
We consider an auxiliary family of
meromorphic $\lambda$-flat bundles
$\nbigpzero\nbige$.
We explain it under the above setting.

Let $\pi_{\gminia,\vecalpha}$ denote the projection
onto $E_{\gminia,\vecalpha}$
in (\ref{eq;08.8.2.10}).
We set 
\[
 g(\lambda):=
 \prod_{\gminia,\vecalpha}
 \exp\Bigl(
 \lambda\cdot \bigl(
 \gminiabar+\sum \alphabar_j\cdot\log|z_j|^2
 \bigr) \Bigr)
 \cdot\pi_{\gminia,\vecalpha}
\]
Let $U(\lambda_0)$ denote 
a small neighbourhood of $\lambda_0\in\cnum$.
We set 
$\nbigx^{(\lambda_0)}:=U(\lambda_0)\times X$,
and $\nbigd^{(\lambda_0)}:=U(\lambda_0)\times D$.
We also set
$\nbigx^{\lambda}:=\{\lambda\}\times X$
and $\nbigd^{\lambda}:=\{\lambda\}\times D$.
Let $p_{\lambda}$ be the projection
of $\nbigx^{(\lambda_0)}-\nbigd^{(\lambda_0)}$
onto $X-D$.
We consider the hermitian metric
\[
 \nbigp^{(\lambda_0)}h
:=g(\lambda-\lambda_0)^{\ast}h
\]
of $p_{\lambda}^{-1}E$
on $\nbigx^{(\lambda_0)}-\nbigd^{(\lambda_0)}$.
Let $\veca\in\real^{\ell}$.
For any open subset $V$ of $\nbigxzero$,
we define
\[
 \nbigpzero_{\veca}\nbige(V):=
 \Bigl\{
 f\in \nbige(V^{\ast})\,\Big|\,
 |f|_{\nbigpzero h}=O\Bigl(
 \prod_{j=1}^{\ell}|z_j|^{-a_j-\epsilon}
 \Bigr),\,\,\forall\epsilon>0
 \Bigr\}
\]
where $V^{\ast}:=V\setminus \nbigdzero$.
Thus, we obtain an increasing sequence 
$\nbigpzero_{\ast}\nbige=\bigl(
 \nbigpzero_{\veca}\nbige\, \big|\,
 \veca\in\real^{\ell}
 \bigr)$
of $\nbigo_{\nbigx^{(\lambda_0)}}$-modules.
We put $\nbigpzero\nbige:=
 \bigcup_{\veca\in\real^{\ell}}
 \nbigpzero_{\veca}\nbige$.
The restrictions to $\nbigx^{\lambda}$
are denoted by
$\nbigpzero_{\ast}\nbigelambda$
and $\nbigpzero\nbigelambda$.

\begin{prop}
\mbox{{}}
\begin{itemize}
\item
(Subsections {\rm 13.1} and {\rm 13.2}
of {\rm \cite{mochi7}})
$(\nbigpzero_{\ast}\nbige,\DD)$
is an unramifiedly good family of
filtered $\lamda$-flat bundles.
The set of irregular values is given by
\[
 \Irr(\nbigpzero\nbige,\DD)
=\bigl\{
 (1+\lambda\lambdabar_0)\cdot\gminia\,\big|\,
 \gminia\in\Irr(\theta)
 \bigr\}.
\]
\item
(Subsection {\rm 13.2.1} of {\rm\cite{mochi7}})
Recall that we have the deformation
mentioned in 
Subsections {\rm \ref{subsubsection;08.8.3.3}}
and {\rm\ref{subsection;08.8.3.4}},
for which 
$(\nbigpzero\nbigelambda,\DDlambda)$
is isomorphic to
$(\nbigp\nbigelambda,\DDlambda)^{T(\lambda)}$
with
$T(\lambda)=
(1+|\lambda|^2)^{-1}
 \cdot(1+\lambda\lambdabar_0)$.
\item
(Subsection {\rm 13.2.3} of {\rm \cite{mochi7}})
Let $U(\lambda_1)\subset U(\lambda_0)$ be small,
and we set 
$\nbigx^{(\lambda_1)}:=U(\lambda_1)\times X$.
Then,
$(\nbigp^{(\lambda_1)}\nbige, \DD)$ 
on $\nbigx^{(\lambda_1)}$
is isomorphic to the deformation
$(\nbigp^{(\lambda_0)}\nbige,\DD)^{
 (T(\lambda_0,\lambda_1))} 
 _{|\nbigx^{(\lambda_1)}}$
with
$T(\lambda_0,\lambda_1)
=(1+\lambda\lambdabar_0)^{-1}
 (1+\lambda\lambdabar_1)$.
\hfill\qed
\end{itemize}
\end{prop}
We should remark that
$\nbigpzero h\neq h$
even in the tame case,
and hence 
$\nbigpzero_{\veca}\nbige$ are different
from $\prolongg{\veca}{\nbige}$
in \cite{mochi2} in the tame case.
We can avoid to use
$\nbigpzero_{\veca}\nbige$
by considering KMS structure
in the tame case.

By the property (D2) of the deformation
(Subsection \ref{subsubsection;08.8.3.3})
and the correspondence between
$\KMS(\nbigp\nbigelambda,i)$
and $\KMS(\nbigp\nbige^0,i)$,
we can show the following.
\begin{lem}
\label{lem;08.8.3.12}
It has the KMS-structure at $\lambda_0$
with the index sets
$\KMS(\nbigp\nbige^0,i)$
$(i=1,\ldots,\ell)$.
\hfill\qed
\end{lem}

\subsubsection{Prolongment
 $\nbigq^{(\lambda_0)}_{\ast}\nbige$
 and $\nbigq\nbige$}
\label{subsection;08.10.4.10}

Applying the deformation procedure
to $(\nbigpzero_{\ast}\nbige,\DD)$ with 
$T=(1+\lambda\lambdabar_0)^{-1}$,
we obtain a family of good filtered 
$\lambda$-flat bundles
$\bigl(
 \nbigq^{(\lambda_0)}_{\ast}\nbige,
 \DD\bigr)$ on 
$(\nbigx^{(\lambda_0)},\nbigd^{(\lambda_0)})$.
Then,
$\nbigq^{(\lambda_0)}_{\veca}\nbige$
is an unramifiedly good lattice of 
$\nbigqzero\nbige$
with the good set of irregular values
$\Irr(\nbigq^{(\lambda_0)}\nbige,\DD)
=\Irr(\theta)$,
i.e.,
\[
 (\nbigqzero_{\veca}\nbige,\DD)_{|\nbigdhat}
=\bigoplus_{\gminia\in \Irr(\theta)}
 \bigl(
 \nbigqzero_{\veca}\nbigehat_{\gminia},
 \DDhat_{\gminia}
\bigr),
\]
such that
$\DDhat_{\gminia}-d\gminia\cdot\id$
has logarithmic singularity for each $\gminia$.
By using the property (D1) of the deformation
explained in Subsection \ref{subsubsection;08.8.3.3},
we obtain the following.
(See Subsection 15.1.1 of \cite{mochi7}
for more details.)

\begin{lem}
\label{lem;08.7.31.5}
$(\nbigqzero\nbige,\DD)_{|\nbigx^{\lambda}}$
is naturally isomorphic to
$\bigl(
 \nbigp\nbigelambda,\DDlambda
 \bigr)^{T_1(\lambda)}$
with $T_1(\lambda)=(1+|\lambda|^2)^{-1}>0$.
\hfill\qed
\end{lem}

By the property (D1) of the deformation,
we have 
$\nbigqzero\nbige_{|\nbigx^{(\lambda_1)}}
=\nbigq^{(\lambda_1)}\nbige$.
Hence, we obtain the global 
family of meromorphic $\lambda$-flat bundles
$(\nbigq\nbige,\DD)$
on $\cnum_{\lambda}\times (X,D)$.
By using the property (D2) of the deformation
and Lemma \ref{lem;08.8.3.12},
we can show the following.
\begin{lem}
For each $\lambda_0$,
$(\nbigq\nbige,\DD)$
has the $KMS$-structure at $\lambda_0$
indexed by
$\KMS\bigl(\nbigp\nbige^0,i\bigr)$
$(i=1,\ldots,\ell)$.
\hfill\qed
\end{lem}

Let $S$ be a small sector in $\{\lambda\}\times (X-D)$.
By Lemma \ref{lem;08.7.31.5},
the Stokes filtrations of
$\nbigq\nbigelambda$ and
$\nbigp\nbigelambda$ in the level $\vecm(0)$
are related as follows:
\[
 \nbigf^S_{\gminia}
 \bigl(\nbigq\nbigelambda_{|\Sbar}\bigr)
=
 \nbigf^S_{(1+|\lambda|^2)\gminia}
 \bigl(\nbigp\nbigelambda_{|\Sbar}\bigr),
\quad
 \gminia\in\Irrbar\bigl(\theta,\vecm(0)\bigr)
\]
Hence, we have the characterization
of the Stokes filtrations of $\nbigq$
in the level $\vecm(0)$,
by growth order of the norms of flat sections
with respect to $h$.
(See Subsection 15.1.1 of \cite{mochi7}
for more details.)

\begin{prop}
\label{prop;08.8.10.20}
Let $f$ be a flat section of $\nbigelambda_{|S}$.
We have 
$f\in \nbigf^{S}_{\gminib}(\nbigq\nbige^{\lambda}_{|S})$
for $\gminib\in\Irrbar(\theta,\vecm(0))$,
if and only if 
\[
\Bigl|
 f\cdot\exp\bigl(
 (\lambda^{-1}+\lambdabar)
 \cdot \gminib
 \bigr)
\Bigr|_h
=O\Bigl(
 \exp\bigl(C\cdot|\vecz^{\vecm(1)}|
  \bigr)
\cdot
 \prod_{k(1)<j\leq\ell}|z_j|^{-N}
 \Bigr)
\]
holds for some $C>0$ and $N>0$,
where $k(1)$ is determined by
$\vecm(1)\in\seisuu_{<0}^{k(1)}\times
 \veczero_{\ell-k(1)}$.
\hfill\qed
\end{prop}

By taking Gr with respect to 
the Stokes filtration $\nbigf^{S}$
in the level $\vecm(0)$
explained in Subsection
\ref{subsubsection;08.8.3.3},
we obtain an unramifiedly good lattice
$\bigl(
 \Gr^{\vecm(0)}_{\gminia}(\nbigq\nbige),
 \DD_{\gminia}
\bigr)$.

\vspace{.1in}

In the case that $D$ is smooth,
we have the following characterization
of the full Stokes filtration $\nbigftilde^S$
(Subsection 15.1.1 of \cite{mochi7}).
\begin{prop}
\label{prop;08.10.19.1}
Let $f$ be a flat section of $\nbigelambda_{|S}$.
We have 
$f\in \nbigftilde^{S}_{\gminib}(\nbigq\nbige^{\lambda}_{|S})$
for $\gminib\in\Irrbar(\theta)$,
if and only if 
\[
\Bigl|
 f\cdot\exp\bigl(
 (\lambda^{-1}+\lambdabar)
 \cdot \gminib
 \bigr)
\Bigr|_h
=O\bigl(
 |z_1|^{-N}
  \bigr)
\]
holds for some $N>0$.
\hfill\qed
\end{prop}

\begin{rem}
We have a characterization of full Stokes filtrations
or more general Stokes filtrations
in the level $\vecm(i)$,
even in the general normal crossing case.
\hfill\qed
\end{rem}

\subsection{Reduction from wild to tame}
\label{subsection;08.8.3.15}

Let $X$, $D$ and $\harmonicbundle$
be as in Subsection \ref{subsection;08.8.3.6}.
By making the same procedure
to $(E,\del_E,\theta^{\dagger},h)$
on $X^{\dagger}-D^{\dagger}$,
we obtain the family of meromorphic 
$\mu$-flat bundles 
$\bigl(
 \nbigq\nbige^{\dagger},\DD^{\dagger}
 \bigr)$
on $\cnum_{\mu}\times (X^{\dagger},D^{\dagger})$.

\begin{lem}
The correspondence
$(a,\alpha)\longleftrightarrow
 (-a,\alphabar)$ induces
a bijection
$\KMS(\nbigp\nbige^0,i)
\simeq
 \KMS(\nbigp\nbige^{\dagger\,0},i)$.
We also have the bijection
$\Irr(\theta)\simeq \Irr(\theta^{\dagger})$
given by $\gminia\longleftrightarrow \gminiabar$.
\end{lem}
\pf
The claim for $\Irr(\theta)$ and $\Irr(\theta^{\dagger})$
is clear.
See Corollary 11.12 of \cite{mochi2}
for the correspondence between
$\KMS(\nbigp\nbige^0,i)$
and 
$\KMS(\nbigp\nbige^{\dagger\,0},i)$.
\hfill\qed

\paragraph{One step reduction I}

Since both the Stokes filtrations of
$(\nbigq\nbige^{\lambda},\DDlambda)$ and 
$\bigl(\nbigq\nbige^{\dagger\,\mu},
 \DD^{\dagger\mu}\bigr)$
are characterized by
growth order of the norms of flat sections
with respect to $h$,
we have the induced isomorphisms of
the associated graded family of flat bundles
for $\gminia\in\Irr(\theta,\vecm(0))$:
\[
 \bigl(
 \Gr^{\vecm(0)}_{\gminia}\nbigq\nbige,
 \DD^{f}_{\gminia}
 \bigr)_{|\cnum_{\lambda}^{\ast}\times(X-D)}
\simeq
 \bigl(
 \Gr^{\vecm(0)}_{\gminiabar}\nbigq\nbige^{\dagger},
 \DD^{\dagger\,f}_{\gminiabar}
 \bigr)_{|\cnum_{\mu}^{\ast}\times(X-D)}
\]
Hence, they give a variation of
$\proj^1$-holomorphic vector bundles
denoted by
$\Gr^{\vecm(0)}_{\gminia}
 (\nbige^{\sankaku},
 \DD^{\sankaku})$
on $\proj^1\times (X-D)$.

We can show that the pairing
$\nbigs:
(\nbige,\DD)\otimes
 \sigma^{\ast}(\nbige^{\sankaku},\DD^{\sankaku})
\lrarr \nbigo_{\nbigx-\nbigd}$
is extended to 
\[
 \nbigq\nbige\otimes
 \sigma^{\ast}\nbigq\nbige^{\dagger}
\lrarr \nbigo_{\cnum_{\lambda}\times X}
 \bigl(\ast (\cnum_{\lambda}\times D)\bigr).
\]
(See Subsection 15.1.3 of \cite{mochi7}.)
By functoriality of Gr
with respect to Stokes structures,
we obtain 
\[
 \Gr_{\gminia}^{\vecm(0)}
 (\nbigq\nbige,\DD)
 \otimes
 \sigma^{\ast}
 \Gr_{\gminiabar}^{\vecm(0)}
 (\nbigq\nbige^{\dagger},
 \DD^{\dagger})
\lrarr \nbigo_{\cnum_{\lambda}\times X}
 \bigl(\ast(\cnum_{\lambda}\times D)\bigr).
\]
Similarly,
we obtain
$\Gr_{\gminiabar}^{\vecm(0)}
 (\nbigq\nbige^{\dagger},
 \DD^{\dagger})
 \otimes
 \sigma^{\ast}
 \Gr_{\gminia}^{\vecm(0)}(\nbigq\nbige,
 \DD)
\lrarr \nbigo_{\cnum_{\mu}\times X^{\dagger}}
 \bigl(\ast(\cnum_{\mu}\times D^{\dagger})\bigr)$.
They give a morphism of variations of 
$\proj^1$-holomorphic vector bundles:
\[
 \Gr^{\vecm(0)}_{\gminia}(\nbigs):
 \Gr^{\vecm(0)}_{\gminia}
 (\nbige^{\sankaku},\DD^{\sankaku})
\otimes
 \sigma^{\ast} \Gr^{\vecm(0)}_{\gminia}
 (\nbige^{\sankaku},\DD^{\sankaku})
\lrarr
 \Tate(0)
\]

One of the main result in the study of
wild harmonic bundles is the following.
(See Subsection 15.2 of \cite{mochi7} for more details.)
\begin{prop}
\mbox{{}}\label{prop;08.8.3.5}
If we shrink $X$ appropriately,
the following holds:
\begin{itemize}
\item
$\Gr^{\vecm(0)}_{\gminia}\bigl(
 \nbige^{\sankaku},
 \DD^{\sankaku},
 \nbigs \bigr)$
 is a variation of pure polarized
 twistor structures.
\item
Let $(E_{\gminia},\delbar_{\gminia},
 h_{\gminia},\theta_{\gminia})$
denote the underlying harmonic bundle
for $\gminia\in\Irr(\theta,\vecm(0))$.
By construction,
the Higgs bundle
$(E_{\gminia},\theta_{\gminia})$ is naturally isomorphic to
\[
 \bigoplus_{
 \substack{\gminib\in\Irr(\theta)\\
 \etabar_{\vecm(0)}(\gminib)=\gminia}}
\bigoplus_{\vecalpha}
 \bigl(E_{\gminib,\vecalpha},\theta_{\gminib,\vecalpha}\bigr)
\]
(Recall the decomposition {\rm(\ref{eq;08.8.2.10})}).
In particular, the harmonic bundle is
unramifiedly good wild.
The set of irregular values is
$\etabar_{\vecm(0)}^{-1}(\gminia)$.
\item
Let $(\nbigq\nbige_{\gminia},\DD_{\gminia})$
be the family of 
meromorphic $\lambda$-flat bundles
on $\cnum_{\lambda}\times (X,D)$
associated to $(E_{\gminia},\delbar_{\gminia},
 h_{\gminia},\theta_{\gminia})$.
Then, we have the natural isomorphism
$(\nbigq\nbige_{\gminia},\DD_{\gminia})
\simeq
 \Gr_{\gminia}^{\vecm(0)}
 (\nbigq\nbige,\DD)$.
\item
Similarly, let 
$(\nbigq\nbige^{\dagger}_{\gminiabar},
 \DD^{\dagger}_{\gminiabar})$
denote the associated family of
meromorphic $\mu$-flat bundles
on $\cnum_{\mu}\times (X,D)$.
Then, we have the natural isomorphism
$(\nbigq\nbige^{\dagger}_{\gminiabar},
 \DD_{\gminiabar}^{\dagger})
\simeq
 \Gr^{\vecm(0)}_{\gminiabar}
 \bigl(
 \nbigq\nbige^{\dagger},
 \DD^{\dagger}
 \bigr)$.
\hfill\qed
\end{itemize}
\end{prop}

\paragraph{One step reduction II}

Let $\Irrbar(\theta,\vecm(j))$
denote the image of $\Irr(\theta)$
via $\etabar_{\vecm(j)}$.
For each $\gminia\in\Irrbar(\theta,\vecm(j))$,
we obtain a variation of $\proj^1$-holomorphic bundles
with a pairing
$\Gr^{\vecm(j)}_{\gminia}
 \bigl(\nbige^{\sankaku},\DD^{\sankaku},\nbigs\bigr)$,
which is naturally isomorphic to
$\Gr^{\vecm(j)}_{\gminia}
 \Gr^{\vecm(j-1)}_{\etabar_{\vecm(j-1)}(\gminia)}
 (\nbige^{\sankaku},\DD^{\sankaku},\nbigs)$.
We explain how to apply Proposition \ref{prop;08.8.3.5}
in this situation.

Let us consider the case in which
$\Irrbar(\theta,\vecm(j-1))$ consists
one element.
We take any $\gminia\in\Irr(\theta)$.
Let $L(-\gminia)$ be the variation of
polarized pure twistor structures
as in Subsection \ref{subsubsection;08.8.3.20}.
The underlying harmonic bundle
is also denoted by $L(-\gminia)$.
We set $(E',\delbar_{E'},\theta',h'):=
 (E,\delbar_E,\theta,h)\otimes L(-\gminia)$.
Note $\Irr(\theta'):=\bigl\{
 \gminia'-\gminia\,\big|\,
 \gminia'\in\Irr(\theta)
 \bigr\}$,
and hence $\vecm(j),\vecm(j+1),\ldots,\vecm(L)$
give an auxiliary sequence for
$\Irr(\theta')$.
We have the natural isomorphisms of
the associated variation of polarized pure twistor structures:
\[
  (\nbige^{\sankaku},\DD^{\sankaku},\nbigs)
\simeq
(\nbige^{\prime\sankaku},
 \DD^{\prime\sankaku},\nbigs')
\otimes L(\gminia)
\]
For each $\gminib\in\Irr(\theta,\vecm(j))$,
we have the natural isomorphism:
\[
 \Gr^{\vecm(j)}_{\gminib}
 (\nbige^{\sankaku},\DD^{\sankaku},\nbigs)
\simeq
 \Gr^{\vecm(j)}_{\gminib-\etabar_{\vecm(j)}(\gminia)}
 (\nbige^{\prime\sankaku},
 \DD^{\prime\sankaku},\nbigs')
\otimes
 L(\gminia)
\]
Hence, 
by shrinking $X$ appropriately,
we obtain that 
$ \Gr^{\vecm(j)}_{\gminib}
 (\nbige^{\sankaku},\DD^{\sankaku},\nbigs)$
is also a variation of pure twistor structures,
due to Proposition \ref{prop;08.8.3.5}.

\paragraph{Full reduction}

Let us consider the general case.
By using the above result inductively,
we obtain that
$\Gr^{\vecm(j)}_{\gminia}
 \bigl(\nbige^{\sankaku},\DD^{\sankaku},\nbigs\bigr)$
are variations of polarized pure twistor structures
for any $\gminia\in\Irr(\theta,\vecm(j))$.
The underlying Higgs field is
\[
 \bigoplus_{\substack{
 \gminib\in\Irr(\theta)\\
 \etabar_{\vecm(j)}(\gminib)=\gminia
 }}
 \bigoplus_{\vecalpha}
 (E_{\gminib,\vecalpha},\theta_{\gminib,\vecalpha})
\]
For any $\gminia\in\Irr(\theta)$,
we set 
$\Gr^{\full}_{\gminia}(\nbige,\DD^{\sankaku},\nbigs)
:=\Gr^{\vecm(L)}_{\gminia}(\nbige,\DD^{\sankaku},\nbigs)$,
which are called the full reductions.
Let $(E_{\gminia},\delbar_{\gminia},h_{\gminia})$
be the underlying harmonic bundles.
Then, $(E_{\gminia},\delbar_{\gminia},h_{\gminia})
 \otimes L(-\gminia)$ are tame.
This procedure is the reduction
from wild harmonic bundles to tame harmonic bundles.

\subsection{Reduction
 from tame to twistor nilpotent orbit}
\label{subsection;08.8.4.10}

Let $X:=\Delta^n$, $D_i=\{z_i=0\}$
and $D:=\bigcup_{i=1}^{\ell}D_i$.
Let $\harmonicbundle$ be a tame harmonic bundle
on $X-D$.
The family of $\lambda$-flat bundles
$(\nbige,\DD)$ is prolonged to
a family of meromorphic $\lambda$-flat
bundle $\bigl(\nbigq\nbige,\DD\bigr)$,
which has the KMS-structure
at $\lambda_0$ 
indexed by
$\KMS(\nbigp\nbige^0,i)$
$(i=1,\ldots,\ell)$
for each $\lambda_0\in\cnum_{\lambda}$.
For later use,
we recall how to obtain 
the limiting mixed twistor structure.
For simplicity,
we assume 
$\KMS\bigl(\nbige^0,i\bigr)\subset 
 \real\times\{0\}$.
See Section 11 of \cite{mochi2}
for the general case.
See also an account
due to Hertling and Sevenheck
in \cite{Hertling-Sevenheck3}
for this case.

In a neighbourhood $U(\lambda_0)$ of $\lambda_0$,
we set
\[
 \nbigg^{(\lambda_0)}_{(\veca,0)}(E):=
 \lefttop{\ellsitabar}\nbigg^{(\lambda_0)}_{(\veca,0)}
 (\nbigqzero\nbige)
_{|U(\lambda_0)\times\{O\}}
\]
for $\veca\in \Par(\nbigp\nbige^0,\ellsitabar)$.
(See (\ref{eq;08.8.3.7})
 for the right hand side.
 In this simpler case,
 we have only to take Gr with respect to
 parabolic filtrations.)
By varying $\lambda_0\in\cnum_{\lambda}$
and gluing them,
we obtain the vector bundle
$\nbigg_{(\veca,0)}(E)$
on $\cnum_{\lambda}$.
It is endowed with the nilpotent maps
$\nbign_i$ $(i=1,\ldots,\ell)$,
which are the nilpotent part of the residues
$\Res_i(\DD)$.
By applying the same procedure
to $(E,\del_E,\theta^{\dagger},h)$
on $X^{\dagger}-D^{\dagger}$,
we obtain the vector bundle
$\nbigg_{(-\veca,0)}^{\dagger}(E)$
on $\cnum_{\mu}$
with nilpotent endomorphisms $\nbign_i^{\dagger}$
induced by residues
$\Res_i(\DD^{\dagger})$.
We would like to glue
$\nbigg_{(\veca,0)}(E)$
and $\nbigg_{(-\veca,0)}^{\dagger}(E)$,
to obtain a vector bundle
$S^{\can}_{(\veca,0)}(E)$ on $\proj^1$.

We have the $\DD$-flat decomposition
$ \nbigq_0\nbige_{|\cnum_{\lambda}^{\ast}\times X}
=\bigoplus_{\veca\in \Par(\nbigp_0\nbige^0,\ellsitabar)}
 \nbigg_{(\veca,0)}\nbige$
with the following property:
\begin{itemize}
\item
Let $M_i$ be the family of the monodromy 
endomorphisms along the path 
$(z_1,\ldots,e^{2\pi\sqrt{-1}\theta}z_i,
 \ldots,z_n)$ ($0\leq \theta\leq 1$)
 with respect to $\DD^f$.
Then, the restriction of $M_i$
to $\nbigg_{(\veca,0)}\nbige$ 
has the unique eigenvalue
$\exp\bigl(2\pi\sqrt{-1} a_i\bigr)$.
\item
$\nbigg_{(\veca,0)}
 \nbige_{|\cnum_{\lambda}^{\ast}\times O}
\simeq\nbigg_{(\veca,0)}(E)$.
\end{itemize}
For $\lambda\neq 0$,
let $H(\nbigelambda)$ be the space
of multivalued flat sections of 
$(\nbigelambda,\DDlambda)$.
We have the holomorphic vector bundle
$\nbigh(E)$ on $\cnum_{\lambda}^{\ast}$
whose fiber over $\lambda$
is $H(\nbigelambda)$.
We have the decomposition
\[
 \nbigh(E)=
 \bigoplus_{
 \veca\in \Par(\nbigp_0\nbige^0,\ellsitabar)}
 \nbigg_{(\veca,0)}\nbigh(E)
\]
such that
(i) it is preserved by the monodromy $M_i$,
(ii) the restriction of $M_i$ to
 $\nbigg_{(\veca,0)}\nbigh(E)$
 has the unique eigenvalue
 $\exp\bigl(2\pi\sqrt{-1} a_i\bigr)$.

Let $U\subset\cnum_{\lambda}^{\ast}$,
and let $s$ be a section of 
$\nbigg_{(\veca,0)}\nbigh(E)$ on $U$.
We regard $s$ as a multi-valued flat section
of $\nbigg_{(\veca,0)}\nbige$.
It is expressed as a finite sum:
\[
 s=\sum f_{\vecm}\cdot 
 \prod_{i=1}^{\ell} \exp\bigl(a_i\log z_i\bigr)
 \cdot\bigl(\log z_i\bigr)^{m_i}
\]
Here, $f_{\vecm}$ are holomorphic sections of
$\nbigg_{(\veca,0)}\nbige_{|U\times X}$.
We set
$\Phi^{\can}_{(\veca,0)}(s)=f_{0|U\times O}$,
and thus we obtain an  isomorphism
\[
 \Phi^{\can}_{(\veca,0)}:
 \nbigg_{(\veca,0)}\nbigh(E)
\lrarr
 \nbigg_{(\veca,0)}
 \nbige_{|\cnum_{\lambda}^{\ast}\times O}
=\nbigg_{(\veca,0)}(E).
\]

Let $\vecdelta=(1,\ldots,1)\in\real^{\ell}$.
We have the $\DD^{\dagger}$-flat decomposition
$ \nbigq_{<\vecdelta}
 \nbige^{\dagger}_{|\cnum_{\lambda}^{\ast}\times X}
=\bigoplus_{\veca\in \Par(\nbigp_0\nbige^0,\ellsitabar)}
 \nbigg_{(-\veca,0)}\nbige^{\dagger}$
with the following property:
\begin{itemize}
\item
The restriction of $M_i^{-1}$
to $\nbigg_{(-\veca,0)}\nbige^{\dagger}$ 
has the unique eigenvalue
$\exp\bigl(-2\pi\sqrt{-1} a_i\bigr)$.
(Because the base space is the complex conjugate
 $X^{\dagger}-D^{\dagger}$, 
 the direction of the loop is reversed.)
\item
$\nbigg_{(-\veca,0)}
 \nbige^{\dagger}_{|\cnum_{\mu}^{\ast}\times O}
\simeq\nbigg_{(-\veca,0)}^{\dagger}(E)$.
\end{itemize}

Similarly,
let $\nbigh^{\dagger}(E)$
be the holomorphic vector bundle on $\cnum_{\mu}^{\ast}$
whose fiber over $\mu$ is
the space of the multi-valued flat sections of
$(\nbige^{\dagger\,\mu},\DD^{\dagger\,\mu})$.
We have the decomposition
\[
 \nbigh^{\dagger}(E)
=\bigoplus_{\veca\in \Par(\nbigp_0\nbige^0,\ellsitabar)}
 \nbigg_{(-\veca,0)}\nbigh^{\dagger}(E)
\]
such that the restriction of $M_i^{-1}$ to 
$\nbigg_{(-\veca,0)} \nbigh^{\dagger}(E)$
has the unique eigenvalue 
$\exp\bigl(-2\pi\sqrt{-1} a_i\bigr)$.
For a section $s$ of 
$\nbigg_{(-\veca,0)}\nbigh^{\dagger}(E)_{|U}$,
we have an expression
\[
 s=\sum f_{\vecm}^{\dagger}\cdot 
 \prod_{i=1}^{\ell}
 \exp\bigl(-a_i\log\zbar_i\bigr)
 \cdot\bigl(\log \zbar_i\bigr)^{m_i},
\]
where 
$f^{\dagger}_{\vecm}$ are sections 
of $\nbigq\nbige^{\dagger}_{|U\times X}$.
We set $\Phi^{\can\dagger}_{(\veca,0)}(s)
=f^{\dagger}_{0|U\times O}$,
and thus we obtain an isomorphism
\[
 \Phi^{\can\,\dagger}_{(\veca,0)}:
 \nbigg_{(-\veca,0)}\nbigh^{\dagger}(E)
\lrarr
 \nbigg_{(-\veca,0)}
 \nbige^{\dagger}_{|\cnum_{\mu}^{\ast}\times O}
=\nbigg_{(-\veca,0)}^{\dagger}(E).
\]
By construction,
we have the natural isomorphism
$\nbigg_{(\veca,0)}\nbigh(E)
\simeq
 \nbigg_{(-\veca,0)}\nbigh^{\dagger}(E)$
under the identification of
$\cnum_{\lambda}^{\ast}=\cnum_{\mu}^{\ast}$
via $\mu=\lambda^{-1}$.
Thus, we obtain the vector bundle
$S^{\can}_{(\veca,0)}(E)$
by gluing 
$\nbigg_{(\veca,0)}(E)$
and
$\nbigg_{(-\veca,0)}^{\dagger}(E)$.
Under the gluing, 
we have the relation
\[
 \lambda^{-1}\nbign_i=-\mu^{-1}\nbign_i^{\dagger}.
\]
Thus, 
$\nbign_i\cdot t_0^{(-1)}$
and $\nbign_{\infty}\cdot t_{\infty}^{(-1)}$
give the morphism
$\nbign^{\sankaku}_i:
 S^{\can}_{(\veca,0)}(E)
\lrarr 
S^{\can}_{(\veca,0)}(E)
\otimes\Tate(-1)$.
The tuple of them is denoted by $\vecN^{\sankaku}$.

The morphism
$\nbigs_0:
 \nbige\otimes\sigma^{\ast}\nbige^{\dagger}
\lrarr \nbigo_{\nbigx-\nbigd}$
is extended to
$\nbigq_{0}\nbige\otimes
 \sigma^{\ast}\nbigq_{<\vecdelta}\nbige^{\dagger}
\lrarr \nbigo_{\nbigx}$.
Similarly, we have 
$\nbigq_{<\vecdelta}\nbige^{\dagger}\otimes
 \sigma^{\ast}\nbigq_0\nbige
\lrarr \nbigo_{\nbigx^{\dagger}}$.
They induce 
\[
 \nbigg_{(\veca,0)}(E)
\otimes
 \sigma^{\ast}\nbigg^{\dagger}_{(-\veca,0)}(E)
\lrarr \nbigo_{\cnum_{\lambda}},
\]
\[
 \nbigg^{\dagger}_{(-\veca,0)}(E)
\otimes
 \sigma^{\ast}\nbigg_{(\veca,0)}(E)
\lrarr \nbigo_{\cnum_{\mu}},
\]
\[
 \nbigg_{(-\veca,0)}\nbigh^{\dagger}(E)
\otimes
 \sigma^{\ast}\nbigg_{(\veca,0)}\nbigh(E)
\lrarr \nbigo_{\cnum_{\lambda}^{\ast}}
\]
They are preserved by the above isomorphisms.
Hence, we obtain
$\nbigs_{(\veca,0)}:
 S^{\can}_{(\veca,0)}(E)
\otimes
 \sigma^{\ast}S^{\can}_{(\veca,0)}(E)
\lrarr
 \Tate(0)$.
Theorem 12.22 of \cite{mochi2}
implies the following.

\begin{prop}
\label{prop;08.8.10.10}
$(S^{\can}_{(\veca,0)}(E),
 \vecN^{\sankaku}, \nbigs)$
is a polarized mixed twistor structure
of weight $0$ in $\ell$-variables.
\hfill\qed
\end{prop}
By Theorem \ref{thm;08.7.26.8},
a polarized mixed twistor structure
induces a nilpotent orbit.
This is the reduction
from tame harmonic bundles
to nilpotent orbits.

\begin{rem}
Although the notation is changed,
the construction explained in this subsection
is the same as that in {\rm\cite{mochi2}}.
In the tame case,
$\nbigq\nbige$ is equal to the sheaf
of holomorphic sections 
whose norms with respect to $h$
are of polynomial growth order.
We also remark the uniqueness 
in Lemma {\rm\ref{lem;07.11.23.5}}.
\hfill\qed
\end{rem}

\paragraph{Family version}

The construction can be done in family
on $D_{\ellsitabar}:=\bigcap_{i=1}^{\ell}D_i$.
As in the construction of
$\nbigg_{(\veca,0)}(E)$,
we obtain the vector bundle 
$\lefttop{\ellsitabar}
\nbigg_{(\veca,0)}(\nbigq\nbige)$
on $\nbigd_{\ellsitabar}:=
 \cnum_{\lambda}\times D_{\ellsitabar}$,
as the gluing of
$\lefttop{\ellsitabar}
\nbiggzero_{(\veca,0)}(\nbigqzero\nbige)$.
They are equipped with the nilpotent maps
$\nbign_i$ $(i=1,\ldots,\ell)$.
By applying the nearby cycle functors for $\nbigr$-modules
along $z_i$ $(i=1,\ldots,\ell)$,
or by a direct consideration
as in Subsection 8.8.3 of \cite{mochi2},
we obtain the induced family of
flat $\lambda$-connections $\DD_{\veca,0}$
of $\lefttop{\ellsitabar}
\nbigg_{(\veca,0)}(\nbigq\nbige)$
for which $\nbign_i$ are flat.
Similarly,
we obtain a family of $\mu$-flat bundles
$\bigl(
\lefttop{\ellsitabar}\nbigg_{(-\veca,0)}
 (\nbigq\nbige^{\dagger}),
 \DD^{\dagger}_{(-\veca,0)}\bigr)$
on $\cnum_{\mu}\times D_{\ellsitabar}^{\dagger}$
with flat nilpotent maps $\nbign_i^{\dagger}$.

Let $q:X-D\lrarr D_{\ellsitabar}$ be the projection.
We naturally obtain a holomorphic vector bundle 
$\nbightilde(E)$ on 
$\cnum_{\lambda}^{\ast}
 \times D_{\ellsitabar}$,
whose fiber over $(\lambda,P)$
is the space of multi-valued flat sections
of $(\nbigelambda,\DDlambda)_{|q^{-1}(P)}$.
It has the generalized eigen decomposition
$\nbightilde(E)
=\bigoplus 
 \lefttop{\ellsitabar}\nbigg_{(\veca,0)}
 \nbightilde(E)$
with respect to the monodromy endomorphisms
around $D_i$. $(i=1,\ldots,\ell)$.
It is naturally equipped with
the family of flat connections $\DD_{\veca,0}^f$.

By using the family of flat bundles
$\bigl(\nbigg_{(\veca,0)}\nbige,
 \DD^f\bigr)$,
we obtain the flat isomorphisms
\[
 \Phi^{\can}_{(\veca,0)}:
 \lefttop{\ellsitabar}
 \nbigg_{(\veca,0)}\nbightilde(E)
\lrarr
\lefttop{\ellsitabar}\nbigg_{(\veca,0)}
 (\nbigq\nbige)_{|\cnum_{\lambda}^{\ast}
 \times D_{\ellsitabar}}.
\]
Similarly, we obtain the flat isomorphisms
$ \Phi^{\can}_{(\veca,0)}:
 \lefttop{\ellsitabar}\nbigg_{(\veca,0)}\nbightilde(E)
\lrarr
 \lefttop{\ellsitabar}\nbigg_{(-\veca,0)}
 (\nbigq\nbige^{\dagger})_{|\cnum_{\mu}^{\ast}
 \times D_{\ellsitabar}^{\dagger}}$.
As the gluing,
we obtain a variation of
$\proj^1$-holomorphic vector bundles
$\bigl(
\lefttop{\ellsitabar}
 \nbige^{\sankaku}_{\veca,0},
\DD^{\sankaku}_{\veca,0}
\bigr)$ with a tuple $\vecN^{\sankaku}$
of flat nilpotent morphisms
\[
 \nbign^{\sankaku}_i:
 \lefttop{\ellsitabar}
 \nbige^{\sankaku}_{\veca,0}
\lrarr
 \lefttop{\ellsitabar}
 \nbige^{\sankaku}_{\veca,0}
\otimes
 \Tate(-1),
\quad
 (i=1,\ldots,\ell)
\]
We also have the induced flat symmetric pairing
$\nbigs:
 \lefttop{\ellsitabar}
 \nbige^{\sankaku}_{\veca,0}
\otimes
 \sigma^{\ast}
  \lefttop{\ellsitabar}
 \nbige^{\sankaku}_{\veca,0}
\lrarr \Tate(0)$.
By Proposition \ref{prop;08.8.10.10},
\[
 \bigl(
 \lefttop{\ellsitabar}
 \nbige^{\sankaku}_{\veca,0},
 \vecN^{\sankaku},
 \DD^{\sankaku}_{\veca,0},
 \nbigs_{\veca,0} \bigr)
\]
is a variation of polarized mixed twistor 
structures of weight $0$
in $\ell$-variables.
(See Subsection
 \ref{subsubsection;08.8.10.40}.)

\section{Prolongation and reductions
 in the integrable case}
\label{section;08.8.20.42}

\subsection{Preliminary Estimate}
\label{subsection;08.8.20.44}

\subsubsection{Statements}

Let $X:=\Delta^n$ and $D:=\{z_1=0\}$.
Let $(E,\delbar_E,\theta,h)$
be an unramifiedly good wild harmonic bundle
on $X-D$.
For simplicity,
we assume that there exists a holomorphic decomposition
\begin{equation}
\label{eq;08.7.26.10}
 (E,\theta)=\bigoplus_{\gminia\in\Irr(\theta)}
 (E_{\gminia},\theta_{\gminia})
\end{equation}
such that 
$\theta_{\gminia}-d\gminia\cdot \pi_{\gminia}$
are tame,
where $\pi_{\gminia}$ denotes the projection
onto $E_{\gminia}$ with respect to the decomposition
(\ref{eq;08.7.26.10}).
\begin{rem}
Since $\harmonicbundle$ is
assumed to be unramifiedly good,
such a decomposition exists 
on a neighbourhood of each point of $D$.
Because we are interested in
the behaviour around $O$,
we may assume such a decomposition
exists globally
by replacing $X$ with a small neighbourhood of $O$.
\hfill\qed
\end{rem}

Let $\nbigu$ be a holomorphic section of
$\End(E)$ on $X-D$
such that $[\theta,\nbigu]=0$.
Let $\nbigq$ be a $C^{\infty}$-section of
$\End(E)$ on $X-D$
such that $\nbigq=\nbigq^{\dagger}$.
We assume the following equations:
\begin{equation}
\label{eq;08.7.26.60}
 \del_E\nbigu-[\theta,\nbigq]+\theta=0
\end{equation}
\begin{equation}
\label{eq;08.7.26.61}
 \del_E\nbigq+[\theta,\nbigu^{\dagger}]=0
\end{equation}
We set
$\nbigutilde:=\nbigu
+\sum_{\gminia\in\Irr(\theta)}
 \gminia\cdot\pi_{\gminia}$.
We will prove the following proposition
in Subsections 
\ref{subsubsection;08.7.26.100}--\ref{subsubsection;08.7.26.101}.

\begin{prop}
\label{prop;08.7.26.55}
$\nbigutilde=O(1)$
and $\nbigq=O\Bigl(\bigl(-\log|z_1|\bigr)^M\Bigr)$
for some $M>0$
with respect to $h$.
\end{prop}

\begin{rem}
Eventually,
we obtain that 
$\nbigq$ is bounded.
(See Corollary {\rm\ref{cor;08.8.19.1}}
and Corollary {\rm\ref{cor;08.8.19.2}}.)
See Corollary {\rm\ref{cor;08.9.11.2}}
for the boundedness of $\nbigutilde$
in the case that $D$ is normal crossing.
\hfill\qed
\end{rem}

We set 
$g_{\irr}(\lambda)=
 \exp\Bigl(
 \bigoplus \lambda\gminiabar\cdot\pi_{\gminia}
 \Bigr)$.
Let $\lambda_0\in\cnum$,
and let $U(\lambda_0)$ be a small neighbourhood
of $\lambda_0$ in $\cnum$.
Let $p_{\lambda}$ be the projection
of $U(\lambda_0)\times (X-D)$ onto $X-D$.
We consider the hermitian metric
\begin{equation}
 \label{eq;08.11.6.12}
\nbigp^{(\lambda_0)}_{\irr}h
:=g_{\irr}(\lambda-\lambda_0)^{\ast}h 
\end{equation}
of $p_{\lambda}^{-1}E$
on $U(\lambda_0)\times (X-D)$.
We regard $\nbigu$ and $\nbigq$
as $C^{\infty}$-sections
of $\End\bigl(p_{\lambda}^{-1}E\bigr)$.
We will prove the following lemma
in Subsection \ref{subsubsection;08.7.26.71}.
\begin{prop}
\label{prop;08.7.26.70}
Assume $U(\lambda_0)$ is sufficiently small.
Then,
$\nbigutilde=O(1)$
and 
$\nbigq=O\bigl((-\log|z_1|)^{M}\bigr)$
with respect to $\nbigp^{(\lambda_0)}_{\irr}h$.
\end{prop}

\subsubsection{Preliminary}
\label{subsubsection;08.7.26.100}

We take orthogonal decompositions
$E=\bigoplus E'_{\gminia,\alpha}
=\bigoplus E'_{\gminia}$
as in Subsection \ref{subsubsection;08.7.30.3}.
For any $f\in \End(E)$,
we have the decompositions:
\[
 f=\sum f'_{\gminia,\gminib},
\quad
 f'_{\gminia,\gminib}
\in \Hom\bigl(E'_{\gminib},E'_{\gminia}\bigr)
\]
\[
 f=\sum f'_{(\gminia,\alpha),(\gminib,\beta)},
\quad
 f'_{(\gminia,\alpha),(\gminib,\beta)}
\in
 \Hom\bigl(
 E'_{\gminib,\beta},
 E'_{\gminia,\alpha}
 \bigr).
\]
We have similar decompositions
for sections of
$\End(E)\otimes\Omega^{p,q}$.
The following lemma is easy to show
by using Proposition \ref{prop;08.7.30.2}.

\begin{lem}
\label{lem;08.7.29.10}
Let $f$ be a $C^{\infty}$-section of $\End(E)$
such that $f$ commutes with $\theta$.
\begin{itemize}
\item
If $\gminia\neq\gminib$,
we have
$|f'_{\gminia,\gminib}|_h
=O\Bigl(
 \exp\big(-\epsilon|z_1|^{\ord(\gminia-\gminib)})
 \Bigr)\cdot |f|_h$ for some $\epsilon>0$.
\item
If $\alpha\neq\beta$,
$\bigl|f'_{(\gminia,\alpha),(\gminia,\beta)}\bigr|_h
=O\bigl(|z_1|^{\epsilon}\bigr)\cdot|f|_h$
 for some $\epsilon>0$.
\hfill\qed
\end{itemize}
\end{lem}

\subsubsection{Step 1}

Let $\theta_1$ denote 
the $dz_1$-component of $\theta$.
\begin{lem}
\label{lem;08.7.26.15}
We have the following estimate with respect to $h$:
\[
 \bigl[\theta_1^{\dagger},\nbigu\bigr]
=O\Bigl(
 \frac{d\zbar_1}{|z_1|\cdot\bigl(-\log|z_1|\bigr)}
 \Bigr)
 \cdot\bigl|\nbigu\bigr|_h
\]
\end{lem}
\pf
In the following,
$\epsilon_i$ denote some positive constants.
We have the decomposition:
\[
 \bigl[
 \theta^{\dagger}_1,\nbigu
 \bigr]
=\sum_{\gminia,\gminib,\gminic}
 \bigl(
 \theta^{\dagger\prime}_{1,\gminia,\gminib}
 \circ \nbigu'_{\gminib,\gminic}
-\nbigu'_{\gminia,\gminib}
 \circ\theta^{\dagger\prime}_{1,\gminib,\gminic}
 \bigr)
\]
By the estimates in Subsection 11.2 of
\cite{mochi7}
(see Subsection \ref{subsubsection;08.7.30.3}),
we have the following estimates for $\gminia\neq\gminib$:
\[
 \theta^{\dagger\prime}_{1,\gminia,\gminib}
=O\Bigl(
 \exp\bigl(-\epsilon_1|z_1|^{\ord(\gminia-\gminib)}\bigr)
 \cdot d\zbar_1
 \Bigr)
\]
Because $\nbigu$ and $\theta$ 
are commutative,
we have the following estimate
for $\gminia\neq\gminib$
due to Lemma \ref{lem;08.7.29.10}:
\[
 \nbigu'_{\gminia,\gminib}
=O\Bigl(\exp\bigl(-\epsilon_2|z_1|^{-1}\bigr)\Bigr)
 \cdot\bigl|\nbigu\bigr|_{h}
\]
Hence, we have the following estimate
with respect to $h$:
\[
 \bigl[\theta_1^{\dagger},\nbigu\bigr]
=\sum_{\gminia}
\bigl[\theta^{\dagger\prime}_{1,\gminia,\gminia},\,
 \nbigu'_{\gminia,\gminia}
\bigr]
+O\Bigl(
 \exp\bigl(-\epsilon_3|z_1|^{-1}\bigr)\cdot d\zbar_1
 \Bigr)\cdot
 \bigl|\nbigu\bigr|_h
\]
Similarly,
we have the following estimates for $\alpha\neq \beta$,
by Theorem 11.12 of \cite{mochi7} and
Lemma \ref{lem;08.7.29.10}:
\[
 \theta^{\dagger\prime}_{1,(\gminia,\alpha),(\gminia,\beta)}
=O\Bigl(
|z_1|^{\epsilon_4}
 \Bigr)\cdot 
 \frac{d\zbar_1}{\zbar_1},
\quad\quad
 \nbigu'_{(\gminia,\alpha),(\gminia,\beta)}
=O\bigl(|z_1|^{\epsilon_4}\bigr)
\]
By Proposition \ref{prop;08.8.3.10},
$\theta^{\dagger\prime}_{(\gminia,\alpha),(\gminia,\alpha)}
-\bigl(d\gminia+\alpha\cdot dz_1/z_1\bigr)\cdot
 \pi'_{\gminia,\alpha}$
is bounded with respect to $h$
and Poincar\'e metric of $X-D$.
Hence, we obtain
\[
 \bigl[\theta_1^{\dagger},\nbigu\bigr]
=\sum_{\gminia,\alpha}
 \bigl[
 \theta^{\dagger\prime}
 _{1,(\gminia,\alpha),(\gminia,\alpha)},\,
 \nbigu'_{(\gminia,\alpha),(\gminia,\alpha)}
 \bigr]
+O\Bigl(|z_1|^{\epsilon_5}\Bigr)
 \cdot\frac{d\zbar_1}{\zbar_1}
\cdot\bigl|\nbigu\bigr|_h
=O\Bigl(
 \frac{d\zbar_1}{|z_1|
 \cdot\bigl(-\log|z_1|\bigr)}
 \Bigr)
\cdot\bigl|\nbigu\bigr|_h.
\]
Thus, we obtain Lemma \ref{lem;08.7.26.15}.
\hfill\qed

\subsubsection{Step 2}
\label{subsubsection;08.7.26.25}

Let $\delbar_1$ denote the $d\zbar_1$-components
of $\delbar_E$ and $\delbar$.
Similarly,
let $\del_{1}$ denote the $dz_1$-component 
of $\del_{E}$ and $\del$.
The following holds:
\[
 \delbar_1\bigl|\nbigu\bigr|^2_h
=\bigl(\nbigu,\del_{1}\nbigu\bigr)_h
=\bigl(\nbigu,[\theta_1,\nbigq]-\theta_1\bigr)_h
=-\tr\Bigl(
 \nbigu\cdot
 \bigl[\theta_1^{\dagger},\nbigq\bigr]
 \Bigr)
-\tr\Bigl(\nbigu\cdot\theta_1^{\dagger}\Bigr)
=-\tr\Bigl(
 \bigl[\nbigu,\theta_1^{\dagger}\bigr]
 \cdot\nbigq
 \Bigr)
-\tr\Bigl(
 \nbigu\cdot\theta_1^{\dagger}
 \Bigr)
\]
Hence, we obtain
\[
 \delbar_1\bigl|\nbigu\bigr|^2_h
=O\Bigl(
 \frac{d\zbar_1}{|z_1|\cdot(-\log|z_1|)}
 \Bigr)
\cdot\bigl|\nbigu\bigr|_h\cdot \bigl|\nbigq\bigr|_h
+O\Bigl(
 \frac{d\zbar_1}{|z_1|^N}
 \Bigr)\cdot
 \bigl|\nbigu\bigr|_h.
\]
We also have
\[
 \delbar_1|\nbigq|^2_h
=-\bigl(
\nbigq,
 [\theta_1,\nbigu^{\dagger}]
 \bigr)_h
+\bigl(
 [\theta_1^{\dagger},\nbigu],
 \nbigq
 \bigr)
=O\Bigl(
 \frac{d\zbar_1}{|z_1|\cdot\bigl(-\log|z_1|\bigr)}
 \Bigr)\cdot
 \bigl|\nbigu\bigr|_h\cdot\bigl|\nbigq\bigr|_h.
\]
Therefore, we obtain
\begin{equation}
 \label{eq;08.7.26.20}
 \delbar_1\Bigl(
 \bigl|\nbigu\bigr|_h^2+\bigl|\nbigq\bigr|_h^2
 \Bigr)
=O\Bigl(
 \frac{d\zbar_1}{|z_1|\bigl(-\log|z_1|\bigr)}
 \Bigr)\cdot
 \bigl|\nbigu\bigr|_h\cdot\bigl|\nbigq\bigr|_h
+O\Bigl(
 \frac{d\zbar_1}{|z_1|^N}
 \Bigr)\cdot\bigl|\nbigu\bigr|_h.
\end{equation}

We set $r:=|z_1|$ and
$F:=\bigl(
 \bigl|\nbigu\bigr|_h^2
+\bigl|\nbigq\bigr|_h^2+1\bigr)^{1/2}$.
We use the polar coordinate
$\bigl(r,\arg(z_1),z_2,\ldots,z_n\bigr)$.
We consider the estimate
on a simply connected region
$Z(\vartheta_0,\vartheta_1):=
 \{\vartheta_0< \arg(z_1)<\vartheta_1\}$
for some fixed $\vartheta_0<\vartheta_1$.
We obtain the following estimate
from (\ref{eq;08.7.26.20}):
\[
 \frac{\del}{\del r}F^2
=G_1\cdot F^2+G_2\cdot F,
\quad\quad
 G_1=O\Bigl(
 \frac{1}{r\cdot(-\log r)}
 \Bigr),
\quad
 G_2=O\bigl(\frac{1}{r^N}\bigr)
\]
We take a solution $H\neq 0$
of the differential equation:
\[
 \frac{\del}{\del r}H=-G_1\cdot H
\]
There exist $C_1>0$  and $M_0>0$
such that
$ C_1^{-1}\cdot(-\log r)^{-M_0}
\leq |H|\leq
 C_1\cdot(-\log r)^{M_0}$.
Since $Z(\vartheta_0,\vartheta_1)$ is simply connected,
we can take $H^{1/2}$.
Then, we have 
\[
 \frac{\del}{\del r}
 \bigl(H\cdot F^2\bigr)
=G_2\cdot H\cdot F
=(G_2\cdot H^{1/2})
\cdot
 (H^{1/2}\cdot F).
\]
Because $G_2\cdot H^{1/2}=O\Bigl( r^{-M_1} \Bigr)$,
we obtain
$H\cdot F^2=O(r^{-M_2})$,
and hence
$F=O\bigl(r^{-M_3}\bigr)$.
Thus, we obtain
the following estimate 
on $Z(\vartheta_0,\vartheta_1)$ for some $M_4>0$:
\begin{equation}
 \label{eq;08.8.19.3}
 \bigl|\nbigu\bigr|_h
=O(r^{-M_4}),
\quad
 \bigl|\nbigq\bigr|_h
=O\bigl(r^{-M_4}\bigr)
\end{equation}
By varying $\theta_0$ and $\theta_1$,
we obtain the estimate (\ref{eq;08.8.19.3})
on $X-D$.
In particular,
we obtain the following estimate
on $X-D$ for $\gminia\neq \gminib$:
\[
 \nbigu'_{\gminia,\gminib}
=O\Bigl(
\exp\bigl(-\epsilon|z_1|^{\ord(\gminia-\gminib)}\bigr)
 \Bigr)
\]

\subsubsection{Step 3}

We have
$\bigl[\theta_1,\nbigu^{\dagger}\bigr]
=\bigl[\theta_1,\nbigutilde^{\dagger}\bigr]
+O\Bigl(
 \exp\bigl(-\epsilon|z_1|^{-1}\bigr)\cdot dz_1
 \Bigr)$
with respect to $h$.
By an argument in the proof of Lemma 
\ref{lem;08.7.26.15},
we obtain the following estimate
with respect to $h$:
\begin{equation}
 \label{eq;08.7.26.50}
 \bigl[\theta_1,\nbigu^{\dagger}\bigr]
=O\Bigl(
 \frac{dz_1}{|z_1|\cdot\bigl(-\log|z_1|\bigr)}
 \Bigr)
 \cdot\bigl|\nbigutilde\bigr|_h
+O\Bigl(
\exp\bigl(-\epsilon|z_1|^{-1}\bigr)\cdot dz_1
 \Bigr)
\end{equation}
According to an estimate in
Subsection 11.5.2 of \cite{mochi7},
we have
\[
 \del_{1}\nbigu
=\del_{1}\nbigutilde
-\sum_{\gminia\in\Irr(\theta)} \del_1\gminia
 \cdot\pi_{\gminia}
+O\Bigl(\exp\bigl(-\epsilon|z_1|^{-1}\bigr)
 \cdot dz_1
 \Bigr).
\]
We set
$\thetatilde:=\theta
-\sum_{\gminia\in\Irr(\theta)}
 d\gminia\cdot\pi_{\gminia}$.
We obtain the following estimates with respect to $h$:
\begin{equation}
 \del_1\nbigutilde-\bigl[\theta_1,\nbigq\bigr]
+\thetatilde_1
=O\Bigl(
 \exp\bigl(-\epsilon|z_1|^{-1}\bigr)
 \Bigr)
\end{equation}
\begin{equation}
 \label{eq;08.7.26.51}
 \del_1\nbigq+\bigl[\theta_1,\nbigutilde^{\dagger}\bigr]
=O\Bigl(\exp\bigl(-\epsilon|z_1|^{-1}\bigr)\Bigr)
\end{equation}
We set
$\Ftilde:=
 \bigl(\bigl|\nbigutilde\bigr|_h^2
+\bigl|\nbigq\bigr|_h^2+1\bigr)^{1/2}$.
As in Step 2,
we consider the estimates on
$Z(\vartheta_0,\vartheta_1)$.
By using an argument in Subsection
\ref{subsubsection;08.7.26.25},
we obtain 
\[
 \frac{\del}{\del r} \Ftilde^2
=\Gtilde_1\cdot\Ftilde^2
+\Gtilde_2\cdot \Ftilde,
\quad\quad
 \Gtilde_1=O\Bigl(
 \frac{1}{r\cdot (-\log r)}
 \Bigr),
\quad
 \Gtilde_2=O\Bigl(
 \frac{1}{r}
 \Bigr).
\]
We take a solution $\Htilde_1\neq 0$ of 
the differential equation:
\[
 \frac{\del}{\del r}\Htilde_1
=-\Gtilde_1\cdot \Htilde_1
\]
Note
$ \log\bigl|\Htilde_1\bigr|
=O\Bigl(
 \log\bigl(-\log r\bigr)
 \Bigr)$.
By choosing $\Htilde_1^{1/2}$,
we obtain
\[
 \frac{\del}{\del r}
 \bigl(\Htilde_1\cdot \Ftilde^2\bigr)
=\bigl(
 \Gtilde_2\cdot \Htilde_1^{1/2}
 \bigr)
 \cdot
 \bigl(\Htilde_1^{1/2}\cdot \Ftilde\bigr).
\]
Because $\Gtilde_2\cdot\Htilde_1^{1/2}
=O\Bigl(
 r^{-1}\cdot\bigl(-\log r\bigr)^{M_5}
 \Bigr)$
for some $M_5>0$,
we obtain 
$ \Htilde_1\cdot\Ftilde^2
=O\Bigl(
 \bigl(-\log r\bigr)^{M_6}
 \Bigr)$
for some $M_6>0$,
and thus
$ \Ftilde=O\Bigl(
 (-\log r)^{M_7}
 \Bigr)$ for some $M_7>0$.
Therefore,
we obtain the following estimates
with respect to $h$:
\begin{equation}
 \label{eq;08.7.26.30}
 \nbigutilde=O\Bigl(
 \bigl(-\log r\bigr)^{M_7}
 \Bigr),
\quad
 \nbigq
=O\Bigl(
 \bigl(-\log r\bigr)^{M_7}
 \Bigr)
\end{equation}

\subsubsection{Step 4}
\label{subsubsection;08.7.26.101}

By (\ref{eq;08.7.26.30}),
$\nbigutilde$ is a holomorphic section of
$\nbigp_0\End(E)$.
Because $[\theta,\nbigutilde]=0$,
we obtain the boundedness of 
$ \bigl|\nbigutilde\bigr|_h$
by an estimate in Subsection 11.7 of \cite{mochi7}.
Thus, the proof of Proposition \ref{prop;08.7.26.55}
is finished.
\hfill\qed

\begin{rem}
From {\rm(\ref{eq;08.7.26.50})}
and {\rm(\ref{eq;08.7.26.51})},
we also have the following estimate:
\[
 \del_1\nbigq=
 O\Bigl(
 \frac{dz_1}{|z_1|\cdot(-\log |z_1|)}
 \Bigr)
\]
Hence, we actually obtain
$\nbigq=O\Bigl(
 \log\bigl(-\log|z_1|\bigr)
 \Bigr)$.
However, we will obtain the boundedness 
later.
\hfill\qed
\end{rem}

\subsubsection{Proof of Proposition 
 \ref{prop;08.7.26.70}}
\label{subsubsection;08.7.26.71}

For an endomorphism $f$ of $E$,
we have the following:
\begin{equation}
 \label{eq;08.8.3.11}
 \bigl|f \bigr|_{\nbigpzero_{\irr}h}
=\bigl|
 g_{\irr}(\lambda-\lambda_0)
 \circ f \circ
 g_{\irr}(\lambda-\lambda_0)^{-1}
 \bigr|_h
\end{equation}
Hence, the claim for $\nbigutilde$ is clear
from $[\nbigutilde,g_{\irr}(\lambda-\lambda_0)]=0$.
We have the decomposition
$\nbigp_0\nbige^0=
\bigoplus\nbigp_0\nbige^0_{\gminia}$
extending $E=\bigoplus E_{\gminia}$.
Let $\vecv=(\vecv_{\gminia})$
be a holomorphic frame of 
$\nbigp_0 \nbige^0$
compatible with the decomposition.
Let $C$ be the matrix-valued function
determined by $\del_1\vecv=\vecv\cdot C\cdot dz_1$.
We have the decomposition 
$C=(C_{\gminia,\gminib})$
corresponding to the decomposition
$\vecv=(\vecv_{\gminia})$.
According to an estimate in Subsection 11.5.2
of \cite{mochi7},
there exists $\epsilon_1>0$
such that the following holds
for $\gminia\neq \gminib$:
\[
 C_{\gminia,\gminib}=
O\Bigl(
 \exp\bigl(-\epsilon_1|z_1|^{\ord(\gminia-\gminib)}\bigr)
 \Bigr)
\]
Let $A$ be the matrix-valued function determined by
$\nbigu\vecv
=\vecv\cdot A$.
Note $A$ is block-diagonal,
i.e.,
$A=\bigoplus A_{\gminia,\gminia}$.
We have
$(\del_1\nbigu)\vecv
=\vecv\cdot 
\bigl(
 \del_1 A+[C,A]\cdot dz_1
\bigr)$.
We set $B\cdot dz_1:=
 \del_1 A+[C,A]\cdot dz_1=
 \bigl(B_{\gminia,\gminib}\cdot dz_1\bigr)$.
Then, there exists $\epsilon_2>0$
such that the following holds for $\gminia\neq\gminib$:
\begin{equation}
\label{eq;08.7.26.105}
 B_{\gminia,\gminib}
=C_{\gminia,\gminib}\cdot A_{\gminib,\gminib}
-A_{\gminia,\gminia}\cdot C_{\gminia,\gminib}
=O\Bigl(
 \exp\bigl(-\epsilon_2|z_1|^{\ord(\gminia-\gminib)}\bigr)
 \Bigr)
\end{equation}

For any section $f$ of $\End(E)\otimes\Omega^{1,0}$,
we have the decomposition
\[
 f=\sum f_{\gminia,\gminib},
\quad
 f_{\gminia,\gminib}\in
 \Hom(E_{\gminib},E_{\gminia})
 \otimes\Omega^{1,0}.
\]
From the relation
$ \del_1\nbigu-[\theta_1,\nbigq]+\theta_1=0$,
we obtain the following:
\[
 (\del_1\nbigu)_{\gminia,\gminib}
-\del_1(\gminia-\gminib)\cdot\nbigu_{\gminia,\gminib}
-(\theta_{1,\gminia}-\del_1\gminia)
 \cdot
 \nbigu_{\gminia,\gminib}
+\nbigu_{\gminia,\gminib}\cdot
 (\theta_{1,\gminib}-\del_1\gminib)=0
\]
Note the following
(see Proposition \ref{prop;08.7.30.2}
and Proposition \ref{prop;08.8.3.10}):
\begin{equation}
\label{eq;08.10.19.10}
\del(\gminia-\gminib)/\del z_1
\sim
 |z_1^{\ord(\gminia-\gminib)-1}|\cdot dz_1,
\quad
 \bigl|
 \theta_{1,\gminia}-\del_1\gminia
 \bigr|_h=O\bigl(dz_1/z_1\bigr),
\quad
 \bigl|
 \theta_{1,\gminib}-\del_1\gminib
 \bigr|_h=O\bigl(dz_1/z_1\bigr)
\end{equation}
The estimate (\ref{eq;08.7.26.105})
implies the following:
\begin{equation}
\label{eq;08.10.19.11}
\bigl|
(\del_1\nbigu)_{\gminia,\gminib}
\bigr|=O\Bigl(
 \exp\bigl(-\epsilon_2|z_1|^{\ord(\gminia-\gminib)}\bigr)
 \Bigr)
\end{equation}
Due to (\ref{eq;08.10.19.10})
and (\ref{eq;08.10.19.11}),
there exists $\epsilon_3>0$
such that the following holds for 
$\gminia\neq\gminib$:
\[
\bigl|
 \nbigq_{\gminia,\gminib}
\bigr|_h
=O\Bigl(
 \exp\bigl(-\epsilon_3|z_1|^{\ord(\gminia-\gminib)}\bigr)
 \Bigr)
\]
By using (\ref{eq;08.8.3.11}),
we obtain the desired estimate for
$\nbigq$ with respect to
$\nbigp^{(\lambda_0)}h$,
if $U(\lambda_0)$ is sufficiently small.
\hfill\qed

\subsubsection{Complement for the normal crossing case}

Let $X:=\Delta^n$
and $D:=\bigcup_{i=1}^{\ell}\{z_i=0\}$.
Let $\harmonicbundle$
be an unramifiedly good wild harmonic bundle
on $X-D$.
Let $\nbigu$ be a holomorphic section
of $\End(E)$ on $X-D$
such that $[\theta,\nbigu]=0$.
Let $\nbigq$ be a $C^{\infty}$-section
of $\End(E)$ on $X-D$
such that $\nbigq^{\dagger}=\nbigq$.
Assume that they satisfy the equations
(\ref{eq;08.7.26.60}) and (\ref{eq;08.7.26.61}).
We also assume that
there exists a holomorphic decomposition
$ (E,\theta)=\bigoplus_{\gminia\in\Irr(\theta)}
 (E_{\gminia},\theta_{\gminia})$
such that 
$\theta_{\gminia}-d\gminia\cdot \pi_{\gminia}$
are tame,
where $\pi_{\gminia}$ denotes the projection
onto $E_{\gminia}$ with respect to the above decomposition.
We set $ \nbigutilde:=\nbigu
+\sum_{\gminia\in\Irr(\theta)}
 \gminia\cdot\pi_{\gminia}$.

\begin{cor}
\label{cor;08.9.11.2}
$\nbigutilde$ is bounded with respect to $h$.
\end{cor}
\pf
It follows from Proposition \ref{prop;08.7.26.55} above
and the estimate in Subsection 11.7 of \cite{mochi7}.
\hfill\qed

\subsection{Prolongation of
 variation of integrable twistor structures}

\subsubsection{Statements}

Let $X$ be a complex manifold,
and let $D$ be a simple normal crossing divisor
of $X$.
Let $(\nbige^{\sankaku},
 \DDtilde^{\sankaku},\nbigs)$
be a variation of pure polarized integrable twistor structures
of weight $0$ on $\proj^1\times(X-D)$.
We have the underlying harmonic bundle
$(E,\delbar_E,\theta,h)$ on $X-D$.

\begin{df}
\mbox{{}}
\begin{itemize}
\item
We say that 
$(\nbige^{\sankaku},\DDtilde^{\sankaku},\nbigs)$
is tame (wild, good wild, unramifiedly good wild),
if $(E,\delbar_E,\theta,h)$ is tame
(wild, good wild, unramifiedly good wild).
\item
If we are given a real structure $\kappa$
of $(\nbige^{\sankaku},\DDtilde^{\sankaku},\nbigs)$,
we say that 
the variation of polarized pure twistor-TERP structures
$(\nbige^{\sankaku},
 \DDtilde^{\sankaku},\nbigs,\kappa,0)$
is tame (wild, good wild, unramifiedly good wild),
if $(\nbige^{\sankaku},\DDtilde^{\sankaku},\nbigs)$
is tame (wild, good wild, unramifiedly good wild).
\end{itemize}
Note that
``wild'' does not imply ``good wild''
as remarked in 
Remark {\rm\ref{rem;08.11.6.11}}.
\hfill\qed
\end{df}

Assume that $(E,\delbar_E,\theta,h)$ is
good wild.
We will show the following proposition later.
(The tame case  was shown in 
 \cite{Hertling-Sevenheck3}.)
\begin{lem}
\label{lem;08.7.26.110}
The sets of 
$\KMS\bigl(\nbigp\nbige^0,i\bigr)$
are contained in
$\real\times\{0\}$.
\end{lem}

We use the notation in 
Subsection \ref{subsubsection;08.7.29.30}.
As explained in Subsection
\ref{subsection;08.8.3.6},
$(\nbige,\DD)$ is prolonged
to the family of meromorphic $\lambda$-flat
bundles $(\nbigq\nbige,\DD)$
on $\cnum_{\lambda}\times (X,D)$,
and 
$(\nbige^{\dagger},\DD^{\dagger})$
is prolonged to the family of meromorphic
$\mu$-flat bundles
$(\nbigq\nbige^{\dagger},\DD^{\dagger})$
on $\cnum_{\mu}\times (X^{\dagger},D^{\dagger})$.

\begin{prop}
\label{prop;08.7.26.111}
\mbox{{}}
\begin{itemize}
\item
$\DDtilde^f$
(resp. $\DDtilde^{\dagger\,f}$)
gives a meromorphic flat connection of
$\nbigq\nbige$
(resp. $\nbigq\nbige^{\dagger}$).
\item
If a real structure $\kappa$ of
$(\nbige^{\sankaku},\DDtilde^{\sankaku},\nbigs)$
is given,
$\kappa_0:
 \gamma^{\ast}\nbige^{\dagger}
\simeq
 \nbige$
is extended to 
the isomorphism
$\gamma^{\ast}\nbigq\nbige^{\dagger}
\simeq
 \nbigq\nbige$.
Similarly,
$\kappa_{\infty}:
 \gamma^{\ast}\nbige\simeq\nbige^{\dagger}$
is extended to
$\gamma^{\ast}\nbigq\nbige
\simeq
 \nbigq\nbige^{\dagger}$.
\end{itemize}
\end{prop}

For the proof of Lemma \ref{lem;08.7.26.110}
and Proposition \ref{prop;08.7.26.111},
we may and will assume
(i) $D$ is smooth, i.e., $\ell=1$,
(ii) $(E,\delbar_E,\theta,h)$ is unramified.

\subsubsection{Meromorphic connection of
$\nbigp^{(\lambda_0)}\nbige$}

Let $\lambda_0\in \cnum_{\lambda}$,
and let $U(\lambda_0)$ be a small
neighbourhood of $\lambda_0$ in $\cnum_{\lambda}$.
We set $\nbigxzero:=U(\lambda_0)\times X$
and $\nbigdzero:=U(\lambda_0)\times D$.
Recall that we have a family of
meromorphic $\lambda$-flat bundles
$(\nbigp^{(\lambda_0)}\nbige,\DD)$
on $(\nbigxzero,\nbigdzero)$,
as explained in Subsection
\ref{subsection;08.8.3.6}.
Note that
$\nbigp^{(\lambda_0)}\nbige$ 
is identified with the sheaf of 
holomorphic sections of $\nbige$
of polynomial order
with respect to $\nbigp^{(\lambda_0)}_{\irr}h$,
because $\nbigpzero_{\irr}h$
and $\nbigpzero h$ are mutually bounded
up to polynomial orders.
(See (\ref{eq;08.11.6.12})
for $\nbigpzero_{\irr}h$.
They are different in general.)

\begin{prop}
\label{prop;08.7.29.30}
$\DDtilde^{f}$
gives a meromorphic flat connection
of $\nbigp^{(\lambda_0)}\nbige$.
\end{prop}
\pf
We have only to show
$\lambda^2\nabla_{\lambda}\bigl(
\del_{\lambda} \bigr)
 \nbigp^{(\lambda_0)}\nbige
\subset
 \nbigp^{(\lambda_0)}\nbige$.
As mentioned in Subsection
\ref{subsubsection;08.7.29.30},
we have the induced holomorphic section
$\nbigu$ of $\End(E)$ on $X-D$
such that $[\theta,\nbigu]=0$,
and the $C^{\infty}$-section $\nbigq$ of $\End(E)$
such that $\nbigq^{\dagger}=\nbigq$,
determined by
\[
 \nabla_{\lambda}
=d_{\lambda}+
 \bigl(
 \lambda^{-1}\nbigu
-\nbigq-\lambda\cdot\nbigu^{\dagger}
 \bigr)\frac{d\lambda}{\lambda},
\]
where $d_{\lambda}$
denote the naturally induced flat connection
of $p_{\lambda}^{-1}E$ along 
the $\lambda$-direction.
They satisfy the equations
(\ref{eq;08.7.26.60}) and (\ref{eq;08.7.26.61}).

Let $\vecv=(\vecv_{\gminia})$ 
be a holomorphic frame of
$\nbigp_0\nbige^0$ compatible with the decomposition
$\nbigp_0\nbige^0=
 \bigoplus_{\gminia}\nbigp_0\nbige^0_{\gminia}$.
Corresponding to the decomposition
$\vecv=(\vecv_{\gminia})$,
the identity matrix is decomposed 
into $\bigoplus_{\gminia\in\Irr(\theta)} I_{\gminia}$.
We regard $\vecv$ as a $C^{\infty}$-frame of
$\nbige_{|\nbigxzero-\nbigdzero}$,
and we set
\[
 \vecvtilde=g_{\irr}(\lambda-\lambda_0)^{-1}\vecv
=\vecv\cdot
\Bigl(
 \bigoplus_{\gminia\in\Irr(\theta)}
 \exp\bigl(-(\lambda-\lambda_0)\cdot\gminiabar\bigr)
 \cdot I_{\gminia}
 \Bigr)
\]
Let 
$H\bigl(\nbigp^{(\lambda_0)}_{\irr}h,\vecvtilde\bigr)$
denote the Hermitian matrix-valued function
whose $(i,j)$-entries are given by
$\nbigp^{(\lambda_0)}_{\irr}h(\vtilde_i,\vtilde_j)$.
Then,
it is clear that
$H\bigl(\nbigp^{(\lambda_0)}_{\irr}h,\vecvtilde\bigr)$
and its inverse
are of polynomial order.
We also have the following relation:
\[
 d_{\lambda}\vecvtilde
=\vecvtilde\cdot A,
\quad\quad
A:=-\bigoplus \gminiabar\cdot d\lambda\cdot I_{\gminia}
\]

Let $\vecw$ be a holomorphic frame of
$\nbigp^{(\lambda_0)}_a\nbige$.
Let 
$H\bigl(\nbigp^{(\lambda_0)}_{\irr}h,\vecw\bigr)$
denote the Hermitian matrix-valued function
whose $(i,j)$-entries are given by
$\nbigp^{(\lambda_0)}_{\irr}h(w_i,w_j)$.
Then, $H\bigl(\nbigp^{(\lambda_0)}_{\irr}h,\vecw\bigr)$
and its inverse are of polynomial order.
(See Subsection 13.1.2 of \cite{mochi7},
 for example.)
Let $G$ be the matrix-valued function
determined by $\vecw=\vecvtilde\cdot G$.
Then,
$G$ and $G^{-1}$
are of polynomial order.
We have
\[
 d_{\lambda}\vecw
=\vecvtilde\cdot 
\Bigl(
 A\cdot G+d_{\lambda}G
\Bigr)
=\vecw\cdot
\Bigl(
 G^{-1}\cdot A\cdot G
+G^{-1}d_{\lambda}G
\Bigr).
\]
Since $\vecvtilde$ and $\vecw$ are $\lambda$-holomorphic,
$G$ is $\lambda$-holomorphic.
Hence, $d_{\lambda}G$ and
$G^{-1}\cdot A\cdot G+G^{-1}d_{\lambda}G$
are of polynomial order

Let $B$ be determined by
$\lambda^2\nabla_{\lambda}(\del_{\lambda})\vecw
=\vecw\cdot B$.
Then, $B$ is of polynomial order,
and hence meromorphic.
Thus, the proof of Proposition \ref{prop;08.7.29.30}
is finished.
\hfill\qed

\vspace{.1in}
We have the irregular decomposition:
\begin{equation}
\label{eq;08.7.26.120}
 (\nbigp^{(\lambda_0)}_a\nbige,\DD)
 _{|\widehat{\nbigd}^{(\lambda_0)}}
=\bigoplus_{\gminia\in\Irr(\theta)}
 (\nbigp^{(\lambda_0)}_a\nbigehat_{\gminia},
 \DDhat_{\gminia})
\end{equation}
\begin{lem}
\mbox{{}} \label{lem;08.7.26.125}
\begin{itemize}
\item
$\lambda^2\nabla_{\lambda}(\del_{\lambda})$
preserves the decomposition
{\rm (\ref{eq;08.7.26.120})}.
\item
Assume $\lambda_0\neq 0$.
Then, {\rm (\ref{eq;08.7.26.120})}
is the irregular decomposition
for 
$\bigl(\nbigp^{(\lambda_0)}\nbige,
 \DDtilde^f
 \bigr)$,
and $\nbigp^{(\lambda_0)}_a\nbige$
is an unramifiedly good lattice
of $\nbigp^{(\lambda_0)}\nbige$.
\end{itemize}
\end{lem}
\pf
Since it can be shown by a standard argument,
we give only an outline.
Let $\vecvhat=(\vecvhat_{\gminia})$
be a frame of
$\nbigp^{(\lambda_0)}_a\nbige_{|\nbigdhat}$
compatible with the decomposition 
(\ref{eq;08.7.26.120}).
Let $A=\sum A_{\gminib,\gminia}$ be determined by
$\lambda^{2}\nabla_{\lambda}(\del_{\lambda})
 \vecvhat=\vecvhat\cdot A$.
For $\gminia\neq\gminib$,
let $F_{\gminib,\gminia}:
 \nbigpzero_a\nbigehat_{\gminia}
 \lrarr
 \nbigpzero_a\nbigehat_{\gminib}$
be given by
$F_{\gminib,\gminia}\vecvhat_{\gminia}
=\vecvhat_{\gminib}\cdot
 A_{\gminib,\gminia}$.
Because
$\bigl[\lambda^2\nabla_{\lambda}(\del_{\lambda}),
 \DD^f\bigr]=0$,
we obtain that $F_{\gminib,\gminia}$ is flat.
However, meromorphic flat section has to be $0$
in the case $\gminib\neq\gminia$.
Thus, we obtain the first claim.

Let us show the second claim.
Let $B_{\gminia}$ be determined by
\[
 \DD^f(z_1\del_1)\vecvhat_{\gminia}
=\vecvhat_{\gminia}
 \cdot\bigl(
 (\lambda^{-1}+\lambdabar_0)\cdot 
  z_1\del_1\gminia
+B_{\gminia}
 \bigr)
\]
Then, $B_{\gminia}$ is regular.
For $\gminia=0$, the following holds:
\[
 \lambda^2\del_{\lambda}B_{0}
+A_{0,0}\cdot B_0-B_{0}\cdot A_{0,0}
-z_1\del_1A_{0,0}=0
\]
We have the expansions
$B_0=\sum_{m\geq 0}
 B_{0;m}\cdot z_1^m$
and 
$A_{0,0}=\sum_{m\geq N}
 A_{0,0;m}\cdot z_1^m$.
We assume 
$N<0$ and $A_{0,0;N}\neq 0$.
We obtain the relation
$\bigl[B_{0;0},A_{0,0;N}\bigr]-NA_{0,0;N}=0$
on $\nbigdzero$.
Note that the eigenvalues of $B_{0;0}$
are of the form
$\lambda^{-1}\eigenmap(\lambda,u)$,
where 
$u\in\KMS(\nbigp\nbige^0)$ and
$a-1<\paramap(\lambda_0,u)\leq a$.
It implies that the difference of
two distinct eigenvalues of $B_{0,0}$
cannot be $N$.
Therefore, we obtain $A_{0,0;N}=0$,
which contradicts with our assumption.
Hence, we obtain $N\geq 0$.

By considering a twist by
a meromorphic flat line bundle
given by 
$\nabla e=e\cdot 
d\bigl((\lambda^{-1}+\lambdabar_0)\gminia\bigr)$,
we obtain that 
$\DDtilde^f=\DD^f+\nabla_{\lambda}$ 
on $\nbigp^{(\lambda_0)}\nbige$ is of the form
\[
 \DDtilde^f
=\bigoplus_{\gminia\in\Irr(\theta)}
 \Bigl(
 d\bigl(
 (\lambda^{-1}+\lambdabar_0)\cdot\gminia
 \bigr)
+\DDtilde^f_{\nbigp^{(\lambda_0)}\nbigehat,\gminia}
 \Bigr),
\]
where
$\DDtilde^f_{\nbigp^{(\lambda_0)}\nbigehat,\gminia}$
are logarithmic with respect to
$\nbigp^{(\lambda_0)}\nbigehat_{\gminia}$.
Thus, the proof of Lemma \ref{lem;08.7.26.125}
is finished.
\hfill\qed

\subsubsection{Proof}

By Lemma \ref{lem;08.7.26.125},
the eigenvalues of
$\Res(\DDtilde^f)$ on 
$\nbigp^{(\lambda_0)}_b\nbige_{|\nbigdzero}$
are constant.
On the other hand,
the eigenvalues of
$\Res(\DD^f)=\Res(\DDtilde^f)$ on 
$\nbigp^{(\lambda_0)}_b\nbige_{|\nbigdzero}$
have to be of the form
$\lambda^{-1}\alpha-a-\lambda\alphabar$
for $(a,\alpha)\in\KMS(\nbigp\nbige^0)$
by Lemma \ref{lem;08.8.3.12}.
Hence, we obtain
$\alpha=0$ for any 
$(a,\alpha)\in\KMS(\nbigp\nbige^0)$,
i.e.,
$\KMS(\nbigp\nbige^0)\subset
\real\times\{0\}$.
Thus, Lemma \ref{lem;08.7.26.110}
is proved.

\vspace{.1in}

Let us show Proposition \ref{prop;08.7.26.111}.
The first claim follows from 
Lemma \ref{lem;08.8.2.5},
Proposition \ref{prop;08.7.29.30}
and the definition of $\nbigq\nbige$
in Subsection \ref{subsection;08.10.4.10}.
To show the second claim,
we remark that $\kappa$ is flat and 
preserves the pluri-harmonic metrics
for $(\nbige^{\sankaku},\DD^{\sankaku},\nbigs)$
and 
$\gamma^{\ast}
 (\nbige^{\sankaku},\DD^{\sankaku},\nbigs)$.
We also remark that we have only to
consider the case in which $D$ is smooth.
We have
$ \Irr(\DD^{\lambda},\nbigq\nbigelambda)=
\Irr(\theta)$
and 
$\Irr\bigl(
 \DD^{\dagger\,\lambdabar},
\nbigq\nbige^{\dagger\,\lambdabar}
 \bigr)=\Irr(\theta^{\dagger})
=\bigl\{\gminiabar\,\big|\,
 \gminia\in\Irr(\theta)\bigr\}$.
Hence, 
we have the natural identification
$\Irr(\DD^{\lambda},\nbigq\nbigelambda)
=\Irr(\gamma^{\ast}\DD^{\dagger\lambdabar},
 \gamma^{\ast}\nbigq\nbige^{\dagger\,\lambdabar})$.
Since the full Stokes filtrations
are characterized by 
growth order of the norms of flat sections
with respect to the pluri-harmonic metrics
(Proposition \ref{prop;08.10.19.1}),
the full Stokes filtrations are 
preserved by $\kappa$.
Thus, the second claim of
Proposition \ref{prop;08.7.26.111}
follows from Lemma \ref{lem;08.9.11.4}.
\hfill\qed

\begin{rem}
Because $\KMS(\nbigp\nbige^0)\subset
 \real\times\{0\}$,
it turns out that
any $\lambda\neq 0$ is generic.
\hfill\qed
\end{rem}

\subsection{Reduction from wild to tame}
\label{subsection;08.8.9.4}

\subsubsection{Construction of the reductions}
\label{subsubsection;08.7.26.200}

Let $X:=\Delta^n$ and
$D:=\bigcup_{i=1}^{\ell}\{z_i=0\}$.
Let $(\nbige^{\sankaku},\DDtilde^{\sankaku},\nbigs)$
be an unramifiedly good wild 
variation of pure polarized integrable twistor structures
of weight $0$ on $\proj^1\times(X-D)$.
We have the underlying harmonic bundle
$(E,\delbar_E,\theta,h)$.
We take an auxiliary sequence
$\nbigm=\bigl(\vecm(0),\vecm(1),\ldots,\vecm(L)\bigr)$
for $\Irr(\theta)$
as in Subsection 3.1.2 of \cite{mochi7}.

\vspace{.1in}
For each $\gminia\in\Irrbar\bigl(\theta,\vecm(0)\bigr)$,
we obtain the variation of
pure polarized twistor structures
$\Gr^{\vecm(0)}_{\gminia}\bigl(
 \nbige^{\sankaku},
 \DD^{\sankaku},
 \nbigs \bigr)$
by taking Gr with respect to Stokes filtrations
in the level $\vecm(0)$,
as explained in Subsection
\ref{subsection;08.8.3.15}.
By Proposition \ref{prop;08.7.26.111}
and Lemma \ref{lem;08.8.3.3},
it is enriched to integrable
$\Gr^{\vecm(0)}_{\gminia}(\nbige^{\sankaku},
 \DDtilde^{\sankaku},\nbigs)$.
If a real structure $\kappa$ of
$(\nbige^{\sankaku},
 \DDtilde^{\sankaku},\nbigs)$ is given,
$\kappa_0$ and $\kappa_{\infty}$
preserve the Stokes filtration
in the level $\vecm(0)$,
which follows from
Proposition \ref{prop;08.7.26.111}
and Lemma \ref{lem;07.12.7.6}.
Hence,
we also have the induced real structure
$\Gr^{\vecm(0)}_{\gminia}(\kappa)$
of $\Gr^{\vecm(0)}_{\gminia}\bigl(
 \nbige^{\sankaku}, 
 \DDtilde^{\sankaku},
 \nbigs\bigr)$,
and  we obtain
a pure polarized variation of 
twistor-TERP structures
$ \Gr^{\vecm(0)}_{\gminia}
 \bigl(
 \nbige^{\sankaku},
 \DDtilde^{\sankaku},
 \nbigs,\kappa,0
 \bigr)$ for each $\gminia\in\Irrbar(\DD,\vecm(0))$.

\vspace{.1in}

Applying the above procedure inductively,
$\Gr^{\vecm(j)}_{\gminia}
 (\nbige^{\sankaku},\DD^{\sankaku},\nbigs)$
are enriched to
integrable 
$\Gr^{\vecm(j)}_{\gminia}
 (\nbige^{\sankaku},\DDtilde^{\sankaku},\nbigs)$
for any $\gminia\in\Irrbar(\theta,\vecm(j))$.
(See the argument in Subsection \ref{subsection;08.8.3.15}.)
If a real structure $\kappa$ is provided,
the reductions are also equipped with
induced real structures,
and we obtain variation of 
twistor-TERP structures
$\Gr^{\vecm(j)}_{\gminia}(
 \nbige^{\sankaku},
 \DDtilde^{\sankaku},
 \nbigs,\kappa,0)$.
In the case $\vecm(L)$,
we use the symbols
$\Gr^{\full}_{\gminia}(\nbige^{\sankaku},
 \DDtilde^{\sankaku},\nbigs)$
and 
$\Gr^{\full}_{\gminia}(\nbige^{\sankaku},
 \DDtilde^{\sankaku},\nbigs,\kappa,0)$.
They are called the full reductions.

For any $\gminia\in\Irr(\theta)$,
we have the harmonic bundles
$L(-\gminia)$ as in Subsection \ref{subsection;08.8.3.15}.
The associated variation of polarized
pure twistor structures is also denoted by 
the same symbol $L(-\gminia)$.
As explained in Subsection \ref{subsubsection;08.8.3.20},
it is naturally enriched to
a variation of pure twistor-TERP structures
of weight $0$.
The underlying harmonic bundle
of $\Gr^{\full}_{\gminia}(\nbige^{\sankaku},
 \DDtilde^{\sankaku},\nbigs)
 \otimes L(-\gminia)$
is tame for each $\gminia\in\Irr(\theta)$.
This procedure is the reduction
``from wild to tame''
in the integrable case.
We have a similar reduction
in the twistor-TERP case.

\subsubsection{Approximating map
 and estimate of the new supersymmetric index}

Let $(\nbige^{\sankaku},\DDtilde^{\sankaku},\nbigs)$
and $\harmonicbundle$ be as above.
Let $\delbar_{\proj^1,\nbige^{\sankaku}}$
denote the $\lambda$-holomorphic structure
of $\nbige^{\sankaku}$.

\paragraph{One step reduction}
By the one step reduction in Subsection 
\ref{subsubsection;08.7.26.200},
we have obtained the unramifiedly good wild variation of
polarized pure integrable twistor structures:
\[
 \bigl(\nbige^{\sankaku}_0,
 \DDtilde^{\sankaku}_0,\nbigs_0\bigr)
:= \bigoplus_{\gminia\in\Irrbar(\theta,\vecm(0))}
 \Gr^{\vecm(0)}_{\gminia}
 \bigl(
 \nbige^{\sankaku},
 \DDtilde^{\sankaku},
 \nbigs
 \bigr)
\]
Let $(E_0,\delbar_{E_0},\theta_0,h_0)$
be the underlying harmonic bundle.
Let $\delbar_{\proj^1,\nbige_0^{\sankaku}}$
denote the $\lambda$-holomorphic structure 
of $\nbige_0^{\sankaku}$.
We fix a hermitian metric $g_{\proj^1}$ of 
$\Omega^{0,1}_{\proj^1}
\oplus 
 \Omega^{1,0}_{\proj^1}
 \bigl(2\{0,\infty\}\bigr)$.
We will prove the following proposition
in Subsection \ref{subsubsection;08.11.2.1}.
\begin{prop}
\label{prop;08.8.4.6}
There exists a $C^{\infty}$-map
$\Phi:\nbige_0^{\sankaku}
\lrarr\nbige^{\sankaku}$
such that the following holds for some $\epsilon>0$
with respect to $h_0$ and $g_{\proj^1}$:
\begin{equation}
 \label{eq;08.8.4.7}
 \Phi^{\ast}\nbigs-\nbigs_0
=O\Bigl(
 \exp\bigl(-\epsilon|\vecz^{\vecm(0)}|\bigr)
 \Bigr),
\quad
 \delbar_{\proj^1,\nbige_0^{\sankaku}}\bigl(
  \Phi^{\ast}\nbigs-\nbigs_0\bigr)
=O\Bigl(
 \exp\bigl(-\epsilon|\vecz^{\vecm(0)}|\bigr)
 \Bigr),
\end{equation}
\begin{equation}
 \label{eq;08.8.4.8}
  \Phi^{\ast}\nabla_{\lambda}-\nabla_{\lambda,0}
=O\Bigl(
 \exp\bigl(-\epsilon|\vecz^{\vecm(0)}|\bigr)
 \Bigr)
\end{equation}
\end{prop}

In fact, the order of the estimates
can be improved as 
$O\bigl(
 \exp(-\epsilon(|\lambda|+|\lambda^{-1}|)
 |\vecz^{\vecm(0)}|)
 \bigr)$.
We give a consequence.
Let $\nbigq_0$ denote the new supersymmetric
index of $(\nbige_0^{\sankaku},
 \DDtilde^{\sankaku}_0,\nbigs_0)$.
\begin{cor}
\label{cor;08.8.19.1}
We have the following estimate
for some $\epsilon>0$ with respect to $h_0$:
\[
 \bigl|
 \Phi^{\ast}h-h_0
 \bigr|_{h_0}
=O\Bigl(\exp\bigl(
 -\epsilon|\vecz^{\vecm(0)}|
 \bigr)\Bigr),
\quad
 \bigl|\Phi^{\ast}\nbigq-\nbigq_0\bigr|_{h_0}
=O\Bigl(
 \exp\bigl(-\epsilon|\vecz^{\vecm(0)}|\bigr)
 \Bigr)
\]
\end{cor}
\pf
It follows from Lemma \ref{lem;08.7.26.11}.
\hfill\qed

\paragraph{Full reduction}

By taking the full reduction in Subsection 
\ref{subsubsection;08.7.26.200},
we have obtained the unramifiedly good wild variation of
polarized pure integrable twistor structures:
\[
 \bigl(\nbige^{\sankaku}_1,
 \DDtilde^{\sankaku}_1,\nbigs_1\bigr)
:= \bigoplus_{\gminia\in\Irrbar(\theta)}
 \Gr^{\full}_{\gminia}
 \bigl(
 \nbige^{\sankaku},
 \DDtilde^{\sankaku},
 \nbigs
 \bigr)
\]
Let $(E_1,\delbar_{E_1},\theta_1,h_1)$ be 
the underlying harmonic bundle,
and let $\nbigq_1$ denote the supersymmetric
index for $(\nbige_1^{\sankaku},\DD^{\sankaku}_1)$.
By applying Proposition \ref{prop;08.8.4.6}
and Corollary \ref{cor;08.8.19.1} inductively
(see Subsection \ref{subsection;08.8.3.15}
 for an inductive use),
we obtain a $C^{\infty}$-map
$\Phi_1:\nbige^{\sankaku}_1\lrarr\nbige^{\sankaku}$
such that the following holds
for some $\epsilon>0$ with respect to $h_1$:
\[
 \bigl|
 \Phi_1^{\ast}h-h_1
\bigr|_{h_1}
=O\Bigl(
 \exp(-\epsilon|\vecz^{\vecm(L)}|)
 \Bigr),
\quad
 \bigl|\Phi^{\ast}_1\nbigq-\nbigq_1
 \bigr|_{h_1}
=O\Bigl(
 \exp(-\epsilon|\vecz^{\vecm(L)}|)
 \Bigr)
\]
Note that 
the new supersymmetric index is unchanged
after taking the tensor product
with $L(-\gminia)$.
(See Subsection \ref{subsubsection;08.8.3.20}.)
Hence, 
the study of asymptotic behaviour of
new supersymmetric index is
reduced to the study in the tame case,
up to decay with exponential orders.

\subsubsection{Construction of
 an approximating map}
\label{subsubsection;08.11.2.1}

We assume that the coordinate
is as in Remark \ref{rem;07.11.15.10}
for the good set $\Irr(\theta)$.
Let $k$ be determined by
$\vecm(0)\in
 \seisuu_{<0}^{k}\times\veczero_{\ell-k}$.
Let $\lambda_0\in\cnum_{\lambda}$.
Let $U(\lambda_0)$ denote a small neighbourhood
of $\lambda_0$.
We set $\nbigxzero:=
 U(\lambda_0)\times X$ and
$\nbigdzero(\leq k):=U(\lambda_0)\times D(\leq k)$.
We also use the symbol
$\nbigdzero_i$ in a similar meaning.
We set $W:=\nbigdzero(\leq k)$ if $\lambda_0\neq 0$,
and $W:=\nbigdzero(\leq k) \cup (\{0\}\times X)$.
Let $\sigma:\cnum_{\lambda}\lrarr\cnum_{\mu}$
be given by $\sigma(\lambda)=-\lambdabar$,
which induces the anti-holomorphic map
$\cnum_{\lambda}\times X
\lrarr \cnum_{\mu}\times X^{\dagger}$.
We set 
 $\nbigx^{\dagger\,(-\lambdabar_0)}:=
 \sigma\bigl(\nbigxzero\bigr)$.

From $\harmonicbundle$,
we obtain the vector bundle
$\nbigpzero_0\nbige$ on $\nbigxzero$
with a meromorphic flat connection
$\DDtilde^f:=\DD^f+\nabla_{\lambda}$.
Similarly we obtain
$\nbigpzero_0\nbige_0$
with $\DDtilde^f_0=\DD_0^f+\nabla_{\lambda,0}$
from $(E_0,\delbar_{E_0},\theta_0,h_0)$.

We also obtain the vector bundle
$\nbigp^{(\mu_0)}_0\nbige^{\dagger}$ with 
the meromorphic flat connection
$\DDtilde^{\dagger\,f}=\DD^{\dagger\,f}+\nabla_{\mu}$
on $\nbigx^{(\mu_0)}$
from $\harmonicbundle$,
and the vector bundle
$\nbigp^{(\mu_0)}_0\nbige_0^{\dagger}$ with 
the meromorphic flat connection
$\DDtilde_0^{\dagger\,f}=
 \DD_0^{\dagger\,f}+\nabla_{\mu,0}$
from $(E_0,\delbar_{E_0},\theta_0,h_0)$.

Let $\DD_{\leq k}$ denote the restriction of $\DD$
to the $(z_1,\ldots,z_k)$-direction.

\paragraph{Preliminary}

Let $S$ be a small multi-sector of
$\nbigxzero-W$.
By Proposition \ref{prop;07.9.30.22},
we take a $\DD_{\leq k}$-flat splitting
\[
\nbigpzero_0\nbige_{|\Sbar}
=\bigoplus_{\gminia\in\Irr(\theta)}
 \nbigpzero_0\nbige_{\gminia,S}
\]
of the Stokes filtration in the level $\vecm(0)$,
such that
the restrictions to $\nbigdzero_j\cap S$ 
$(j=k+1,\ldots,\ell)$
are compatible with $\Res_j(\DD)$
and the filtrations $\lefttop{j}\Fzero$.
If $\lambda_0\neq 0$,
we may assume that it is $\DD^f$-flat
by Proposition \ref{prop;07.9.30.5}.
(Note that the $\DD^f$-flatness implies
 the compatibility with the residues
 and the parabolic filtrations.)
By construction of $\Gr^{\vecm(0)}$,
it induces the isomorphism
$\bigl(
 \nbigpzero_0\nbige_0,
 \DD_{0,\leq k}
 \bigr)_{|\Sbar}
\simeq
 \bigl(\nbigpzero_0\nbige,
 \DD_{\leq k}\bigr)_{|\Sbar}$.
Let $\Phi_S^p$ $(p=0,\ldots,m)$
be such isomorphisms.
Let $a_p$ $(p=0,\ldots,m)$ be
non-negative $C^{\infty}$-functions on $S$
such that (i) $\sum a_p=1$,
(ii) $\del_ia_p$ and $\del_{\lambda}a_p$
are $O\bigl(
 |\lambda|^{-C}\cdot\prod_{i=1}^k |z_i|^{-C}
 \bigr)$
for some $C>0$.
We set 
$\Phi_S:=\sum a_p\cdot \Phi^p_S$.
We also set
$G:=(\Phi^0_S)^{-1}\circ\Phi_S$
and
$G^p:=(\Phi^0_S)^{-1}\circ\Phi^p_S$.

\begin{lem}
We have the following estimates
with respect to $h_0$
for some $\epsilon>0$:
\begin{equation}
\label{eq;08.11.2.10}
G^p-\id
=O\Bigl(
 \exp\bigl(
 -\epsilon|\lambda^{-1}\vecz^{\vecm(0)}|
 \bigr)
 \Bigr)
\end{equation}
\begin{equation}
\label{eq;08.11.2.11}
 (\Phi_S^0)^{-1}\circ \bigl(
\lambda^2\nabla_{\lambda}(\del_{\lambda})
\bigr)\circ\Phi_S^0
-\lambda^2\nabla_{\lambda,0}(\del_{\lambda})
=O\bigl(
 \exp(-\epsilon|\lambda^{-1}\vecz^{\vecm(0)}|)
 \bigr)
\end{equation}
\end{lem}
\pf
Let $\nbigg$ be the left hand side of
(\ref{eq;08.11.2.10}) or (\ref{eq;08.11.2.11}).
It is flat with respect to $\DD_{0\leq k}$,
and strictly decreases the Stokes filtration
in the level $\vecm(0)$.
Moreover,
$\nbigg_{|\nbigdzero_i\cap S}$
preserves the filtrations $\lefttop{i}\Fzero$
and the residues $\Res_j(\DD)$
for $j=k+1,\ldots,\ell$.
Then, we obtain the desired estimate
by using the estimate in Subsection 13.3 
of \cite{mochi7}.
(It is easy to show it directly.)
\hfill\qed

\vspace{.1in}

Hence, we have
$\bigl| G-\id\bigr|_{h_0}
=O\bigl(
 \exp\bigl(
 -\epsilon|\lambda^{-1}\vecz^{\vecm(0)}|
 \bigr)
 \bigr)$.
We set
$\Phi_S^{\ast}\nabla_{\lambda}(\del_{\lambda}):=
 \Phi_S^{-1}\circ
 \bigl(\nabla_{\lambda}(\del_{\lambda})\bigr)\circ
 \Phi_S$.
We use the symbol
$(\Phi_S^0)^{\ast}\nabla_{\lambda}(\del_{\lambda})$
in a similar meaning.
By the previous lemma,
we have the following estimate for some $\epsilon>0$
with respect to $h_0$:
\[
 (\Phi_S^0)^{\ast}\nabla_{\lambda}(\del_{\lambda})
-\nabla_{\lambda,0}(\del_{\lambda})
=O\bigl(
 \exp\bigl(-\epsilon|\lambda^{-1}\vecz^{\vecm(0)}|\bigr)
 \bigr)
\]

\begin{lem}
\label{lem;08.8.4.1}
The following estimate holds for some $\epsilon>0$
with respect to $h_0$:
\[
 \Phi_S^{\ast}
 \nabla_{\lambda}(\del_{\lambda})
-\nabla_{\lambda,0}(\del_{\lambda})
=O\Bigl(
 \exp\bigl(
 -\epsilon|\lambda^{-1}\vecz^{\vecm(0)}|
 \bigr) \Bigr)
\]
\end{lem}
\pf
We have the following equalities:
\begin{multline}
 \Phi_S^{\ast}
 \nabla_{\lambda}(\del_{\lambda})
-\nabla_{\lambda,0}(\del_{\lambda})
=\bigl(\Phi_S^{-1}\circ\Phi^0_S\bigr)
 \circ
 (\Phi_S^0)^{\ast}
\nabla_{\lambda}(\del_{\lambda})
 \circ\bigl((\Phi_S^0)^{-1}\circ\Phi_S\bigr)
-\nabla_{\lambda,0}(\del_{\lambda})
 \\
=G^{-1}\circ
 \bigl(
 (\Phi_S^0)^{\ast}\nabla_{\lambda}(\del_{\lambda})
-\nabla_{\lambda,0}(\del_{\lambda})
\bigr)
  \circ G
+G^{-1}\circ\nabla_{\lambda,0}(\del_{\lambda})
 \circ G
-\nabla_{\lambda,0}(\del_{\lambda})\\
=
G^{-1}\circ
 \bigl(
 (\Phi_S^0)^{\ast}\nabla_{\lambda}(\del_{\lambda})
-\nabla_{\lambda,0}(\del_{\lambda})
\bigr)
  \circ G
+G^{-1}\cdot\bigl(
 \nabla_{\lambda,0}(\del_{\lambda})
 G \bigr)
\end{multline}
We have the following:
\[
 \nabla_{\lambda,0}(\del_{\lambda})G
=\sum
 \frac{\del a_p}{\del \lambda}
 \cdot G^p
=\sum
 \frac{\del a_p}{\del \lambda}
 \cdot (G^p-\id)
=O\Bigl(
 \exp\bigl(-\epsilon
 |\lambda^{-1}\vecz^{\vecm(0)}|\bigr)
 \Bigr)
\]
Thus, we obtain Lemma 
\ref{lem;08.8.4.1}.
\hfill\qed

\vspace{.1in}

Assume we are also given 
morphisms on sectors $\sigma(S)$
of $\nbigx^{\dagger(-\lambdabar_0)}-W^{\dagger}$
\[
 \Phi^{\dagger\,q}_{\sigma(S)}: 
 \bigl(
 \nbigp^{(-\lambdabar_0)}
 \nbige^{\dagger}_{0},
 \DD^{\dagger}_0\bigr)
 _{|\sigma(\Sbar)}
\lrarr
 \bigl(
 \nbigp^{(-\lambdabar_0)}
 \nbige^{\dagger},
 \DD^{\dagger}\bigr)_{|\sigma(\Sbar)},
\quad
 (q=0,\ldots,m'),
\]
induced by $\DD_{\leq k}^{\dagger}$-flat
of the Stokes filtration in the level $\vecm(0)$
such that the restriction to
$\sigma(S)\cap \nbigd^{\dagger(-\lambdabar_0)}_j$
$(j=k+1,\ldots,\ell)$
are compatible with the residue
$\Res_j(\DD^{\dagger})$
and the filtration $\lefttop{j}F^{(-\lambdabar_0)}$.
If $\lambda_0\neq 0$,
we may assume that the splittings are
$\DD^{\dagger}$-flat.
Let $b_q$ $(q=0,\ldots,m')$ be non-negative
$C^{\infty}$-functions on $\sigma(S)$
satisfying similar conditions for $a_p$.
We set
$\Phi_{\sigma(S)}^{\dagger}:=
 \sum b_q\cdot \Phi_{\sigma(S)}^{\dagger\,q}$.

\begin{lem}
\label{lem;08.8.4.5}
We set 
$H:=
 \nbigs\circ
 \bigl(\Phi_S\otimes
 \sigma^{\ast}\Phi^{\dagger}_{\sigma(S)}\bigr)
-\nbigs_0$.
Then, we have the following estimate
with respect to $h_0$
for some $\epsilon>0$:
\[
 H=O\Bigl(
 \exp\bigl(
 -\epsilon|\lambda^{-1}\vecz^{\vecm(0)}|
 \bigr)\Bigr),
\quad
 \delbar_{\nbige^{\sankaku}_0,\proj^1}H=
O\Bigl(
 \exp\bigl(
 -\epsilon|\lambda^{-1}\vecz^{\vecm(0)}|
 \bigr)\Bigr)
\]
\end{lem}
\pf
We set 
$H_{p,q}:=
 \nbigs\circ
 \bigl(\Phi_S^p\otimes
 \sigma^{\ast}\Phi^{\dagger\,q}_{\sigma(S)}\bigr)
-\nbigs_0$.
According to an estimate
in Subsection 15.3.2 of \cite{mochi7},
we have
\[
H_{p,q}=O\Bigl(
 \exp\bigl(-\epsilon|\lambda^{-1}\vecz^{\vecm(0)}|\bigr)
 \Bigr) 
\]
with respect to $h_0$ for some $\epsilon>0$.
We also have
$\delbar_{\nbige_0^{\sankaku},\proj^1}
 H_{p,q}=0$.
Then, the claim of Lemma \ref{lem;08.8.4.5}
follows.
\hfill\qed

\paragraph{Construction}

We take a compact region $\nbigk$
of $\cnum_{\lambda}$
such that the union of
the interior parts of $\nbigk$
and $\sigma(\nbigk)$ cover
$\proj^1$.
We take a covering of
\[
\bigl(\nbigk\times X\bigr)
-\Bigl(
 \bigl( \nbigk\times D(\leq k) \bigr)
\cup
 \bigl(\{0\}\times X\bigr)
\Bigr)
\]
by multi-sectors $S_i$ $(i=1,\ldots,N)$
such that $S_i$ are sufficiently small
as in {\bf Preliminary} above.
Then, 
we have 
$\proj^1=\bigcup S_i\cup\bigcup \sigma(S_i)$.
We take a partition 
$\bigl( 
\chi_{S_i},\chi_{\sigma(S_i)}\,\big|\,
 i=1,\ldots,N
 \bigr)$
of unity on $\proj^1$ 
subordinated to the covering.
We assume that
$\del_j \chi_{S_i}$ and
$\del_{\lambda}\chi_{S_i}$
are $O\bigl(
 |\lambda|^{-C}\cdot\prod_{i=1}^k |z_i|^{-C}
 \bigr)$
for some $C>0$.
We assume similar conditions for
$\del_j \chi_{\sigma(S_i)}$ and
$\del_{\mu}\chi_{\sigma(S_i)}$.

For each $S_i\subset\nbigxzero-W$,
we take isomorphisms:
\[
 \Phi_{S_i}:
 \bigl(
 \nbigpzero_0\nbige_0,
 \DD_{0}
 \bigr)_{|\Sbar_i}
\simeq
 \bigl(\nbigpzero_0\nbige,
 \DD\bigr)_{|\Sbar_i},
\quad
 \Phi^{\dagger}_{\sigma(S_i)}: 
 \bigl(
 \nbigp_0^{(-\lambdabar_0)}
 \nbige^{\dagger}_{0},
 \DD^{\dagger}_0\bigr)
 _{|\sigma(\Sbar_i)}
\lrarr
 \bigl(
 \nbigp_0^{(-\lambdabar_0)}
 \nbige^{\dagger},
 \DD^{\dagger}\bigr)
 _{|\sigma(\Sbar_i)}
\]
induced by $\DD_{\leq k}$-flat or
$\DD^{\dagger}_{\leq k}$-flat splittings
of Stokes filtrations as above.
If $\lambda_0\neq 0$,
we assume that $\DD^f$-flatness
and $\DD^{\dagger\,f}$-flatness.
We set
\[
  \Phi:=\sum_{i=1}^N
 \chi_{S_i}\cdot\Phi_{S_i}
+\sum_{i=1}^N
 \chi_{\sigma(S_i)}\cdot\Phi_{\sigma(S_i)}.
\]
It is easy to check
that $\Phi$ satisfies the desired estimates
(\ref{eq;08.8.4.7})
and (\ref{eq;08.8.4.8}),
by using Lemma \ref{lem;08.8.4.1}
and Lemma \ref{lem;08.8.4.5}.
Note that a $\DD$-flat splitting of
the Stokes filtration of $\nbigpzero\nbige_{|\Sbar}$
in the level $\vecm(0)$
naturally gives a $\DD^{\dagger}$-flat splitting
of the Stokes filtration of 
$\nbigp^{(\lambda_0^{-1})}
 \nbige^{\dagger}_{|\Sbar'}$
in the level $\vecm(0)$,
where $S'$ is the multi-sector of
$\nbigx^{\dagger(\lambda_0^{-1})}
 \setminus  W^{\dagger}$,
which follows from the characterization
of the Stokes filtrations by
the growth order of the norms
of flat sections.
Thus, we obtain 
Proposition \ref{prop;08.8.4.6}.
\hfill\qed

\subsection{Reduction from tame to twistor nilpotent orbit}
\label{subsection;08.8.9.3}

\subsubsection{Reduction}

Let $X:=\Delta^n$,
$D_i:=\{z_i=0\}$,
$D:=\bigcup_{i=1}^{\ell}D_i$
and $D_{\ellsitabar}=\bigcap_{i=1}^{\ell}D_i$.
Let $(\nbige^{\sankaku},
 \DDtilde^{\sankaku},\nbigs)$
be a tame variation of 
pure polarized integrable twistor structures
of weight $0$ on $\proj^1\times(X-D)$.
We have the underlying harmonic bundle
$(E,\delbar_E,\theta,h)$.
As explained in Subsection
\ref{subsection;08.8.4.10},
we have the limiting polarized 
mixed twistor structure
$\bigl(
 S^{\can}_{\veca,0}(E),
 \vecN,\nbigs_{\veca,0}
 \bigr)$
associated to $(E,\delbar_E,\theta,h)$.
We also have the variation of
polarized mixed twistor structures 
$\bigl(
 \lefttop{\ellsitabar}
 \nbige^{\sankaku}_{\veca,0},
 \vecN^{\sankaku},
 \DD^{\sankaku}_{\veca,0},
 \nbigs_{\veca,0} \bigr)$
of weight $0$ in $\ell$-variables.
Hertling and Sevenheck observed
the following (see \cite{Hertling-Sevenheck3}).
\begin{prop}
$\bigl(
 S^{\can}_{\veca,0}(E),\vecN,\nbigs_{\veca,0}
 \bigr)$ is naturally enriched to 
a polarized mixed integrable twistor structure
$\bigl(
 S^{\can}_{\veca,0}(E),\nabla,
 \vecN,\nbigs_{\veca,0}
 \bigr)$.
Similarly,
$\bigl(
 \lefttop{\ellsitabar}
 \nbige^{\sankaku}_{\veca,0},
 \vecN^{\sankaku},
 \DD^{\sankaku}_{\veca,0},
 \nbigs_{\veca,0} \bigr)$
is naturally enriched to
a variation of polarized mixed 
integrable twistor structures
$\bigl(
 \lefttop{\ellsitabar}
 \nbige^{\sankaku}_{\veca,0},
 \vecN^{\sankaku},
 \DDtilde^{\sankaku}_{\veca,0},
 \nbigs_{\veca,0} \bigr)$.

If $(\nbige^{\sankaku},
 \DDtilde^{\sankaku},\nbigs)$
has a real structure $\kappa$,
they are also equipped with 
induced real structures.
\hfill\qed
\end{prop}

\subsubsection{Approximating maps}

For $0<R<1$,
we set $X^{\ast}(R):=
\bigl\{
 (z_1,\ldots,z_n)\,\big|\,
 0<|z_i|<R,\,\,i=1,\ldots,n
 \bigr\}$
and
$D_{\ellsitabar}:=\bigl\{
 (z_{\ell+1},\ldots,z_n)\,\big|\,
 |z_i|<R
 \bigr\}$.
By the natural projection
$X^{\ast}(R)\lrarr D_{\ellsitabar}(R)$,
we regard $X^{\ast}(R)$ as
$D_{\ellsitabar}(R)\times
 \bigl\{(z_1,\ldots,z_{\ell})\,\big|\,
 0<|z_i|<R\bigr\}$.
Due to Theorem \ref{thm;08.7.26.8},
we have the integrable twistor nilpotent orbit 
$\TNIL(\nbige^{\sankaku}_{\veca,0},
 \DDtilde^{\sankaku}_{\veca,0}, 
 \vecN,\nbigs_{\veca,0})$
on $X^{\ast}(R)$ for some $R$.
Thus, we obtain 
a tame variation of pure polarized integrable
twistor structures:
\[
 \bigl(\nbige_0^{\sankaku},
 \DDtilde^{\sankaku}_0,
 \nbigs_0\bigr)
:=\bigoplus
 _{\veca\in\Par(\nbigp_0\nbige^0,\ellsitabar)}
 \TNIL(\nbige^{\sankaku}_{\veca,0},
 \DDtilde^{\sankaku}_{\veca,0},
 \vecN,\nbigs_{\veca,0})
 \otimes L(\veca)
\]
(See Subsection 
 \ref{subsubsection;08.8.10.30}
 for $L(\veca)$.)
We have the underlying
tame harmonic bundle
\[
 \bigl(E_0,\delbar_{E_0},\theta_0,h_0\bigr)
=\bigoplus
 \bigl(E_{\veca},\delbar_{\veca},
 \theta_{\veca},h_{\veca}\bigr).
\]
We would like to explain 
that we can approximate
the original 
$(\nbige^{\sankaku},\DDtilde^{\sankaku}, \nbigs)$
with 
$(\nbige^{\sankaku}_0,\DDtilde^{\sankaku}_0,
 \nbigs_0)$.

Let $\delbar_{\proj^1,\nbige_0^{\sankaku}}$
denote the $\lambda$-holomorphic structure 
of $\nbige_0^{\sankaku}$.
We fix a hermitian metric $g_{\proj^1}$ of 
$\Omega^{0,1}_{\proj^1}
\oplus 
 \Omega^{1,0}_{\proj^1}
 \bigl(2\{0,\infty\}\bigr)$.
For a permutation $\sigma$
of $\{1,\ldots,\ell\}$
and for $C>0$,
we set 
\[
 Z(\sigma,C):=\bigl\{
 (z_1,\ldots,z_n)
 \in X^{\ast}(R)
 \,\big|\,
 |z_{\sigma(i-1)}|^C<|z_{\sigma(i)}|,\,\,
 i=1,\ldots,\ell-1
 \bigr\}
\]
If we take a sufficiently large $C>0$,
we have 
$X^{\ast}(R)=\bigcup_{\sigma}Z(\sigma,C)$.
For any $\epsilon>0$,
we set 
$\Lambda_0(\epsilon):=
 \sum_{i=1}^{\ell}|z_i|^{\epsilon}$.
We will prove the following proposition
in Subsection \ref{subsubsection;08.8.5.31}.

\begin{prop}
\label{prop;08.8.5.2}
There exists a $C^{\infty}$-map
$\Phi_{\sigma}:\nbige_0^{\sankaku}
\lrarr\nbige^{\sankaku}$
such that the following estimate holds
for some $\epsilon>0$
with respect to $h_0$ and $g_{\proj^1}$
on $\proj^1\times Z(\sigma,C)$:
\begin{equation}
\label{eq;08.8.5.30}
 \Phi_{\sigma}^{\ast}\nbigs-\nbigs_0
=O\bigl(\Lambda_0(\epsilon)
 \bigr),
\quad
 \delbar_{\proj^1,\nbige_0^{\sankaku}}
\bigl(
 \Phi_{\sigma}^{\ast}\nbigs-\nbigs_0
\bigr)
=O\bigl(
 \Lambda_0(\epsilon)
 \bigr),
\quad
  \Phi_{\sigma}^{\ast}\nabla_{\lambda}
 -\nabla_{\lambda,0}
=O\bigl(
 \Lambda_0(\epsilon)
 \bigr)
\end{equation}
\end{prop}

Before going into the proof,
we give a consequence.
Let $\nbigq_0$ and $\nbigq$
denote the new supersymmetric
indices of $(\nbige_0^{\sankaku},
 \DDtilde^{\sankaku}_0)$
and $(\nbige^{\sankaku},\DDtilde^{\sankaku})$.
By using Lemma \ref{lem;08.7.26.11},
we obtain the following estimates on $Z(\sigma,C)$
for some $\epsilon>0$ with respect to $h_0$:
\begin{equation}
\label{eq;08.8.5.1}
 \bigl|
 \Phi_{\sigma}^{\ast}h-h_0
 \bigr|_{h_0}
=O\bigl(
 \Lambda_0(\epsilon)
 \bigr),
\quad
 \bigl|
 \Phi_{\sigma}^{\ast}\nbigq
-\nbigq_0\bigr|_{h_0}
=O\bigl(
 \Lambda_0(\epsilon)
 \bigr)
\end{equation}

\begin{cor}
\label{cor;08.8.19.2}
The eigenvalues of
$\nbigq$ and $\nbigq_{0}$
are the same up to
$O\bigl(\Lambda_0(\epsilon)\bigr)$
for some $\epsilon>0$.
\end{cor}
\pf
By using (\ref{eq;08.8.5.1}),
we obtain the estimate on 
$Z(\sigma,C)$.
Because $X^{\ast}(R)=\bigcup Z(\sigma,C)$,
the claim of the corollary follows.
\hfill\qed

\vspace{.1in}

We also give a more rough but global estimate,
for which the proof is much simpler.
For $M>0$ and $\epsilon>0$,
we set
\[
 \Lambda(M,\epsilon):=
 \prod_{i=1}^{\ell}(-\log|z_i|)^M
 \sum_{i=1}^{\ell}|z_i|^{\epsilon}. 
\]
\begin{prop}
\label{prop;08.8.4.10}
There exists a $C^{\infty}$-map
$\Phi:\nbige_0^{\sankaku}
\lrarr\nbige^{\sankaku}$
such that the following holds for
some $\epsilon>0$ and $M>0$
with respect to $h_0$ and $g_{\proj^1}$:
\begin{equation}
 \label{eq;08.8.4.21}
 \Phi^{\ast}\nbigs-\nbigs_0
=O\bigl(
\Lambda(M,\epsilon)
 \bigr),
\quad
 \delbar_{\proj^1,\nbige_0^{\sankaku}}
\bigl(
 \Phi^{\ast}\nbigs-\nbigs_0
\bigr)
=O\bigl(
 \Lambda(M,\epsilon)
 \bigr),
\quad
  \Phi^{\ast}\nabla_{\lambda}-\nabla_{\lambda,0}
=O\bigl(
 \Lambda(M,\epsilon)
 \bigr)
\end{equation}
\end{prop}

Note that $\Phi^{\ast}h$
and $h_0$ are mutually bounded
up to log order,
which follows from the weak norm estimate
for acceptable bundles.
(See Lemma \ref{lem;08.11.3.1} below.)
Hence, we obtain the following estimate
for some $M'>0$ and $\epsilon'>0$
by using Lemma \ref{lem;08.10.27.15}:
\[
\bigl|\Phi^{\ast}\nbigq-\nbigq_0\bigr|_{h_0}
=O\bigl(
 \Lambda(M',\epsilon')
 \bigr)
\]

In the one dimensional case,
the estimates in the two propositions
are not so different.
We also remark that
$\Phi_{\sigma}$ in Proposition \ref{prop;08.8.5.2}
also satisfies the estimates (\ref{eq;08.8.4.21}).

\subsubsection{Proof of Proposition
 \ref{prop;08.8.5.2}}
\label{subsubsection;08.8.5.31}

For the proof of Proposition
\ref{prop;08.8.5.2},
we have only to consider the case
that $\sigma$ is the identity.
We use the symbol $Z(C)$
instead of $Z(\id,C)$.
Instead of considering $X^{\ast}(R)$,
we will shrink $X$ around the origin.

\paragraph{Decomposition}
For any subset $I\subset\ellsitabar$,
let $m(I)$ be determined by the condition
$ m(I):=\min\bigl\{
 m\in I\,|\,m+1\not\in I
 \bigr\}$,
in other words,
$\{1,\ldots,m(I)\}\subset I$
but $m(I)+1\not\in I$.
Let 
$q_I:\Par\bigl(\nbigp_0\nbige^0,\ellsitabar\bigr)
\lrarr \Par(\nbigp_0\nbige^0,I)$
and $r_{m(I)}:\seisuu^{\ell}\lrarr \seisuu^{m(I)}$
be the natural projections.
Let $\lambda_0\in\cnum_{\lambda}$.
Let $\nbigk$ denote a small neighbourhood of 
$\lambda_0$ in $\cnum_{\lambda}$.
We set $\nbigx:=\nbigk\times X$.
We use the symbols
$\nbigd_i$, $\nbigd_I$, $\nbigd$, etc.,
in similar meanings.

We have the induced filtrations
$\lefttop{i}F$ $(i\in I)$
of $\nbigq_0\nbige_{|\nbigd_I}$.
For any $i\in I$,
we have the residue endomorphisms
$\Res_i(\DD)$
on $\lefttop{I}\Gr_{\vecb}(\nbigq_0\nbige_{|\nbigd_I})$,
which have the unique eigenvalues $-b_i\cdot \lambda$.
Hence, the nilpotent part $\nbign_i$ is well defined.
For $i\leq m(I)$,
we set 
$\nbign(\ibar):=
 \sum_{j\leq i}\nbign_j$.
Recall that the conjugacy classes of
$\nbign(\ibar)_{|(\lambda,P)}$ 
are independent of 
$(\lambda,P)\in\nbigd_I$
(Lemma 12.47 of \cite{mochi2}).
By considering the weight filtration of $\nbign(\ibar)$,
we obtain the filtration
$W(\ibar)$ of 
$\lefttop{I}\Gr_{\vecb}(\nbigq_0\nbige_{|\nbigd_I})$
indexed by $\seisuu$
in the category of vector bundles on $\nbigd_I$.

\begin{lem}
\label{lem;08.8.5.3}
We have a decomposition
\begin{equation}
 \label{eq;08.8.5.5}
 \nbigq_0\nbige_{|\nbigx}
=\bigoplus
 _{\substack{ 
 \veca\in\Par(\nbigp_0\nbige^0,\ellsitabar)\\
 \veck\in\seisuu^{\ell} }}
 U_{\veca,\veck}
\end{equation}
with the following property:
\begin{itemize}
\item
 For any subset $I\subset\ellsitabar$,
 $\vecb\in \Par(\nbigp_0\nbige^0,I)$ and
 $\vech\in\seisuu^{m(I)}$,
 we put 
\[
 \lefttop{I}U_{\vecb,\vech}
 =\bigoplus_{\substack{
 \veca\in q_I^{-1}(\vecb)\\
 \veck\in r_{m(I)}^{-1}(\vech) }}
 U_{\veca,\veck}
\quad
 \mbox{\rm and }\,\,\,
 \lefttop{I}U_{\vecb}
=\bigoplus_{\vech\in\seisuu^{m(I)}}
 \lefttop{I}U_{\vecb,\vech}
\]
Then, the following holds
for any $\vecc\in\real^{I}$:
\begin{equation}
 \label{eq;08.8.5.2}
 \bigoplus_{\vecb\leq\vecc}
 \lefttop{I}U_{\vecb|\nbigd_I}
=\bigcap_{i\in I}
 \lefttop{i}F_{c_i}\bigl(
 \nbigq_0\nbige_{|\nbigd_I}
 \bigr)
\end{equation}
Moreover, 
the following holds 
for any $\vecn\in\seisuu^{m(I)}$
under the identification
$\lefttop{I}U_{\vecb|\nbigd_I}
\simeq
 \lefttop{I}\Gr_{\vecb}(\nbigq_0\nbige)$
induced by {\rm(\ref{eq;08.8.5.2})}:
\[
 \bigoplus_{\vech\leq\vecn}
 \lefttop{I}U_{\vecb,\vech|\nbigd_I}
=\bigcap_{1\leq i\leq m(I)}
 W_{n_i}(\ibar)
 \bigl(
 \lefttop{I}\Gr_{\vecb}(\nbigq_0\nbige_{|\nbigd_I})
 \bigr)
\]
\end{itemize}
\end{lem}
\pf
Although this is essentially
Corollary 4.47 of \cite{mochi2},
we recall an outline for later use.
The theorems and the definitions 
referred in this proof
are given in \cite{mochi2}.
By Theorem 12.43,
the tuple $\bigl(\lefttop{i}F,\nbign(\jbar)\,\big|\,
 i\in\ellsitabar, j\in\ellsitabar \bigr)$
is sequentially compatible in the sense of
Definition 4.43.
Hence, 
$\bigl(
 \lefttop{i}F,W(\jbar)\,\big|\,
 i\in\ellsitabar,j\in\ellsitabar
 \bigr)$ is compatible
in the sense of Definition 4.39,
as remarked in Lemma 4.44.
By Proposition 4.41,
there exists a splitting of
$\bigl(
 \lefttop{i}F,W(\jbar)\,\big|\,
 i\in\ellsitabar,j\in\ellsitabar
 \bigr)$
in the sense of Definition 4.40.
By applying Lemma 2.16,
we can take a frame compatible
with splittings.
It is easy to take a decomposition
as in the claim of Lemma \ref{lem;08.8.5.3},
by using such a compatible frame.
\hfill\qed

\vspace{.1in}

Let $(\nbigq\nbige_{\veca},\DD)$ be
the prolongment of
$(E_{\veca},\delbar_{\veca},\theta_{\veca},h_{\veca})$.
Similarly, we have a decomposition 
\begin{equation}
 \label{eq;08.8.5.6}
 \nbigq_0\nbige_{\veca|\nbigx}
=\bigoplus_{\veck\in\seisuu^{\ell}}
 U_{0,\veca,\veck}
\end{equation}
satisfying a similar condition.
By our construction of
$(\nbige_0^{\sankaku},
 \DDtilde^{\sankaku}_0,\nbigs_0)$,
we are given the isomorphism
for each $\veca\in\Par(\nbigp_0\nbige^0,\ellsitabar)$:
\[
 \nu_{\veca}:
 \lefttop{\ellsitabar}\Gr_{\veca}(\nbigq_0\nbige)
\simeq
 \nbigq_0\nbige_{\veca|\nbigd_{\ellsitabar}}.
\]
\begin{lem}
\label{lem;08.8.5.4}
We may assume that 
$\nu_{\veca}$ preserves the decompositions
$\bigoplus_{\veck} U_{0,\veca,\veck|\nbigd_{\ellsitabar}}$
and 
$\bigoplus_{\veck} U_{\veca,\veck|\nbigd_{\ellsitabar}}$.
\end{lem}
\pf
In Proposition 4.41 of \cite{mochi2},
the construction of a splitting is given
in a descending inductive way,
and we can take any splitting 
of $\lefttop{\ellsitabar}\Gr_{\veca}(\nbigq\nbige)$
of the filtrations $W(\jbar)$ $(j=1,\ldots,\ell)$
in the beginning.
Thus, we obtain Lemma \ref{lem;08.8.5.4}.
\hfill\qed

\vspace{.1in}
Let $\nu_{\veca,\veck}$
denote the induced map
$U_{0,\veca,\veck|\nbigd_{\ellsitabar}}
 \simeq
U_{\veca,\veck|\nbigd_{\ellsitabar}}$.

\paragraph{Norm estimate}

We recall the norm estimate
for tame harmonic bundles.
We take a $C^{\infty}$-frame
$h'_{\veca,\veck}$ of
$U_{\veca,\veck}$ in (\ref{eq;08.8.5.5}).
We set 
\[
 h^{(1)}_{\veca,\veck}:=
 h'_{\veca,\veck}\cdot
 \prod_{j=1}^{\ell}
 |z_j|^{-2a_j}
 \bigl(-\log|z_j|\bigr)^{k_j-k_{j-1}}
=h'_{\veca,\veck}\cdot
 \prod_{j=1}^{\ell}|z_j|^{-2a_j}
 \prod_{j=1}^{\ell-1}
 \left(
 \frac{-\log|z_j|}{-\log|z_{j+1}|}
 \right)^{k_j}
\cdot
 \bigl(-\log|z_{\ell}|\bigr)^{k_{\ell}}
\]
(We formally set $k_0:=0$.)
We obtain a $C^{\infty}$-hermitian metric
$h^{(1)}=\bigoplus 
 h^{(1)}_{\veca,\veck}$ of
$\nbigq\nbige_{|\nbigx-\nbigd}$.
Theorem 13.25 of \cite{mochi2}
implies the following lemma.
\begin{lem}
\label{lem;08.8.5.10}
$h$ and $h^{(1)}$ are mutually bounded
on $\nbigk\times Z(C)$.
\hfill\qed
\end{lem}

\paragraph{An estimate}

\begin{lem}
\label{lem;08.8.5.20}
Let $f$ be a holomorphic endomorphism
of $\nbigq_0\nbige_0$ satisfying the following
conditions:
\begin{itemize}
\item
 It preserves the filtrations
 $\lefttop{i}F$ $(i=1,\ldots,\ell)$.
\item
For each $\vecb\in \real^I$, 
 the induced endomorphism 
 $\lefttop{I}\Gr^F_{\vecb}(f)$ of 
 $\bigoplus_{q_I(\veca)=\vecb}
 \nbigq_0\nbige_{\veca|\nbigd_I}$
preserves the weight filtrations
$W(\jbar)$ $(j=1,\ldots,m(I))$.
\item
 For each $\veca\in\real^{\ell}$,
 the induced endomorphism 
 $\lefttop{\ellsitabar}\Gr^F_{\veca}(f)$ of
 $\nbigq_0\nbige_{\veca|\nbigd_{\ellsitabar}}$
 is $0$.
\end{itemize}
Then, we have 
$|f|_{h_0}=O\bigl(
 \Lambda_0(\epsilon)
 \bigr)$
for some $\epsilon>0$
on $\nbigk\times Z(C)$.
\end{lem}
\pf
We take decompositions (\ref{eq;08.8.5.6}).
Applying Lemma \ref{lem;08.8.5.10}
to $(E_{\veca},\delbar_{\veca},\theta_{\veca},h_{\veca})$
with the decomposition (\ref{eq;08.8.5.6}),
we take a $C^{\infty}$-hermitian metric
$h^{(1)}_{0,\veca}
=\bigoplus h^{(1)}_{0,\veca,\veck}$ 
of $\nbigq_0\nbige_{\veca|\nbigx-\nbigd}$
and $h^{(1)}_0:=
 \bigoplus h^{(1)}_{0,\veca}$
of $\nbigq_0\nbige_{0|\nbigx-\nbigd}$
as above.
We have the decomposition:
\[
f=
 \sum f_{(\veca,\veck),(\veca',\veck')},
\quad\quad
 f_{(\veca,\veck),(\veca',\veck')}
\in
 \Hom\bigl(
 U_{0,\veca',\veck'},\,
 U_{0,\veca,\veck}
 \bigr)
\]
We have only to show 
\begin{equation}
\label{eq;08.8.5.9}
\bigl|
 f_{(\veca,\veck),(\veca',\veck')}
\bigr|_{h^{(1)}_0}
=O\bigl(\Lambda_0(\epsilon)\bigr)
\end{equation}
for any $(\veca,\veck)$ and $(\veca',\veck')$
on $\nbigk\times Z(C)$.
Note that 
the induced metrics on 
$ \Hom\bigl(
 U_{0,\veca',\veck'},\,
 U_{0,\veca,\veck}
 \bigr)_{|\nbigx-\nbigd}$ are of the form
\begin{equation}
 \label{eq;08.8.5.11}
 g_{(\veca,\veck),(\veca',\veck')}
 \cdot
 \prod_{j=1}^{\ell}
 |z_j|^{2(-a_j+a_j')}
\cdot
 \prod_{j=1}^{\ell-1}
\left(
 \frac{-\log|z_j|}{-\log|z_{j+1}|}
\right)^{k_j-k_j'}
\cdot
\bigl(-\log|z_{\ell}|\bigr)^{k_{\ell}-k_{\ell}'},
\end{equation}
where $g_{(\veca,\veck),(\veca',\veck')}$
are $C^{\infty}$-metrics
of $ \Hom\bigl(
 U_{0,\veca',\veck'},\,
 U_{0,\veca,\veck}
 \bigr)$ on $\nbigx$.

\vspace{.1in}
{\bf (I)} Let us consider the case
$\veca\neq\veca'$.
We define
\[
 I_+:=\bigl\{
 i\,\big|\,a_i>a_i'
 \bigr\},
\quad
 I_-:=\bigl\{
 i\,\big|\,a_i<a_i'
 \bigr\},
\quad
 I_0:=\bigl\{
 i\,\big|\,a_i=a_i'
 \bigr\}.
\]
Let $m$ be the number determined by
$\{1,\ldots,m\}\subset I_0$
and $m+1\not\in I_0$.
Since the parabolic filtrations are preserved,
we have
$f_{(\veca,\veck),(\veca',\veck')|\nbigd_{i}}=0$
for any $i\in I_+$.
Hence, 
there exists a holomorphic section
$f'_{(\veca,\veck),(\veca',\veck')}$
of $ \Hom\bigl(
 U_{0,\veca',\veck'},\,
 U_{0,\veca,\veck}
 \bigr)$
such that
\begin{equation}
 \label{eq;08.8.5.12}
 f_{(\veca,\veck),(\veca',\veck')}
=f_{(\veca,\veck),(\veca',\veck')}'\cdot
 \prod_{i\in I_+}z_i
\end{equation}
We have the following inequality
for some $\epsilon>0$:
\begin{equation}
 \label{eq;08.8.5.13}
 \prod_{i\in I_+}
 |z_i|^{1-a_i+a_i'}
\cdot
 \prod_{i\in I_-}
 |z_i|^{-a_i+a_i'}
\leq
 \prod_{i\in I_+\cup I_-}
 |z_i|^{\epsilon}
\leq
 \bigl|z_{m+1}\bigr|^{\epsilon}.
\end{equation}
Let us consider the set
$S=\bigl\{
 p\leq m\,\big|\,
 k_p>k_p'
 \bigr\}$.
If $S$ is not empty,
let $p$ be the minimum.
Note that 
$k_t\leq k_t'$ for any $t<p$
and $k_p>k_p'$ by our choice.
Since the weight filtrations 
$W(\jbar)$ $(j=1,\ldots,p)$
are preserved
on $\lefttop{\pbar}\Gr^{F}$,
we have 
$f'_{(\veca,\veck),(\veca',\veck')|\nbigd_{\pbar}}=0$.
Hence, there exist holomorphic sections
$f''_{t,(\veca,\veck),(\veca',\veck')}$ $(t=1,\ldots,p)$
of $ \Hom\bigl(
 U_{0,\veca',\veck'},\,
 U_{0,\veca,\veck}
 \bigr)$ such that
\begin{equation}
 \label{eq;08.8.5.14}
 f'_{(\veca,\veck),(\veca',\veck')}
=\sum_{t=1}^p
 z_t\cdot
 f''_{t,(\veca,\veck),(\veca',\veck')}.
\end{equation}
We remark the following for any $t\leq p$:
\begin{multline}
 \label{eq;08.8.5.15}
  |z_t|\cdot
 \prod_{j=1}^{\ell-1}
 \left(
 \frac{-\log|z_j|}{-\log|z_{j+1}|}
 \right)^{k_j-k_j'}
 \!\!\!
 \bigl(-\log|z_{\ell}|\bigr)^{k_{\ell}-k_{\ell}'}
\leq
 |z_t|\cdot
\prod_{j=1}^{t-1}C^{k_j-k_{j}'}
\prod_{j=t}^{\ell-1}
 \left(
 \frac{-\log|z_j|}{-\log|z_{j+1}|}
 \right)^{k_j-k_j'}
 \!\!\!
 \bigl(-\log|z_{\ell}|\bigr)^{k_{\ell}-k_{\ell}'}
 \\
=O\bigl(|z_t|^{1/2}\bigr)
\end{multline}
By using (\ref{eq;08.8.5.11}), (\ref{eq;08.8.5.12}),
(\ref{eq;08.8.5.13}),
(\ref{eq;08.8.5.14}) and (\ref{eq;08.8.5.15}),
we obtain 
$\bigl|
 f_{(\veca,\veck),(\veca',\veck')}
\bigr|_{h^{(1)}_0}
=\sum_{t=1}^p
 O\bigl(|z_t|^{1/2}\bigr)
=O\bigl(\Lambda_0(1/2)\bigr)$.

If $S$ is empty,
we have $k_{j}\leq k_{j}'$ for $j=1,\ldots,m$.
Hence, we have the following:
\begin{multline}
 \label{eq;08.8.5.16}
  |z_{m+1}|^{\epsilon}\cdot
 \prod_{j=1}^{\ell-1}
 \left(
 \frac{-\log|z_j|}{-\log|z_{j+1}|}
 \right)^{k_j-k_j'}
 \!\!\!
 \bigl(-\log|z_{\ell}|\bigr)^{k_{\ell}-k_{\ell}'} 
 \\
\leq
 |z_{m+1}|^{\epsilon}\cdot
 \prod_{j=1}^{m}
 C^{k_j-k_j'}
\prod_{j=m+1}^{\ell-1}
 \left(
 \frac{-\log|z_j|}{-\log|z_{j+1}|}
 \right)^{k_j-k_j'}
 \!\!\!
 \bigl(-\log|z_{\ell}|\bigr)^{k_{\ell}-k_{\ell}'}
=O\bigl(|z_{m+1}|^{\epsilon/2}\bigr)
\end{multline}
By using (\ref{eq;08.8.5.11}),
(\ref{eq;08.8.5.12}) 
(\ref{eq;08.8.5.13})
and (\ref{eq;08.8.5.16}),
we obtain 
(\ref{eq;08.8.5.9}).

\vspace{.1in}
{\bf (II)}
Let us consider the case $\veca=\veca'$.
Because
$\Gr^F_{\veca}\Phi_{\nbigk}
=\Gr^F_{\veca}\Phi'_{\nbigk}
=\nu_{\veca}$,
there exist holomorphic sections
$f_{i,\veca,(\veck,\veck')}$
of $\Hom\bigl(
 U_{0,\veca,\veck'},\,
 U_{0,\veca,\veck}
 \bigr)$ such that
\begin{equation}
 \label{eq;08.8.5.17}
 f_{(\veca,\veck),(\veca,\veck')}
=\sum_{i=1}^{\ell}
 z_i\cdot f_{i,\veca,(\veck,\veck')}
\end{equation}
Let us consider the set
$S=\bigl\{
 p\,\big|\,
 k_p>k_p'
 \bigr\}$.
If $S$ is not empty,
let $p$ be the minimum.
Note that 
$k_t\leq k_t'$ for any $t<p$
and $k_p>k_p'$ by our choice.
Since the weight filtrations 
$W(\jbar)$ $(j=1,\ldots,p)$
are preserved
on $\lefttop{\pbar}\Gr^{F}$,
we have
$f_{i,\veca,(\veck,\veck')|\nbigd_{\pbar}}=0$.
Hence, there exist holomorphic sections
$f'_{t,i,\veca,(\veck,\veck')}$ $(t=1,\ldots,p)$
of $ \Hom\bigl(
 U_{0,\veca,\veck'},\,
 U_{0,\veca,\veck}
 \bigr)$ such that
\begin{equation}
\label{eq;08.8.5.18}
 f_{i,\veca,(\veck,\veck')}
=\sum_{t=1}^p
 z_t\cdot
 f'_{t,i,\veca,(\veck,\veck')}.
\end{equation}
By using (\ref{eq;08.8.5.15}) and
(\ref{eq;08.8.5.18}),
we obtain 
$\bigl|
 f_{i,\veca,(\veck,\veck')}
\bigr|_{h^{(1)}_0}
=O\bigl(\Lambda_0(1/2)\bigr)$.

If $S$ is empty,
we have $k_{j}\leq k_{j}'$ for $j=1,\ldots,\ell$.
Hence, we have the following:
\begin{equation}
\label{eq;08.8.5.19}
 \prod_{j=1}^{\ell-1}
 \left(
 \frac{-\log|z_j|}{-\log|z_{j+1}|}
 \right)^{k_j-k_j'}
 \!\!\!
 \bigl(-\log|z_{\ell}|\bigr)^{k_{\ell}-k_{\ell}'} 
=O(1).
\end{equation}
Hence, 
we obtain 
$|f_{i,\veca,(\veck,\veck')}|_{h_0^{(1)}}=O(1)$.
By using (\ref{eq;08.8.5.17}),
we obtain 
(\ref{eq;08.8.5.9}).
Thus, we obtain
Lemma \ref{lem;08.8.5.20}.
\hfill\qed

\paragraph{Local isomorphism
 with a nice property}

\begin{lem}
\label{lem;08.8.5.7}
There exists a holomorphic isomorphism
$\Phi_{\nbigk}:
 \nbigq_0\nbige_{0|\nbigx}
\lrarr
 \nbigq_0\nbige_{|\nbigx}$
with the following property:
\begin{itemize}
\item
 It preserves the filtrations
 $\lefttop{i}F$ $(i=1,\ldots,\ell)$.
\item
 For each $\vecb\in \real^I$, 
 the induced map 
 $\bigoplus_{q_I(\veca)=\vecb}
 \nbigq_0\nbige_{\veca|\nbigd_I}
\lrarr
 \lefttop{I}\Gr^F_{\vecb}(\nbigq_0\nbige_{|\nbigd_I})$
preserves the weight filtrations
$W(\jbar)$ $(j=1,\ldots,m(I))$.
\item
For each $\veca\in\real^{\ell}$,
the induced map
$\nbigq_0\nbige_{\veca|\nbigd_{\ellsitabar}}
\lrarr
 \lefttop{\ellsitabar}
 \Gr^{F}_{\vecb}(\nbigq_0\nbige_{|\nbigd_{\ellsitabar}})$ 
 is equal to $\nu_{\veca}$.
\end{itemize}
\end{lem}
\pf
We take decompositions
(\ref{eq;08.8.5.5}) and (\ref{eq;08.8.5.6})
as in Lemma \ref{lem;08.8.5.4}.
We take an isomorphism
$\nutilde_{\veca,\veck}:
 U_{0,\veca,\veck}\simeq
 U_{\veca,\veck}$
such that
$\nutilde_{\veca,\veck|\nbigd_{\ellsitabar}}
 =\nu_{\veca,\veck}$.
We set
$\Phi_{\nbigk}:=
 \sum \nutilde_{\veca,\veck}$.
It is easy to check that
$\Phi_{\nbigk}$ has the desired property.
Thus, we obtain Lemma \ref{lem;08.8.5.7}.
\hfill\qed

\vspace{.1in}

By the norm estimate (Lemma \ref{lem;08.8.5.10}),
$\Phi_{\nbigk}$ and $\Phi_{\nbigk}^{-1}$
are bounded on $\nbigk\times Z(C)$.

\begin{lem}
\label{lem;08.8.5.27}
We have the following estimate
for some $\epsilon>0$
with respect to $h_0$
on $\nbigk\times Z(C)$:
\begin{equation}
 \label{eq;08.8.5.26}
 \Phi_{\nbigk}^{\ast}\nabla_{\lambda}
-\nabla_{\lambda,0}
=O\bigl(
 \Lambda_{0}(\epsilon)
 \bigr)
\end{equation}
\end{lem}
\pf
Let $F$ denote the left hand side of 
(\ref{eq;08.8.5.26}).
It is easy to observe that
$F$ satisfies the conditions
in Lemma \ref{lem;08.8.5.20}.
Hence, Lemma \ref{lem;08.8.5.27}
follows from Lemma \ref{lem;08.8.5.20}.
\hfill\qed

\vspace{.1in}

Let $\Phi_{\nbigk}$
and $\Phi_{\nbigk}'$ be morphisms
as in Lemma \ref{lem;08.8.5.7}.
We set $G:=\Phi_{\nbigk}^{-1}\circ\Phi_{\nbigk}'$.
\begin{lem}
\label{lem;08.8.5.25}
We have the following estimates
for some $\epsilon>0$ 
on $\nbigk\times Z(C)$:
\[
 \bigl| G-\id \bigr|_{h_0}
=O\bigl(\Lambda_0(\epsilon)\bigr),
\quad
 \bigl|
 \nabla_{\lambda,0}(\lambda^2\del_{\lambda})G
 \bigr|_{h_0}
=O\bigl(
 \Lambda_0(\epsilon)
 \bigr)
\]
\end{lem}
\pf
We have only to apply 
Lemma \ref{lem;08.8.5.20}
to $G-\id$ and 
$\nabla_{\lambda,0}(\lambda^2\del_{\lambda})G$.
\hfill\qed

\vspace{.1in}

Let $\sigma:\cnum_{\lambda}\lrarr\cnum_{\mu}$
given by $\sigma(\lambda)=-\lambdabar$.
The induced map
$\cnum_{\lambda}\times X
\lrarr
 \cnum_{\mu}\times X^{\dagger}$
is also denoted by $\sigma$.
\begin{lem}
\label{lem;08.8.5.21}
We can take a holomorphic isomorphism
$\Phi^{\dagger}_{\sigma(\nbigk)}:
 \nbigq_{<\vecdelta}
 \nbige^{\dagger}_{0|\sigma(\nbigx)}
\lrarr
 \nbigq_{<\vecdelta}
 \nbige^{\dagger}_{|\sigma(\nbigx)}$
satisfying the conditions
(i) it preserves the filtrations
$\lefttop{i}F$ $(i=1,\ldots,\ell)$,
(ii) the induced morphism
 on $\lefttop{I}\Gr^{F}_{-\veca}$
 preserves the weight filtrations
 $W(\jbar)$ $(j=1,\ldots, m(I))$,
(iii) the induced morphism on
 $\lefttop{\ellsitabar}\Gr^{F}_{-\veca}$
 is equal to the given one.
\end{lem}
\pf
It can be shown by the argument
in the proof of Lemma \ref{lem;08.8.5.7}.
More directly,
we have the isomorphisms
$\nbigq_{<\vecdelta}
 \nbige^{\dagger}_{0|\sigma(\nbigx)}
\simeq
 \sigma^{\ast}
 \bigl(
 \nbigq_0\nbige_{0|\nbigx}
 \bigr)^{\lor}$
and
$\nbigq_{<\vecdelta}
 \nbige^{\dagger}_{|\sigma(\nbigx)}
\simeq
 \sigma^{\ast}\bigl(
  \nbigq_{0}
 \nbige_{|\nbigx}
 \bigr)^{\lor}$,
and $\sigma^{\ast}(\Phi_{\nbigk})^{\lor}$
satisfies the conditions.
\hfill\qed

\begin{lem}
\label{lem;08.8.5.22}
Let $\Phi_{\nbigk}$ and $\Phi_{\sigma(\nbigk)}^{\dagger}$
satisfy the above conditions.
We set
\[
 H:=\nbigs_0-\nbigs\bigl(
 \Phi_{\nbigk}\otimes\sigma^{\ast}\Phi_{\sigma(\nbigk)}^{\dagger}
 \bigr):
 \nbigq_0\nbige_{0|\nbigx}
\otimes
 \sigma^{\ast}\bigl(
 \nbigq_{<\vecdelta}\nbige^{\dagger}_{0|\sigma(\nbigx)}
 \bigr)
\lrarr \nbigo_{\nbigx}
\]
Then,
$H=O\bigl(\Lambda_0(\epsilon)\bigr)$
with respect to $h_0$
for some $\epsilon>0$
on $\nbigk\times Z(C)$.
\end{lem}
\pf
If $\Phi^{\dagger}_{\sigma(\nbigk)}$
is given by $\sigma^{\ast}\Phi_{\nbigk}^{\lor}$,
we have $H=0$.
Hence, we have only to show that
the property is independent of the choice 
of $\Phi^{\dagger}_{\sigma(\nbigk)}$.

Let $\Phi^{\dagger}_{i,\sigma(\nbigk)}$
$(i=1,2)$ be as in Lemma \ref{lem;08.8.5.21}.
Note that $h$ and $h_0$ are mutually bounded
through $\Phi^{\dagger}_{1,\sigma(\nbigk)}$
on $\sigma(\nbigk)\times Z(C)$.
By using Lemma \ref{lem;08.8.5.22},
we obtain 
$\Phi^{\dagger}_{1,\sigma(\nbigk)}
-\Phi^{\dagger}_{2,\sigma(\nbigk)}
=O\bigl(\Lambda_0(\epsilon)\bigr)$
for some $\epsilon>0$
with respect to $h$ and $h_0$.
Then, we obtain 
$ \nbigs\circ\Bigl(
 \Phi_{\nbigk}\otimes\sigma^{\ast}
 \bigl(
 \Phi^{\dagger}_{1,\sigma(\nbigk)}
-\Phi^{\dagger}_{2,\sigma(\nbigk)}
 \bigr)
 \Bigr)
=O\bigl( \Lambda_0(\epsilon) \bigr)$
with respect to $h_0$.
Thus, the proof of Lemma \ref{lem;08.8.5.22}
is finished.
\hfill\qed

\paragraph{Local $C^{\infty}$-isomorphisms}

Let $\Phi_{\nbigk}^p$ $(p=0,\ldots,m)$
be as in Lemma \ref{lem;08.8.5.7},
and let $a_p$ $(p=0,\ldots,m)$ be non-negative 
$C^{\infty}$-functions on $\nbigk$
such that $\sum a_p=1$.
We set $\Phi_{\nbigk}:=
 \sum_{p=0}^m a_p\cdot \Phi_{\nbigk}^p$.
We also set 
$G:=(\Phi_{\nbigk}^0)^{-1}\circ\Phi_{\nbigk}$
and 
$G^p:=(\Phi_{\nbigk}^0)^{-1}\circ\Phi_{\nbigk}^p$.
By Lemma \ref{lem;08.8.5.25},
$|G^p-\id|_{h_0}=O\bigl(
 \Lambda_0(\epsilon)
 \bigr)$,
and hence
$|G-\id|_{h_0}=O\bigl(
 \Lambda_0(\epsilon)
 \bigr)$
for some $\epsilon>0$
on $\nbigk\times Z(C)$.

\begin{lem}
\label{lem;08.8.5.28}
The following estimate holds for
some $\epsilon>0$ with respect to $h_0$
on $\nbigk\times Z(C)$:
\[
 \Phi_{\nbigk}^{-1}\circ
 \nabla_{\lambda}(\lambda^2\del_{\lambda})
 \circ \Phi_{\nbigk}
-\nabla_{\lambda,0}(\lambda^2\del_{\lambda})
=O\bigl(
 \Lambda_0(\epsilon)
 \bigr)
\]
\end{lem}
\pf
We have the following equalities:
\begin{multline}
 \Phi_{\nbigk}^{-1} \circ
 \nabla_{\lambda,0}(\del_{\lambda})
 \circ\Phi_{\nbigk}
-\nabla_{\lambda,0}(\del_{\lambda})
=\bigl(\Phi_{\nbigk}^{-1}\circ\Phi^0_{\nbigk}\bigr)
 \circ
 (\Phi_{\nbigk}^0)^{-1}\circ
\nabla_{\lambda}(\del_{\lambda})
\circ\bigl(\Phi_{\nbigk}^0\bigr)
 \circ\bigl((\Phi_{\nbigk}^0)^{-1}\circ\Phi_{\nbigk}\bigr)
-\nabla_{\lambda,0}(\del_{\lambda})
 \\
=
 G^{-1}\circ\bigl(
 (\Phi_{\nbigk}^0)^{\ast}\nabla_{\lambda}(\del_{\lambda})
-\nabla_{\lambda,0}(\del_{\lambda})
 \bigr)\circ G
+G^{-1}\cdot
 \nabla_{\lambda,0}(\del_{\lambda})G
\end{multline}
By Lemma \ref{lem;08.8.5.27},
we have
$(\Phi^0_{\nbigk})^{\ast}
 \nabla_{\lambda}(\lambda^2\del_{\lambda})
-\nabla_{\lambda,0}(\lambda^2\del_{\lambda})=
O\bigl( \Lambda_0(\epsilon) \bigr)$.
We also have 
\[
 \nabla_{\lambda,0}(\lambda^2\del_{\lambda})G
=\sum
 \lambda^2\frac{\del a_p}{\del \lambda}
 \cdot (G^p-\id)
=O\bigl(
 \Lambda_0(\epsilon)
 \bigr)
\]
Thus, we obtain Lemma \ref{lem;08.8.5.28}.
\hfill\qed

\vspace{.1in}

Let $\Phi^{\dagger\,q}_{\sigma(\nbigk)}$
$(q=0,1,\ldots,m')$
be as in Lemma \ref{lem;08.8.5.21},
and let $b_q$ be non-negative 
$C^{\infty}$-functions on $\sigma(\nbigk)$
such that $\sum b_q=1$.
We set 
$\Phi^{\dagger}_{\sigma(\nbigk)}:=
\sum b_q\cdot \Phi^{\dagger\,q}_{\sigma(S)}$.

\begin{lem}
\label{lem;08.8.5.29}
We set $H:=\nbigs\bigl(
 \Phi_{\nbigk}\otimes
 \sigma^{\ast}(\Phi^{\dagger}_{\sigma(\nbigk)})
 \bigr)-\nbigs_0$.
Then, we have the following estimate
on $\nbigk\times Z(C)$
with respect to $h_0$
for some $\epsilon>0$:
\[
 H=O\bigl( 
 \Lambda_0(\epsilon)
 \bigr),
\quad
 \delbar_{\nbige^{\sankaku}_0,\proj^1}H
=O\bigl(
 \Lambda_0(\epsilon)
 \bigr)
\]
\end{lem}
\pf
It follows from Lemma \ref{lem;08.8.5.22}.
\hfill\qed

\paragraph{Construction 
 of an approximating map}

We take $0<R_1<R_2<1$.
We set 
$\nbigk_1:=\bigl\{\lambda\,\big|\,
 |\lambda|\leq R_2 \bigr\}$
and 
$\nbigk_2:=\bigl\{\lambda\,\big|\,
 R_1\leq|\lambda|\leq R_1^{-1}
 \bigr\}$.
We take a partition of unity 
$\bigl(\chi_{\nbigk_1},
 \chi_{\nbigk_2},
 \chi_{\sigma(\nbigk_1)}\bigr)$
on $\proj^1$
which subordinates to
$\{\nbigk_1,\nbigk_2,\sigma(\nbigk_1)\}$.

We take a holomorphic isomorphism
$\Phi_{\nbigk_1}:\nbigq_0\nbige_{0|\nbigk_1\times X}
 \lrarr \nbigq_0\nbige_{|\nbigk_1\times X}$
as in Lemma \ref{lem;08.8.5.7}.
Similarly,
we take a holomorphic isomorphism
$\Phi^{\dagger}_{\sigma(\nbigk_1)}:
 \nbigq_{<\vecdelta}
 \nbige^{\dagger}_{0|
 \sigma(\nbigk)\times X^{\dagger}}
 \lrarr \nbigq_{<\vecdelta}
 \nbige^{\dagger}_{|
 \sigma(\nbigk)\times X^{\dagger}}$
as in Lemma \ref{lem;08.8.5.21}.

We can take a flat isomorphism
$\Phi_{\nbigk_2}:
 \bigl(
  \nbige_0,\DDtilde_0^f
 \bigr)_{|\nbigk_2\times (X-D)}
\lrarr
 \bigl(\nbige,\DDtilde^f\bigr)
 _{|\nbigk_2\times (X-D)}$.
We may assume that $\Phi_{\nbigk_2}$
is extended to the isomorphisms
$\nbigq_0\nbige_{0|\nbigk_2\times X}
\simeq
 \nbigq_0\nbige_{|\nbigk_2\times X}$
and 
$\nbigq_{<\vecdelta}
 \nbige^{\dagger}_{0|\nbigk_2\times X^{\dagger}}
\simeq
 \nbigq_{<\vecdelta}
 \nbige^{\dagger}_{|\nbigk_2\times X^{\dagger}}$
equipped with the property in
Lemmas \ref{lem;08.8.5.7} and \ref{lem;08.8.5.21}.
We set 
\[
 \Phi:=
 \chi_{\nbigk_1}\cdot \Phi_{\nbigk_1}
+\chi_{\nbigk_2}\cdot\Phi_{\nbigk_2}
+\chi_{\sigma(\nbigk_1)}\cdot 
 \Phi^{\dagger}_{\sigma(\nbigk_1)}.
\]
By using Lemmas \ref{lem;08.8.5.28}
and \ref{lem;08.8.5.29},
we can check that
$\Phi$ satisfies the estimates 
in (\ref{eq;08.8.5.30}).
Thus, the proof of Proposition \ref{prop;08.8.5.2}
is finished.
\hfill\qed

\subsubsection{Proof of Proposition
 \ref{prop;08.8.4.10}}

\paragraph{Decomposition}

We have a decomposition
\begin{equation}
\label{eq;08.8.5.40}
 \nbigq_0\nbige_{\nbigx}
=\bigoplus
 _{\veca\in\Par(\nbigp_0\nbige^0,\ellsitabar)}
 U_{\veca}
\end{equation}
with the following property:
\begin{itemize}
\item
 For any subset $I\subset\ellsitabar$
 and $\vecb\in \Par(\nbigp_0\nbige^0,I)$,
 we put 
$\lefttop{I}U_{\vecb}
 =\bigoplus_{
 \veca\in q_I^{-1}(\vecb)}
 U_{\veca}$.
Then, the following holds
for any $\vecc\in\real^{I}$:
\begin{equation}
 \bigoplus_{\vecb\leq\vecc}
 \lefttop{I}U_{\vecb|\nbigd_I}
=\bigcap_{i\in I}
 \lefttop{i}F_{c_i}\bigl(
 \nbigq_0\nbige_{|\nbigd_I}
 \bigr)
\end{equation}
\end{itemize}

\paragraph{Weak norm estimate}

We take a $C^{\infty}$-frame
$h'_{\veca}$ of
$U_{\veca}$ in (\ref{eq;08.8.5.40}).
We set 
$h^{(2)}_{\veca}:=
 h'_{\veca}\cdot
 \prod_{j=1}^{\ell}
 |z_j|^{-2a_j}$.
We obtain a $C^{\infty}$-hermitian metric
$h^{(2)}=\bigoplus 
 h^{(2)}_{\veca}$ of
$\nbigq\nbige_{|\nbigx-\nbigd}$.
Proposition 8.70 of \cite{mochi2}
implies the following lemma.
\begin{lem}
\label{lem;08.11.3.1}
$h$ and $h^{(2)}$ are mutually bounded
up to log order, namely,
\[
 h^{(2)}\cdot
  C^{-1}\cdot
 \Bigl(\sum_{i=1}^{\ell} -\log|z_i|\Bigr)^{-N}
\leq h\leq 
 h^{(2)}\cdot
  C \cdot\Bigl(\sum_{i=1}^{\ell} -\log|z_i|\Bigr)^{N}
\]
holds for some $C>0$ and $N>0$.
\hfill\qed
\end{lem}

\paragraph{An estimate}

\begin{lem}
\label{lem;08.8.19.30}
Let $f$ be a holomorphic endomorphism
of $\nbigq_0\nbige_0$ satisfying the following
conditions:
\begin{itemize}
\item
 It preserves the filtrations
 $\lefttop{i}F$ $(i=1,\ldots,\ell)$.
\item
 For each $\veca\in\real^{\ell}$,
 the induced endomorphism 
 $\lefttop{\ellsitabar}\Gr^F_{\veca}(f)$ of
 $\nbigq\nbige_{\veca|\nbigd_{\ellsitabar}}$
 is $0$.
\end{itemize}
Then, we have 
$|f|_{h_0}=O\bigl(
 \Lambda(M,\epsilon)
 \bigr)$
for some $M>0$ and $\epsilon>0$.
\end{lem}
\pf
We take a decomposition of $\nbigq_0\nbige_0$
like (\ref{eq;08.8.5.40}).
Applying the weak norm estimate
to $(E_{\veca},\delbar_{\veca},\theta_{\veca},h_{\veca})$
with the decomposition (\ref{eq;08.8.5.6}),
we take a $C^{\infty}$-hermitian metric
$h^{(2)}_{\veca}$
of $\nbigq\nbige_{\veca|\nbigx-\nbigd}$,
and $h^{(2)}_0=\bigoplus h^{(2)}_{\veca}$
of $\nbigq\nbige_{0|\nbigx-\nbigd}$.
We have the decomposition:
\[
f=
 \sum f_{\veca,\veca'},
\quad\quad
 f_{\veca,\veca'}
\in
 \Hom\bigl(
 U_{0,\veca'},\,
 U_{0,\veca}
 \bigr)
\]
We have only to show 
$\bigl|
 f_{\veca,\veca'}
\bigr|_{h^{(2)}_0}
=O\bigl(\Lambda(M,\epsilon)\bigr)$
for any $\veca$ and $\veca'$.
Assume $\veca\neq\veca'$.
We define
\[
 I_+:=\bigl\{
 i\,\big|\,a_i>a_i'
 \bigr\},
\quad
 I_-:=\bigl\{
 i\,\big|\,a_i<a_i'
 \bigr\},
\quad
 I_0:=\bigl\{
 i\,\big|\,a_i=a_i'
 \bigr\}.
\]
Since the parabolic filtrations are preserved,
we have
$f_{\veca,\veca',|\nbigd_{i}}=0$
for any $i\in I_+$.
Hence, 
there exists a holomorphic section
$f'_{\veca,\veca'}$
such that
$f_{\veca,\veca'}
=f_{\veca,\veca'}'\cdot
 \prod_{i\in I_+}z_i$.
We have the inequality
as in (\ref{eq;08.8.5.13}).
Then, we obtain the desired estimate
for $f_{\veca,\veca'}$ in the case
$\veca\neq\veca'$.

If $\veca=\veca'$,
$f_{\veca,\veca|\nbigd_{\ellsitabar}}=0$.
Hence,
there are holomorphic sections
$f_{t,\veca}$
of $\Hom\bigl(\nbigq_0\nbige_{\veca},
 \nbigq_0\nbige_{\veca}\bigr)$
such that
$f_{\veca,\veca}=
 \sum z_t\cdot f_{t,\veca}$.
Because
$|f_{t,\veca}|_{h_0}=
O\Bigl(\bigl(
 \sum_{i=1}^{\ell}-\log|z_i| \bigr)^N
 \Bigr)$,
we obtain the desired estimate.
\hfill\qed

\paragraph{Local isomorphism with a nice property}

We can show the following lemma
by the argument in the proof of
Lemma \ref{lem;08.8.5.7}.

\begin{lem}
\label{lem;08.8.5.45}
There exists a holomorphic isomorphism
$\Phi_{\nbigk}:
 \nbigq\nbige_{0|\nbigx}
\lrarr
 \nbigq\nbige_{|\nbigx}$
such that
(i) it preserves the filtrations
 $\lefttop{i}F$ $(i=1,\ldots,\ell)$,
(ii) for each $\veca\in\real^{\ell}$,
the induced map
$\nbigq\nbige_{\veca|\nbigd_{\ellsitabar}}
\lrarr
 \lefttop{\ellsitabar}
 \Gr^{F}_{\vecb}(\nbigq\nbige_{|\nbigd_{\ellsitabar}})$ 
 is equal to $\nu_{\veca}$.

Similarly,
we can take a holomorphic isomorphism
$\Phi^{\dagger}_{\sigma(\nbigk)}:
 \nbigq_{<\vecdelta}
 \nbige^{\dagger}_{0|\sigma(\nbigx)}
\lrarr
 \nbigq_{<\vecdelta}
 \nbige^{\dagger}_{|\sigma(\nbigx)}$
satisfying the conditions
(i) it preserves the filtrations
$\lefttop{i}F$ $(i=1,\ldots,\ell)$,
(iii) the induced morphism on
 $\lefttop{\ellsitabar}\Gr^{F}_{-\veca}$
 is equal to the given one.
\hfill\qed
\end{lem}

By the weak norm estimate,
$\Phi_{\nbigk}$ and $\Phi_{\nbigk}^{-1}$
are bounded up to log order.
We can show the following lemma
by using Lemma \ref{lem;08.8.19.30}.

\begin{lem}
\label{lem;08.8.5.47}
We have
$\Phi_{\nbigk}^{\ast}\nabla_{\lambda}
-\nabla_{\lambda,0}
=O\bigl(
 \Lambda(M,\epsilon)
 \bigr)$
for some $\epsilon>0$ and $M>0$
with respect to $h_0$.
\hfill\qed
\end{lem}

Let $\Phi_{\nbigk}$
and $\Phi_{\nbigk}'$ be morphisms
as in Lemma \ref{lem;08.8.5.45}.
We set $G:=\Phi_{\nbigk}^{-1}\circ\Phi_{\nbigk}'$.
\begin{lem}
\label{lem;08.8.5.46}
We have the following estimates
for some $\epsilon>0$ and $M>0$:
\[
 \bigl| G-\id \bigr|_{h_0}
=O\bigl(\Lambda(M,\epsilon)\bigr),
\quad
 \bigl|
 \nabla_{\lambda,0}(\lambda^2\del_{\lambda})G
 \bigr|_{h_0}
=O\bigl(
 \Lambda(M,\epsilon)
 \bigr)
\]
\end{lem}
\pf
It follows from 
Lemma \ref{lem;08.8.19.30}.
\hfill\qed

\begin{lem}
\label{lem;08.8.5.48}
Let $\Phi_{\nbigk}$ and $\Phi_{\sigma(\nbigk)}^{\dagger}$
satisfy the above conditions.
We set
\[
 H:=\nbigs_0-\nbigs\bigl(
 \Phi_{\nbigk}\otimes\sigma^{\ast}\Phi_{\sigma(\nbigk)}^{\dagger}
 \bigr):
 \nbigq_0\nbige_{0|\nbigx}
\otimes
 \sigma^{\ast}\bigl(
 \nbigq_{<\vecdelta}\nbige^{\dagger}_{0|\sigma(\nbigx)}
 \bigr)
\lrarr \nbigo_{\nbigx}
\]
Then,
$H=O\bigl(\Lambda(M,\epsilon)\bigr)$
with respect to $h_0$
for some $\epsilon>0$ and $M>0$.
\end{lem}
\pf
It can be shown by the argument
in the proof of Lemma \ref{lem;08.8.5.22}.
\hfill\qed

\paragraph{Local $C^{\infty}$-isomorphisms}

Let $\Phi_{\nbigk}^p$ $(p=0,\ldots,m)$
be as in Lemma \ref{lem;08.8.5.45},
and let $a_p$ $(p=0,\ldots,m)$ be non-negative 
$C^{\infty}$-functions on $\nbigk$
such that $\sum a_p=1$.
We set $\Phi_{\nbigk}:=
 \sum_{p=0}^m a_p\cdot \Phi_{\nbigk}^p$.
We also set 
$G:=(\Phi_{\nbigk}^0)^{-1}\circ\Phi_{\nbigk}$
and 
$G^p:=(\Phi_{\nbigk}^0)^{-1}\circ\Phi_{\nbigk}^p$.
By Lemma \ref{lem;08.8.5.46},
$|G^p-\id|_{h_0}=O\bigl(
 \Lambda_0(\epsilon)
 \bigr)$,
and hence
$|G-\id|_{h_0}=O\bigl(
 \Lambda_0(\epsilon)
 \bigr)$
for some $\epsilon>0$ and $M>0$.

We can show the following estimate
by using an argument in the proof of
Lemma \ref{lem;08.8.5.28}
with Lemma \ref{lem;08.8.5.47}:
\begin{equation}
 \label{eq;08.8.5.50}
 \Phi_{\nbigk}^{-1}\circ
 \nabla_{\lambda}(\lambda^2\del_{\lambda})
 \circ \Phi_{\nbigk}
-\nabla_{\lambda,0}(\lambda^2\del_{\lambda})
=O\bigl(
 \Lambda(M,\epsilon)
 \bigr)
\end{equation}
Let $\Phi^{\dagger\,q}_{\sigma(\nbigk)}$
$(q=0,1,\ldots,m')$
be as in Lemma \ref{lem;08.8.5.21},
and let $b_q$ be non-negative 
$C^{\infty}$-functions on $\sigma(\nbigk)$
such that $\sum b_q=1$.
We set 
$\Phi^{\dagger}_{\sigma(\nbigk)}:=
\sum b_q\cdot \Phi^{\dagger\,q}_{\sigma(S)}$.
We set $H:=\nbigs\bigl(
 \Phi_{\nbigk}\otimes
 \sigma^{\ast}(\Phi^{\dagger}_{\sigma(\nbigk)})
 \bigr)-\nbigs_0$.
Then, we can show the following estimate
with respect to $h_0$
for some $\epsilon>0$ and $M>0$,
by using Lemma \ref{lem;08.8.5.48}:
\begin{equation}
 \label{eq;08.8.5.51}
 H=O\bigl( 
 \Lambda(M,\epsilon)
 \bigr),
\quad
 \delbar_{\nbige^{\sankaku}_0,\proj^1}H
=O\bigl(
 \Lambda(M,\epsilon)
 \bigr)
\end{equation}

\paragraph{Construction}

We take $0<R_1<R_2<1$.
We set 
$\nbigk_1:=\bigl\{\lambda\,\big|\,
 |\lambda|\leq R_2 \bigr\}$
and 
$\nbigk_2:=\bigl\{\lambda\,\big|\,
 R_1\leq|\lambda|\leq R_1^{-1}
 \bigr\}$.
We take a partition of unity 
$\bigl(\chi_{\nbigk_1},
 \chi_{\nbigk_2},
 \chi_{\sigma(\nbigk_1)}\bigr)$
on $\proj^1$
which subordinates to
$\{\nbigk_1,\nbigk_2,\sigma(\nbigk_1)\}$.

We take a holomorphic isomorphism
$\Phi_{\nbigk_1}:\nbigq\nbige_{0|\nbigk\times X}
 \lrarr \nbigq\nbige_{|\nbigk\times X}$
as in Lemma \ref{lem;08.8.5.45}.
Similarly,
we take a holomorphic isomorphism
$\Phi^{\dagger}_{\sigma(\nbigk_1)}:
 \nbigq_{<\vecdelta}
 \nbige^{\dagger}_{0|
 \sigma(\nbigk)\times X^{\dagger}}
 \lrarr \nbigq_{<\vecdelta}
 \nbige^{\dagger}_{|
 \sigma(\nbigk)\times X^{\dagger}}$
as in Lemma \ref{lem;08.8.5.45}.

We can take a flat isomorphism
$\Phi_{\nbigk_2}:
 \bigl(
  \nbige_0,\DDtilde_0^f
 \bigr)_{|\nbigk_2\times (X-D)}
\lrarr
 \bigl(\nbige,\DDtilde^f\bigr)
 _{|\nbigk_2\times (X-D)}$.
We set 
\[
 \Phi:=
 \chi_{\nbigk_1}\cdot \Phi_{\nbigk_1}
+\chi_{\nbigk_2}\cdot\Phi_{\nbigk_2}
+\chi_{\sigma(\nbigk_1)}\cdot 
 \Phi^{\dagger}_{\sigma(\nbigk_1)}.
\]
By using (\ref{eq;08.8.5.50})
and (\ref{eq;08.8.5.51}),
we can check that
$\Phi$ satisfies the estimates in (\ref{eq;08.8.4.21}).
Thus, the proof of Proposition \ref{prop;08.8.4.10}
is finished.
\hfill\qed

\section{An application to HS-orbit}
\label{section;08.8.20.45}
\subsection{Preliminary}

\subsubsection{Compatibility of
real structure and Stokes structure}

Let $X$ be a complex manifold.
We set $\nbigx:=\cnum_{\lambda}\times X$
and $\nbigx^0:=\{0\}\times X$.
Let $(H,H'_{\real},\nabla)$ be a TER-structure
on $\nbigx$.
We say that $H$ is unramifiedly pseudo-good
if the following holds:
\begin{itemize}
\item
We are given a 
good set of irregular values
$\Irr(\nabla)\subset
 M(\nbigx,\nbigx^0)/H(\nbigx)$
in the level $-1$.
Namely,
(i) any elements $\gminia$ of $\Irr(\nabla)$
are of the form
$\gminia=\lambda^{-1}\gminia'$
for some holomorphic functions $\gminia'$
on $X$,
(ii) $\gminia'-\gminib'$
are nowhere vanishing for distinct 
$\lambda^{-1}\gminia',
 \lambda^{-1}\gminib'\in \Irr(\nabla)$.
\item
$H$ has the formal decomposition
\[
 (H,\nabla)_{|\nbigxhat^0}
=\bigoplus_{\gminia\in\Irr(\nabla)}
 (\Hhat_{\gminia},\nablahat_{\gminia}),
\]
such that 
$\nablahat_{\gminia}-d\gminia$ is regular.
Note that they are not assumed to be logarithmic.
\end{itemize}
(See also Subsection \ref{subsubsection;08.8.12.4}.)
If $X$ is a point,
it means that $H$ requires no ramification
in the sense of \cite{Hertling-Sevenheck}.

By a classical theory
(see also Subsection \ref{subsubsection;08.8.12.4}),
we have the Stokes filtration $\nbigf^S$
indexed by $\bigl(\Irr(\nabla),\leq_S\bigr)$
for each small sector $S$ of $\nbigx-\nbigx^0$.
We say that 
the real structure and the Stokes structure
are compatible,
if the Stokes filtrations on any small sectors $S$
come from a flat filtration of $H'_{\real|S}$.
(See \cite{Katzarkov-Kontsevich-Pantev}.)

By taking Gr of $(H,\nabla)$
with respect to the Stokes filtrations,
we obtain a TE-structure
$\Gr_{\gminia}(H,\nabla)$
for $\gminia\in\Irr(\nabla)$.
As observed in \cite{Hertling-Sevenheck},
if the real structure and the Stokes structure
are compatible,
$\Gr_{\gminia}(H,\nabla)$
is enriched to a TER-structure
denoted by
$\Gr_{\gminia}(H,H'_{\real},\nabla)$.
If $(H,H'_{\real},\nabla)$
is enriched to
a TERP-structure
$(H,H'_{\real},\nabla,P,w)$,
$\Gr_{\gminia}(H,H'_{\real},\nabla)$
is also naturally enriched to
a TERP-structure denoted by
$\Gr_{\gminia}(H,H'_{\real},\nabla,P,w)$.

\paragraph{Another formulation}

In \cite{Hertling-Sevenheck},
a compatibility of real structure
and Stokes structure is 
formulated in a slightly different way.
Let us check that it is equivalent to the above.
For simplicity, we consider the case 
in which $X$ is a point.

Let $H$ be a vector bundle on $\cnum_{\lambda}$
with a meromorphic flat connection
$\nabla:
 H\lrarr 
 H\otimes \Omega^{1}_{\cnum_{\lambda}}(\ast 0)$
such that 
$H$ requires no ramification
 with the good set of irregular values 
 $\Irr(\nabla)\subset \lambda^{-1}\cdot\cnum$.
Take $\theta_0\in\real$ such that
$\Re(\gminia-\gminib)
 (r\cdot e^{\sqrt{-1}\theta_0})\neq 0$
for any distinct $\gminia,\gminib\in \Irr(\nabla)$.
Take a sufficiently small $\epsilon>0$,
and let us consider the sector
\[
 \nbigs:=\bigl\{
 r\cdot e^{\sqrt{-1}\theta}\,\big|\,
 \theta_0-\epsilon
\leq \theta\leq \theta_0+\pi+\epsilon
 \bigr\}
\]
Let $\nbigsbar$ denote the closure of
$\nbigs$ in the real blow up
$\cnumtilde_{\lambda}(0)\lrarr \cnum_{\lambda}$
along $0$.
Let $\nbigz:=\nbigsbar\cap \pi^{-1}(0)$.
As a version of Hukuhara-Turrittin theorem,
it is well known that
we have a {\em unique} flat decomposition
\begin{equation}
 \label{eq;08.8.7.1}
 (H,\nabla)_{|\nbigsbar}
=\bigoplus_{\gminia\in \Irr(\nabla)}
 \bigl(
 H_{\gminia,\nbigs},
 \nabla_{\gminia,\nbigs}
 \bigr)
\end{equation}
such that the restriction of (\ref{eq;08.8.7.1})
to $\widehat{\nbigz}$
is the same as the pull back 
of the irregular decomposition of 
$H_{|\widehat{0}}$.

Assume that the flat bundle
$(H,\nabla)_{|\cnum_{\lambda}^{\ast}}$
is equipped with a real structure,
i.e.,
a $\cnum$-anti-linear flat involution
$\kappa:H\lrarr H$.
In other words,
$(H,\nabla,\kappa)$ is a TER-structure.
In Section {\rm 8} of \cite{Hertling-Sevenheck},
the real structure and 
the Stokes structure are defined to be compatible,
if $\kappa(H_{\gminia,\nbigs})=H_{\gminia,\nbigs}$
for any $\gminia\in \Irr(\nabla)$ and 
any $\nbigs$ as above.

If a small sector $S$ is contained in $\nbigs$,
the restriction of (\ref{eq;08.8.7.1}) to $S$
gives a splitting of $\nbigf^S$.
Hence,
if $H_{\gminia,\nbigs}$ are preserved by $\kappa$
for any $\gminia$,
the filtration $\nbigf^S$ is also preserved by $\kappa$.
Let $S_1$ and $S_2$ be small sectors
containing the rays
$\{r\cdot e^{\sqrt{-1}\theta_0}\,|\,r> 0\}$ 
and $\{-r\cdot e^{\sqrt{-1}\theta_0} |\,r> 0\}$,
respectively.
Then, $\gminia\leq_{S_1}\gminib$
if and only if $\gminia\geq_{S_2}\gminib$.
By the parallel transform on $\nbigs$,
the flat bundle $H_{|\nbigs}$ is trivialized,
and we can observe that
$H_{\gminia,\nbigs}
=\nbigf^{S_1}_{\gminia}\cap
 \nbigf^{S_2}_{\gminia}$.
Hence, if $\nbigf_{\gminia}^{S_i}$ $(i=1,2)$
are preserved by $\kappa$,
$H_{\gminia,\nbigs}$ is also preserved by $\kappa$.
The equivalence of
two notions of compatibilities
follows from these considerations.

\subsubsection{Two Stokes filtrations of 
 integrable twistor structures}

Let $(V,\DDtilde^{\sankaku})$
be a variation of
integrable twistor structures over
$\proj^1\times X$.
It is obtained as the gluing of $TE$-structure
$(V_0,\DDtilde_0^f)$ on 
$\nbigx:=\cnum_{\lambda}\times X$
and $\Ttilde E$-structure
$(V_{\infty},\DDtilde_{\infty}^{\dagger\,f})$
on $\nbigx^{\dagger}:=\cnum_{\mu}\times X^{\dagger}$.
We set $\nbigx^0:=\{0\}\times X\subset\nbigx$
and $\nbigx^{\dagger\,0}:=\{0\}\times X^{\dagger}
\subset\cnum_{\mu}\times X^{\dagger}$.

\begin{df}
We say that
$(V,\DDtilde^{\sankaku})$ is 
unramifiedly pseudo-good,
if both $(V_0,\DDtilde_0^f)$
and $(V_{\infty},\DDtilde_{\infty}^{\dagger\,f})$
is unramifiedly pseudo-good.
In that case,
let $\Irr(\DDtilde^f_0)$ and
$\Irr(\DDtilde^{\dagger\,f}_{\infty})$
denote the sets of irregular values
of $\DDtilde^f_0$ and $\DDtilde^{\dagger\,f}_{\infty}$,
respectively.

If $X$ is a point,
it is also said that
$(V,\DDtilde^{\sankaku})$ requires
no ramification.
\hfill\qed
\end{df}

\begin{df}
\label{df;08.8.7.11}
Assume $(V,\DDtilde^{\sankaku})$
is unramifiedly pseudo-good.
\begin{itemize}
\item
We say that the sets of the irregular values of
$(V,\DDtilde^{\sankaku})$
are compatible,
 if $\Irr\bigl(\DDtilde_0^f\bigr)$
 and $\Irr\bigl(\DDtilde_{\infty}^f\bigr)$
 bijectively correspond
 by $\gminia
 \longleftrightarrow
 \overline{\gamma^{\ast}\gminia}$.
\item
We say that
$(V,\DDtilde^{\sankaku})$ has
compatible Stokes structures,
if the following holds:
\begin{itemize}
\item
 The sets of irregular values of $(V,\DDtilde^{\sankaku})$
 are compatible.
\item
 For a small sector $S$ of $\nbigx-\nbigx^0$,
 we have the Stokes filtration $\nbigf^S$
 of $(V_0,\DDtilde_0^f)$.
 We also have the Stokes filtration $\nbigf^{\gamma(S)}$
 of $(V_{\infty},\DDtilde_{\infty}^f)$,
 where we regard $\gamma(S)$
 as a small sector of 
 $\nbigx^{\dagger}-\nbigx^{\dagger\,0}$.
 Then,
 $\nbigf^S$ and $\nbigf^{\gamma(S)}$
 are the same
 under the parallel transform 
 along any rays
 connecting $S$ and $\gamma(S)$.
\hfill\qed
\end{itemize}
\end{itemize}
\end{df}

\begin{rem}
In the above definition,
a ray means a line
$\bigl\{(t\cdot e^{\sqrt{-1}\varphi},P)\,\big|\,
 0<t<\infty \bigr\}$
in $\cnum_{\lambda}^{\ast}\times \{P\}
\subset \cnum_{\lambda}^{\ast}\times X$.
We say that
it connects $S$ and $\gamma(S)$,
if (i) $(t\cdot e^{\sqrt{-1}\varphi},P)$
is contained in $S$ for any sufficiently small $t$,
(ii) $(t\cdot e^{\sqrt{-1}\varphi},P)$
is contained in $\gamma(S)$ for any sufficiently large $t$.
\hfill\qed
\end{rem}

\begin{lem}
If $(V,\DDtilde^{\sankaku})$ is equipped 
with either a real structure $\kappa$
or a perfect pairing $\nbigs$ of weight $w$,
then the irregular values of
$\DDtilde_0^f$ and $\DDtilde_{\infty}^f$
are compatible.
\end{lem}
\pf
We have
$\Irr\bigl(
 \gamma^{\ast}\DDtilde^f_{\infty}
 \bigr)
=\bigl\{
 \overline{\gamma^{\ast}\gminia}\,
 \big|\,
 \gminia\in\Irr(\DDtilde^f_{\infty})
 \bigr\}$.
If $(V,\DDtilde^{\sankaku})$
is equipped with a real structure,
$\gamma^{\ast}
 (V_{\infty},\DDtilde^f_{\infty})
\simeq
(V_0,\DDtilde^f_0)$.
Hence,
the irregular values of
$\DDtilde^f_0$ and $\DDtilde^f_{\infty}$
are compatible.

We have
$\Irr\bigl(
 \sigma^{\ast}\DDtilde^f_{\infty}
 \bigr)
=\bigl\{
 \overline{\sigma^{\ast}\gminia}\,\big|\,
 \gminia\in\Irr(\DDtilde^f_{\infty})
 \bigr\}$.
Note that $\gminia\in \Irr(\DDtilde^f_{\infty})$
are of the form $\mu^{-1}\gminia'$,
where $\gminia'$ are 
holomorphic functions on $X^{\dagger}$.
Hence, 
$\overline{\sigma^{\ast}\gminia}
=-\overline{\gamma^{\ast}\gminia}$.
If $(V,\DDtilde^{\sankaku}_f)$ is equipped 
with a perfect pairing,
$(V_0,\DDtilde^f_0)$  is isomorphic to
the dual of 
$\sigma^{\ast}(V_{\infty},\DDtilde^f_{\infty})$.
Therefore, the irregular values of
$\DDtilde^f_0$ and $\DDtilde^f_{\infty}$
are compatible.
\hfill\qed

\vspace{.1in}

If $(V,\DDtilde^{\sankaku})$ is
unramifiedly pseudo-good,
we obtain $TE$-structure
$\Gr_{\gminia}(V_0,\DDtilde_0^f)$ 
on $\nbigx$
for $\gminia\in \Irr(\DDtilde_0^f)$,
and $\Ttilde E$-structure
$\Gr_{\gminib}(V_{\infty},\DDtilde_{\infty}^f)$
on $\nbigx^{\dagger}$
for $\gminib\in\Irr(\DDtilde_{\infty}^f)$,
by taking Gr with respect to the Stokes filtrations.
If $(V,\DDtilde^{\sankaku})$ 
has compatible Stokes structures,
we have the natural isomorphism
\[
 \Gr_{\gminia}(V_0,\DDtilde_0^f)_{|
 \nbigx-\nbigx^0
 }
\simeq
 \Gr_{\overline{\gamma^{\ast}\gminia}}
 \bigl(V_{\infty},\DDtilde_{\infty}^f\bigr)_{|
 \nbigx^{\dagger}-\nbigx^{\dagger\,0}}.
\]
Hence, 
we obtain a variation of integrable twistor structures
$\Gr_{\gminia}(V,\DDtilde^{\sankaku})$
for each $\gminia\in\Irr(\DDtilde^f_0)$
as the gluing of them.
We have the following functoriality
(Lemma \ref{lem;08.8.12.10}).

\begin{lem}
Let $(V^{(a)},\DDtilde^{(a)\sankaku})$
be unramifiedly pseudo-good.
Assume (i) $(V^{(a)},\DDtilde^{(a)\sankaku})$
$(a=1,2)$
have compatible Stokes filtrations,
(ii) the union 
$\nbigi:=\Irr(\DDtilde^{(1)\,f}_{0})\cup
 \Irr(\DDtilde^{(2)\,f}_{0})$
is good.
Then, a morphism 
$(V^{(1)},\DDtilde^{(1)\,\sankaku})
\lrarr 
 (V^{(2)},\DDtilde^{(2)\,\sankaku})$
induces 
$\Gr_{\gminia}\bigl(V^{(1)},
 \DDtilde^{(1)\,\sankaku}
 \bigr)
\lrarr
 \Gr_{\gminia}\bigl(V^{(2)},
 \DDtilde^{(2)\sankaku}\bigr)$
for each $\gminia\in \nbigi$.
\hfill\qed
\end{lem}

We have the natural isomorphisms
\[
 \gamma^{\ast}\Gr_{\gminia}(V,\DDtilde^{\sankaku})
\simeq
 \Gr_{\gminia}\bigl(
 \gamma^{\ast}(V,\DDtilde^{\sankaku})
 \bigr),
\quad
 \sigma^{\ast}\Gr_{\gminia}(V,\DDtilde^{\sankaku})
\simeq
 \Gr_{-\gminia}\bigl(
 \sigma^{\ast}(V,\DDtilde^{\sankaku})
 \bigr).
\]
The following lemma follows from functoriality.
\begin{lem}
Assume $(V,\DDtilde^{\sankaku})$
has compatible Stokes structures.
If $(V,\DDtilde^{\sankaku})$
is equipped with a real structure,
(resp. a perfect pairing of weight $w$),
each $\Gr_{\gminia}(V,\DDtilde^{\sankaku})$
is also equipped with
an induced real structure
(resp. an induced perfect pairing of weight $w$).
\hfill\qed
\end{lem}

\begin{lem}
\label{lem;08.8.8.2}
Let $(H,H'_{\real},\nabla,P',-w)$
be a variation of TERP-structures,
and let $(V,\DDtilde^{\sankaku},\nbigs,\kappa,-w)$
be the corresponding variation of 
twistor-TERP structures.
(See Subsection 
{\rm\ref{subsubsection;08.9.11.30}}
for the correspondence.)
Assume that
 $(H,H'_{\real},\nabla,P',-w)$
 is unramifiedly pseudo-good,
or equivalently,
 $(V,\DDtilde^{\sankaku},\nbigs,\kappa,-w)$
 is unramifiedly pseudo-good.
\begin{itemize}
\item
 The real structure and the Stokes structures of
 $(H,\nabla)$ are compatible,
 if and only if
 $(V,\DDtilde^{\sankaku})$ has compatible
  Stokes structures.
\item
 If the real structures and the Stokes structures are compatible,
 $\Gr_{\gminia}(V,\DDtilde^{\sankaku},\nbigs,\kappa,-w)$
 is the variation of twistor-TERP structures
 corresponding to
 $\Gr_{\gminia}(H,H'_{\real},\nabla,P',-w)$.
\end{itemize}
\end{lem}
\pf
Note that the Stokes filtrations of 
$\gamma^{\ast}(H,\nabla)$ on $\gamma^{\ast}(S)$
is given by the composite of 
the conjugate with respect to $H_{\real}'$
and the parallel transport along the rays
connecting $S$ and $\gamma(S)$,
with the change of the index sets from 
$\Irr(\nabla)$ to $\bigl\{
 \overline{\gamma^{\ast}\gminia}\,\big|\,
 \gminia\in\Irr(\nabla) \bigr\}$.
Then, the first claim follows.

Let us consider the second claim.
We have only to consider the case $w=0$.
We may assume that
$(H,H'_{\real},\nabla,P')$ is obtained from
 $(V,\nabla,\nbigs,\kappa)$
by the procedure explained in
Subsection \ref{subsubsection;08.9.11.30}.
By construction,
we have
$\Gr_{\gminia}(H,\nabla)
=\Gr_{\gminia}(V_0,\nabla_0)$.
For comparison of
induced real structures and pairings,
we have only to consider the case
in which $X$ is a point.

Let us compare the induced real structures.
The flat real structure of $H'$ is obtained
as the composite:
\[
\begin{CD}
 \overline{H}_{|\lambda}
@>{\rm parallel\,\, transform}>>
 \overline{H}_{|\lambdabar^{-1}}
@>{\kappa_{|\lambda}}>>
H_{|\lambda}
\end{CD}
\]
Hence, we have the following factorization
of the real structure on
$\Gr_{\gminia}(H)_{|\lambda}$
obtained as Gr of the Stokes filtration:
\[
 \begin{CD}
 \overline{\Gr_{\gminia}(H)}_{|\lambda}
@>{\rm parallel\,\, transform}>>
 \overline{\Gr_{\gminia}(H)}_{|\lambdabar^{-1}}
@>{\Gr_{\gminia}(\kappa)_{|\lambda}}>>
 \Gr_{\gminia}(H)_{|\lambda}
\end{CD}
\]
It is the same as the real structure
induced by $\Gr_{\gminia}(\kappa)$
on $\Gr_{\gminia}(V,\nabla)$.

Let $P:H\otimes j^{\ast}H\lrarr 
 \nbigo_{\cnum_{\lambda}}$ be 
the pairing induced by $\kappa$ and $\nbigs$
as in (\ref{eq;08.7.21.12}),
whose restriction to $H'$ is $P'$.
Let $S$ be a small sector in
$\cnum_{\lambda}^{\ast}$.
We have the following factorization of $P_{|S}$:
\[
\begin{CD}
 \nbigf^S_{\gminia}(H)
\otimes
 j^{\ast}\nbigf^{j(S)}_{\gminib}(H)
=
\nbigf^S_{\gminia}(V_0)
\otimes
 \sigma^{\ast}
 \gamma^{\ast}\nbigf^{j(S)}_{\gminib}(V_0)
@>{1\otimes\sigma^{\ast}\kappa}>>
 \nbigf^S_{\gminia}(V_0)
\otimes
 \sigma^{\ast}
 \nbigf^{\sigma(S)}_{\overline{\gamma^{\ast}(\gminib)}}
(V_{\infty})
@>{\nbigs}>>
 \nbigo_{S}
\end{CD}
\]
The restriction to 
 $\nbigf^S_{\gminia}(H)\otimes
 j^{\ast}\nbigf^{j(S)}_{\gminib}(H)$
is $0$
unless $\gminia-\gminib\geq_S 0$.
The induced pairing 
$P_{\gminia}$
for $\Gr_{\gminia}(V_0)$
is factorized as follows:
\[
 \begin{CD}
 \Gr_{\gminia}(V_0)_{|S}
\otimes
 j^{\ast}\Gr_{\gminia}(V_0)_{|j(S)}
@>{1\otimes\sigma^{\ast}\Gr_{\gminia}\kappa}>>
 \Gr_{\gminia}(V_0)_{|S}
\otimes
 \sigma^{\ast}
 \Gr_{\overline{\gamma^{\ast}(\gminia)}}
 (V_{\infty})_{|\sigma(S)}
@>{\Gr_{\gminia}\nbigs}>>
 \nbigo_{S}
\end{CD}
\]
Hence, it is the same as
the pairing induced by
$\Gr_{\gminia}(V,\nabla,\nbigs,\kappa)$.
Thus, the proof of Lemma \ref{lem;08.8.8.2}
is finished.
\hfill\qed

\subsubsection{Preliminary for pull back}

We set $X:=\cnum_z$,
$D=\bigl\{0\bigr\}$,
$\nbigx:=\cnum_{\lambda}\times X$,
$\nbigd:=\cnum_{\lambda}\times D$
and $W:=\nbigd\cup(\{0\}\times X)$.
Let $\pi:\nbigxtilde(W)\lrarr \nbigx$ be
a real blow up of $\nbigx$ along $W$.
Let $\pi_{1}:
 \cnumtilde_{\lambda}(0)\lrarr\cnum_{\lambda}$
be the real blow up of $\cnum_{\lambda}$
along $\{0\}$.
Let $\phi_0:\nbigx\lrarr \cnum_{\lambda}$ be
given by $\phi_0(\lambda,z)=\lambda\cdot z$.
It induces the map
$\phitilde_0:\nbigxtilde(W)\lrarr
 \cnumtilde_{\lambda}(0)$.

Let $H$ be a vector bundle on
$\cnum_{\lambda}$
with a meromorphic flat connection
$\nabla:
 H\lrarr 
 H\otimes \Omega^{1}_{\cnum_{\lambda}}(\ast 0)$
such that 
$(H,\nabla)$ requires no ramification
 with the good set of irregular values 
 $\nbigi\subset\cnum\cdot \lambda^{-1}$.
Let $\gbigv$ denote the flat bundle on 
$\cnumtilde_{\lambda}(0)$
associated to $H_{|\cnum_{\lambda}^{\ast}}$.
For each $Q\in \pi_{1}^{-1}(0)$,
we have the Stokes filtration
$\nbigf^Q$ of $\gbigv_{|Q}$
for the meromorphic prolongment $H$.
(See Subsection {\rm\ref{subsubsection;08.8.8.3}})
We can naturally regard
$\phitilde_0^{\ast}\gbigv$
as the flat bundle on $\nbigxtilde(W)$
associated to $(\phi_0^{\ast}H)_{|\nbigx-W}$.

\begin{lem}
\label{lem;08.8.7.20}
The following holds:
\begin{itemize}
\item
 $\phi_0^{\ast}(H,\nabla)$ is 
 unramifiedly pseudo-good in the level 
 $\vecm=(-1,-1)$.
 (See Subsection {\rm\ref{subsubsection;08.8.12.4}}.)
 The set of irregular values is given by
 $\phi_0^{\ast}\nbigi:=
 \bigl\{
 \phi_0^{\ast}\gminia\,\big|\,
 \gminia\in\nbigi
 \bigr\}$.
\item
 For each $P\in \pi^{-1}(W)$,
 the Stokes filtration $\nbigf^P$ of
 $\phitilde_0^{\ast}(\gbigv)_{|P}$
 for $\phi_0^{\ast}H$
 is the pull back of the Stokes filtration of
 $\gbigv_{|\phitilde_0(P)}$ .
\item
 We have the natural isomorphism
 $\phi_0^{\ast}\Gr_{\gminia}(H)
\simeq
 \Gr_{\phi_0^{\ast}\gminia}\bigl(
 \phi_0^{\ast}H \bigr)$.
\end{itemize}
\end{lem}
\pf
We have the decomposition
$(H,\nabla)_{|\widehat{0}}
=\bigoplus_{\gminia\in\nbigi}
 (H_{\gminia},\nablahat_{\gminia})$,
where $\nablahat_{\gminia}-d\gminia$
are regular.
It induces the decomposition of
$\phi_0^{\ast}(H,\nabla)_{|\What}$.
Hence, the first claim is clear.

We set $Q:=\phitilde_0(P)$.
Note that the orders $\leq_Q$
and $\leq_P$ are the same
under the identification
$\nbigi\simeq\phi_0^{\ast}\nbigi$.
Let $H_1\supset H$ be an unramifiedly good lattice.
Then, $\phi_0^{\ast}H_1$ is an unramifiedly good lattice.
We take a small sector
$S_{Q}\in
 \Multisector(Q,\cnum_{\lambda}^{\ast},\nbigi)$
such that there exists the Stokes filtration
$\nbigf^{S_Q}$ of $H_{1|\Sbar_Q}$.
We take a small multi-sector
$S_P\in
 \Multisector(P,\nbigx-W,\phi^{\ast}_0\nbigi)$
such that
$\phi_0(S_P)\subset S_Q$.
Then, we obtain the filtration
$\phitilde_0^{\ast}\nbigf^{S_Q}$
of $\phi_0^{\ast}(H_1)_{|\Sbar_P}$
indexed by
$\bigl(\phi_0^{\ast}\nbigi,\leq_P\bigr)$.
It gives the Stokes filtration
of $\phi_0^{\ast}(H_1)_{|\Sbar_P}$,
which follows from the characterization
in Proposition \ref{prop;07.6.16.6}.
Since the filtration of 
$\phitilde_0^{\ast}(\gbigv)_{|P}$
induced by $\phitilde_0^{\ast}\nbigf^{S_Q}$
is the same as the pull back of 
$\nbigf^Q$ on $\gbigv_{|Q}$,
we obtain the second claim.
Note that we also obtain that
the Stokes filtration of 
$\phi_0^{\ast}(H)_{|\Sbar_P}$
is given by the pull back of
the Stokes filtration of $H_{|\Sbar_Q}$.

Let $S_P$ be a small multi-sector as above.
By the above compatibility of the Stokes filtrations
and Lemma \ref{lem;08.8.20.1},
we obtain the natural isomorphisms
\begin{equation}
 \label{eq;08.8.7.15}
 \phi_0^{\ast}\bigl(
 \Gr_{\gminia}(H)\bigr)_{|\Sbar_P}
\simeq 
 \Gr_{\phi_0^{\ast}\gminia}(\phi_0^{\ast}H)
 _{|\Sbar_P}.
\end{equation}
By varying $S_P$ and gluing them,
we obtain
$\phi_0^{\ast}\bigl(
 \Gr_{\gminia}(H) \bigr)_{|\nbigutilde(W)}
\simeq
 \Gr_{\phi_0^{\ast}\gminia}
 (\phi_0^{\ast}H)_{|\nbigutilde(W)}$,
 where $\nbigu$ is a neighbourhood of $W$,
and $\nbigutilde(W)$
denote the real blow up of $\nbigu$
along $W$.
By using the flatness,
it is extended to 
$\phi_0^{\ast}\bigl(
 \Gr_{\gminia}(H) \bigr)_{|\nbigxtilde(W)}
\simeq
 \Gr_{\phi_0^{\ast}\gminia}
 (\phi_0^{\ast}H)_{|\nbigxtilde(W)}$.
Hence, we obtain an isomorphism on $\nbigx$.
\hfill\qed

\subsubsection{Rescaling and HS-orbit}
\label{subsubsection;08.8.20.2}

We recall a rescaling construction
in \cite{Hertling} and \cite{Hertling-Sevenheck}.
See also \cite{sabbah10}.
We set $X:=\cnum_z$,
$D=\bigl\{0\bigr\}$
and $X^{\ast}:=X-D$.
For $R>0$,
we set $X(R):=\bigl\{
 z\in X\,\big|\,|z|<R
 \bigr\}$
and $X^{\ast}(R):=X(R)\cap X^{\ast}$.
We set $\nbigx:=\cnum_{\lambda}\times X$.
We use the symbols
$\nbigx^{\ast}$, $\nbigd$,
$\nbigx(R)$ and $\nbigx^{\ast}(R)$
in similar meanings.
Let $\phi_0:\nbigx\lrarr \cnum_{\lambda}$ be
given by $\phi_0(\lambda,z)=\lambda\cdot z$.
The restriction to $\nbigx^{\ast}$
is denoted by $\psi_{0}$.

\paragraph{TERP-structure}

We consider only TERP-structures of weight $0$.
Hence, we omit to specify weights.
Let $(H,H'_{\real},\nabla,P)$ be a TERP-structure.
Hertling and Sevenheck studied
the variation of TERP-structures
$\psi_0^{\ast}(H,H'_{\real},\nabla,P)$
on $X^{\ast}$.
If there exists an $R>0$
such that 
$\psi_0^{\ast}(H,H'_{\real},\nabla,P)_{|X^{\ast}(R)}$
is pure and polarized,
the variation is called an HS-orbit 
(Hertling-Sevenheck orbit),
and we say in this paper
that $(H,H'_{\real},\nabla,P)$ induces an HS-orbit.
\begin{rem}
An HS-orbit is called a ``nilpotent orbit''
in {\rm\cite{Hertling-Sevenheck}}.
We use ``HS-orbit''
for distinction from twistor nilpotent orbit.
It matches their terminology ``Sabbah-orbit''.
\hfill\qed
\end{rem}

\begin{lem}
We assume (i) $(H,\nabla)$ requires no ramification,
(ii) the Stokes structure
and the real structure of $(H,H'_{\real},\nabla)$
are compatible.
Then, the following holds:
\begin{itemize}
\item
 $\psi_0^{\ast}(H,\nabla)$ 
 is unramifiedly pseudo-good.
 The set of irregular values is given by
 $\bigl\{
 \psi_0^{\ast}\gminia\,\big|\,
 \gminia\in\Irr(\nabla)
 \bigr\}$.
\item
The real structure and the Stokes structure of
$\psi^{\ast}_0(H,\nabla)$ are compatible.
\item
We have the natural isomorphism
$\psi_0^{\ast}\Gr_{\gminia}
 \bigl(H,H'_{\real},\nabla,P\bigr)
\simeq
 \Gr_{\psi_0^{\ast}\gminia}\psi_0^{\ast}
 \bigl(H,H'_{\real},\nabla,P\bigr)$.
\end{itemize}
\end{lem}
\pf
The first two claims follow from
Lemma \ref{lem;08.8.7.20}.
To show the third claim,
we have only to compare the induced flat pairings.
It can be done directly,
or by considering the restriction
to $\cnum_{\lambda}\times \{1\}$.
\hfill\qed

\paragraph{Integrable twistor structure}

We set 
$\nbigx^{\dagger}:=\cnum_{\mu}\times X^{\dagger}$,
$\nbigd^{\dagger}:=\cnum_{\mu}\times D^{\dagger}$,
$\nbigx^{\ast\,\dagger}
 :=\nbigx^{\dagger}-\nbigd^{\dagger}$
and $W^{\dagger}:=
 \nbigd^{\dagger}\cup\bigl(\{0\}\times X^{\dagger}\bigr)$.
Let $\phi_{\infty}:\nbigx^{\dagger}
\lrarr \cnum_{\mu}$ be given by
$\phi_{\infty}(\mu,z)=\mu\cdot\zbar$.
The restriction to $\nbigx^{\ast\, \dagger}$
is denoted by $\psi_{\infty}$.

Let $(V,\nabla)$ be an integrable twistor
structure on $\proj^1$
which requires no ramification.
It is obtained as the gluing of
$(V_0,\nabla_0)$ and $(V_{\infty},\nabla_{\infty})$.
The gluing is denoted by
$g:V_{0|\cnum_{\lambda}^{\ast}}
\simeq
 V_{\infty|\cnum_{\mu}^{\ast}}$,
which is flat with respect to $\nabla$.

We set $\HS(V)_0:=\psi_0^{\ast}(V_0)$
and $\HS(V)_{\infty}:=\psi_{\infty}^{\ast}(V_{\infty})$.
They are naturally equipped with 
$TE$-structure $\HS(\nabla)_0$
and $\Ttilde E$-structure 
$\HS(\nabla)_{\infty}$.
Note that
$\HS(V,\nabla)_0$ and
$\HS(V,\nabla)_{\infty}$
are unramifiedly pseudo good.
Let us construct
a flat isomorphism $\Phi$ between
$\HS(V,\nabla)_{0|
 \cnum_{\lambda}^{\ast}\times X^{\ast}}$
and 
$\HS(V,\nabla)_{\infty|
 \cnum_{\mu}^{\ast}\times X^{\dagger\ast}}$.
The fibers $\HS(V)_{0|(\lambda,z)}$
and $\HS(V)_{\infty|(\mu,z)}$
are naturally identified with
$V_{0|\lambda\cdot z}$
and $V_{\infty|\mu\cdot \zbar}$,
respectively.
If $\lambda=\mu^{-1}$,
we have 
$(\lambda\cdot z)^{-1}
=\mu\cdot \zbar\cdot |z|^{-2}$.
Hence, we have an isomorphism
$\Phi_{(\lambda,z)}:
 H(V)_{0|(\lambda,z)}
\simeq
 H(V)_{\infty|(\lambda^{-1},z)}$
induced by the gluing $g$ with
the parallel transform along the segments 
connecting $\lambda^{-1}\cdot \zbar$
and $\lambda^{-1}\cdot \zbar\cdot |z|^{-2}$.
Thus, we obtain the isomorphism $\Phi$
as desired.

Let $\HS(V,\nabla)$ denote 
the variation of integrable twistor structures
obtained as the gluing of
$\HS(V,\nabla)_0$ and $\HS(V,\nabla)_{\infty}$.
The following lemma is clear from the construction
and the functoriality (Lemma \ref{lem;08.8.12.10}).
\begin{lem}
\mbox{{}}
\begin{itemize}
\item
Let $F:(V^{(1)},\nabla^{(1)})\lrarr
 (V^{(2)},\nabla^{(2)})$
be a morphism of integrable pure twistor structures.
Then, we have the induced morphisms
$\HS(F):\HS(V^{(1)},\nabla^{(1)})
\lrarr \HS(V^{(2)},\nabla^{(2)})$.
\item
Let $f$ be $\gamma$ or $\sigma$.
Then,
$\HS\circ f^{\ast}(V,\nabla)$
is naturally isomorphic to
$f^{\ast}\HS(V,\nabla)$.
\hfill\qed
\end{itemize}
\end{lem}

By the above lemma,
a real structure $\kappa$ of $(V,\nabla)$
induces a real structure 
$\HS(\kappa)$ of $\HS(V,\nabla)$.
Since 
we have the natural isomorphism
$\HS\bigl(\Tate(0)\bigr)
\simeq
 \Tate(0)_{X^{\ast}}$,
a paring $\nbigs$ of $(V,\nabla)$ with weight $0$
induces 
a pairing $\HS(\nbigs)$
of $\HS(V,\nabla)$ with weight $0$.
Hence,
an integrable twistor structure with a pairing
$(V,\nabla,\nbigs)$
induces 
$\HS(V,\nabla,\nbigs)$
on $\proj^1\times X^{\ast}$,
and if $(V,\nabla,\nbigs)$ is
 equipped with a real structure,
$\HS(V,\nabla,\nbigs)$ is also 
equipped with a naturally induced
real structure.

\begin{lem}
\label{lem;08.8.20.5}
Assume that 
$(V,\nabla)$ has compatible Stokes structures.
Then, 
$\HS(V,\nabla)$ also has 
compatible Stokes structures,
and we have the natural isomorphism
\begin{equation}
 \label{eq;08.8.7.25}
\HS\Gr_{\gminia}(V,\nabla)
\simeq
\Gr_{\psi_0^{\ast}\gminia}\HS(V,\nabla)
\end{equation}
If $(V,\nabla)$ is equipped
with a pairing of weight $0$
(resp. a real structure),
{\rm(\ref{eq;08.8.7.25})}
preserves the induced pairings
(resp. real structures).
\end{lem}
\pf
It follows from Lemma \ref{lem;08.8.7.20}.
\hfill\qed

\begin{lem}
\label{lem;08.8.20.3}
Let $(H,H'_{\real},\nabla,P')$ be
a TERP-structure,
and let $(V,\nabla,\nbigs,\kappa)$
be the corresponding twistor-TERP structure.
Then, 
$\HS(V,\nabla,\nbigs,\kappa)$
is the variation of twistor-TERP structure
corresponding to 
$\psi_0^{\ast}
 (H,H'_{\real},\nabla,P')$.
\end{lem}
\pf
By construction,
we have the natural isomorphism
$\HS(V,\nabla)_0\simeq (H,\nabla)$.
We have only to compare 
the induced real structures and pairings
on them.
Since they are flat,
we have only to compare them
on the fiber over $z=1$.
Then, the claim is clear.
\hfill\qed

\vspace{.1in}

If there exists an $R>0$
such that
$\HS(V,\nabla,\nbigs)_{|\proj^1\times X^{\ast}(R)}$ 
is pure and polarized,
it is called a twistor HS-orbit,
and we say that
$(V,\nabla,\nbigs)$ induces a twistor HS-orbit.

\subsection{Reduction of wild HS-orbit}

\subsubsection{Statement}

We use the notation in 
Subsection \ref{subsubsection;08.8.20.2}.
Let $(V,\nabla)$ be 
an integrable twistor structure
with a perfect pairing $\nbigs$
of weight $0$,
which requires no ramification.
Assume that
$(V,\nabla,\nbigs)$ induces
a twistor HS-orbit
on $\proj^1\times X^{\ast}(R)$ for some $R>0$.
We obtain the underlying 
unramifiedly good wild harmonic bundle
$(E,\delbar_E,\theta,h)$ on $X^{\ast}(R)$
of $\HS(V,\nabla,\nbigs)_{|\proj^1\times X^{\ast}(R)}$,
which is unramifiedly good.
Let $\nbigi$ denote the set of irregular values
of $(V,\nabla)$ at $0$.
It is easy to see 
\[
 \Irr(\theta)=\bigl\{
 \gminia(z)\,\big|\,
 \gminia(\lambda)\in \nbigi
\bigr\}\simeq\nbigi.
\]
We will not distinguish them in the following.

Let $(\nbige^{\sankaku},\DD^{\sankaku},\nbigs_E)$
denote the variation of polarized pure twistor structure
associated to $(E,\delbar_E,\theta,h)$.
It is enriched to integrable one
$(\nbige^{\sankaku},\DDtilde^{\sankaku},\nbigs_E)$.
Although it is naturally isomorphic to $\HS(V,\nabla,\nbigs)$,
it is non-trivial that
the natural meromorphic extensions
$\nbigq\nbige_0$ and 
$\phi_0^{\ast}(V_0)\otimes\nbigo_{\nbigd}(\ast \nbigd)$
are isomorphic.
Hence, we use the symbol
$(\nbige^{\sankaku},\DDtilde^{\sankaku},\nbigs_E)$
for distinction.
By applying the construction 
in Subsection \ref{subsubsection;08.7.26.200}
to $(\nbige^{\sankaku},\DDtilde^{\sankaku},\nbigs_E)$,
we obtain a wild variation of
pure polarized integrable twistor structures
$\Gr_{\gminia}(\nbige^{\sankaku},
 \DDtilde^{\sankaku},
 \nbigs_E)$
for each $\gminia\in\nbigi$.
We will prove the following theorem
in Subsection \ref{subsubsection;08.8.7.31}.

\begin{thm}
\mbox{{}}\label{thm;08.8.7.30}
\begin{itemize}
\item
 $(V,\nabla)$ has compatible
 Stokes structures.
\item
 $\HS\Gr_{\gminia}(V,\nabla,\nbigs)$
 is naturally isomorphic to
 $\Gr_{\gminia}(\nbige^{\sankaku},
 \DDtilde^{\sankaku},
 \nbigs_E)$ for each $\gminia\in\nbigi$.
In particular,
$\Gr_{\gminia}(V,\nabla,\nbigs)$
induces a twistor HS-orbit.
\end{itemize}
\end{thm}

Before going into the proof,
we give a consequence.

\begin{cor}
\label{cor;08.9.12.5}
Let $(H,H'_{\real},\nabla,P,0)$ be a TERP
structure which requires no ramification.
If $(H,H'_{\real},\nabla,P,0)$ induces an HS-orbit,
it is a mixed-TERP structure in the sense of 
Definition {\rm 9.1} of {\rm \cite{Hertling-Sevenheck}}.
\end{cor}
\pf
Thanks to Theorem 9.3 of \cite{Hertling-Sevenheck},
$(H,H'_{\real},\nabla,P,0)$ is a mixed-TERP structure,
if and only if
(i) the real structure and the Stokes structure of
 $(H,H'_{\real},\nabla)$ are compatible,
(ii) $\Gr_{\gminia}(H,H'_{\real},\nabla,P,0)$
induces an HS-orbit for each $\gminia\in\Irr(\nabla)$.
Hence, this corollary
follows from Theorem \ref{thm;08.8.7.30},
Lemma \ref{lem;08.8.8.2}
and Lemma \ref{lem;08.8.20.3}.
\hfill\qed

\vspace{.1in}
The claim of the corollary was established
by Hertling and Sevenheck
{\rm\cite{Hertling-Sevenheck}}
in the case that $(H,\nabla)$
has regular singularity.
They also showed the converse of the claim
 in general.

\begin{rem}
In their study of the case 
that $(H,\nabla)$ has regular singularity,
Hertling and Sevenheck
closely investigated the limiting object.
In particular,
they showed that the limiting TERP-structure
is generated by elementary sections,
for which the eigenvalues of
the new supersymmetric index
can be described in terms of
the Hodge filtrations of the corresponding
mixed Hodge structure.

Even in the irregular case,
the limiting object can be obtained 
from the reduced regular one.
Hence, the limit of the eigenvalues
of the new supersymmetric index
of $\phi_0^{\ast}(H,H'_{\real},\nabla,P)$
can be described in terms of 
their mixed Hodge structures.
\hfill\qed
\end{rem}

\subsubsection{Proof of Theorem 
\ref{thm;08.8.7.30}}
\label{subsubsection;08.8.7.31}

We have the natural identifications
$ \HS(V,\nabla)_0
\simeq\bigl(
 \nbige,\DDtilde^f
 \bigr)$
and 
$ \HS(V,\nabla)_{\infty}
\simeq\bigl(
 \nbige^{\dagger},\DDtilde^{\dagger\,f}
 \bigr)$.
We have the following
locally free $\nbigo_{\nbigx}(\ast \nbigd)$-modules
\[
 \HSbar(V)_0:=
 \phi_0^{\ast}(V_0)
 \otimes \nbigo_{\nbigx}(\ast\nbigd).
\]
We also have the following
locally free 
$\nbigo_{\nbigx^{\dagger}}
 (\ast \nbigd^{\dagger})$-modules
\[
 \HSbar(V)_{\infty}:=
 \phi_{\infty}^{\ast}(V_{\infty})
 \otimes \nbigo_{\nbigx^{\dagger}}(\ast\nbigd^{\dagger}).
\]

\paragraph{Comparison of
 $\nbigq\nbige$ and $\HSbar(V)_0$}

We would like to show that
$\nbigq\nbige$ and $\HSbar(V)_0$
are naturally isomorphic.
We set $W:=\nbigd
 \cup
\bigl(\{0\}\times X\bigr)$.
Let $\pi:\nbigxtilde(W)\lrarr \nbigx$ be 
the real blow up of $\nbigx$ along $W$.
Let $\gbigv$ be the flat bundle
on $\nbigxtilde(W)$ associated to
$(\nbige,\DDtilde^f)_{|\nbigx-W}$.
We set $\phi_0^{\ast}\nbigi:=\bigl\{
 \phi_0^{\ast}\gminia\,\big|\,
 \gminia\in\nbigi\bigr\}$.

As remarked in Lemma \ref{lem;08.8.7.20},
$\phi_0^{\ast}H$ is a pseudo-good lattice
of $\HSbar(V)_0\otimes \nbigo_{\nbigx}(\ast W)$
in the level $\vecm=(-1,-1)$.
\begin{lem}
$\nbigq_0\nbige$ is a good lattice
of $\nbigq\nbige\otimes\nbigo_{\nbigx}(\ast W)$
in the level $\vecm=(-1,-1)$
around $(\lambda,z)=(0,0)$.
\end{lem}
\pf
We have the decomposition
$(\nbigq_0\nbige,\DD)_{|\What}
=\bigoplus
 _{\gminia\in\nbigi}
 (\nbigq_0\nbigehat_{\gminia},\DDhat_{\gminia})$
such that
$\DDhat_{\gminia}-d_X\gminia$ is logarithmic.
We have the corresponding decomposition
$\DDtilde^f
=\bigoplus \DDtilde^f_{\gminia}$.
Let us show that
\[
 \Bigl(
 \DDtilde^f_{\gminia}-d\bigl(\gminia(z)/\lambda\bigr)
\Bigr)
 \nbigq\nbigehat_{\gminia}
\subset
 \lambda^{-1}\cdot
 \nbigq\nbigehat_{\gminia}
\otimes
 \Omega_{\nbigx}^{1,0}(\log W)
\]
Since $\DD_{\gminia}-d_X\gminia$ is logarithmic,
it is satisfied for the derivatives along the $X$-direction.
Since 
we have already known that
the restriction of $\nbigq_0\nbige$ to
$\cnum_{\lambda}^{\ast}\times X$
is unramifiedly good
by Lemma \ref{lem;08.7.26.125},
it is also satisfied for the derivatives
along the $\lambda$-direction.
\hfill\qed

\vspace{.1in}

Let $P\in \pi^{-1}(0,0)$.
We have the Stokes filtration
$\nbigf^{P}_1$ of $\gbigv_{|P}$
corresponding to the meromorphic prolongment 
$\nbigq\nbige\otimes\nbigo_{\nbigx}(\ast W)$,
and the Stokes filtration $\nbigf^{P}_2$ of
$\gbigv_{|P}$
corresponding to the meromorphic prolongment
$\HSbar(V)_0\otimes\nbigo_{\nbigx}(\ast W)$.
(See Subsection \ref{subsubsection;08.8.8.3}
for such filtrations
in the pseudo-good case.)

\begin{lem}
\label{lem;08.8.7.40}
$\nbigf^{P}_1=\nbigf^P_2$.
\end{lem}
\pf
Let $S_P\in \Multisector(P,\nbigx-W,\phi_0^{\ast}\nbigi)$ 
be a small sector such that 
there exist the Stokes filtrations
$\nbigf^{S_P}_1$ of
$\nbigq\nbige_{|\Sbar_P}$
and $\nbigf^{S_P}_2$
of $\phi_0^{\ast}(H)_{|\Sbar_P}$.
We can take 
$Q\in \Sbar_P\cap \pi^{-1}(W\setminus\nbigd)$
such that the orders
$\leq_Q$ and $\leq_P$ on $\nbigi$
are the same.
We have the filtrations
$\nbigf^Q_i$ $(i=1,2)$ of $\gbigv_{|Q}$
corresponding to the meromorphic prolongments
$\nbigq\nbige(\ast W)$
and $\HSbar(V)_0(\ast W)$.
Because $ \HS(V,\nabla)_0
\simeq\bigl( \nbige,\DDtilde^f \bigr)$,
we have $\nbigf^Q_1=\nbigf^{Q}_2$.

Let us show that $\nbigf^{P}_i$ is obtained
as the parallel transport of $\nbigf^Q_i$,
which implies $\nbigf^P_1=\nbigf^P_2$.
We take
$S_Q\in\Multisector(Q,\nbigx-W,\phi_0^{\ast}\nbigi)$
such that there exist
the Stokes filtrations
$\nbigf^{S_Q}_1$ of
$\nbigq\nbige_{|\Sbar_Q}$
and $\nbigf^{S_Q}_2$
of $\phi_0^{\ast}(H)_{|\Sbar_Q}$.
By using the characterization in
Proposition \ref{prop;08.8.12.5},
we obtain 
$\bigl(\nbigf^{S_P}_i\bigr)_{|\Sbar_Q}
=\nbigf^{S_Q}_i$.
Hence, we can conclude that
$\nbigf^P_i$ are obtained as 
the parallel transport of 
$\nbigf^Q_i$.
\hfill\qed

\begin{lem}
The isomorphism
$\nbige\simeq \HS(V)_0$ on $\nbigx-\nbigd$
is extended to the isomorphism
$\nbigq\nbige\simeq \HSbar(V)_0$
on $\nbigx$.
\end{lem}
\pf
Let $P\in\pi^{-1}(0,0)$.
We take a small multi-sector 
$S_P\in\Multisector(P,\nbigx-W,\nbigi)$
such that
we have the Stokes filtrations $\nbigf^{S_P}$
for $\nbigq_0\nbige_{|\Sbar_P}$ 
and $\phi_0^{\ast}(H)_{|\Sbar_P}$.
By Lemma \ref{lem;08.8.7.40},
the restrictions of them to $S_P$ are the same.
We take a flat splitting 
$\nbige_{|S_P}=
 \bigoplus_{\gminia\in\nbigi}
 \nbige_{\gminia,S}$,
which is extended to the decompositions:
\[
 \nbigq\nbige_{|\Sbar_P}
=\bigoplus_{\gminia\in\nbigi} \nbigq\nbige_{\gminia,S},
\quad
 \phi_0^{\ast}(H)_{|\Sbar_P}
=\bigoplus
 \phi_0^{\ast}(H)_{\gminia,S}
\]
Let $\nbigltilde(-\gminia)$ be a line bundle
$\nbigo_{\nbigx}(\ast W)\cdot e$
with $\nabla e=
 e\cdot\bigl(-d(\lambda^{-1}\gminia)\bigr)$.
We remark that 
$\Gr_{\gminia}(\phi_0^{\ast}H)
 \otimes\nbigltilde(-\gminia)$
and 
$\Gr_{\gminia}(\nbigq\nbige)
 \otimes\nbigltilde(-\gminia)$
have regular singularity along $W$.
Hence, the isomorphism on $\nbigx-W$
is naturally extended to the isomorphism
$\Gr_{\gminia}(\phi_0^{\ast}H)
 \otimes\nbigltilde(-\gminia)
\simeq
\Gr_{\gminia}(\nbigq\nbige)
 \otimes\nbigltilde(-\gminia)$.
Since the restrictions of 
$\Gr_{\gminia}(\phi_0^{\ast}H)
 \otimes \nbigo(\ast \nbigd)$
and 
$\Gr_{\gminia}(\nbigq\nbige)$
to $\nbigx-\nbigd$ are naturally isomorphic,
we obtain the isomorphism
$\Gr_{\gminia}(\phi_0^{\ast}H)
 \otimes \nbigo(\ast \nbigd)
\simeq\Gr_{\gminia}(\nbigq\nbige)$.

Let $\vecw_{\gminia}$
and $\vecv_{\gminia}$ be frames of
$\Gr_{\gminia}(\phi_0^{\ast}H)
 \otimes\nbigo(\ast \nbigd)$ and 
$\Gr_{\gminia}(\nbigq\nbige)$,
respectively.
We have the relation
$\vecw_{\gminia}
=\vecv_{\gminia}\cdot A_{\gminia}$,
where $A_{\gminia}$ are meromorphic
along $\nbigd$.
We take lifts $\vecw_{\gminia,S}$
and $\vecv_{\gminia,S}$ to 
$\nbige_{\gminia,S}$
by using the above splittings
for any small sectors in $\nbigx-W$.
We have the relation
$\vecw_{\gminia,S}
=\vecv_{\gminia,S}\cdot A_{\gminia}$.
Then, the isomorphism
$\nbige\simeq \HS(V)_0$
is extended to
$\nbigq\nbige_{|\nbigu}
\simeq
 \HSbar(V)_{0|\nbigu}$
on some small neighbourhood 
$\nbigu$ of $(0,0)$,
which follows from Proposition 
\ref{prop;08.8.12.11}.
(We may apply Lemma \ref{lem;08.11.6.15}.
 But, since $\phi_0^{\ast}V_0$
 may not be a good lattice,
 we replace $V_0$ with an unramifiedly good lattice,
 or we use a variant of Lemma \ref{lem;08.11.6.15}
 for a pseudo-good lattice.)
Then, it is easy to observe
that the isomorphism is extended to
$\nbigq\nbige
\simeq
 \HSbar(V)_0$
by using Hartogs theorem.
(Sabbah also independently 
obtained an argument
to extend such isomorphisms
in this kind of situation.)
\hfill\qed

\vspace{.1in}

Similarly,
$\nbige^{\dagger}\simeq \HS(V)_{\infty}$
on $\nbigx^{\dagger}-\nbigd^{\dagger}$
is extended to the isomorphism 
$ \nbigq\nbige^{\dagger}\simeq \HSbar(V)_{\infty}$
on $\nbigx^{\dagger}$.

\paragraph{Proof of the first claim}

Let $\nbigxtilde^{\sankaku}$
denote the real blow up
of $\proj^1\times X$ along 
$(\proj^1\times D)\cup
  (\{0\}\times X)\cup
 (\{\infty\}\times X)$.
Let $\pi^{\dagger}:
 \nbigxtilde^{\dagger}(W^{\dagger})
\lrarr
 \nbigx^{\dagger}$
denote the real blow up 
of $\cnum_{\mu}\times X^{\dagger}$ along 
$W^{\dagger}=\nbigd^{\dagger}
 \cup
 (\{0\}\times X^{\dagger})$.
We have
\[
 \nbigxtilde^{\sankaku}
=\nbigxtilde(W)\cup
 \nbigxtilde^{\dagger}(W^{\dagger}).
\]
Let $\gbigv^{\sankaku}$
denote the flat bundle on $\nbigxtilde^{\sankaku}$
associated to the flat bundle
$(\nbige,\DDtilde^f)_{|
 \cnum_{\lambda}^{\ast}\times (X-D)}$.

We have the $C^{\infty}$-map
$\nbigx-W\lrarr 
(\real_{\geq\,0}\times S^1)^2$
given by
\[
 (\lambda,z)\longmapsto
 \Bigl(
 \bigl(|\lambda|,\,
 \lambda/|\lambda|
 \bigr),\,\,
  \bigl(|z|,\, z/|z| \bigr)
 \Bigr)
\]
It induces the natural identification
$\nbigxtilde(W)\simeq
 (\real_{\geq\,0}\times S^1)^2$.
We set 
\[
 P_0=\Bigl(
 \bigl(0,\exp(\sqrt{-1}\varphi)
 \bigr),
 (1,1) \Bigr)
\in 
\pi^{-1}\bigl((0,1)\bigr)
 \subset \nbigxtilde(W).
\]

Similarly,
we identify
$\nbigxtilde^{\dagger}(W^{\dagger})$
with 
$(\real_{\geq\, 0}\times S^1)^2$
via the map induced by
\[
 (\mu,z)\longmapsto
 \Bigl(
 \bigl(|\mu|,\,\mu/|\mu|\bigr),\,\,
 \bigl(|z|,\,z/|z|\bigr)
 \Bigr).
\]
We set
$ Q_0:=\Bigl(
 \bigl(0,\exp(-\sqrt{-1}\varphi)\bigr),
 (1,1)
\Bigr)
\in (\pi^{\dagger})^{-1}\bigl((0,1)\bigr)
\subset
 \nbigxtilde^{\dagger}(W^{\dagger})$.
Note that we can identify
$(V,\nabla)$ with
$\HS(V,\nabla)_{|\proj^1\times\{1\}}$.
Hence, we have only to compare the Stokes filtrations
$\nbigf^{P_0}(\gbigv^{\sankaku}_{|P_0})$
and $\nbigf^{Q_0}(\gbigv^{\sankaku}_{|Q_0})$
under the parallel transport
along the ray 
$\Bigl(
 \bigl(s,\exp(\sqrt{-1}\varphi)\bigr),
 (1,1)
\Bigr)$ $(s\in \real_{\geq\,0}\cup\{+\infty\})$
connecting
$P_0$ and $Q_0$.
(Note that the signature of the arguments
are reversed by the coordinate change
$\lambda^{-1}=\mu$.)

\vspace{.1in}

Let us consider the map
$G:[0,1]\times[0,1]\lrarr \nbigxtilde(W)$
given by
\[
 G(s,t)=\Bigl(
 \bigl(s,\exp(\sqrt{-1}\varphi)\bigr),
 (t,1)
 \Bigr).
\]
Note $G(0,1)=P_0$.
We set $P_1:=G(1,0)$ and $P_2:=G(1,1)$.
The image of
$\Gamma_0:=
 \bigl([0,1]\times\{0\}\bigr)
\cup\bigl(\{0\}\times [0,1]\bigr)$
is contained in $\pi^{-1}(W)$.
The orders $\leq_P$ 
are independent of 
$P\in G(\Gamma_0)$.
Hence, the Stokes filtrations
are unchanged along $G(\Gamma_0)$.

Similarly, let us consider the map
$G^{\dagger}:
 [0,1]\times[0,1]\lrarr 
 \nbigxtilde^{\dagger}(W^{\dagger})$
given by
\[
 G^{\dagger}(s,t)=\Bigl(
 \bigl(s,\exp(-\sqrt{-1}\varphi)\bigr),
 (t,1)
 \Bigr).
\]
Note $G^{\dagger}(0,1)=Q_0$.
We set $Q_1:=G^{\dagger}(1,0)$ 
and $Q_2:=G^{\dagger}(1,1)$.
The image of
$\Gamma_{\infty}:=
 \bigl([0,1]\times\{0\}\bigr)
\cup\bigl(\{0\}\times [0,1]\bigr)$
is contained in 
$(\pi^{\dagger})^{-1}(W^{\dagger})$.
The orders $\leq_Q$ are independent
of the choice of
$Q\in G^{\dagger}(\Gamma_{\infty})$.
Hence, the Stokes filtrations
are unchanged along 
$G^{\dagger}(\Gamma_{\infty})$.

Under the identification
$\nbigx-W=\nbigx^{\dagger}-W^{\dagger}$,
we have $P_2=Q_2$,
and the union of the paths
$G([0,1]\times\{1\})$
and $G^{\dagger}([0,1]\times\{1\})$
is the ray connecting $P_0$ and $Q_0$.
Hence, 
for the comparison of
$\nbigf^{P_0}$ and $\nbigf^{Q_0}$,
we have only to show that
$\nbigf^{P_1}(\gbigv_{0|P_1})$ of $\gbigv_{|P_1}$
and $\nbigf^{Q_1}(\gbigv_{\infty|Q_1})$ of $\gbigv_{|Q_1}$
are the same.
It follows from the characterization
of the Stokes filtrations of
$(\nbigq\nbige,\DD)$
and $(\nbigq\nbige^{\dagger},\DD^{\dagger})$
by growth order of the norms of flat sections
with respect to the metric $h$.
(See Subsection \ref{subsection;08.8.3.6}.)
Thus, we obtain the first claim 
of Theorem \ref{thm;08.8.7.30}.

\paragraph{Proof of the second claim}

By using Corollary \ref{cor;08.8.7.5}
and Lemma \ref{lem;08.8.20.5},
we obtain the isomorphisms
on $\proj^1\times X^{\ast}(R)$
for some $R>0$:
\[
 \Gr_{\gminia}
 (\nbige^{\sankaku},
 \DDtilde^{\sankaku},
 \nbigs_E)
\simeq
 \Gr_{\psi_0^{\ast}\gminia}\HS(V,\nabla,\nbigs)
\simeq
 \HS\Gr_{\gminia}(V,\nabla,\nbigs)
\]
Thus, the second claim is also proved.
\hfill\qed

\noindent
{\em Address\\
Research Institute for Mathematical Sciences,
Kyoto University,
Kyoto 606-8502, Japan\\
takuro@kurims.kyoto-u.ac.jp
}

\end{document}